\documentclass[11pt]{article} 
\usepackage{amsmath,amssymb,amsfonts,amsthm,amscd,graphicx,psfrag,epsfig}
\usepackage{color}
\definecolor{blue}{rgb}{0,0,0.7}
\definecolor{red}{rgb}{0.75, 0, 0}
\usepackage{stmaryrd}
\usepackage[titletoc,toc]{appendix}
\usepackage{comment}
\usepackage[all]{xy} 

\newtheorem{theorem}{Theorem}[section]
\newtheorem{lemma}[theorem]{Lemma}
\newtheorem{proposition}[theorem]{Proposition}
\newtheorem{corollary}[theorem]{Corollary}
\newtheorem{conjecture}[theorem]{Conjecture}
\newtheorem{definition}[theorem]{Definition}

\newtheorem{example}[theorem]{Example}

\begin{document}
\newcommand{\bs}{\begin{split}}
\newcommand{\es}{\end{split}}
\newcommand{\bpr}{\begin{proof}}
\newcommand{\epr}{\end{proof}}
\newcommand{\be}{\begin{equation}}
\newcommand{\ee}{\end{equation}}
\newcommand{\bt}{\begin{theorem}}
\newcommand{\et}{\end{theorem}}
\newcommand{\bd}{\begin{definition}}
\newcommand{\ed}{\end{definition}}
\newcommand{\bp}{\begin{proposition}}
\newcommand{\ep}{\end{proposition}}
\newcommand{\bl}{\begin{lemma}}
\newcommand{\el}{\end{lemma}}
\newcommand{\bc}{\begin{corollary}}
\newcommand{\ec}{\end{corollary}}
\newcommand{\bcon}{\begin{conjecture}}
\newcommand{\econ}{\end{conjecture}}
\newcommand{\la}{\label}
\newcommand{\Z}{{\Bbb Z}}
\newcommand{\R}{{\Bbb R}}
\newcommand{\Q}{{\Bbb Q}}
\newcommand{\C}{{\Bbb C}}
\newcommand{\hra}{\hookrightarrow}
\newcommand{\lra}{\longrightarrow}
\newcommand{\lms}{\longmapsto}
\newcommand{\bex}{\begin{example}}
\newcommand{\eex}{\end{example}}

\begin{titlepage}
\title{Hodge correlators}
\author{A.B. Goncharov }

 \date{\it To Alexander Beilinson for his 50th birthday}

\end{titlepage}
\maketitle

\tableofcontents 
\begin{abstract}\begin{footnotesize}
Hodge correlators are complex numbers given by certain integrals assigned to a smooth complex curve. 
We show that they are correlators of a Feynman integral, and describe the real mixed Hodge structure 
on the pronilpotent completion of the fundamental group of the curve. We introduce motivic correlators, which 
are elements of the motivic Lie algebra and whose periods are the Hodge correlators. They describe the motivic 
fundamental group of the curve. We describe variations of real mixed Hodge structures on a variety by certain connections 
on the product of the variety by an affine line. We call them twistor connections. Generalising this, we suggest a DG 
enhancement of the subcategory of Saito's Hodge complexes with smooth cohomology. We show that when the curve varies, 
the Hodge correlators are the coefficients of the twistor connection describing the corresponding variation of real MHS. 
Examples of the Hodge correlators include classical and elliptic polylogarithms, and their generalizations. 
The simplest Hodge correlators on the modular curves are the Rankin-Selberg integrals. Examples of the motivic 
correlators include Beilinson's elements in the motivic cohomology, e.g. the ones delivering the Beilinson - 
Kato Euler system on modular curves. 
\end{footnotesize}

\end{abstract}

\section{Introduction}

\subsection{Summary} 

Let $X$ be a smooth compact complex curve, $v_0$  a non-zero 
tangent vector 
at a point $s_0$ of $X$, and $S^* = \{s_1, \ldots , s_m\}$ 
a collection of distinct points of $X $ 
different from $s_0$. We introduce {\it Hodge correlators} 
related to this datum. They are complex numbers, given by 
integrals of certain differential forms over products of copies of $X $. 
Let $S = S^* \cup \{s_0\}$. When the data $(X, S, v_0)$ 
varies, Hodge correlators satisfy 
a system of non-linear quadratic differential equations. 
Using this, 
we show that Hodge correlators 
encode a variation of real mixed Hodge structures, called below mixed $\R$-Hodge structures. We prove that 
it coincides 
with the standard mixed $\R$-Hodge structure
 on the pronilpotent completion $\pi^{\rm nil}_1(X -S, v_0)$ 
of the fundamental group of 
$X -S$. The latter was defined, by 
different methods, by J. Morgan \cite{M},  R. Hain \cite{H} and A.A. Beilinson (cf. 
\cite{DG}). 
The real periods of $\pi^{\rm nil}_1(X -S, v_0)$ 
have a well known description via Chen's iterated integrals. 
 Hodge correlators give a completely different
 way to describe them.
\vskip 3mm

Mixed Hodge structures are relatives 
of $l$-adic Galois representations. 
There are two equivalent definitions of mixed $\R$-Hodge structures: 
as vector spaces with the weight and Hodge filtrations 
satisfying some conditions \cite{D}, 
and as representations of the {\it Hodge Galois group} \cite{D2}. 

The true analogs of $l$-adic representations are  
Hodge Galois group modules. However 
usually we describe mixed $\R$-Hodge structures 
arising from geometry  by constructing the weight and Hodge filtrations, 
getting the Hodge Galois group modules only {\it a posteriori}. 
The Hodge correlators 
describe the mixed 
$\R$-Hodge structure on $\pi_1^{\rm nil}(X-S, v_0)$ directly as a module over 
the Hodge Galois group.

\vskip 3mm

We introduce a Feynman integral related to $X$. It does not have 
a rigorous mathematical meaning. However
the standard perturbative series 
expansion procedure 
 provides a collection of its {\it correlators} assigned to  
the data $(X,S, v_0)$, 
which turned out to be 
convergent finite dimensional integrals. We show that they  
coincide with the Hodge correlators, 
thus explaining the name of the latter. Moreover: 

\begin{figure}[ht]
\centerline{\epsfbox{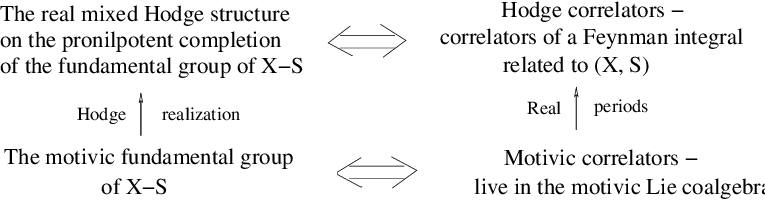}}
\label{hc1}
\end{figure}

\vskip 3mm
We define
{\it motivic correlators}.  Their periods are the Hodge correlators. 
Motivic correlators lie in the {\it motivic Lie coalgebra}, and describe the motivic fundamental group of $X-S$.  
The coproduct in the {motivic Lie coalgebra} is a new feature, 
which is missing when we work just with numbers. 
We derive a simple explicit formula for the coproduct 
of motivic correlators. 
It allows us to perform {\it arithmetic analysis} of Hodge correlators. 
This is one of the essential advantages of the Hodge correlator description 
of the real periods of $\pi_1^{\rm nil}(X-S)$, which has a lot of arithmetic 
applications.  It was available before 
only for the rational curve case \cite{G7}.

\vskip 3mm

The Lie algebra of the unipotent part of the Hodge Galois group 
is a free graded 
Lie algebra. Choosing a set of its generators we arrive at  
a collection of periods of a real MHS. 
 We introduce new generators of the Hodge Galois group, 
which differ from Deligne's generators \cite{D2}. 
The periods of variations mixed $\R$-Hodge structures 
corresponding to these generators satisfy non-linear quadratic 
Maurer-Cartan type 
differential equations. For the subcategory of Hodge-Tate structures 
they were defined by A. Levin \cite{L}. The Hodge correlators 
are the periods for this set of generators.

We introduce a DG Lie coalgebra ${\cal L}_{{\cal H}, X}^*$. 
The category of ${\cal L}_{{\cal H}, X}^*$-comodules 
  is supposed to be a DG-enhancement 
 of  the category of {\it smooth $\R$-Hodge sheaves}, i.e.  
 the subcategory of Saito's mixed $\R$-Hodge sheaves whose cohomology are 
variations of mixed $\R$-Hodge structures. We show that 
 the category of comodules  over the Lie coalgebra $H^0{\cal L}_{{\cal H}, X}^*$ 
is equivalent to the category of variations of mixed $\R$-Hodge structures. 

\vskip 3mm

The simplest Hodge correlators for 
the rational and elliptic curves deliver  
 single-valued versions of the classical polylogarithms 
and their elliptic counterparts, the classical 
Eisenstein-Kronecker series \cite{We}. The latter were 
interpreted by A.A. Beilinson and A. Levin \cite{BL} as periods of 
variations of mixed $\R$-Hodge structures. 
More generally,  when $X= \C{\Bbb P}^1$ the Hodge correlators are real periods of 
variations of mixed Hodge structures related to multiple polylogarithms. 
For an elliptic curve $E$ they deliver 
the multiple Eisenstein-Kronecker series defined in \cite{G1}.

When $X$ is a modular curve and $S$ is the set of its cusps,  
Hodge correlators 
generalize the Rankin-Selberg integrals. 
Indeed,  one of the simplest of them  is the Rankin-Selberg 
 convolution of a pair $f_1, f_2$ of weight two cuspidal Hecke eigenforms   
with an Eisenstein series. It 
computes, up to certain 
 constants,  the special value $L(f_1\times f_2, 2)$. 
A similar Hodge correlator gives 
 the Rankin-Selberg convolution 
of a weight two cuspidal Hecke eigenform $f$ 
and two Eisenstein series, and computes, up to certain 
 constants, $L(f,2)$. 
 (A generalization to higher weight modular forms will appear elsewhere).

The simplest motivic correlators on modular curves deliver  
Beilinson's elements in motivic cohomology, e.g. 
the Beilinson-Kato Euler system in $K_2$. 
We use motivic correlators to define  {\it motivic multiple L-values} 
related to: Dirichlet characters of $\Q$ 
(for $X= \C^*-\mu_N$, where $\mu_N$ is the group of $N$-th roots of unity); 
 Hecke Gr\"ossencharacters imaginary quadratic fields 
(for CM elliptic curves minus torsion points); 
Jacobi Gr\"ossencharacters of cyclotomic fields (for affine Fermat curves);  
the weight two modular forms. 

\vskip 3mm
The Hodge correlators considered in this paper admit 
a  generalization when $X-S$ is replaced by 
an arbitrary regular complex projective variety. 
We prove that they describe the mixed $\R$-Hodge structure on the rational 
motivic homotopy type of $X$,  defined by a different method 
in \cite{M}. We describe them in a separate paper since   
the case of curves is more transparent, and so far has more applications\footnote{The case of a smooth compact Kahler manifold was treated in \cite{GII}}. 
\vskip 3mm
Key constructions of this paper were outlined in Sections 8-9 of 
\cite{G1}, 
which may serve as an introduction. 

\subsection{Arithmetic motivation I: special values of $L$-functions} 
Beilinson's conjectures \cite{B} imply that special values of 
$L$-functions of motives over $\Q$ 
at the integral points to the left of the critical line are  periods of mixed motives over $\Q$. 
This means that they can be written as integrals 
$$
\int_{\Delta_B}\omega_A
$$
where $A$ and $B$ are normal crossing divisors over $\Q$ 
in a smooth $n$-dimensional projective variety $X$ over $\Q$, $\omega_A\in \Omega_{\rm log}^n(X-A)$ is an $n$-form with logarithmic singularities at $A$, and 
$\Delta_B$ is an $n$-chain with the boundary on $B(\C)$. 
Moreover, the special values are periods of their motivic avatars -- 
the motivic $\zeta$-elements -- which play a crucial role in arithmetic applications. 

The classical example is given by the Leibniz formula for the special values of the Riemann $\zeta$-function:
$$
\zeta(n) = \int_{0 \leq t_1 \leq \ldots \leq t_n \leq 1}\frac{dt_1}{1-t_1}\wedge \frac{dt_2}{t_2}.
\wedge \ldots \wedge \frac{dt_n}{t_n}.
$$
The corresponding motivic $\zeta$-elements were extensively studied (\cite{BK}, \cite{BD}, \cite{HW}).  
Beilinson's conjectures \cite{B} are known 
 for the special values of the Dedekind $\zeta$-function of a number field thanks to the work of Borel 
\cite{Bo1}, \cite{Bo2}. 

\vskip 3mm
Here is another crucial example. Let $f(z)$ be a weight two cuspidal Hecke eigenform. 
So $f(z)dz$ is a holomorphic  $1$-form on a compactified modular curve $\overline M$. 
Recall that by the Manin-Drinfeld theorem the image in the Jacobian of 
any degree zero cuspidal divisor $a$ 
on $\overline M$ is torsion. So there is  a modular unit  
$g_a \in {\cal O}^*(M)$ such that  ${\rm div}(g_a)$ is a multiple of 
the divisor $a$. According to Bloch and Beilinson \cite{B} we get an element 
\be \la{BEE}
\{g_a, g_b\}\in K_2(\overline M)\otimes \Q.
\ee
Applying the regulator map to this element, and evaluating it on the 
 $1$-form $f(z)dz$, we get an integral
\be \la{BEERS}
\int_{\overline M(\C)}\log|g_a|\overline \partial \log|g_b| f(z)dz.
\ee
It is a Rankin-Selberg convolution integral. Indeed, $\log|g_a|$ is the value of a
 non-holomorphic Eisenstein series at $s=1$, and $\partial \log|g_b|$ is a holomorphic 
weight two Eisenstein series. Therefore according to the Rankin - Selberg method, 
the integral is proportional to the product of special values of $L(f, s)$ at $s=1$ and $s=2$. 
Moreover (see \cite{SS} for a  detailed account) one can 
find cuspidal divisors $a$ and $b$ such that the proportionality coefficient is non-zero, 
i.e. 
 \be \la{BEERS213}
\int_{\overline M(\C)}\log|g_a|\overline \partial \log|g_b| f(z)dz \sim L(f,2).
\ee
Using the functional equation for $L(f,s)$, 
one can easily deduce that 
$L'(f,0)$ is a period. 
Finally, suitable modifications 
of Beilinson's $\zeta$-elements (\ref{BEE})
give rise to Kato's Euler system \cite{Ka}. 

There are similar results, due to Beilinson \cite{B2} in the weight two case, and  
Beilinson (unpublished) and, independently,  Scholl (cf. \cite{DS}) for the special values of L-functions of 
cuspidal Hecke eigenforms of arbitrary weight $w\geq 2$ at any integral point to the 
left of the critical strip, that is at $s \in \Z_{\leq 0}$. 

\vskip 3mm
This picture, especially the 
motivic $\zeta$-elements, seem to come out of the blue. 
Furthermore,  already for $L''({\rm Sym}^2f, 0)$, related by the functional equation to 
$L({\rm Sym}^2f, 3)$ 
we do not know in general how to prove that it is a period. 

One may ask whether there is a general framework,  
which delivers naturally 
both Rankin-Selberg integrals (\ref{BEERS}) and Beilinson's $\zeta$-elements (\ref{BEE}), 
and tells where to look for  generalizations related to non-critical 
special values $L({\rm Sym}^mf, n)$. 

We suggest that one should look at the motivic fundamental group of the 
universal modular curve. We show that the simplest Hodge correlator 
for the modular curve coincides with the Rankin-Selberg integral (\ref{BEERS}), and 
the corresponding motivic correlator is Beilinson's motivic $\zeta$-element (\ref{BEE}) -- 
see Section \ref{mccurves} for an elaborate discussion. 

We will show elsewhere that the Rankin-Selberg integrals related to the 
values at non-critical 
special points of $L$-functions 
of arbitrary cuspidal Hecke eigenforms and the corresponding motivic $\zeta$-elements 
appear naturally as the Hodge and motivic correlators related 
to the standard local systems on the modular curve: 
the present paper deals with the trivial local system.

\subsection{Arithmetic motivation II: arithmetic analysis of periods}

The goal of {\it arithmetic analysis} is to investigate periods without actually computing them, using instead 
arithmetic theory of mixed motives.\footnote{The latter so far mostly is  a conjectural theory, so 
we often arrive to conjectures rather then theorems.}  Here is  
the simplest non-trivial example, provided by the values of the dilogarithm at rational numbers. 
Recall that the classical dilogarithm 
$$
{\rm Li}_2(z):= -\int_0^z\log (1-t) \frac{dt}{t}
$$
is a multivalued function on $\C{\Bbb P}^1-\{0, 1, \infty\}$. Its value at $z$ depends on the homotopy class of the path 
from $0$ to $z$ on $\C{\Bbb P}^1-\{0, 1, \infty\}$ used to define the integral. Computing the monodromy 
of the dilogarithm, we see that ${\rm Li}_2(z)$ is well defined modulo the subgroup of $\C$ generated by 
 $(2\pi i)^{2}$ and $2\pi i\log(z)$. 
So the values  ${\rm Li}_2(z)$ at $z\in \Q$ are well defined modulo the subgroup 
$2\pi i\log \Q^*$  
spanned by  $\Z(2):= (2\pi i)^{2}\Z$ and  
$2\pi i\log q$, $q\in \Q^*$ -- the latter is well defined modulo $\Z(2)$.

Let $\Q^*$ denotes the multiplicative group of the field $\Q$, and 
$\Q^*\otimes\Q^*$ the tensor product of two abelian groups $\Q^*$. 
Consider a map of abelian groups
\be \la{3.27.10.1}
\Delta: \Z[\Q^*-\{1\}] \lra \Q^*\otimes \Q^*, \qquad \{z\} \lms -(1-z) \otimes z. 
\ee
\bcon \la{3.24.01.10}
 Let $z_i \in \Q, a_j, b_j\in \Q^*$ and $n_i \in \Z$. Then one has 
\be \la{3.27.10.1a}
\sum_i n_i{\rm Li}_2(z_i) + \sum_j \log a_j \cdot \log b_j = 0 ~ \mbox{\rm mod}~(2\pi i\log \Q^*)
\ee
if and only if 
\be \la{3.27.10.1b}
\Delta\Bigl(\sum_in_i\{z_i\}\Bigr) + \sum_j( a_j \otimes  b_j + b_j\otimes a_j)  =0 ~~
\mbox{\rm in}~~ (\Q^*\otimes \Q^*)\otimes_{\Z}\Q.
\ee  
\econ
This gives a complete conjectural description 
of the $\Q$-linear relations between the 
values of the dilogarithm and the products of two logarithms at rational numbers.

\bt \la{3.25.10.1}
The condition (\ref{3.27.10.1b}) implies  (\ref{3.27.10.1a}).
\et
This tells the identities between integrals without actually computing the integrals.

\bex {\em One has}  \eex
\be \la{5teq}
{\rm Li}_2(\frac{1}{3}) +{\rm Li}_2(-\frac{1}{2}) + \frac{1}{2}(\log \frac{3}{2})^2 
  =0 ~~\mbox{\rm mod}~(2\pi i\log \Q^*).
\ee
Indeed, modulo $2$-torsion in $\Q^* \otimes \Q^*$ (notice that $3/2 \otimes (-1)$ is a $2$-torsion) one has: 
\be \la{5teq1}
\Delta\Bigl(\{\frac{1}{3}\}+ \{- \frac{1}{2}\}\Bigr) = -(1-\frac{1}{3})\otimes \frac{1}{3} -(1+\frac{1}{2})\otimes (-\frac{1}{2}) 
= -\frac{2}{3}\otimes \frac{1}{3} -
\frac{3}{2}\otimes \frac{1}{2} = -\frac{3}{2} \otimes \frac{3}{2}.  
\ee

\paragraph{Motivic avatars of the logarithm and the dilogarithm.} 
The values of the dilogarithm 
and the product of two logarithms at rational arguments are 
weight two mixed Tate periods over $\Q$. In fact they describe all such periods. 
Here is the precise meaning of this. 
\vskip 2mm
The category of mixed Tate motives over a number field $F$ 
is canonically equivalent to the category of graded comodules over  
a Lie coalgebra ${\cal L}_{\bullet}(F)$ over $\Q$, graded by positive integers, the weights (cf. \cite{DG}). 
The graded Lie coalgebra ${\cal L}_{\bullet}(F)$ is free, with the space of generators in degree $n$ 
isomorphic to $K_{2n-1}(F)\otimes \Q$. Equivalently, the kernel of the coproduct map 
$$
\delta: {\cal L}_{\bullet}(F) \lra \Lambda^2 {\cal L}_{\bullet}(F)
$$
is isomorphic to the graded space $\oplus_{n>0}K_{2n-1}(F)\otimes \Q$, 
where $K_{2n-1}(F)\otimes \Q$ is in degree $n$. 
One has $K_1(F) =F^*$, and the group $K_{2n-1}(F)$ is of finite rank, which is 
$r_1+r_2$ for $n$ odd and $r_2$ for $n$ even due to 
result of Borel \cite{Bo2}. Here  $r_1$ and $r_2$ 
are the numbers of real and complex places of $F$, so that $[F:\Q]=2r_2+r_1$.  
It follows that  ${\cal L}_{1}(F) = F^* \otimes_{\Z}\Q$. Let us describe the $\Q$-vector space ${\cal L}_{2}(F)$.

For any field $K$, consider a version of the map  (\ref{3.27.10.1}):
$$
\delta: \Z[K^*-\{1\}] \lra K^*\wedge K^*, \qquad \{z\} \lms - (1-z) \wedge z. 
$$
Let $r(*,*,*,*)$ be the cross-ratio of four points on ${\Bbb P}^1$, normalized by 
$r(\infty, 0, 1, x) = x$. Let $R_2(K)$ be the 
subgroup of $\Z[K^*-\{1\}]$ generated by the elements
\be \la{5trel}
\sum_{i=1}^5 (-1)^i\{r(x_1, ..., \widehat x_i, ..., x_5)\}, \qquad x_i \in P^1(F), 
\quad x_i \not = x_j. 
\ee
Then $\delta (R_2(K)) =0$. Let $B_2(K):= \Z[K^*-\{1\}]/R_2(K)$ be the Bloch group of $K$. 
The following proposition is deduced  in Section 1 of \cite{G10} from a theorem of Suslin \cite{S}.

\bp \la{3.25.10.2}
Let $F$ be a number field. Then one has ${\cal L}_{2}(F) = B_2(F) \otimes_{\Z}\Q$. The same 
should be true for any field $F$. However the Lie coalgebra ${\cal L}_\bullet(F)$ in that case is a conjectural object. 
\ep

The dual of the universal enveloping algebra of the Lie algebra graded dual to the 
Lie coalgebra ${\cal L}_{\bullet}(F)$ is a commutative 
graded Hopf algebra.\footnote{We want to note that although it is 
traditionally graded by the nonnegative integers, the Tate weights, 
there is no sign rule involved here. Perhaps we would be better off by multiplyinbg
 the weights by the factor of two.} 
 Its weight two component is 
\be \la{Motavat}
{\cal L}_{2}(F) \oplus {\rm Sym}^2{\cal L}_{1}(F) = (B_2(F) \oplus {\rm Sym}^2F^* )\otimes_{\Z}\Q. 
\ee
The coproduct is provided by the map 
\be \la{Motavat1}
\Delta_2: B_2(F) \oplus {\rm Sym}^2F^* \lra (F^*\otimes F^*)\otimes_\Z \Q
\ee
$$ \{z\}_2 \lms 
\frac{1}{2}\Bigl(z \otimes (1-z)-(1-z) \otimes z \Bigr) , \quad 
a\cdot b \lms a \otimes b + b  \otimes a. 
$$

Let us return to our example, when $F=\Q$. The period homomorphism 
$$
{\cal L}_{1}(\Q) \lra \C/\Z(1), \qquad q \lms \log q
$$ is known to be injective. 
So $q\in \Q^*$ is the motivic avatar of $\log q$.

\bp \la{3.25.10.3}
There is a weight two period map
\be \la{3.25.10.4}
{\cal L}_{2}(\Q) \oplus {\rm Sym}^2{\cal L}_{1}(\Q) \lra \C/2\pi i\log \Q^*, \qquad 
\ee
$$
\{q\}_2 \lms {\rm Li}_2(q) + \frac{1}{2}\log(1-q)\log(q), \quad a\otimes b + b \otimes a \lms \log a \log b.
$$
\ep 
The point is that the period map kills the 
subgroup $R_2(\C)$, which is equivalent to Abel's five term relation for the dilogarithm. 

The period map (\ref{3.25.10.4}) is conjectured to be injective: 
this is a special case of a Grothendieck type conjecture on periods. 
Finally, ${\rm Ker}\Delta_2 = K_3(\Q)\otimes \Q$, 
which is zero by Borel's theorem. The injectivity of $\Delta_2$ 
plus Propositions \ref{3.25.10.2} and \ref{3.25.10.3} imply Theorem \ref{3.25.10.1}, 
and, modulo the injectivity of the period map (\ref{3.25.10.4}), Conjecture \ref{3.24.01.10}. 

Although relation (\ref{5teq}) follows from the five 
term relation for the dilogarithm, there is no algorithm which allows to write the element $\{\frac{1}{3}\} + 
\{-\frac{1}{2}\}$ as a sum of the five term relations (\ref{5trel}). 
The only effective way to know that they exist is 
to calculate the coproduct, as we did in (\ref{5teq1}). 

Let us now discuss the case when $F$ is an arbitrary number field. 
\vskip 3mm
{\it Real periods.} There is a single valued version of the dilogarithm, the Bloch-Wigner function 
$$
{\cal L}_2(z):= -{\rm Im}\Bigl(\int_0^z\log (1-t) \frac{dt}{t} + \int_0^z\frac{dt}{1-t}\cdot \log|z|\Bigr).
$$
where both integrals are defined by using the same integration 
path from $0$ to $z$. It satisfies the five term relation, and thus 
provides a group homomorphism
$$
B_2(\C) \lra \R, \qquad \{z\}_2 \to {\cal L}_2(z).
$$
Using the (conjectural) isomorphism ${\cal L}_2(\C) = B_2(\C)$, it can be interpreted as 
the {\it real  period map} ${\cal L}_2(\C)\to \R$. Notice that 
${\cal L}_2(z) =0$ for $z\in \R$. So the real period map looses a lot of 
information about the motivic Lie algebra. It relates, however, 
the kernel of the coproduct map with the special values of $L$-functions. 
Here is how it works in our running example. 

Let $F$ be a number field. 
Then by Borel's theorem \cite{Bo2} the  real period map provides an 
injective regulator map
$$
K_3(F)\otimes \Q  \stackrel{\sim}{=} {\rm Ker}\delta \subset B_2(F)\otimes \Q \lra \R^{r_2}, 
\qquad \{z\}_2 
\lms 
\Bigl({\cal L}_2(\sigma_1(z)), ..., {\cal L}_2(\sigma_{r_2}(z))\Bigr).
$$
Its image is a rational lattice -- that is, a lattice tensor $\Q$ -- 
 in $\R^{r_2}$, whose covolume, well defined up to  $\Q^*$, 
 is  a $\Q^*$-multiple of $\zeta_F(-1)$. 

\vskip 2mm
So  to perform the arithmetic analysis of the values of the dilogarithm / products of two logarithms 
at rational arguments we 
upgrade them to their motivic avatars lying in (\ref{Motavat}), 
and determine their coproduct there. The kernel of the coproduct is captured by the 
regulator map. 

The weight $m$ periods of mixed Tate motives over 
a number field $F$ are studied similarly using the iterated coproduct in $\otimes^mF^*$, 
see Section 4 of \cite{G7}. 

\vskip 2mm
In this paper we develope a similar picture for the periods of 
the pronilpotent completions of fundamental groups of curves. Namely,  we introduce motivic correlators, 
which are the motivic avatars of the periods of fundamental groups of curves. 
They span a Lie subcoalgebra in the motivic Lie coalgebra. The key point is that the coproduct 
of a motivic correlator is given explicitly in terms of the motivic correlators. We prove that 
their real periods 
are given by the Hodge correlators. 

\subsection{Mixed motives and the motivic Lie algebra}\la{secaap}
Let us recall the conjectural motivic setting. 
For any field $F$ there is a conjectural {\it abelian} category of 
pure Grothendieck motives over $F$. It is expected to be a semisimple abelian tensor category. 
Moreover, there is a conjectural abelian category of mixed motivic sheaves over any base $X$, 
whose properties were discussed by Beilinson in \cite{B1}. 
In particular it contains the Tate motives over $X$, denoted by $\Q(n)_X$. 
Let us asume that $X$ is regular. 
The Ext groups between them in this category are related by Beilinson's formula \cite{B1} 
to the Adams graded quotients of the algebraic K-groups of the base $X$:
\be \la{Beilf}
{\rm Ext}^m(\Q(0)_X, \Q(n)_X) = {\rm gr}^n_\gamma K_{2n-m}(X)\otimes \Q. 
\ee

A useful tool investigate the Ext groups in general is the hypothetical 
motivic Leray spectral sequence. Let us assume that $X$ is a regular projective variety 
over a field $F$. Then using (\ref{Beilf}) and the pair of adjoint functors $p^*, p_*$ where 
$p: X \to {\rm Spec}(F)$, we have 
\be \la{BMLSS1}
{\rm gr}^n_\gamma K_{2n-m}(X)\otimes \Q = {\rm Ext}^m(\Q(0)_X, \Q(n)_X) = 
{\rm Ext}^m(p^*\Q(0)_F, \Q(n)_X) = {\rm Ext}^m(\Q(0)_F, p_*\Q(n)_X).  
\ee
Since $p_*\Q(n)_X$ is expected to have a filtration with the associate graded are pure motives 
$h^iX(n)$, we expect a spectral sequence with the  $E^2$-term
\be \la{BMLSS}
{\rm Ext}^j(\Q(0)_F, h^iX(n)) => {\rm gr}^n_\gamma K_{2n-m}(X)\otimes \Q.
\ee
The trianguilated category of motives is available, together with Beilinson's formula (\ref{Beilf}), 
thatnks to the works of V. Voevodsky \cite{V} and M. Levine \cite{Le}. The abelian category 
of mixed motives is not available yet. However if we restrict to 
the case when $X = {\rm Spec}({\cal O}_{F, S})$ is the spectrum of the ring of $S$-integers in a number field 
$F$, and consider the subcategory of mixed Tate motives, the corresponding category with all necessary properties is 
available \cite{Le2} for $X= {\rm Spec}(F)$ and \cite{DG}. 

Let us assume that $X= {\rm Spec}(F)$ where $F$ is a field. 
The Tannakian category of mixed motives over $F$ is expected to be 
canonically equivalent to the category of 
modules over a Lie algebra ${\rm L}_{\rm Mot/F}$ in the category of 
pure motives over ${\rm Spec}(F)$. We also use notation ${\rm L}_{\rm Mot}$. 
The equivalence is provided by the fiber functor 
$M \lms {\rm gr}^WM$. 
Denote by ${\cal L}_{\rm Mot/F}$ the dual Lie coalgebra. 
So there is an isomorphism
\be \la{biso1}
H^i({\cal L}_{\rm Mot/F}) = \oplus_{[M]}{\rm Ext}_{\rm Mot/F}^i(\Q(0), M)\otimes M^\vee.
\ee
where the sum is over the isomorphism classes of pure motives. 
In particular, one has
\be \la{biso}
{\rm Ker} \Bigl({\cal L}_{\rm Mot/F} \stackrel{\delta}{\lra} 
\Lambda^2 {\cal L}_{\rm Mot/F} \Bigr) = \oplus_{[M]}{\rm Ext}_{\rm Mot/F}^1(\Q(0), M)\otimes M^\vee.
\ee

If $F$ is a number field,  
then, conjecturally, $H^i({\cal L}_{\rm Mot/F})=0$ for $i>1$. So 
the motivic Lie algebra ${\rm L}_{\rm Mot/F}$ is a free Lie algebra in the category of pure motives. 
So (\ref{biso}) implies that  
 the $M^\vee$-isotipical component of the space of 
generators is isomorphic to ${\rm Ext}_{\rm Mot/F}^1(\Q(0), M)$. 

Finally, the size of ${\rm Ext}_{\rm Mot/F}^1(\Q(0), M)$ 
is predicted by Beilinson's conjecture \cite{B} on regulators. 

We conclude that in the case when $F$ is  a number field, the 
structure of the motivic Lie algebra ${\rm L}_{\rm Mot/F}$ is completely predicted by Beilinson's conjectures.

\subsection{Pronilpotent completions of 
 fundamental groups of curves} \la{sec1.2a}
 
The  fundamental 
group $\pi_1 = \pi_1(X  - S, v_0)$ is a free group 
with generators provided by  loops around the 
rest of the punctures $s_i \not = s_0$ and loops 
generating $H_1(X )$. 
Let ${\cal I}:= {\rm Ker}(\Q[\pi]\to \Q)$  be the augmentation ideal of the group algebra of  
$\pi_1$. Then there is  
a complete cocommutative Hopf algebra over $\Q$
\begin{equation} \label{10/25/04/1}
{\rm A}^{\rm Betti}(X - S, v_0):= 
\lim_{{\longleftarrow}}\Q[\pi_1]/{\cal I}^n.
\end{equation} 
Its coproduct is induced by the map $g \lms g \otimes g$, $g \in \pi_1$.  
It is called  the {\it fundamental Hopf algebra of $X  - S$} with the tangential base point $v_0$. 
The subset of its primitive elements  is a free pronilpotent Lie 
algebra over $\Q$, the Maltsev completion of $\pi_1$. It is denoted 
$\pi_1^{\rm nil}(X-S, v_0)$ 
and called the {\it fundamental Lie algebra} 
of $X -S$. 
The Hopf algebra (\ref{10/25/04/1}) 
is its universal enveloping algebra.

Denote by ${\rm T}(V)$ the tensor algebra of a vector space $V$. 
The associated graded of the Hopf algebra 
(\ref{10/25/04/1}) for the ${\cal I}$-adic filtration is isomorphic to the 
tensor  algebra  
of $H_1(X -S, \Q)$: 
$$
{\rm gr}^{\cal I}{\rm A}^{\rm Betti}(X - S, v_0) = 
\oplus_{n=0}^{\infty} {\cal I}^n/{\cal I}^{n+1} = {\rm T}\Bigl(H_1(X -S, \Q)\Bigr).
$$

There is a weight filtration on $H_1(X -S, \Q)$, given by the  
extension 
\begin{equation} \label{6.11.06.9}
0 \lra \Q(1)^{|S|-1} \lra H_1(X -S, \Q) \lra H_1(X , \Q) \lra 0,
\end{equation}
where the right arrow is provided by the embedding 
$
X -S \hookrightarrow X . 
$  
The Hopf algebra (\ref{10/25/04/1}) is equipped with a weight filtration $W$ compatible with 
the Hopf algebra structure. 
The  
associated graded for the weight filtration 
of the Hopf algebra (\ref{10/25/04/1}) is isomorphic to the tensor algebra 
of ${\rm gr}^WH_1(X -S, \Q)$:
\begin{equation} \label{6.11.06.10}
{\rm gr}^{W}{\rm A}^{\rm Betti}(X - S, v_0) = {\rm T}\left( {\rm gr}^WH_1(X -S, \Q)
 \right).
\end{equation}

\subsection{Hodge correlators} \la{sec1.5} 
Recall $S^*= S-\{s_0\}$. Set
\be \la{9.25.13.1}
{\rm V}_{X, S^*} := H_1(X, \C) \oplus \C[S^*], \qquad {\rm A}_{X, S^*} := 
{\rm T}({\rm V}_{X, S^*}).
\ee
One sees from (\ref{6.11.06.9}), (\ref{6.11.06.10}) that there are canonical isomorphisms
 of vector spaces
\be \la{9.25.13.2}
{\rm V}_{X, S^*} = {\rm gr}^{W}H_1(X -S, \C), \quad {\rm A}_{X, S^*} = 
{\rm gr}^{W}{\rm A}^{\rm Betti}(X - S, v_0)\otimes \C.
\ee
We also need the dual objects: 
$$
{\rm V}^{\vee}_{X, S^*} := H^1(X, \C) \oplus \C[S^*], 
\qquad {\rm A}^{\vee}_{X, S^*} := {\rm T}({\rm V}^{\vee}_{X, S^*}).
$$

For an associative algebra $A$ over $k$,  let $A^+:= {\rm Ker}(A \to k)$ be the augmentation ideal. Let $[A^+, A^+]$ be the subspace (not the ideal) of $A^+$ generated by 
commutators in $A^+$. The quotient space ${\cal C}(A):= A^+/[A^+,A^+]$ 
is called the {\it cyclic envelope} 
of $A$. If $A$ is 
freely generated by a set ${\cal S}$ then the vector space 
${\cal C}(A)$ has a basis parametrised by cyclic words in  ${\cal S}$.

We use the shorthand ${\cal C}{\rm T}(V)$ for 
${\cal C}({\rm T}(V))$. 
We define the subspace of {\it shuffle relations} in 
${\cal C}{\rm T}({\rm V})$ as the subspace generated by the 
elements
\be \la{shshrel}
\sum_{\sigma \in \Sigma_{p,q}}(v_0 \otimes v_{\sigma(1)} \otimes \ldots \otimes v_{\sigma(p+q)}), 
\qquad p, q \geq 1, 
\ee
where the sum is over all $(p,q)$-shuffles. 
Set
\be \la{mvlca}
{\cal C}^{\vee}_{X, S^*}:= {\cal C}({\rm A}^{\vee}_{X, S^*}), \qquad 
{\cal C}{{\cal L}ie}^{\vee}_{X, S^*}:= 
\frac{ {\cal C}^{\vee}_{X, S^*}}{\mbox{Shuffle relations}};
\ee
\be \la{mvla}
{\cal C}_{X, S^*}:= {\cal C}({\rm A}_{X, S^*})
\qquad 
{\cal C}{{\cal L}ie}_{X, S^*}:= \mbox{the dual of 
${\cal C}{{\cal L}ie}^{\vee}_{X, S^*}$}.
\ee
We show that ${\cal C}^{\vee}_{X, S^*}$ is a Lie coalgebra, 
${\cal C}{{\cal L}ie}^{\vee}_{X, S^*}$ is its quotient Lie coalgebra. 
Equivalently,  
${\cal C}_{X, S^*}$ is a Lie algebra,  
${\cal C}{{\cal L}ie}_{X, S^*}$ is its Lie subalgebra.

\begin{figure}[ht]
\centerline{\epsfbox{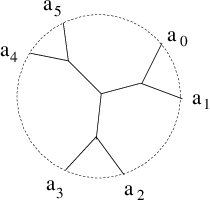}}
\caption{A plane trivalent tree decorated by ${\cal C}(\{a_0\} \otimes \ldots \otimes \{a_5\})$.}
\label{feyn36}
\end{figure}

\vskip 3mm
Given a mass one volume form $\mu$ on $X $, 
our Feynman integral construction in Section \ref{hc2sec}
 provides a linear map, called the 
{\it Hodge correlator map}
\be \la{1.15.08.1}
{\rm Cor}_{{\cal H}, \mu}: {\cal C}{{\cal L}ie}^{\vee}_{X, S^*} \lra \C.
\ee
\bex {\em For a cyclic word} \eex 
\begin{equation} \label{11:42}
W = {\cal C}\Bigl(
\{a_0\} \otimes \{a_{1}\} \otimes ... \otimes \{a_{n}\} \Bigr), 
\qquad a_i \in S^*, 
\end{equation}
 the Hodge correlator is given by a 
sum over all plane trivalent trees $T$ whose external edges are decorated by  
elements $a_0, \ldots , a_n$, see Fig \ref{feyn36}. To define the integral 
corresponding to such a tree $T$ we proceed as follows. 
The volume form $\mu$ provides a Green function $G(x,y)$ 
on $X ^2$. Each edge $E$ of $T$ 
contributes a Green function on $X ^{\{\mbox{\rm vertices of $E$}\}}$, 
which we lift to a function 
on  $X ^{\{\mbox{\rm vertices of $T$}\}}$. Further, given 
any $m+1$ smooth functions on a complex manifold $M$, there is a canonical 
linear map $\omega_m: \Lambda^{m+1}{\cal A}^0_M \to {\cal A}^m_M$, where 
${\cal A}^k_M$ is the space of smooth $k$-forms on $M$ (Section 2.2). 
Applying it to the Green functions assigned to  
the edges of $T$, we get a differential form 
of the top degree 
on $X ^{\{\mbox{\rm internal vertices of $T$}\}}$. Integrating it, we 
get the integral assigned to $T$. Taking the sum over all trees $T$, we get 
the number ${\rm Cor}_{{\cal H}, \mu}(W)$. The shuffle relations result  
from taking the sum over all trees with a given decoration. 

Dualising the Hodge correlator map (\ref{1.15.08.1}), we get a {\it Green operator} 
\begin{equation} \label{5.22.06.1}
{\bf G}_{\mu} \in {\cal C}{{\cal L}ie}_{X, S^*}. 
\end{equation}

Let 
${\rm L}_{X,S^*}$ be the free Lie algebra generated by the space 
${\rm V}_{X,S^*}$. There is a canonical isomorphism 
\be \la{varl}
{\rm L}_{X,S^*} = {\rm gr}^W\pi_1^{\rm nil}(X-S, v_0)\otimes \C.
\ee
 The Green operator can be viewed as a derivation of the Lie algebra ${\rm L}_{X,S^*}$ as follows.

Let $(p_i, q_i)$ be a symplectic basis of $H_1(X)$. Denote by $X_s$ the generator of 
${\rm L}_{X,S^*}$ assigned to  $s\in S^*$. There is a  canonical generator 
of ${\rm gr}^W\pi_1^{\rm nil}(X-S, v_0)$ corresponding to a loop around $s_0$.
 Its projection to ${\rm L}_{X,S^*}$ is written as 
\begin{equation} \label{6.15.06.1}
X_{s_0}:= -\sum_{s\in S^*}X_{s}+ \sum[p_i, q_i].
\end{equation} 

 Say that a derivation of 
the Lie algebra ${\rm L}_{X,S^*}$ is {\it special} if it kills 
the generator (\ref{6.15.06.1}) and 
 preserves the conjugacy classes 
of the generators $X_{s}$, where $s\in S^*$. 

This definition is motivated by 
the $l$-adic picture (Section \ref{hc8sec}.2): If $X$ is a curve 
over a field $F$ and $S \subset X(F)$, 
the Galois group ${\rm Gal}(\overline F/F(\mu_{l^\infty}))$ 
acts by special automorphisms 
of ${\rm L}_{X,S^*}\otimes \Q_l$. Indeed, by the comparison theorem 
the latter is identified 
with the associate graded 
for the weight filtration on the Lie algebra of 
the pro-$l$-completion $\pi_1^{(l)}(X-S, v_0)$ of the fundamental group.  

  Denote the Lie algebra of special derivations by 
${\rm Der}^S{\rm L}_{X,S^*}$. 
We show that the Lie algebra 
${\cal C}{{\cal L}ie}_{X,S^*}$ acts by special derivations on
 the Lie algebra ${\rm L}_{X,S^*}$. This generalizes constructions of Drinfeld \cite{Dr1} 
(when $X$ is of genus zero) 
and Kontsevich \cite{K} (when $S$ is empty).  
Namely, let $F$ be a cyclic polynomial in non-commuting variables 
$Y_i$, i.e. $F \in {\cal C}{\rm T}(Y)$ 
where $Y$ is a vector space with a  basis $\{Y_i\}$.  
Then there are ``partial derivatives'' maps 
$$ 
\partial/\partial Y_i: {\cal C}{\rm T}(Y) \lra {\rm T}(Y), \qquad F \lms {\partial F}/{\partial Y_i}, 
$$
defined by deleting one of $Y_i$'s from $F$, getting as a result a non-commutative polynomial, 
and taking the sum over all  possibilities. For instance if $F = {\cal C}(Y_1Y_2Y_1Y_3)$ then 
${\partial F}/{\partial Y_1} = Y_2Y_1Y_3+ Y_3Y_1Y_2$. 

Let $F$ be a cyclic polynomial in non-commuting variables 
$X_{s}, p_i, q_i$, where $s\in S^*$. 
Then the derivation $\kappa_F$ assigned to $F$ acts on the generators by 
$$
p_i \lms -\frac{\partial F}{\partial q_i}, \quad q_i \lms \frac{\partial F}{\partial p_i}, \quad 
X_s \lms [X_s, \frac{\partial F}{\partial X_s}], \quad s\in S^*. 
$$
{\it A priori} it is a derivation of the associative algebra 
${\rm A}_{X, S^*}$. However if $F$ annihilates the shuffle relations,  
it is a derivation of the Lie algebra ${\rm L}_{X, S^*}$. 
We show that we get an isomorphism
\begin{equation} \label{9.17.05.5}
{\cal C}{{\cal L}ie}_{X,S^*} \stackrel{\sim}{\lra} {\rm Der}^S{\rm L}_{X,S^*}, \quad F \lms \kappa_F.
\end{equation}

Combining (\ref{5.22.06.1}) and (\ref{9.17.05.5}), we conclude that 
the Green 
operator can be viewed as an element 
\begin{equation} \label{5.22.06.1qa}
{\bf G}_{\mu} \in {\rm Der}^S{\rm L}_{X,S^*}. 
\end{equation}
We show that it encodes a mixed $\R$-Hodge structure using a 
construction from Section \ref{1.5}.

\subsection{A Feynman integral for Hodge correlators} \la{1.4}
Let $\varphi$ be a smooth function 
on a complex curve $X(\C)$ with values in $N\times N$ complex matrices. We say that 
the space 
$\{\varphi\}$ of such functions is our space of  fields. 
Given $N$, consider the following correlator corresponding to $W\in {\cal C}^{\vee}_{X,S}$, 
formally defined  via 
a Feynman integral over the space of fields $\{\varphi\}$
\begin{equation} \label{3.12.05.3asd}
{\rm Cor}_{X, N, h}(W):= 
\int {\cal F}_W(\varphi)
e^{iS(\varphi)}{\cal D}\varphi,
\end{equation}
where 
$$
S(\varphi):= \frac{1}{2\pi i}\int_{X } {\rm Tr}
\Bigl( \frac{1}{2}\partial  \varphi \wedge \overline \partial\varphi + 
\frac{1}{6}\hbar  \cdot 
\varphi [\partial \varphi, \overline \partial \varphi]  \Bigr).
$$
and ${\cal F}_W(\varphi)$ is a function on the space of fields, see 
Section \ref{hc12sec}.
For example, for the cyclic word (\ref{11:42}) 
we set 
\begin{equation} \label{7.28.06.1}
{\cal F}_W(\varphi):= {\rm Tr}\Bigl(\varphi(a_0) \ldots \varphi(a_n) \Bigr). 
\end{equation}
Formula (\ref{3.12.05.3asd}) does not have a precise mathematical meaning. 
We understand it by postulating the perturbative series 
expansion with respect to a small parameter $h$, 
using the standard Feynman rules, and then taking the leading term in the 
asymptotic expansion as $N \to \infty, \hbar=N^{-1/2}$. This way we get 
a sum over a finite number of Feynman diagrams given by plane trivalent trees 
decorated by the factors of $W$. One needs to specify 
a volume form $\mu$ on $X $ to determine a measure ${\cal D}\varphi$ on the space of fields.  
Mathematically, we need $\mu$ to specify the Green function which is used 
to write the perturbative series 
expansion.

We prove in Theorem \ref{1.10.05.1} that these Feynman intregral correlators
coincide with the suitably normalised Hodge correlators.

\paragraph{Problem.} The Hodge correlators provide only the leading term of the $N \to \infty$ 
asymptotics of the Feynman integral correlators. 
The next terms in the asymptotics are given by the correlators corresponding to higher genus 
ribbon graphs. They could be divergent. How to renormalise them? 
What is their role in the Hodge theory?

\subsection{Variations of mixed $\R$-Hodge structures by twistor connections} \la{1.5}
\la{hc1.5}

\paragraph{Real Hodge structures.} A {\it real Hodge structure} is a real vector space $V$ 
whose complexification is equipped with a bigrading $V_\C = \oplus V^{p,q}$ such that 
$\overline V^{p,q} = V^{q,p}$.  A real Hodge structure is {\it pure} if $p+q$ 
is a given number, called the weight. 
So a real Hodge structure is a direct sum of 
pure ones of different weights. 

The category ${\rm HS}_{/ \R}$of real 
Hodge structures is equivalent to the category of representations of the 
algebraic group ${{\bf G}_m}_{\C/\R}$. The group of its real points 
is  the group $\C^*$.  The group of its complex points is the group 
$\C^*\times \C^*$ with the action of the complex conjugation interchanging the factors. 
Abusing notation, we denote it by $\C^*_{\C/\R}$.

\paragraph{Real mixed Hodge structures.} According to P. Deligne \cite{D}, a {\it mixed $\R$-Hodge structure} 
is a real vector space $V$  equipped with a weight filtration $W_{\bullet}$, 
 and a Hodge filtration $F^{\bullet}$ of its complexification $V_\C$, 
satisfying the following condition. The filtration $F^{\bullet}$ 
and its conjugate $\overline F^{\bullet}$ induce on ${\rm gr}_n^WV$ a pure 
weight $n$ real Hodge structure: 
$$
{\rm gr}^WV_\C = \oplus_{p+q=n}F_{(n)}^{p}\cap \overline F_{(n)}^{q}.
$$
Here $F_{(n)}^{\bullet}$ is the filtration on ${\rm gr}_n^WV_\C$ induced by 
$F^{\bullet}$, and similarly $\overline F_{(n)}^{\bullet}$.

\paragraph{The real Hodge Galois group.} The category ${\rm MHS}_{/\R}$ of mixed $\R$-Hodge structures is a 
Tannakian category
 with 
a canonical 
fiber functor $\omega_{\rm Hod}$ to the category ${\rm Vect}_\R$ of $\R$-vector spaces:
$$
\omega_{\rm Hod}: {\rm MHS}_{/\R} \lra {\rm Vect}_\R.
$$
It assigns to a real mixed Hodge structure $(V,W_{\bullet}V,  F^{\bullet}V_\C)$  the 
underlying real vector space $V$. 
The {\it Hodge Galois group} $G_{\rm Hod}$ is defined as the group of 
automorphisms of the fiber functor respecting the tensor structure:
$$
G_{\rm Hod}:= {\rm Aut}^\otimes\omega_{\rm Hod}.
$$
It is a  pro-algebraic group over $\R$. 
Thanks to the Tannakian formalism, 
the functor $\omega_{\rm Hod}$ provides a canonicall equivalence of 
the tensor category of real mixed Hodge structures with the tensor 
 category of representations of the proalgebraic group $G_{\rm Hod}$:
$$
\omega_{\rm Hod}: {\rm MHS}_{/\R} \stackrel{\sim}{\lra} G_{\rm Hod}-{\rm modules}.
$$ 
There are two canonical functors 
$$
{\rm gr}^W: {\rm MHS_{/\R}} \lra {\rm HS}_{/\R}, ~~ V \lms {\rm gr}^WV, 
~~~~ i: {\rm HS}_{/\R}\hra {\rm MHS}_{/\R}, ~~ V \lms V, ~~~~
{\rm gr}^W\circ i = {\rm Id}. 
$$
Here the functor $i$ is the canonical embedding. Therefore there are two homomorphisms  
$$
s: \C^*_{\C/\R}\lra G_{\rm Hod}, ~~~~p: G_{\rm Hod}\lra \C^*_{\C/\R}, ~~~~p\circ s={\rm Id}.
$$
These homomorphisms provide the  group $G_{\rm Hod}$ with a structure of  a 
semidirect product of a real prounipotent algebraic group ${U}_{\rm Hod}:={\rm Ker}~p$ and $\C^*_{\C/\R}$: 
$$
0\lra {U}_{\rm Hod} \lra G_{\rm Hod} \stackrel{p}{\lra} \C^*_{\C/\R}\lra 0.
$$ 
Let ${\rm L}_{\rm Hod}$ be the Lie algebra of the group ${U}_{\rm Hod}$.  
The action of  $s(\C^*_{\C/\R})$ provides the Lie algebra ${\rm L}_{\rm Hod}$ with 
 a structure of a Lie algebra 
in the category of $\R$-Hodge structures. 
The category of mixed $\R$-Hodge structures  is equivalent to the category of  
representations of the Lie algebra ${\rm L}_{\rm Hod}$ in the category of $\R$-Hodge structures. 
The Lie algebra ${\rm L}_{\rm Hod}$ is a free Lie algebra in the category of $\R$-Hodge structures. 
Thus the Lie algebra ${\rm L}_{\rm Hod}$ is generated by the $\R$-Hodge structure
$$
\oplus_{[H]}{\rm Ext}_{{\rm MHS}_{/ \R}}^1(\R(0), H)^{\vee}\otimes H,
$$
where the sum is over the set of isomorphism classes of simple $\R$-Hodge structures $H$. 
The ones with non-zero ${\rm Ext}^1$  are parametrised by pairs of 
negative integers $(-p,-q)$ up to a permutation (the Hodge degrees of $H$),
and one has  ${\rm dim}_\R{\rm Ext}_{{\rm MHS}_{/ \R}}^1(\R(0), H)=1$ 
for those $H$. 

It follows that the complex Lie algebra ${\rm L}_{\rm Hod}\otimes_\R\C$ is generated by 
certain elements  
$n_{p, q}$, $p,q \geq 1$, of bidegree $(-p, -q)$, 
with the only relation $\overline n_{p, q}= -n_{q, p}$. 
Such a generators were defined by Deligne \cite{D2}. We call them Deligne's generators.

This just means that 
a mixed $\R$-Hodge structure can be described by a pair $(V, g)$, where 
$V$ is a real vector space whose complexification is equipped with a bigrading $V_\C = \oplus V^{p,q}$ such that 
$\overline V^{p,q} = V^{q,p}$, and 
$g = \sum_{p,q\geq 1}g_{p,q}$ is an imaginary operator on $V_\C$, so 
that $g_{p,q}$ is of bidegree $(-p,-q)$ and  
$\overline g_{p,q} = -g_{q,p}$. The operators $g_{p,q}$ are the images of 
Deligne's  
generators  
$n_{p,q}$ in the representation $V$. 
One can describe variations of real Hodge structures 
using the operators $g_{p,q}$ fiberwise, but the Griffiths transversality condition 
is encoded by so  complicated nonlinear 
differential equations on $g_{p,q}$ that it was not even written down. 

A.A. Beilinson emphasized that 
Deligne's generators may not be 
the most natural ones, and asked whether one should exist 
a canonical choice of the generators. 
For the subcategory of mixed $\R$-Hodge-Tate  structures 
(i.e. when the Hodge numbers are zero unless $p=q$) 
a different set of generators $n_p$ 
was suggested by 
A. Levin \cite{L}.

\vskip 3mm
Let us decompose the Hodge correlator ${\bf G} = \sum 
{\bf G}_{p,q}$ according to the Hodge bigrading. 
We proved that the  ${\bf G}_{p,q}$ 
are images of certain generators of the Lie algebra ${\rm L}_{\rm Hod}$,
thus providing the mixed $\R$-Hodge structure on $\pi_1^{\rm nil}(X-S, v_0)$.  
We proved in Section \ref{hc6sec}
 that, when the data $(X, S, v_0)$ varies, 
 the  $\{{\bf G}_{p,q}\}$ satisfy a Maurer-Cartan 
nonlinear quadratic  differential equations.

Motivated by this, 
we show that a variation of mixed $\R$-Hodge structures can be 
described by a direct sum of variations of pure Hodge structures 
plus a {\it Green datum} $\{G_{p,q}, \nu\}$ consisting of 
operator-valued functions $G_{p,q}$, $p,q \geq 1$ and a $1$-form $\nu$,  
satisfying a Maurer-Cartan system of quadratic nonlinear differential equations. 
The Green operators $G_{p,q}$ are 
 nonlinear modifications 
of Deligne's operators $g_{p,q}$, 
providing a canonical set of generators of the Hodge 
Lie algebra ${\rm L}_{\rm Hod}$.  
For the Hodge-Tate structures we recover 
the Green operators  of \cite{L}, although even in this case our 
treatment is somewhat different, 
 based on the twistor transform defined in Section \ref{hc5sec}.

\vskip 3mm
Let us  describe our construction. 
Let ${\cal V}$ be a variation of real Hodge structures on $X $. 
It is a real local system ${\cal V}$, whose complexification ${\cal V}_{\C}$ has a 
canonical decomposition
\begin{equation} \label{8.20.05.1i}
{\cal V}_{\C}  = \oplus_{p,q}{\cal V}^{p,q}.  
\end{equation}
Tensoring ${\cal V}$ with the sheaf of smooth complex functions on $X $, we get 
a $C^{\infty}$-bundle ${\cal V}_{\infty}$ with a flat 
connection ${\bf d}$ satisfying the 
Griffiths transversality condition.
Pick a real closed $1$-form 
\be \la{formnu}
g_{0,0}\in \Omega^1\otimes 
 {\rm End}^{-1, 0}{\cal V}_\infty \oplus \overline \Omega^1 \otimes 
{\rm End}^{0, -1}{\cal V}_\infty, \qquad {\bf d}g_{0,0} =0, \quad \overline g_{0,0} = g_{0,0}, 
\ee
and  a collection of imaginary operators 
\be \la{functionG}
G_{p,q}\in  {\rm End}^{-p, -q}{\cal V}_{\infty}, \quad p,q \geq 1, \qquad \overline G_{p,q} = - G_{p,q}.
\ee
Let ${\bf d} = \partial' + \partial''$ be the decomposition 
into the holomorphic and antiholomorphic components. 

Let $\varphi^{a,b}_{p,q}$ be an ${\rm End}^{-p, -q}{\cal V}$-valued 
$(a,b)$-form.  
We define its {\it Hodge bigrading} $(s,t)$ by setting $s:= a-p$ and $t:=b-q$.  
Denote by $\varphi_{s,t}$ the 
 component of $\varphi$ 
of the  Hodge bidegree $(-s, -t)$. 

\paragraph{Twistor connections.} Let us introduce a {\it twistor plane} 
$\C^2$ with coordinates $(z,w)$. 
Consider 
 the product  $X\times \C^2$, and the  projection
$\pi: X \times \C^2 \to X$. 

Given a datum $\{G_{p,q}, g_{0,0}\}$, we 
introduce a {connection} $\nabla_{\cal G}$ on $\pi^*{\cal V}_{\infty}$: 
\be \la{TCONH}
\nabla_{\cal G}:= {\bf d}+ g_{0,0}+ 
\sum_{s,t\geq 0}z^{s}w^{t}\Bigl((s+t+1)G_{s+1,t+1}
(zdw-wdz)
 + (w\partial' - z \partial '')G_{s+1, t+1} \Bigr). 
\ee
The group $\C^* \times \C^*$ 
acts on the twistor plane: $(z,w) \lms (\lambda_1z, \lambda_2w)$, as well as 
on the Dolbeaux complex of the local system ${\rm End}{\cal V}$ on $X$, 
providing the Hodge bigrading: $\varphi_{s,t} \lms 
\lambda_1^{-s}\lambda_2^{-t}\varphi_{s,t}$.

There is an antiholomorphic 
involution 
$\sigma: (z,w)  \lms (\overline w, \overline z)$ 
of the twistor plane $\C^2$.

The connection $\nabla_{\cal G}$ has the following two basic properties: 

\begin{itemize}

\item 
The  connection $\nabla_{\cal G}$ is invariant under the action of the group $\C^* \times \C^*$.

\item The connection $\nabla_{\cal G}$ 
is invariant under the composition 
of  $\sigma$ with the 
 complex conjugation $c$ if and only if 
$G_{p,q}$ are imaginary (\ref{functionG}) and $g_{0,0}$ is real (\ref{formnu}):   
$$
(c\circ \sigma)^*\nabla_{\cal G} = \nabla_{\cal G} ~~
\mbox{if and only if} ~~ \overline G_{p,q} = -G_{p,q}, \quad \overline g_{0,0} = g_{0,0}.  
$$
\end{itemize}

Restricting 
the connection $\nabla_{\cal G}$ to  the {\it twistor line} $z+w=1$, parametrised by 
$z=1-u, w=u$,  we get a {connection}
\be \la{TCON}
{\rm Res}_{z+w=1}\nabla_{\cal G}:= {\bf d}+ g_{0,0}+ 
\sum_{s,t\geq 0}(1-u)^{s}u^{t}\Bigl((s+t+1)G_{s+1,t+1}du+ (u\partial'  - (1-u)\partial'')G_{s+1, t+1}\Bigr). 
\ee

\bd 
A connection $\nabla_{\cal G}$ as in (\ref{TCONH}) is called a twistor connection if 
its restriction to the twistor line  ${\rm Res}_{z+w=1}\nabla_{\cal G}$  is flat. 
\ed

It is easy to see that:
\begin{itemize}
\item The connection ${\rm Res}_{z+w=1}\nabla_{\cal G}$  is flat if and only if the 
datum $\{G_{p,q}, g_{0,0}\}$ 
satisfies  a Maurer-Cartan  system of 
 differential equations -- we call  such a   
datum $\{G_{p,q}, g_{0,0}\}$  a {\it Green datum}. 
\end{itemize}

The twistor connections
form an abelian tensor category in an obvious way.

\bt \la{2.22.08.1}
The abelian tensor 
category of twistor connections is canonically equivalent to 
the category of variations of mixed $\R$-Hodge 
structures on $X$. 
\et

A  twistor connection gives rise to 
a variation of mixed $\R$-Hodge structures on $X$ as follows.  
Restricting $\nabla_{\cal G}$ to $X  \times \{\frac{1}{2}, \frac{1}{2}\}$ 
we get  a flat connection 
$\nabla^{\frac{1}{2}}_{\cal G}$ on 
${\cal V}_{\infty}$.  We equip  the local system 
$({\cal V}_{\infty}, \nabla^{\frac{1}{2}}_{\cal G})$ with 
a structure of a variation of mixed $\R$-Hodge structures.  
The weight filtration 
is the standard weight filtration on a bigraded object: 
$$
W_n{\cal V}:= \oplus_{p,q\leq n}{\cal V}^{p,q}. 
$$
Let  $P$ be the operator of parallel transport 
for the connection $\nabla_{\cal G}$ along the line segment in the twistor plane  
from $X \times \{0,1\}$ to $X \times \{\frac{1}{2}, \frac{1}{2}\}$. 
Take the standard Hodge filtration on the bigraded object 
$F^p_{\rm st}{\cal V}_\C:= \oplus_{i\geq p}{\cal V}^{i, *}$ 
on the restriction of $p^*{\cal V}$ to $X \times \{0,1\}$.  
We define the Hodge filtration $F^{\bullet}$ on ${\cal V}_{\C}$ by applying 
 the operator of parallel transport $P$ to the standard Hodge filtration  on $p^*{\cal V}_{\C}$ at $\{0,1\}$: 
$$
F^{p}{\cal V}_\C := P(F^p_{\rm st}{\cal V}_\C).
$$ 
We prove that $({\cal V}_{\infty}, \nabla^{\frac{1}{2}}_{\cal G}, W_{\bullet}, F^{\bullet})$ 
is a variation of real 
mixed Hodge structures on $X $, and that any variation 
is obtained this way. 


\vskip 3mm
 Just recently M. Kapranov \cite{Ka} gave a neat interpretation of the 
category of mixed $\R$-Hodge structures as the 
category of $\C^*_{\C/\R}$-equivariant  (not necessarily flat) connections 
on $\C$. His description is equivalent to 
the one by twistor connections when $X$ is a point. Precisely,   twistor 
connections (\ref{TCONH}), restricted to $\C = (\C^2)^{\sigma}$, 
are natural representatives of the 
gauge equivalence classes of $\C^*_{\C/\R}$-equivariant  connections 
on $\C$.

\subsection{The twistor transform, Hodge DG algebra, 
and variations of real MHS} \la{hc1.6}
Theorem \ref{2.22.08.1} describes the category of 
variations of mixed $\R$-Hodge structures via twistor connections. 
The next step would be to describe complexes of 
variations of mixed $\R$-Hodge structures. Precisely, 
we want an alternative 
description of  the subcategory ${\rm Sh}^{\rm sm}_{{\rm Hod}}(X)$
of Saito's (derived) category ${\rm Sh}_{{\rm Hod}}(X)$ of mixed $\R$-Hodge sheaves on $X$ 
consisting of complexes 
whose cohomology are smooth, i.e. are 
 variations of mixed $\R$-Hodge structures. 
We are going to define a DG Lie coalgebra, 
and conjecture that the DG-category of comodules over it 
is a DG-enhancement of the triangulated category of smooth mixed $\R$-Hodge sheaves on a compact complex manifold $X$. 
To include the case of a smooth open complex variety $X$, one needs to impose 
conditions on the behavior of the forms forming the DG Lie coalgebra at infinity.

\paragraph{The semi-simple tensor category ${\rm Hod}_X$.} 
Let $X$ be a regular complex projective variety. 
Denote by ${\rm Hod}_X$ the 
category of variations of real Hodge structures on $X(\C)$. 
The category 
${\rm Hod}_X$ 
is a semi-simple abelian  
tensor 
category.  Denote by $w({\cal L})$ the weight of a pure variation 
${\cal L}$. 
The commutativity morphism is given by 
$
{\cal L}_1 \otimes {\cal L}_2 \lra  
(-1)^{w({\cal L}_1)w({\cal L}_2)}{\cal L}_2 \otimes {\cal L}_1. 
$

\paragraph{Hodge complexes.} 
Given a variation ${\cal L}$ of real Hodge structures on  
$X$, 
we introduce a complex ${\cal C}^{\bullet}_{{\cal H}_\R}({\cal L})$ and show that it 
calculates ${\rm RHom}_{{\rm Sh}_{{\rm Hod}}(X)}(\R(0), {\cal L})$ 
in the category of 
mixed $\R$-Hodge sheaves on $X$. 
We call it the {\it Hodge complex} of ${\cal L}$. 
When  ${\cal L} = \R(n)$ it is 
the complex calculating the weight $n$ Beilinson-Deligne 
cohomology of $X$. 

\vskip 3mm
The smooth de Rham complexes of variations of real 
Hodge structures on $X$ can be organized into a commutative 
DG algebra in the category 
${\rm Hod}_X$, called the de Rham DGA:
$$
{\cal A}_{X}:= \oplus_{{\cal L} } {\cal A}^{\bullet}({\cal L})\bigotimes 
{\cal L}^{\vee}, 
$$
where ${\cal L}$ runs through the isomorphism
 classes of simple objects in the category 
${\rm Hod}_X$. The product is given by the wedge 
product of differential forms and tensor product of variations.

We define a  commutative 
DG algebra ${\cal D}_{X}$ in the 
category ${\rm Hod}_X$. As an object,  it is 
a direct sum 
of complexes  
\be \la{complex}
{\cal D}_{X}:= \oplus_{{\cal L} } {\cal C}^{\bullet}_{\cal H}({\cal L})\bigotimes 
{\cal L}^{\vee}, 
\ee
where ${\cal L}$ runs through the isomorphism
 classes of simple objects in the category 
${\rm Hod}_X$. To define a DGA structure on ${\cal D}_{X}$, we 
recall the projection $p:X \times \R \to X$, and  introduce
our main hero:  

\paragraph{The twistor transform.}
It is a linear map 
$$
\gamma: \mbox{The Hodge complex of ${\cal V}$ on $X$} \lra 
\mbox{The de Rham complex of $p^*{\cal V}$ on $X\times \R$}, 
$$
Its definition is very  similar to the formula 
(\ref{TCON}) / (\ref{TCONH}) for 
the twistor connection -- 
see Definition \ref{12.09.ups.15}. 
The twistor transform is evidently injective, and gives 
rise to an injective linear map 
\be \la{12.16.ups.1}
\gamma: {\cal D}_X \hra {\cal A}_{X\times \R}.
\ee
Our key result in Section \ref{hc3sec} is the following theorem. 

\bt \la{12.8.ups.1s} Let $X$ be a complex manifold. Then the image of the 
twistor transform (\ref{12.16.ups.1}) is 
closed under  the differential and product, and thus   
is a DG subalgebra of the de Rham DGA. 
\et

Theorem \ref{12.8.ups.1s} provides ${\cal D}_X$ with 
a structure of a commutative DGA in the category ${\rm Hod}_X$. 
We call it the {\it Hodge DGA} 
of the complex variety $X$. 
The 
differential $\delta$ on ${\cal D}_X$  is obtained 
by conjugation by a degree-like operator $\mu$, see (\ref{mu}), of the 
natural differential  on (\ref{complex}). 

It is surprising that one can 
realize the Hodge DGA inside of the de Rham DGA on 
$X\times \R$.

\vskip 3mm
Given a semi-simple abelian tensor category ${\cal P}$, denote by 
 ${DGCom}_{\cal P}$ and $DGCoLie_{\cal P}$ the categories of DG commutative 
and Lie coalgebras in the category ${\cal P}$. 
Recall the two standard functors
$$
{\cal B}: DGCom_{\cal P} \lra DGCoLie_{\cal P}, \qquad {\cal C}: 
DGCoLie_{\cal P} \lra DGCom_{\cal P}.  
$$
The functor ${\cal B}$ is the bar construction followed 
by projection to the indecomposables. The functor ${\cal C}$ is 
given by the Chevalley standard complex: 
$$
({\cal G}^{\bullet}, d) \lms {\cal C}({\cal G}):= ({\rm Sym}^*({\cal G}^{\bullet}[1]), \delta), \qquad 
\delta:= d_{Ch} + d.
$$

Applying the functor ${\cal B}$ to the commutative DGA  
${\cal D}_{X}$ in the category 
${\rm Hod}_X$ we get a DG Lie coalgebra ${\cal L}_{{\cal H}; X}^*:= 
{\cal B}({\cal D}_{X})$ in the same 
category. Let ${\cal L}_{{\cal H}; X}:= H^0({\cal L}_{{\cal H}; X}^*)$ be the Lie coalgebra 
given by its zero cohomology. 
We show that
 Theorem \ref{2.22.08.1}  is 
equivalent to
\bt\la{12.09.ups.20}
The category of comodules over the Lie coalgebra ${\cal L}_{\cal H; X}$ in the category 
${\rm Hod}_X$ is canonically equivalent to the 
category of variations of mixed $\R$-Hodge structures on $X$. 
The equivalence is given by the functor 
${\rm gr}^W$ of 
the associate graded for the weight filtration.
\et

When $X$ is a point it provides a 
canonical set of generators for the Lie algebra of the Hodge Galois group. 
The DGA ${\cal D}_X$ contains as a 
sub DGA the Hodge-Tate algebra of A.Levin \cite{L}, and our 
description of variations of Hodge-Tate structures essentially 
coincides with the one 
in {\it loc. cit}.

\bcon \la{12.1.ups.7}
The category ${\rm Sh}^{\rm sm}_{\rm Hod}(X)$ of smooth 
complexes of 
real Hodge sheaves on $X$ 
is equivalent to the DG-category of DG-modules over the DG Lie coalgebra 
${\cal L}^*_{{\cal H}; X}$.
\econ

\subsection{Mixed $\R$-Hodge structure on $\pi_1^{\rm nil}$  via Hodge correlators} \la{sec1.9}

Recall that ${\rm L}_{X,S^*}$, see (\ref{varl}), 
 is a free Lie algebra generated by 
${\rm gr}^WH^1(X-S)$. When the data $(X, S, v_0)$ varies, the latter 
forms a
 variation ${\cal L}$ of real Hodge structures over the 
enhanced moduli space ${\cal M}'_{g, n}$, $n=|S|$, 
(we add a tangent vector $v_0$ to the standard data). 
We are going to equip it with a Green datum, thus getting a variation 
mixed $\R$-Hodge structures. 

\vskip 3mm
A point $s_0\in X $ determines a $\delta$-current 
given by the evaluation of a test 
function at $s_0$. Viewed as 
a generalized volume form $\mu$, it provides a Green function $G(x,y)$ 
up to a constant. We use the tangent vector $v_{0}$ at $s_0$ 
to specify the constant. Namely, let 
$t$ be a local parameter at $s_0$ with $\langle dt, v_0\rangle=1$. 
We normalize the Green function so that 
$G(s_0, s) - \frac{1}{2\pi i} \log |t|$ vanishes at $s=s_0$. 

Let ${\bf G}_{v_0}$ be the corresponding Green 
operator (\ref{5.22.06.1qa}). It 
 is an endomorphism of the smooth bundle ${\cal L}_{\infty}$. 
Thus the construction described in Section \ref{1.5} provides  
a collection of mixed $\R$-Hodge structures in the fibers of ${\cal L}$. 

Furthermore, the Hodge correlators in families 
deliver in addition to the Green operator ${\bf G}_{v_0}$
 an operator valued $1$-form $g_{0,0}$. 
We prove in Section \ref{hc6sec}   
that the pair ${\bf G}_{v_0}:= ({\bf G}_{v_0}, g_{0,0})$ 
satisfies the Green data differential equations. Therefore 
we get 
a variation of mixed $\R$-Hodge structures. 

\begin{theorem} \label{9.16.05.2} When the data 
$(X, S, v_0)$ varies, 
the 
variation of mixed $\R$-Hodge structures given by the Green datum 
${\bf G}_{v_0}$ 
is isomorphic to the 
variation formed by the standard   
mixed $\R$-Hodge structures on $\pi_1^{\rm nil}(X -S, v_0)$. 
\end{theorem}

\vskip 3mm
The standard MHS on $\pi_1^{\rm nil}(X -S, v_0)$ is described by 
 using Chen's theory of iterated integrals \cite{Ch}. 
Our approach is different: it is
given by integrals of non-holomorphic differential forms over products
of copies of $X $. Theorem \ref{9.16.05.2} implies that the two sets
of periods obtained in these two descriptions coincide. For example
the periods of the real MHS on $\pi_1^{\rm nil}({\Bbb P}^1-\{0,1,\infty\}, v_{0})$, 
where $v_0 = \partial/\partial t$ at $t=0$,  x
are given by the multiple $\zeta$-values \cite{DG}. So Theorem
\ref{9.16.05.2} implies that 
the Hodge correlators in this case are $\Q$-linear combinations of the
multiple $\zeta$-values, which is far from being obvious from their
definition. 

Here is another benefit  of our approach: for a modular curve $X$, 
it makes obvious  that the Rankin-Selberg convolutions are periods of the 
motivic fundamental group of $X- \{\mbox{the cusps}\}$.

\subsection{Motivic correlators on curves} \la{mccurves}
Let $X$ be a regular projective curve over a field $F$, $S$  a 
non-empty 
subset of $X(F)$, and $v_0$  a tangent vector at a point $s_0 \in S$ defined 
over $F$. 
The {\it motivic fundamental group} of $X-S$ with the tangential base point 
$v_0$  
is supposed to be   a pro-unipotent Lie algebra ${\rm L}(X-S, v_0)$ 
in the hypothetical abelian category of mixed motives over $F$. 

At the moment 
we have it either in realizations, $l$-adic or 
Hodge \cite{D1}, or when $X = {\Bbb P}^1$ and $F$ 
is a number field \cite{DG}. We work in 
one of these settings, or assume the 
existence of the abelian category of mixed motives. 
A precise description of possible settings see in Section 
\ref{hc8.2sec}.

Since ${\rm L}(X-S, v_0)$ is a pro-object in the category of mixed motives, 
the Lie algebra ${\rm L}_{\rm Mot/F}$, see Section \ref{secaap},  acts by special derivations 
on its associate graded for the weight filtration:
\begin{equation} \label{7.18.06.2}
{\rm L}_{\rm Mot/F} \lra {\rm Der}^S\Bigl({\rm gr}^W{\rm L}(X-S, v_0)\Bigr). 
\end{equation} 
Starting from the pure motive ${\rm gr}^W{H}_1(X-S)$ 
rather then from its Betti realization, and following (\ref{mvla}), we define 
the motivic version 
${\cal C}{{\cal L}ie}_{{\cal M}; X, S^*}$ 
of ${\cal C}{{\cal L}ie}_{X,S^*}$. It is a Lie algebra in the semisimple 
abelian category of pure motives. Let 
${\cal C}{{\cal L}ie}^\vee_{{\cal M}; X, S^*}$ be the dual Lie coalgebra. 
Thanks to 
(\ref{9.17.05.5}), there is a
Lie algebra  isomorphism (which gets a Tate twist in the motivic setting) 
$$ 
{\rm Der}^S\Bigl({\rm gr}^W{\rm L}(X-S, v_0)\Bigr) = 
{\cal C}{{\cal L}ie}_{{\cal M}; X, S^*}(-1).
$$
So dualising (\ref{7.18.06.2}), 
we arrive at a  map of Lie coalgebras, called the {\it motivic correlator map}:
\begin{equation} \label{7.18.06.1}
{\rm Cor}_{\rm Mot}: {\cal C}{{\cal L}ie}^\vee_{{\cal M}; X, S^*}(1)
 \lra {\cal L}_{\rm Mot/F}.
\end{equation} 
The left hand side 
is decomposed into a direct sum of cyclic tensor 
products of simple  
pure motives. Their images under the motivic correlator map 
are called {\it motivic correlators}.

\vskip 2mm
Our next goal is to relate the motivic and Hodge correlators.

\paragraph{The canonical period map.}  
\bl \la{4.5.10.1}
A choice of generators of the Lie algebra ${\rm L}_{\rm Hod}$ provides a 
period map 
\be \la{periio}
p: {\cal L}ie_{\rm Hod} \lra i\R.
\ee 
\el

{\bf Proof}. A choice of generators of the Lie algebra ${\rm L}_{\rm Hod}$ 
is just the same thing as a projection 
\be \la{splitmap}
{\cal L}ie_{\rm Hod} \lra \oplus_{[H]}{\rm Ext}_{\R-{\rm MHS}}^1(\R(0), H)\otimes H^\vee,
\ee
where the sum is over the set of isomorphism classes of simple objects in the category of 
$\R$-Hodge structures. Next, for any simple object $H$ 
there is a canonical map 
\be \la{4.6.10.5}
p_{[H]}: {\rm Ext}_{\R-{\rm MHS}}^1(\R(0), H)\otimes H^\vee\lra i\R.
\ee
Recall that given an $\R$-Hodge structure $H$, we have the formula
\be \la{11.19.ups.1}
{\rm Ext}^1_{\R-{\rm MHS}}(\R(0), H) = {\rm CoKer}\Bigl(W_0H \oplus F^0(W_0H)_\C \lra (W_0H)_\C\Bigr).
\ee
Thus for a pure $\R$-Hodge structure $H$ of Hodge degrees $(-p,-q) + (-q,-p)$, 
where $p,q \geq 1$, we have  ${\rm Ext}^1_{\R-{\rm MHS}}(\R(0), H) = (H \otimes \C)/H = H \otimes_\R i\R$. 
So we arrive to a canonical pairing (\ref{4.6.10.5}).

We define a period map (\ref{periio}) as the composition of the map (\ref{splitmap}) 
with the sum of the maps $p_{[H]}$. 
The Lemma is proved. 

\vskip 2mm
Therefore the 
 canonical generators of the Lie algebra ${\rm L}_{\rm Hod}$, discussed in Section 1.7, 
give rise to  a canonical  period map.

\paragraph{Relating motivic and Hodge correlators.} 
Let $X$ be a smooth complex curve. Just like in the case of the motivic correlators, 
the mixed $\R$-Hodge structure on ${\rm L}(X-S, v_0)$ plus the isomorphism (\ref{9.17.05.5}) lead to 
the Hodge version of the motivic correlator map (\ref{7.18.06.1}):
 \be \la{4.5.10.2}
{\rm Cor}_{\rm Hod}: {\cal C}{{\cal L}ie}^\vee_{{\cal H}; X, S^*}(1) \lra
 {\cal L}ie_{\rm Hod}.
\ee
Here ${\cal C}{{\cal L}ie}^\vee_{{\cal H}; X, S^*}$ 
is the cyclic Lie coalgebra assigned to the $\R$-Hodge structure ${\rm gr}^WH^1(X-S;\R)$. 
Forgetting the Hodge bigrading, there is an isomorphism of Lie algebras over $\C$
$$
{\cal C}{{\cal L}ie}^\vee_{{\cal H}; X, S^*}\otimes \C \lra {\cal C}{{\cal L}ie}^\vee_{X, S^*}.
$$
Using this, and combining (\ref{4.5.10.2}) with the canonical period map (\ref{periio}) we arrive at a map
\be \la{ccoommppoo}
 {\cal C}{{\cal L}ie}^\vee_{X, S^*}(1) \stackrel{{\rm Cor}_{\rm Hod}} {\lra}
 {\cal L}ie_{\rm Hod}\otimes \C \stackrel{p}{ \lra}  \C.
\ee 
The following result is an immediate corollary of Theorem \ref{9.16.05.2} and 
the definitions. 
\bt\la{4.6.10.2}
The  composition (\ref{ccoommppoo}) coincides with the Hodge correlator map. 
\et
Let us assume now the motivic formalism. Then the $\R$-Hodge realization functor 
provides a homomorphism of Lie coalgebras
$$
r_{\rm Hod}: {\cal L}_{\rm Mot/\C} \lra {\cal L}_{\rm Hod}. 
$$
\bc \la{4.6.10.1} Assuming the motivic formalism, 
the complexification of the  composition
$$
 {\cal C}{{\cal L}ie}^\vee_{{\cal M}; X, S^*}(1) \stackrel{{\rm Cor}_{\rm Mot}}{\lra}
 {\cal L}_{\rm Mot/\C} \stackrel{r_{\rm Hod}}{ \lra} {\cal L}ie_{\rm Hod} \stackrel{p}{ \lra} i\R.
$$ 
coincides  with the Hodge correlator map. In other words, we get a commutative diagram
$$
\begin{array}{ccc}
{\cal C}{{\cal L}ie}^\vee_{{\cal M}; X, S^*}(1)& \stackrel{{\rm Cor}_{\rm Mot}}{\lra}&{\cal L}_{\rm Mot/\C} \\
&&\downarrow r_{\rm Hod}\\
&\searrow{\rm Cor}_{\rm Hod}  &{\cal L}ie_{\rm Hod}\\
&& \downarrow p\\
&&i\R
\end{array}
$$
\ec

\paragraph{Conclusion.} The source of the motivic correlator map (\ref{7.18.06.1}) 
is an explicitly defined Lie coalgebra. 
Therefore motivic correlators come together  with an explicit formula for their coproduct. 
The real period of the motivic correlator is the Hodge correlator. 
This, togerther with basic isomorphism (\ref{biso}), is all we need for the arithmetic analysis 
of the Hodge correlators.

\begin{figure}[ht]
\centerline{\epsfbox{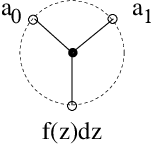}}
\caption{The Feynman diagram  for Beilinson's element.}
\label{feyn-10}
\end{figure}

\paragraph{Example: Rankin-Selberg integrals and Beilinson's elements as correlators.} 
Here is how the story discussed in Section 1.2 fits in the correlator framework. 
Recall that $M$ is 
a modular curve, $\overline M$ its compactification, $a,b$
degree zero cuspidal divisors, and $g_a, g_b$ are invertible functions  on $M$ 
whose divisors are integral multiples of $a$ and $b$. 
The Green function $G(a,t)$ of the divisor $a$ 
is an integral multiple of $\log|g_a(t)|^2$, and similarly for $G(b,t)$. Therefore 
integral (\ref{BEERS213}) equals, up to a rational multiple, to the Hodge correlator integral 
provided by the cyclic word ${\cal C}(\{a_0\}\otimes \{a_1\} \otimes f(z)dz)$, see Fig \ref{feyn-10}: 
\be \la{3.28.10.2}
{\rm Cor}_{\cal H}\Bigl((\{a_0\}\otimes \{a_1\} \otimes f(z)dz)\Bigr) = 
\int_{M(\C)} G(a_0, t) d^\C G(a_1, t) \wedge f(z)dz. 
\ee
The latter coincides with a Rankin-Selberg integral. 

Let $M_f$ be the dual to the pure weight two motive corresponding to the Hecke eigenform $f(z)$. 
It is a direct summand of the motive $H_1(\overline M)$.  
The motivic correlator is an element
\be \la{2.28.10.1}
{\rm Cor}_{\rm Mot}\Bigl((\{a_0\}\otimes \{a_1\} \otimes M^\vee_f)(1)\Bigr) \in {\cal L}_{\rm Mot}.
\ee
It lies in the $M_f(1)^\vee$-isotipical component of ${\cal L}_{\rm Mot}$. 
Since $a_i$ are torsion classes in the Jacobian of $\overline M$,  we prove  (Lemma \ref{7.6.06.1}) 
that (\ref{2.28.10.1}) is in the kernel of the coproduct $\delta$. 
So, combining with the isomorphism (\ref{biso}),  we get an element
$$
{\rm Cor}_{\rm Mot}\Bigl((\{a_0\}\otimes \{a_1\} \otimes M_f^\vee)(1)\Bigr) \in 
{\rm Ext}_{\rm Mot/\Q}^1(\Q(0), M_f(1))\otimes M_f(1)^\vee.
$$
Using Beilinson's conjecture \cite{B}, and the conjectural motivic Leray spectral sequence, 
we have
\be \la{EXTK2}
{\rm Ext}_{\rm Mot/\overline \Q}^1(\Q(0), H_1(\overline M)(1)) = K_2(\overline M)\otimes \Q. 
\ee

Indeed, for any regular projective curve $X$ over a number field 
$F$ with the structure map $p: X \to {\rm Spec}(F)$, 
we should have using (\ref{Beilf}), (\ref{BMLSS1}), and (\ref{BMLSS}):
$$
K_2(X)\otimes \Q =
{\rm Ext}_{\rm Mot/F}^2(\Q(0), p_*\Q(2)_X) = 
{\rm Ext}_{\rm Mot/F}^1(\Q(0), H^1(X)(2)).
$$
To get the last isomorphism notice that 
${\rm Ext}_{\rm Mot/F}^i(\Q(0), M)=0$ for $i\not =1$ for any pure simple motif $M$ over 
a number field $F$ 
non-isomorphic to $\Q$. 
By Poincare duality, 
$H^1(X)=H_1(X)(-1)$. 

This implies a similar claim for the direct summands of the motif $H^1(\overline M)$. 

Using  (\ref{EXTK2}),  the motivic correlator is related to Beilinson's element (\ref{BEE}) by the formula
$$
{\rm Cor}_{\rm Mot}\Bigl((\{a_0\}\otimes \{a_1\} \otimes M_f^\vee)(1)\Bigr)  = 
\{g_a, g_b\} \otimes M_f(1)^\vee \in K_2(\overline M)\otimes M_f(1)^\vee.
$$
By Corollary \ref{4.6.10.1}, 
the real period map on this element is given by the Rankin-Selberg integral. 

\paragraph{A conjecture expressing $L({\rm Sym}^2M_f, 3)$ via Hodge correlators on modular curves.} 
Let us assume, for simplicity only, that the coefficients of the Hecke form $f(z)$ are  in $\Q$, i.e. 
the motive $M_f$ is defined over $\Q$. Then 
Beilinson's conjecture \cite{B} predicts that 
$$
{\rm dim}_\Q {\rm Ext}_{\rm Mot/\Q}^1(\Q(0), {\rm Sym}^2M_f(1))=2.
$$

\bcon \la{9.30.13.1} The space 
 ${\rm Ext}_{\rm Mot/\Q}^1(\Q(0), {\rm Sym}^2M_f(1))$ is generated by 
the motivic correlators 
\be \la{4.3.10.1}
{\rm Cor}_{\rm Mot}\Bigl((\{a_1\}\otimes \{a_2\}\otimes M^\vee_f \otimes M^\vee_f)(1) 
 + 
(\{b_1\}\otimes M^\vee_f \otimes \{b_2\} \otimes M^\vee_f)(1)\Bigr)
\ee
for certain degree zero cuspidal divisors $a_i$ and $b_i$  
on the universal modular curve, which satisfy an explicit  ``coproduct zero'' condition. 
\econ
The two terms in (\ref{4.3.10.1}) reflect the two different kinds 
of the cyclic tensor product of the following four pure motives: 
$\Q(-1), \Q(-1), M^\vee_f, M^\vee_f$. 

Precisely, an element (\ref{4.3.10.1}) gives rise to an element of 
${\rm Ext}_{\rm Mot/\Q}^1(\Q(0), {\rm Sym}^2M_f(1))$ if and only if it is killed by the 
coproduct in the motivic Lie coalgebra. We elaborate in Example 2 of Section 11.3 
 the ``coproduct zero'' condition for any pair of cuspidal divisors $a_i, b_i$ 
entering element (\ref{4.3.10.1}).  
Finally, one should 
relate $L({\rm Sym}^2M_f, 3)$ to the determinant of the $2 \times 2$  
matrix, whose entries are given by the Hodge correlators of elements (\ref{4.3.10.1}).

\subsection{Motivic multiple $L$-values} Here are
 some arithmetic applications obtained by picking certain specific
 pairs $(X,S)$, as well as certain specific motivic correlators.  

The image ${\cal G}(X-S, v_0)$ of the map (\ref{7.18.06.1}) is called the {\it 
motivic Galois 
Lie algebra of $X-S$}. In the $l$-adic setting it is the 
Lie algebra of the image of the Galois group 
${\rm Gal}(\overline F/F(\mu_{l^\infty}))$ 
acting on the pro-$l$ completion 
$\pi^{(l)}_1(X-S, v_0)$ of the fundamental group. 

\vskip 3mm
Specializing to the case 
when $X$ is one of the following:  

(i) ${\Bbb G}_m$, ~~(ii) an elliptic CM curve, 
~~(iii) a Fermat curve, ~~(iv) a modular curve, 

\noindent
and choosing the subset $S$ 
appropriately, we arrive at the definition of 
{\it motivic multiple $L$-values} 
corresponding to certain algebraic Hecke characters of 
$\Q$, $\Q(\sqrt{-d})$, $\Q(\mu_N)$, or weight two modular 
forms, respectively. 
They span a 
 Galois Lie coalgebra, called the  {\it multiple $L$-values 
Lie coalgebra}.  

In the $l$-adic setting there is one more 
interesting case: ~(v) a Drinfeld modular curve. 

\vskip 3mm
Given triples  $(X',S',v_0')$ and $(X,S,v_0)$, 
a map $p:X'\to X$ such that $p(S') = S$ and $p(v_0')=v_0$ 
 induces a map of the motivic fundamental Lie algebras $p: {\rm L}(X'-S', v'_0)\to 
{\rm L}(X-S, v_0)$, and hence 
 a map of the Galois Lie algebras $p: {\cal G}(X'-S', v'_0)\to {\cal G}(X-S, v_0)$, 
and the dual map of the Galois Lie coalgebras. 

In each of the  cases (i)-(v),  
the open curves $X-S$ form a  tower. 
For instance in the case (iv) it is the modular tower. 
Therefore the  corresponding  Galois 
Lie coalgebras  form an inductive system. 
In many cases there is a canonical  
``averaged tangential base point'' (Section \ref{hc8.4sec}), 
so below we ignore the base point. 
There is an adelic description of the limiting Galois 
 Lie coalgebras. It is given in terms of 
 the corresponding base field in the first three cases, and via the 
automorphic representations attached to the weight two modular forms in the 
modular curve case. In the latter case there is a 
generalization dealing with modular  forms of arbitrary integral 
weight $k \geq 2$, which we will elaborate elsewhere.

Here is a description of the 
pairs $(X,S)$  and the corresponding towers. 

\begin{enumerate}
\item  Let $X = {\Bbb P}^1$, and $S$ is an arbitrary subset of ${\Bbb P}^1$. 
This is the situation studied in \cite{G7}, \cite{G9}. 

\subitem a) A very interesting case is when 
$S = \{0\} \cup  \{\infty\} \cup \mu_N$. 
It was studied in \cite{G4}, \cite{DG}. 

{\it The tower}. The curves ${\Bbb G}_m - \mu_N$, 
parametrized by positive integers 
$N$, form a ``tower'' for the isogenies 
${\Bbb G}_m - \mu_{MN} ~\to~ {\Bbb G}_m - \mu_N$. 

\item  The curve  is an elliptic curve $E$, and $S$ is any subset. 
An interesting special case is when $S= E[N]$ is 
the subgroup of $N$-torsion points of 
$E$. {\it The tower} in this case 
is formed by the isogenies   
$E-E[MN] ~\to~ E-E[N]$.

There are the following ramifications: 

\subitem a) The  universal elliptic 
curve ${\cal E}$ with the level $N$ structure over the modular curve. 

\subitem s) The curve  $E$ is a CM curve, ${\cal N}$ an 
ideal in the endomorphism ring ${\rm End}(E)$, and 
$S$ is the subgroup $E[{\cal N}]$ of the ${\cal N}$-torsion points.

\item  The curve  is the Fermat curve ${\Bbb F}_N$, given in ${\Bbb P}^2$ by the equation
$
x^N + y^N = z^N
$, 
and $S$ is the intersection of ${\Bbb F}_N$ with 
the coordinate triangle $x=0$, $y=0$, $z= 0$ in ${\Bbb P}^2$. 

{\it The tower}. It is given by  the  
maps ${\Bbb F}_{MN} \to{\Bbb F}_{N}$, $(x,y,z) \lms (x^M, y^M, z^M)$. 

\item  The curve $X = X(N)$ is the level $N$ modular curve, $S$ is the set of the cusps. 
There is a tower of modular curves, and especially interesting situation appears at the limit, 
i.e. when 
$
{\cal M} := \lim_{\longleftarrow}X(N) 
$  
is the universal modular curve and $S$ is the set of its cusps.

\item  Let $A$ be the ring of rational functions on a regular projective curve over a finite field, which are 
regular at a chosen point $\infty$ of the curve. Let $I$ be an ideal in $A$. 
The Drinfeld modular curve 
is the moduli space of the rank $2$ elliptic modules 
with the $I$-level structure. Take it as our curve $X$, and let  $S$ is the set of the cusps. 
There is a tower of the Drinfeld modular curves. 
In  the limit we get the universal 
 Drinfeld modular curve.

\end{enumerate}
 
\vskip 3mm
In each of these cases we get a supply of 
elements of the motivic Lie coalgebras. These are

\begin{enumerate}
\item ``Cyclic'' versions of motivic multiple  polylogarithms. 

\subitem a) Motivic multiple  $L$-values for the Dirichlet 
characters of $\Q$. 

\item  Motivic multiple 
 elliptic polylogarithms, and their restrictions to the torsion points.


\subitem a) Higher (motivic) analogs of modular units. 

\subitem b) Motivic multiple $L$-values for Hecke Gr\"ossencharacters 
of imaginary quadratic fields.

\item Motivic multiple $L$-values for the Jacobi sums Gr\"ossencharacters 
of $\Q(\mu_N)$. 

\item  Motivic multiple $L$-values for the modular forms of weight $2$. 

\item  Similar objects in the function field case. 
\end{enumerate}

We prove in Section 11 that the real periods of the motivic multiple 
 elliptic polylogarithms are the generalized Eisenstein-Kronecker series, 
defined in \cite{G1} as integrals of the classical Eisenstein-Kronecker series, and 
obtain similar results in the rational case.

\paragraph{What are the multiple $L$-values?} 
In each of the above cases the motivic multiple $L$-values are certain 
natural collections of elements of the motivic Lie coalgebra, 
closed under the Lie cobracket. Their periods  are numbers, called 
the multiple $L$-values. 

The ``classical'' motivic $L$-values 
are the elements killed by the Lie cobracket. 
Moreover, they are precisely 
the cogenerators of the corresponding motivic 
Lie coalgebra, providing the name for the latter. 
Their periods are 
special values of $L$-functions. 

 The ``classical'' motivic $L$-values are 
 motivic cohomology classes. In general the motivic multiple $L$-values are 
no longer motivic cohomology classes.

\paragraph{Example.} The first multiple $L$-values discovered 
were Euler's multiple $\zeta$-numbers
$$
\zeta(n_1, ..., n_m) = \sum_{0<k_1< \dots < k_m}
\frac{1}{k_1^{n_1} \ldots k_m^{n_m}}. 
$$ 
 The motivic $\zeta$-values are the elements $\zeta_{\cal M}(2n-1) 
\in {\rm Ext}_{\rm Mot/\Z}^1(\Q(0), \Q(n)$. The latter 
${\rm Ext}$-group is one-dimensional, and is generated by $\zeta_{\cal M}(2n-1)$. 
The motivic multiple zeta values  
$\zeta_{\cal M}(n_1, ..., n_m)$ are the elements of the 
motivic Lie coalgebra of the category of mixed Tate motives over 
${\rm Spec}(\Z)$. The conjecture that they span the latter was recently
 proved by Brown \cite{Br}.  

\vskip 3mm
Similarly, each motivic multiple $L$-values Lie coalgebra  
is tied up to a certain subcategory of the category of mixed motives. 
A very interesting question is how far  it is from the motivic Lie coalgebra 
which governs that mixed category. 
See \cite{G4} for the case 1a) 
 and \cite{G10} for the case 2b), 
where the gap in each case was 
related to the geometry of certain modular varieties.

In all above cases the motivic cohomology group 
responsible for the $L$-values were one-dimensional.

\paragraph{Question.} Can one define multiple $L$-values in a 
more general situation?

\subsection{Coda: 
Feynman integrals and  motivic correlators}

We suggested in Section 8 of \cite{G7} the following picture relating 
Feynman integrals and mixed motives. 
If correlators of a Feynman integral 
are periods, they should be upgraded to motivic correlators, which 
lie in the motivic Lie coalgebra, and whose periods are the Feynman correlators. 
The obtained motivic correlators should be closed under the coproduct in the 
motivic Lie coalgebra 
- otherwise the original set of Feynman correlators was not complete. 
Finally,  the coproduct of the motivic correlators should be 
calculated as follows. There should be a 
combinatorially defined Lie coalgebra -- the renormalization Lie coalgebra -- 
and the motivic correlators should provide 
a homomorphism 
\be \la{RenMot}
\mbox{the renormalization Lie coalgebra} \lra \mbox{the motivic Lie coalgebra.} 
\ee

This paper gives an example of realization 
of this program. 
Indeed, the Hodge correlators are the tree level correlators 
of the  Feynman integral from Section \ref{1.4}.  
We upgrade them 
to motivic correlators. 
They span a Lie coalgebra, 
which is nothing else but the Galois Lie coalgebra ${\cal G}(X-S, v_0)$. 
The Lie coalgebra ${\cal C}{{\cal L}ie}_{X,S^*}^{\vee}$ 
plays the role of the renormalization Lie coalgebra:
 There is a canonical Lie coalgebra homomorphism 
${\cal C}{{\cal L}ie}_{X,S^*}^{\vee} \to {\cal G}(X-S, v_0)$ 
describing the coproduct of  
motivic correlators, which we think of as the homomorphism (\ref{RenMot}).

\paragraph{Problem.} Develop a similar picture 
for Feynman integrals whose correlators are periods.

\subsection{The structure of the paper} The paper is build from four 
parts. Their architecture is described as follows. 

\begin{enumerate}

\item  In Section \ref{hc2sec} we introduce a key tool to define Hodge correlators: a polydifferential 
operator $\omega_m$. Its properties become evident in its twistor version, $\widehat \omega_m$. 

Using the operator $\omega_m$, and the sum of the plane decorated trees construction, 
we define in Section  \ref{formersec2} Hodge correlators on punctured curves. 

We show in Section \ref{hc12sec} that the Hodge correlators on punctured curves
indeed are correlators of a Feynman integral, i.e. can be defined via a 
Feynman integral procedure. 
However we do not use this anywhere else in the paper.

\item 
In Section \ref{hc3sec} we define a Hodge complex of a variation of Hodge structures.  
We introduce a commutative DGA structure on the sum of Hodge complexes over 
isomorphism classes of simple variations  of Hodge structures on a complex manifold $X$.  
We call it the Hodge DGA of $X$. 
We define a twistor transform and show that it embeds 
the Hodge DGA of $X$ into the De Rham DGA of $X \times \R$. 
This elucidates the definition of the Hodge DGA.

In Section \ref{hc4sec}   we 
show that variations of mixed Hodge structures are described by twistor connections. 
The flatness of the twistor connection encodes the Griffith transversality. 
It amounts to Maurer-Cartan type differential equations on the Green data describing 
a twistor connection. 

\item In Section \ref{hc5sec}, given a symplectic vector space $H$ and a set $S$, 
 we define a Lie algebra ${\cal C}_{H, S}$. It 
generalizes the Lie algebras defined in \cite{G4}, \cite{G9} 
in the case $H=0$. We introduce an $L_\infty$-algebra    
of plane trees   decorated by $H \oplus\Q[S]$.   
It is a resolution of the Lie algebra ${\cal C}_{H, S}$. 

The Lie algebra ${\cal C}_{H, S}$ has a Lie subalgebra ${\cal C}{\cal L}ie_{H, S}$.

Our basic example is assigned to a smooth family of compact curves $p: X\to B$ 
with a smooth divisor $S \to B$, and a base point $s_0\in S$. Let $S^*:= S-s_0$.  
The Lie algebra ${\cal C}{\cal L}ie_{H, S^*}$ assigned to $H= R_1p_*(X, \R)$ 
is denoted by ${\cal C}{\cal L}ie_{X, S^*}$. 

In Section \ref{hc6sec}, given such a family  $p: X\to B$, 
 we define the Hodge correlator twistor connection $\nabla_{\bf G}$ on $B \times \C^2$. 
It takes values in a variation of Lie algebras  ${\cal C}{\cal L}ie_{X, S^*}$. 
The Hodge correlators defined in Section \ref{formersec2} are just the coefficients 
of the twistor connection $\nabla_{\bf G}$. The 
differential equations for the 
Hodge correlators express the fact that the restriction of the twistor connection 
to $B\times \C$, where $\C\subset \C^2$ is the twistor line $z+w=1$, is flat. 

We stress that our construction in Section  \ref{hc6sec} delivers the 
Hodge correlator twistor connection $\nabla_{\bf G}$ directly, using the twistor 
operator $\widehat \omega_m$ and the sum over all plane trivalent trees construction. 
So the very notion of twistor connection is naturally build in 
the definition of the Hodge correlators. 

\item In Section \ref{hc7sec} we identify the  
Lie algebra ${\cal C}{{\cal L}ie}_{X,S^*}$ with the Lie algebra of special derivations 
of ${\rm gr}^W\pi^{\rm nil}_1(X-S, v_0)$, where $v_0$ is a tangent vector at $s_0$. 
Since the Hodge correlator connection $\nabla_{\bf G}$ lies in 
 ${\cal C}{{\cal L}ie}_{X,S^*}$, it gives rise to  
a special derivation ${\bf G}$
of the Lie algebra ${\rm gr}^W\pi^{\rm nil}_1(X-S, v_0)$. 
By the Section \ref{hc4sec} construction, it 
provides a 
real mixed Hodge structure on the variation $\pi^{\rm nil}_1(X-S, v_0)$ over  $B$. 
In Section \ref{hc9sec} we prove that it is the standard one.

In Section \ref{hc8sec} we 
translate  results of Section \ref{hc7sec} into the motivic framework. 
Using this we introduce {\it motivic correlators}. They 
 are canonical elements 
in  the motivic Lie coalgebra. 
The Hodge correlators are the canonical 
real periods of the motivic correlators. 

In Section \ref{hc10sec} we consider examples. We show how the 
classical and elliptic polylogarithms appear  
as Hodge correlators for simple Feynman diagrams. 
We prove that for an elliptic curve the Hodge correlators 
can be expressed by the multiple Eisenstein-Kronecker series. 

In Section \ref{hc11sec} we show that Hodge / motivic  
correlators on modular curves 
generalize the Rankin-Selberg integrals / Beilinson's elements in $K_2$ 
of modular curves. 
Going to the limit in the tower of modular curves we get an 
automorphic adelic description of the correlators on modular curves. 
\end{enumerate}

 \paragraph{Acknowledgments.} I am very grateful to A.A. Beilinson for 
many fruitful discussions, and especially for
 encouragement over the years to work on Hodge correlators. 

Many ideas of this work 
 were worked out, and several Sections written,  
during my stays at the MPI(Bonn) and IHES. 
A part of this work was written during my stay at the 
Okayama University (Japan) in May 2005. 
I am indebted to Hiraoki Nakamura for the hospitality there. 
I am grateful to these institutions for the hospitality and support. 
I was supported by the NSF grants DMS-0400449, DMS-0653721,  DMS-1059129 and 
DMS-1301776. 
I am very grateful to referees for the extraordinary job. 
I am grateful to C. Hertling for pointing out 
some errors, and sending me notes of his seminar 
talks in 2008-2009 on the paper.  

                
\section{Polydifferential operator $\omega_{m}$} \la{hc2sec}

In Section \ref{hc2sec} we introduce 
a polydifferential operator $\omega_{m}(\varphi_1, ..., 
\varphi_{m})$. 

In Section \ref{hc2.1.sec} we present a down to earth definition of the 
 polydifferential operator $\omega_{m+1}$ on the Dolbeaut complex of a complex manifold.
The definition extends naturally to a {\it cohomological Dolbeaut bicomplex}, axiomatizing 
the Dolbeaut bicomplex of a variation of Hodge structures. 

A much better definition of a twistor cousin $\widehat \omega_{m}$  of the 
 operator $\omega_{m}$ 
is given in Section \ref{hc2.2sec}. 
The twistor approach is more  conceptual, and allows to give transparent proofs of the main properties of the 
operator $\omega_{m}$. However since it looks a bit out of the bleu  
we include both.  

The twistor approach  makes evident the two main features of the operator $\omega_{m}$: 

\begin{itemize}

\item The differential equation  the operator $\omega_{m}$ satisfies. 

\item The multiplicativity property of the operator $\omega_{m}$.

\end{itemize} 

Namely, let ${\cal A}^{\ast, \ast}(M)$ be the Dolbeaut bicomplex of a complex manifold $M$. 
Then the two operators fit into the following commutative diagram, where 
$\int$ is the integration over the line segment connecting the points $(1,0)$ and $(0,1)$ in $\C^2$:

\begin{displaymath}
    \xymatrix{
      && {\cal A}^{\ast, \ast}(M) \otimes \Omega^{\leq 1}_{\C^2}  \ar[d]_{\int}  \\
                                 & S^m({\cal A}^{\ast, \ast}(M)[-1]) \ar[ur]^{\widehat \omega_m} \ar[r]^{\omega_m}&
{\cal A}^{\ast, \ast}(M)[-1]\\
                                 }
\end{displaymath}

The map $\widehat \omega_m$ is a graded algebra map. 
It is easy to prove a differential equation it satisfies. 

The integration over a segment 
almost commutes with the the differential: there are additional boundary terms. 
This implies the differential equation for the operator  $\omega_m$.

\subsection{The form $\omega_{m}(\varphi_1, ..., 
\varphi_{m})$ and its properties} \la{hc2.1.sec}

Given a function $F(f_1, \ldots , f_m)$ of certain 
variables  $f_i$ of degrees $|f_i|$, set  
\begin{equation} \label{5.25.06.1}
{\rm Sym}_{m}F(f_1, \ldots , f_m):= 
\sum_{\sigma \in \Sigma_m}{\rm sgn}_{\sigma; f_1, ..., f_m}F(f_{\sigma(1)}, \ldots , f_{\sigma(m)}), \qquad 
\end{equation}
where the sign of the transposition of $f_1, f_2$ is $ (-1)^{(|f_1|+1)(|f_2|+1)}$,
 and the sign of a permutation 
written as a product of transpositions is 
the product of the signs of the transpositions. 
It is the symmetrisation of elements $\overline f_i$ 
obtained by shifting the degree of $f_i$ by one: ${\rm deg}(\overline f_i) = {\rm deg}(f_i)+1$.

\bd Let $\varphi_0$, ..., 
$\varphi_{m}$ be any $m+1$ smooth forms on a  complex manifold $M$. Set
\begin{equation} \label{12.29.04.1}
\omega_{m+1}(\varphi_0, ..., \varphi_{m}) := 
\frac{1}{(m+1)!}{\rm Sym}_{m+1}
\Bigl( \sum_{k=0}^{m}(-1)^{k} \varphi_0 \wedge \partial' \varphi_1 \wedge ... \wedge \partial' \varphi_k \wedge 
\partial'' \varphi_{k+1} \wedge ... \wedge \partial'' \varphi_{m}\Bigr).
\end{equation}
\ed

Here we write the deRham differential $d$ on  $M$ as a sum of its holomorphic and antiholomorphic components 
$d=\partial'+\partial''$. We apply the symmetrisation (\ref{5.25.06.1}) to a function of  
$\varphi_0, ..., \varphi_{m}$. 

Let ${\cal A}^k_M$ be the space of smooth $k$-forms on $M$. We use the notation $\circ$ for the 
symmetric product. 
The map (\ref{12.29.04.1}) gives rise to a degree zero linear map 
$$
\omega_{m+1}: S^{m+1} ({\cal A}^{\bullet}_M[-1]) \lra {\cal A}^{\bullet}_M[-1], 
\qquad 
\omega_{m+1}: \varphi_0 \circ ... \circ
\varphi_{m} \lms \omega_{m+1}(\varphi_0, ... , \varphi_{m}).
$$

\bex {\em One has} $\omega_1(\varphi_0) = \varphi_0$, 
$$
\omega_2(\varphi_0, \varphi_{1}) = 
\frac{1}{2}\Bigl(\varphi_0\wedge \partial'' \varphi_1  - \varphi_0\wedge \partial'\varphi_1  + 
(-1)^{(|\varphi_0|+1)( |\varphi_1|+1)}
(\varphi_1\wedge \partial'' \varphi_0 -
\varphi_1\wedge \partial'\varphi_0)   \Bigr).
$$
\eex 
The first key property of this map 
is the differential equation it satisfies:

\begin{lemma} \label{5-20.1zxc}
\begin{equation} \label{5-20.1}
\begin{split}
&d \omega_{m+1}(\varphi_0, \ldots, \varphi_{m}) \quad = \quad \frac{1}{m!}{\rm Sym}_{m+1}\Bigl((-1)^{|\varphi_0|}\partial'' \partial' \varphi_{0} \wedge 
\omega_{m}(\varphi_1, \ldots, \varphi_{m})\Bigr) +\\
&
(-1)^m \partial'\varphi_0 \wedge \ldots \wedge \partial'\varphi_{m} + 
\partial'' \varphi_0 \wedge \ldots \wedge \partial''\varphi_{m}.\\
\end{split}
\end{equation}
\end{lemma}
(The first term can be written as 
$
\sum_{j=0}^{m} (-1)^j \pm \partial'' \partial'\varphi_{j} \wedge 
\omega_{m}(\varphi_0, \ldots, \widehat \varphi_j, \ldots, \varphi_{m})
$).
\vskip 3mm

Below we give a direct proof of Lemma \ref{5-20.1zxc}. A conceptual proof, 
which requires almost no calculations, is presented in Section \ref{hc2.2sec}. 
The reader can skip the direct proof. 

\begin{proof} There are the following three identities, which allow to move 
differential operators under the alternation sign, without moving the forms $\varphi_i$:
\begin{equation} \label{5.25.06.23}
{\rm Sym}_{2}(\partial'' \varphi_1 \wedge\partial'\varphi_2) = 
{\rm Sym}_{2}(\partial'\varphi_1 \wedge\partial'' \varphi_2). 
\end{equation}
\begin{equation} \label{5.25.06.21}
 {\rm Sym}_{2}\Bigl((-1)^{|\varphi_{1}|}\varphi_{1}\wedge \partial'\partial'' \varphi_{2}\Bigr) = 
-{\rm Sym}_{2}\Bigl((-1)^{|\varphi_{1}|}\partial'\partial'' \varphi_{1}\wedge \varphi_{2}\Bigr). 
\end{equation}
\begin{equation} \label{5.25.06.24}
 {\rm Sym}_{2}\Bigl(\partial'\varphi_{1}\wedge \partial'' \partial'\varphi_{2}\Bigr) = 
{\rm Sym}_{2}\Bigl((-1)^{|\varphi_{1}|+1}\partial'' \partial'\varphi_{1}\wedge 
\partial'\varphi_{2}\Bigr). 
\end{equation}
Indeed, identity (\ref{5.25.06.23}) is proved as follows: 
$$
{\rm Sym}_{2}(\partial'' \varphi_1 \wedge\partial'\varphi_2) = 
(-1)^{(|\varphi_1|+1)(|\varphi_2|+1)}{\rm Sym}_{2}(\partial'' \varphi_2 \wedge\partial'\varphi_1)= 
{\rm Sym}_{2}(\partial'\varphi_1 \wedge\partial'' \varphi_2). 
$$
Identity (\ref{5.25.06.21}) is proved by 
$$
{\rm Sym}_{2}\Bigl((-1)^{|\varphi_{1}|}\varphi_{1}\wedge \partial'\partial'' \varphi_{2}\Bigr) = 
{\rm Sym}_{2}\Bigl((-1)^{|\varphi_{1}|+ |\varphi_1||\varphi_2| }
\partial'\partial'' \varphi_{2}\wedge \varphi_{1}\Bigr) 
$$
$$
= {\rm Sym}_{2}\Bigl((-1)^{|\varphi_{2}|+ |\varphi_2||\varphi_1| + 
(|\varphi_{1}|+1)(|\varphi_{2}|+1)}\partial'\partial'' \varphi_{1}\wedge \varphi_{2}\Bigr) = 
-{\rm Sym}_{2}\Bigl((-1)^{|\varphi_{1}|}\partial'\partial'' \varphi_{1}\wedge \varphi_{2}\Bigr). 
$$
The third identity is checked in a similar way. These formulas are easily generalised to ${\rm Sym}_{m}$. 

We begin by computing $d \omega_{m+1}(\varphi_0, \ldots, \varphi_{m})$ as 
\be
\begin{split}
&d \omega_{m+1}(\varphi_0, \ldots, \varphi_{m}) = \\
&\frac{1}{(m+1)!}{\rm Sym}_{m+1}
\Bigl( \sum_{k=0}^{m}(-1)^{k} \cdot (\partial'\varphi_0 + \partial''\varphi_0)\wedge 
\partial'\varphi_1 \wedge \ldots \wedge \partial'\varphi_k \wedge 
\partial'' \varphi_{k+1} \wedge \ldots \wedge \partial'' \varphi_{m}+ \\
&\sum_{k=1}^{m}(-1)^{k+|\varphi_0|} k \cdot  \varphi_0 \wedge \partial'' \partial'\varphi_1 \wedge \partial'\varphi_2 
\wedge \ldots \wedge \partial'\varphi_k \wedge 
\partial'' \varphi_{k+1} \wedge \ldots \wedge \partial'' \varphi_{m}+ \\
&\sum_{k=0}^{m-1} (-1)^{|\varphi_0|+ \ldots + |\varphi_k|} (m-k) \cdot  \varphi_0 \wedge \partial'\varphi_1 \wedge \ldots \wedge \partial'\varphi_k \wedge 
\partial'\partial'' \varphi_{k+1} \wedge\partial'' \varphi_{k+2} \wedge \ldots \wedge \partial'' \varphi_{m} \Bigr).\\
\end{split}
\ee
Using (\ref{5.25.06.23}), the first term contributes to the first term in the right hand side of (\ref{5-20.1}).

Using (\ref{5.25.06.24}), the third term is written as 
$$
\frac{1}{(m+1)!}{\rm Sym}_{m+1}\Bigl(\sum_{k=0}^{m-1} (m-k) (-1)^{k+ |\varphi_0|} \cdot  
\varphi_0 \wedge \partial'\partial'' \varphi_{1} \wedge 
\partial'\varphi_2 \wedge \ldots \wedge \partial'\varphi_{k+1} 
\wedge\partial'' \varphi_{k+2} \wedge \ldots \wedge \partial'' \varphi_{m} \Bigr).
$$
Using (\ref{5.25.06.21}) and $\partial'\partial'' = - \partial'' \partial'$,  
the  last two terms in the expression for $d \omega_{m+1}(\varphi_0, \ldots, \varphi_{m}) $ give
$$
\frac{1}{(m+1)!}{\rm Sym}_{m+1}
\Bigl( \sum_{k=1}^{m} (-1)^{k-1+|\varphi_0|} k \cdot  \partial'' \partial'
\varphi_0 \wedge   \varphi_1 \wedge \partial'\varphi_2 
\wedge \ldots \wedge \partial'\varphi_k \wedge 
\partial'' \varphi_{k+1} \wedge \ldots \wedge \partial'' \varphi_{m}
$$
$$
+ \sum_{k=0}^{m-1} (-1)^{k+|\varphi_0|} (m-k) \cdot  \partial'' \partial'
\varphi_0 \wedge \varphi_1 \wedge \partial'\varphi_2 \wedge \ldots \wedge 
 \partial'\varphi_{k+1} \wedge\partial'' \varphi_{k+2} 
\wedge \ldots \wedge \partial'' \varphi_{m}\Bigr).
$$
Changing the summation in the second term from $0 \leq k\leq m-1$ to $1 \leq k\leq m$, and using 
$k+(m+1-k) =m+1$, we see that it coincides with the term containing Laplacians in (\ref{5-20.1}). 
\end{proof}

{\bf Remark}. Let $f_i$ be rational functions on 
a complex algebraic variety. Let $\varphi_i := \log|f_i|$. Then the form 
(\ref{12.29.04.1}) is a part of a cocycle 
representing the product in real Deligne cohomology of
1-cocycles $(\log|f_i|, d\log f_i)$ representing classes in $H^1_{\cal D}(X, \R(1))$. 
The form  (\ref{12.29.04.1}) in this set up was used for the definition of 
Chow polylogarithms  and, more generally, in the 
construction of the canonical regulator map on motivic complexes \cite{G2}, \cite{G8}.

\subsection{Twistor definition of the operator $\omega_{m}$. } \la{hc2.2sec}

\bd \la{12.28.15.1} A cohomological Dolbeaut complex $({\cal A}^{\ast, \ast}, \partial', \partial'')$ is given by  the following data:

\begin{itemize}

\item A 
vector space ${\cal A}^{\ast, \ast}$ 
with a   {\em Hodge bigrading}.  

\item A {\em cohomological grading} and a differential 
 ${\bf d}$ of cohomological  degree $1$,  written as $$
{\bf d}=\partial' +\partial'',
$$ where 
$\partial'$  and 
 $\partial''$ have the cohomological degree $1$, and the Hodge bidegrees  $(1,0)$ and 
  $(0,1)$. 

\item   $({\cal A}^{\ast, \ast}, \partial',  \partial'')$
satisfies the 
$\partial'\partial''$-lemma.
\end{itemize}
 \ed

\noindent Decomposing  ${\bf d}^2=0$ into the components of Hodge bidegrees 
$(2,0)$, $(1,1)$ and $(0,2)$ we get 
\be \la{5.2.20.1}
\partial'^2 = \partial''^2= \partial' \partial''+ \partial''\partial'=0.
\ee
We set 
$$
{\bf d}^\C:= \partial' - \partial''. 
$$

An example of a cohomological Dolbeaut bicomplex is given by the Dolbeaut bicomplex of a complex manifold.
A more general example is given by the Dolbeaut bicomplex of a variation of 
Hodge structures, discussed in Section \ref{hc3sec}.

Although the $\partial'\partial''$-lemma is essential, we do not use it in Section 2.

\paragraph{The twistor polydifferential operator $\widehat \omega_m$.} 
Let $\C^2$ be the twistor plane with canonical coordinates $(z,w)$. 
The 
 group $\C^*\times \C^*$ acts on the twistor plane and 
on the bicomplex ${\cal A}^{\ast, \ast}$ by  
$$
(\lambda, \mu): (z, w)  \lms (\lambda z, \mu w), ~~(\lambda, \mu): f_{s,t}\in {\cal A}^{s, t} \lms \lambda^{s} \mu^{t}f_{s,t}. 
$$
So the Hodge bidegree can be defined as  the character of the action of the group 
$\C^*\times \C^*$. 

Take the complex vector space $\langle \partial', \partial''\rangle_\C$ spanned by the 
operators $\partial'$ and $\partial''$. Then the  map 
$$
\langle \partial', \partial''\rangle_\C \otimes {\cal A}^{\ast, \ast} \lra {\cal A}^{\ast, \ast}, ~~~~\partial'\otimes f \lms 
\partial'f, ~~\partial''\otimes f \lms 
\partial''f
$$
is  $\C^*\times \C^*$-equivariant for the following action of the group $\C^*\times \C^*$ on the operators $\partial'$, $\partial''$: 
$$
(\lambda, \mu): (\partial', \partial'') \lms (\lambda \partial', \mu \partial''). 
$$

The  space 
$\langle \partial', \partial''\rangle_\C$ is identified canonically, as $\C^*\times \C^*$-module, 
with the space of linear functions $(z,w)$ on 
the twistor plane,  
``explaning''  the role of the latter.

Consider the tensor product of the complex ${\cal A}^{\ast, \ast}$ with the 
algebraic de Rham complex of  $\C^2$;
$$
{\cal A}^{\ast, \ast}\otimes_\C \Omega^{\bullet}_{\C^2}, ~~{\bf d} + {\rm d}_{\C^2}. 
$$
Let us introduce the twistor version of the  operator ${\bf d}^\C$: 
\be \la{DCF}
{\rm D}^\C:= w\partial'-z\partial''+ (zdw-wdz).
\ee
It has cohomological degree $+1$ and Hodge bidegree $(1,1)$. 

Denote by $\C(-1)[-1]$ a one dimensional vector space 
of Hodge bidegree $(1,1)$ and cohomological degree $+1$. 
So ${\cal A}^{\ast, \ast}(-1)[-1]$ is ${\cal A}^{\ast, \ast}$ with each of its three degrees 
 shifted by one.

\bd The twistor operator $\widehat \omega_{m}$  
 is the following map: 
\be \la{4.13.14.10}
\begin{split} 
&\widehat \omega_{m}:~ S^{m}
\Bigl({\cal A}^{\ast, \ast}(-1)[-1]\Bigr)\lra 
{\cal A}^{\ast, \ast}\otimes \Omega^{\leq 1}_{\C^2};\\
&\widehat \omega_{m}(\varphi_1, \ldots , \varphi_m;z,w):= {\rm D}^\C \varphi_1 \wedge \ldots \wedge {\rm D}^\C \varphi_m. \\
\end{split}
\ee
\ed

\bl The twistor operator $\widehat \omega_{m}$ in (\ref{4.13.14.10})
  preserves each of the three degrees. 

It 
is a graded algebra homomorphism: 
$$
\widehat \omega_{m+n}(\varphi_1, \ldots , \varphi_{m+n};z,w) = \widehat \omega_{m}(\varphi_1, \ldots , \varphi_{m};z,w) 
\wedge \widehat \omega_{n}(\varphi_{m+1}, \ldots , \varphi_{m+n};z,w). 
$$
\el

\begin{proof}
Since the  ${\rm D}^\C$ has cohomological degree $+1$ and Hodge bidegree $(1,1)$, 
the operator $\widehat \omega_{m}$ preserves each of the three degrees. The second claim is clear. 
\end{proof}

The next question is how the differential acts on the map $\widehat \omega_{m}$.  
The answer is most naturally given by using the restriction to the twistor line.

\paragraph{Restriction to the twistor line and the differential.} Let us consider the twistor line: 
$$
\C:= \{z,w \in \C^2~|~ z+w=1\}. 
$$
It is parametrised  by 
$$
z=1-u,~~w=u. 
$$
We denote by $(\Omega^{\bullet}_\C, {\rm d}_{ u})$ the algebraic de Rham complex of the twistor line. Observe that 
$$
{\rm Res}_{z+w=1}~(zdw-wdz) = du. 
$$
Denote by $\widehat \omega_{m}(\varphi_1, \ldots , \varphi_m;u)$ the restriction 
of $\widehat \omega_{m}(\varphi_1, \ldots , \varphi_m; z,w)$ to the twistor line. 
Let ${\bf d} + {\rm d}_{ u}$ be the total differential on ${\cal A}^{\ast, \ast}\otimes \Omega^{\bullet}_\C$. 
\bl One has 
\be \la{5.22.15.2}
\begin{split}
&({\bf d} + {\rm d}_{ u})\widehat \omega_{m}(\varphi_1, \ldots , \varphi_m;u)= 
 \partial'' \partial' \varphi_1\wedge 
{\rm D}^\C \varphi_2 \wedge \ldots \wedge {\rm D}^\C \varphi_m + \ldots \\
&= \sum_{k=1}^m\pm\partial'  \partial'' \varphi_k\wedge \widehat \omega_{m-1}(\varphi_1, \ldots , \widehat \varphi_k, 
\ldots , \varphi_m;u). 
\end{split}
\ee
Here the signs are calculated in the standard way using 
$|{\rm D}^\C \varphi_i|=|\varphi_i|+1$. 
\el
  
\begin{proof} Denote by ${\rm D}^\C_u$ the differential ${\rm D}^\C$ restricted to the twistor line 
${\C}$. 

Then 
one has the following simple but fundamental equality: 
\be \la{SFE}
\begin{split}
&({\bf d} + {\rm d}_{ u}){\rm D}^\C_u =  ({\bf d} + {\rm d}_{u})\Bigl(u\partial' - (1-u)
 \partial'' + du\Bigr)= \\
&du (\partial' + \partial'') + u\partial''\partial' - (1-u)\partial' \partial'' -
 du (\partial' +  \partial'')) =  \partial'' \partial'.\\
\end{split}
\ee\end{proof}

\paragraph{A conceptual proof of Lemma \ref{5-20.1zxc}.} 
First, we 
 relate the operators $\widehat \omega_{m+1}$ and $\omega_{m+1}$. 

Let $\Delta$ be the oriented line segment in the twistor plane with the ends $(1,0), (0,1)$:
\be \la{SDELTA}
\Delta = \{(1-u, u)~|~ 0 \leq u \leq 1\} \subset \C^2.
\ee

\bl \la{5.5.14.1sa} The integral of the $\widehat \omega_{m+1}$ over the 
segment $\Delta$ is the operator $(-1)^m\omega_{m+1}$:
\be \la{stateL}
\int_{\Delta}\widehat \omega_{m+1}(\varphi_0, ..., \varphi_{m}; z,w) = (-1)^m\omega_{m+1}(\varphi_0, ..., \varphi_{m}).
\ee
\el

\begin{proof}  A typical term of $\widehat 
\omega_{m+1}(\varphi_0, ..., \varphi_{m}; z,w)$ containing a differential $zdw-wdz$ is given by 
\be \la{typtxia}
(-1)^{m-k} {(zdw-wdz)} \wedge \varphi_0 
 \wedge {w}\partial'\varphi_1  
\wedge \ldots \wedge {w}\partial'\varphi_k \wedge 
{z}\partial'' \varphi_{k+1}  
\wedge \ldots \wedge {z}\partial'' \varphi_m .
\ee 

Let us restrict it to the twistor line and integrate over $0 \leq u \leq 1$. Observe that 
\be \la{2.19.15.1}
\int^{1}_{0}(1-u)^{s}u^{t}
du  = \frac{\Gamma(s+1)\Gamma(t+1)}{\Gamma(s+t+2)} = 
\frac{1}{s+t+1}{s+t \choose s}^{-1}.
\ee
Therefore after the integration  we get 
$$
(-1)^{m-k}\frac{k!(m-k)!}{(m+1)!}\varphi_0 
 \wedge \partial'\varphi_1  
\wedge \ldots \wedge \partial'\varphi_k \wedge 
\partial'' \varphi_{k+1}  
\wedge \ldots \wedge \partial'' \varphi_m. 
$$
Comparing with (\ref{12.29.04.1}), we get the right hand side of (\ref{stateL}). 
\end{proof}

To prove Lemma \ref{5-20.1zxc} notice that we integrate in Lemma \ref{5.5.14.1sa} along the segment 
$\Delta$ on the twistor plane $\C^2$, which has boundary. So 
the integration operator 
does not quite commute with the de Rham differential  on the twistor plane: one gets 
additional boundary terms 
provided by the  ends. The boundary terms match the two terms in (\ref{5-20.1}), 
corresponding to the evaluation of the monomials $z^sw^t$ at $z=1$ when $(s,t) = (m+1, 0)$, and 
at $w=1$ when $(s,t) = (0, m+1)$.

\paragraph{Decomposition of $\widehat \omega_m$.} 
Since $\Omega^{\leq 1}_{\C^2}=\Omega^0_{\C^2} \oplus \Omega^1_{\C^2}$, the map $\widehat \omega_m$ 
 has two components: 
\be 
\begin{split}
&\widehat \omega_{m} = \widehat \eta_{m} + \widehat \xi_{m}. \\
&\widehat \omega^0_{m}: S^{m}\Bigl({\cal A}^{\ast, \ast}[-1]\Bigr) \lra 
{\cal A}^{\ast, \ast}\otimes \Omega^0_{\C^2}.\\
&\widehat \omega^1_{m}: S^{m}\Bigl({\cal A}^{\ast, \ast}[-1]\Bigr) \lra 
{\cal A}^{\ast, \ast} \otimes \Omega^1_{\C^2},
\end{split}
\ee
 Explicitly, we have, using the notation $\ldots$ for the terms obtained by the alternation: 
$$
\widehat \omega^1_{m+1}(\varphi_0, \ldots, \varphi_{m}; {z,w}) := 
$$
\be \la{typtxi}
\sum_{k=0}^{m} (-1)^{m-k} \Bigl({(zdw-wdz)}\wedge \varphi_0 
 \wedge {w}\partial'\varphi_1  
\wedge \ldots \wedge {w}\partial'\varphi_k \wedge 
{z}\partial'' \varphi_{k+1}  
\wedge \ldots \wedge {z}\partial'' \varphi_m + ... \Bigr).
\ee

$$
\widehat \omega^0_{m+1}(\varphi_0, \ldots, \varphi_{m}; {z,w}) := \sum_{k=0}^{m} (-1)^{m-k} \Bigl(
{w}\partial'\varphi_0  
\wedge \ldots \wedge{w}\partial'\varphi_k \wedge 
{z}\partial'' \varphi_{k+1}  
\wedge \ldots \wedge {z}\partial'' \varphi_m + ... \Bigr).
$$

\section{Hodge correlators for curves} \la{formersec2}

\subsection{Green functions on a Riemann surface} \la{2.2}
Let $X$ be a smooth compact  Riemann surface of genus $g$. 
There is a canonical hermitian structure on the space $\Omega^1(X)$ of 
holomorphic differentials on $X$ given by 
\begin{equation} \label{1.8.05.6}
<\alpha_{k}, \alpha_{l}>:= \frac{i}{2}\int_{X(\C)} \alpha_{k}
 \wedge \overline \alpha_{l}. 
\end{equation}
Let $\alpha_1, .., \alpha_g$ be an 
orthonormal basis with respect to this form. 

Recall that a $p$-current on a manifold is a differential $p$-form with generalised coefficients. 
Let $\mu$  be a $2$-current on $X$ of volume one: 
$\int_{X}\mu =1$.  The  corresponding 
Green current  $G_{\mu}(x,y)$ is a $0$-current on 
$X \times X$ satisfying a differential equation 
\begin{equation} \label{7.3.00.1}
\overline \partial \partial G_{\mu}(x,y) = \delta_{\Delta_X} - 
\Bigl( p_1^*\mu + p_2^*\mu  - \frac{i}{2}\sum_{k=1}^g 
 (p_1^*\alpha_{k} \wedge p_2^*\overline \alpha_{k} + 
p_2^*\alpha_{k} \wedge p_1^*\overline \alpha_{k} )\Bigr).
\end{equation}

\paragraph{Remark.} Our Green function is $(2\pi i)^{-1}$ times the 
traditional one, defined  using $(2\pi i)^{-1}\overline \partial \partial$ in 
(\ref{7.3.00.1}). For example, the Green function 
for the projective line in our normalisation is 
$$
G(x,y) = (2\pi i)^{-1}\log |x-y|. 
$$
 \vskip 2mm

There exists a unique up to a constant  solution of this equation. 
Indeed, the kernel of the 
Laplacian $\overline \partial \partial $ on $X\times X$ 
consists of constants. 
One has a symmetry relation:  
$$
G_{\mu}(x,y) = G_{\mu}(y, x).
$$ A point $a \in X$ provides  a volume one $2$-current $\delta_a$ on $X$. 
So letting 
$\mu:= \delta_a$, we get a Green function 
denoted $G_{a}(x,y)$. 
Then the 
current on the right hand side of (\ref{7.3.00.1}) 
is smooth 
outside of the divisor $\Delta \cup (\{a\}\times X) \cup (X \times \{a\})$. 
Since the Laplacian is  an elliptic operator, 
the restriction of the 
Green current  to the complement 
of this divisor is represented by a smooth function, 
called the Green function, also 
 denoted by 
$G_{a}(x,y)$. The  Green function $G_{a}(x,y)$ near the 
singularity divisor  looks like $\log r$, where 
$r$ is a distance to the divisor, plus a 
smooth function.  So it has an integrable 
singularity, 
and thus provides 
a current on $X \times X$, which coincides with the Green current, see \cite{La}.

\paragraph{The Arakelov Green function.} 
A classical choice of $\mu$ is given by a smooth volume form 
of total mass $1$ on $X$. In particular a 
metric on $X$ provides such a volume form. 
There is a canonical {\it Arakelov volume form} on $X$:
$$
{\rm vol}_X:= \frac{i}{2g}\sum_{a=1}^g 
\omega_{a} \wedge \overline \omega_{a}, 
\qquad \int_{X}{\rm vol}_X = 1.
$$ 
The corresponding Green function is 
denoted by $G_{\rm Ar}(x,y)$. 
It can be  normalized so that 
\begin{equation} \label{1.8.05.1}
p_{2*}(G_{\rm Ar}(x,y)p_1^*{\rm vol}_X) = p_{1*}(G_{\rm Ar}(x,y)p_2^*{\rm vol}_X) =
0.
\end{equation}

The Green function $G_a(x,y)$ 
is expressed via the Arakelov Green function: 
\begin{equation} \label{7.3.00.1ds}
G_a(x,y) = G_{\rm Ar}(x,y) - G_{\rm Ar}(a,y) - G_{\rm Ar}(x,a) +C. 
\end{equation}
Indeed, the right hand side satisfies differential equation 
(\ref{7.3.00.1}).

\paragraph{A normalization of $G_{a}(x,y)$.} 
Let us choose a non-zero tangent vector $v$ at $a$.  Let $t$ be a 
local parameter at $a$ such that $dt(v)=1$. Then there exists a unique 
 solution $G_{v}(x,y)$ 
of (\ref{7.3.00.1}) with $\mu = \delta_a$  
such that $G_{v}(x,y) - (2\pi i)^{-1}\log|t|$ vanishes at $x=a$. 
Indeed, there is a unique normalized Arakelov Green function 
$G_{{\rm Ar}, v}(x,y)$ such that $G_{{\rm Ar}, v}(x,a)- (2\pi i)^{-1}\log|t|$ vanishes at $x=a$. 
Thanks to the symmetry, the same is true for $G_{\rm Ar}(a,y)- (2\pi i)^{-1}\log|t|$ 
which vanishes at $y=a$. It remains to use (\ref{7.3.00.1ds}) with $C=0$ there.

\paragraph{Specialization.}  
Let $F$ be a function defined in a punctured neighborhood 
of a point $a$ on a complex curve $X$.
Let $v$  be a tangent vector 
at $a$. Choose a local parameter 
$t$ in the neighborhood of $a$ such that $dt(v) =1$. 
Suppose that the function $F(t)$ admits an asymptotic expansion 
$$
F(t) = \sum_{n \geq 0}
F_{n}(t)\log^n|t|
$$
where $F_n(t)$ are smooth functions near $t=0$. 
We define the specialization 
${\rm Sp}^{t\to 0}_{v}F(t):= F_0(0)$. 
We skip sometimes the superscript $t\to 0$. 
It follows immediately from the definitions that
\be \la{1.11.08.2}
{\rm Sp}^{x\to a}_{v}G_{v}(x,y) = 
{\rm Sp}^{y\to a}_{v}G_{v}(x,y) = 0. 
\ee

We extend  the Green function by linearity to the group of 
divisors on $X$ by setting
  $G_a(\sum n_ix_i, \sum m_jy_j) = \sum_{i,j}n_im_jG_a(x_i, y_j)$. 

\begin{lemma} \label{3.18.05.5} If $D_0$ is a degree zero divisor, then 
$G_{\mu}(D_0, y) - G_{\mu'}(D_0, y)$ does not depend on $y$. 
\end{lemma}

\begin{proof} One has 
$
\overline \partial_y \partial_y 
\Bigl(G_{\mu}(D_0, y) - G_{\mu'}(D_0, y)\Bigr) 
= {\rm deg}(D_0)p_2^*(\mu-\mu') =0.  
$ 
\end{proof}

\subsection{Construction of Hodge correlators for curves} \la{sec2.2x}

\paragraph{Decorated plane trivalent trees.} 
Let us recall some terminology.  
A {\it tree} is a connected graph without loops.  
The {\it external vertices} of a tree are the ones of valency $1$. 
The rest of the vertices are {\it internal vertices}. 
The edges of a tree consist of {\it external edges}, i.e. the ones containing the external vertices, 
and {\it internal edges}.  
A {\it trivalent tree} is a tree whose internal vertices are of valency $3$. 
We allow a trivalent tree without internal vertices. It is just one edge 
with two (external) vertices. A trivalent tree with $m+1$ external vertices has $2m-1$ edges. 
A {\it plane 
tree} is a tree 
without self intersections  
located on the plane.

\begin{figure}[ht]
\centerline{\epsfbox{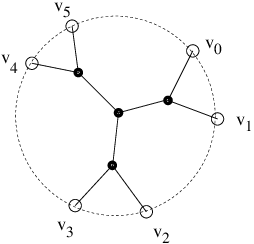}}
\caption{A decorated plane  trivalent tree.}
\label{feyn9}
\end{figure}

\begin{definition} \label{1.7.05.1}
A  {\it decoration} of a tree $T$ by elements of a set $R$ is a map 
$\{\mbox{external vertices of $T$}\}\to R$. 
\end{definition}
An example of  a decorated plane trivalent tree is given on Figure \ref{feyn9}. 

\vskip 3mm
Recall (Section 1.3) that $S^* \subset X$ and 
$$
{\rm V}^{\vee}_{X, S^*} = \C[S^*] \oplus (\Omega_X^1 \oplus \overline \Omega_X^1).
$$ 
Further,  ${\cal C}^{\vee}_{X,S^*}$ 
is the cyclic envelope of the tensor algebra of  ${\rm V}^{\vee}_{X, S^*}$. 
Its elements  are cyclic words in ${\rm V}^{\vee}_{X,S^*}$: for instance, 
one has ${\cal C}(abc) = 
{\cal C}(cab)$.

A $2$-current $\mu$ is {\it admissible} if it is 
either a smooth volume form 
of total mass $1$ on $X$, 
or the $\delta$-current $\delta_a$ for some point $a \in X$.  
Given an admissible $2$-current $\mu$,  we are going to define a linear map, 
the {\it  Hodge correlator map}:
$$
{\rm Cor}_{{\cal H}, \mu}: {\cal C}^{\vee}_{X,S^*}  \lra \C.
$$
To define it, take a degree $m+1$ 
decomposable cyclic element  
$$
W = {\cal C}\Bigl(v_0 \otimes v_1 \otimes \ldots \otimes v_m\Bigr), 
\quad v_i \in {\rm V}^{\vee}_{X,S^*}.
$$ 
The external vertices of a plane tree have a cyclic order 
provided by a chosen 
(say, clockwise) orientation of the plane. 
 We say that a plane trivalent tree  $T$ is  decorated by $W$ if the external 
vertices of the tree are decorated by the elements $v_0, \ldots, v_m$, 
respecting the cyclic order of the vertices, as on 
Figure \ref{feyn9}. Although 
the decoration depends on the presentation of $W$ as a cyclic tensor product of the vectors $v_i$, 
constructions below depend only on $W$ and not on such presentation.
We are going to assign to $W$ a  
$2m$-current $\kappa_W$ on 
\begin{equation} \label{7.1.00.1q} 
X^{\mbox{\{vertices of $T$\}}}. 
\end{equation}
We assume that each 
of the vectors $v_i$ in $W$ is either a 
generator $\{s_i\}\in \C[S^*]$, or a $1$-form 
$\omega$ on $X$. In Section 2.2 we assume 
in addition that if $T$ has just one edge, 
both its vertices are $S$-decorated. 
Then 
a  decoration of $T$ provides 
a decomposition of the set of external edges into  
the $S^*$-{\it decorated} and {\it special} edges: 
the latter are decorated by $1$-forms, the former by points of $S^*$.

Denote by  $\delta_X$ the $\delta$-current of the diagonal in $X\times X$, and 
by  $\delta_X\omega$ its product by a smooth form $p^*\omega$, where $\omega$ is a form on $X$ 
and $p: X\times X\to X$ is the projection to one of the factors. 
Given an edge $E$ of $T$, we assign to it a current  
$\widetilde G_E$ on $X\times X$:  
\begin{equation} \label{11.48.04.11df}
\widetilde G_E:= 
\left\{ \begin{array}{ll}
\delta_X \omega& \quad \mbox{if $E$ is a special external edge,} \\
G_\mu(x_1, x_2)
& \quad \mbox{otherwise.} 
\end{array} \right. 
\end{equation}

For each edge $E$ of $T$, there is a projection 
of the set of vertices of $T$ onto the set of vertices of $E$. It gives rise to 
a canonical projection 
$$
p_E: X^{\mbox{\{vertices of $T$\}}} \quad \lra \quad 
X^{\mbox{\{vertices of E\}}} = X\times X.
$$
We set 
$$
G_E:= p_E^*\widetilde G_E.
$$
Since the map $p_E$ is transversal to the wave front of the current 
$\widetilde G_E$, the inverse image is well defined. 

Given a finite set ${\cal X}$, there is a $\Z/2\Z$-torsor ${\rm or}_{\cal X}$, 
called the {\it orientation torsor of 
${\cal X}$}. Its elements are expressions $x_1 \wedge \ldots \wedge x_{|{\cal X}|}$, 
where $\{x_1, \ldots, x_{|{\cal X}|}\} = {\cal X}$; 
interchanging two neighbors we change the sign of the expression. 
The {\it orientation torsor ${\rm or}_T$} of a graph $T$ is 
the orientation torsor of the set of edges of $T$. 
An orientation of the plane induces an orientation of a  plane trivalent tree:  
Indeed, the orientation of 
the plane provides  orientations of links of each of the vertices. 

Let us choose an element 
$$
(E_0 \wedge \ldots \wedge E_r) \wedge (E_{r+1} \wedge \ldots \wedge E_{2m}) 
\in {\rm or}_T
$$
of the orientation torsor of our tree $T$, 
such that  the edges $E_0, \ldots, E_r$ are internal or $S$-decorated, 
and the others are decorated by $1$-forms. 
Set\begin{equation} \label{3.8.05.15}
\kappa_T(W):= {\rm sgn}(E_0\wedge \ldots \wedge E_{2m})\omega_r(G_{E_0}\wedge \ldots \wedge 
G_{E_r})
\wedge G_{E_{r+1}}\wedge \ldots \wedge G_{E_{2m}},
\end{equation} 
where 
$
{\rm sign}(E_0 \wedge \ldots \wedge E_{2m}) \in \{\pm 1\}
$  
is the difference between the 
element $E_0 \wedge \ldots \wedge E_{2m}$ and the 
canonical generator of ${\rm or}_T$. 

\begin{lemma}
$\kappa_T(W)$ is a $2m$-current on (\ref{7.1.00.1q}).
\end{lemma}

\begin{proof} 
Consider  
the smooth differential $r$-form 
\begin{equation} \label{3.7.04.1}
\omega_{r}(G_{E_0} \wedge \ldots \wedge G_{E_{r}})
\end{equation} 
on the complement to the diagonals in (\ref{7.1.00.1q}). We claim that 
it 
has integrable singularities, i.e.  provides an $r$-current on 
(\ref{7.1.00.1q}). Indeed, the Green function is smooth outside of the diagonal.    
Its singularity near the diagonal is of type $\log |t|$, where $t$ is a local equation 
of the diagonal
i.e. $G(x,y) - \log|t|$ is smooth near the diagonal. 
Since $\omega_{m+1}$ is a linear map,  it is sufficient to prove the claim for 
a specific collection of functions with such singularities. 
Since the problem is local, we can choose local equations $f_i=0$ of all the diagonals involved, 
and consider the functions $\log|f_i|$ as such a specific collection. 
Then the claim follows from Theorem 2.4 in \cite{G8}. 
\end{proof}

Consider  the projection map
\begin{equation} \label{7.1.00.1qq}
p_T: X^{\{\mbox{vertices of T}\}} \lra X^{\{\mbox{$S^*$-decorated 
vertices of T}\}}. 
\end{equation}
\begin{definition} Let $\mu$ be an admissible  $2$-current on $X$, 
and  $W \in {\cal C}^{\vee}_{X,S^*}$. 
Then the  correlator ${\rm Cor}_{{\cal H}}(W)$ is a $0$-current on 
the right hand side of (\ref{7.1.00.1qq}), 
given by 
$$
{\rm Cor}_{{\cal H}}(W):= {p_T}_*\Bigl(\sum_{T}\kappa_T(W)\Bigr),
$$
where the sum is over all plane trivalent trees $T$ 
decorated by the cyclic word $W$. 
\end{definition}
The restriction of the current ${\rm Cor}_{{\cal H}}(W)$ to
 the complement of the  diagonals 
is a smooth function.

\paragraph{Remarks.} 1. The direct image ${p_T}_*$ is given by  
integration over 
$ 
X^{\{\mbox{internal vertices of $T$}\}}.
$  
Indeed, $\kappa_T(W)$ is the exterior product of (\ref{3.7.04.1}) 
and a smooth differential form on the product of copies of $X$ 
corresponding to the internal and $S^*$-decorated vertices, 
extended by the $\delta$-function to (\ref{7.1.00.1q}).

2. 
We may assume that the $S$-decorated vertices are decorated by 
divisors supported on $S$ by extending the Hodge correlator 
map by linearity. 

3. For a given measure $\mu$, the corresponding Green function $G_\mu$ 
is defined up to adding a constant. The Hodge 
correlators depend on it. However this dependence is minimal. 
If $W$ has at least one degree zero divisor $D$ as 
a factor, the corresponding Hodge correlator does not depend on 
the choice of the constant. Indeed, we may assume 
that in the form $\omega_{m+1}$ 
we do not differentiate the Green function $G_E$ assigned to 
the edge $E$ decorated by $D$.

\subsection{Shuffle relations for Hodge correlators} \la{sec2.3x}

Our goal is the following shuffle relations for Hodge correlators:

\begin{proposition} \label{7.3.06.4} Let $v_i \in {\rm V}^{\vee}_{X, S^*}$. 
Then for any $p,q \geq 1$ one has 
\begin{equation} \label{7.3.06.5}
\sum_{\sigma \in \Sigma_{p,q}}{\rm Cor}_{\cal H}{\cal C}(v_0 \otimes v_{\sigma(1)} \otimes \ldots \otimes v_{\sigma(p+q)}) = 0, 
\end{equation}
where the sum is over all $(p,q)$-shuffles $\sigma \in \Sigma_{p,q}$. 
\end{proposition}

To prove this we interpret the ``sum over plane trivalent trees'' 
construction 
using cyclic operads. See \cite{LV} for the background on operads. 

\vskip 3mm 
{\it Trees and the cyclic Lie co-operad.}
Let ${\cal C}{\cal A}ss_{\bullet}$ be the cyclic associative operad. 
This means, in particular, the following. 
Let ${\cal A}ss_m \{v_1,...,v_m\}$ 
be the space of all associative words formed by the letters $v_1,...,v_m$, each used once. 
Then ${\cal C}{\cal A}ss_{m+1}$ is a vector space generated by the expressions 
$$
(A, v_0), \quad \mbox{where} \quad A \in {\cal A}ss_m \{v_1,...,v_m\}; \quad 
(xy, z) = (x, yz), \quad (x,y) =(y,x),
$$ 
i.e. $(*,*)$ is an invariant scalar product on an associative  algebra. 
So as a vector space ${\cal C}{\cal A}ss_{m+1}$ is isomorphic to the space 
of $m$-ary operations ${\cal A}ss_{m}$ in the associative operad, or simply speaking to 
${\cal A}ss_m \{v_1,...,v_m\}$. 

Similarly let ${\cal C}{\cal L}ie_{\bullet}$ be the cyclic Lie operad. Let 
 ${\cal L}ie_m\{v_1,...,v_m\}$ be the space of all Lie words formed by the letters $v_1,...,v_m$, each used once. 
 Then  as a vector space ${\cal C}{\cal L}ie_{m+1}$ is isomorphic to the space of 
$m$-ary operations in the Lie operad. We think about it as of 
the space generated by expressions 
$$
(L, v_0), \quad \mbox{where} \quad L \in {\cal L}ie_m \{v_1,...,v_m\}; \quad 
([x,y],z) = (x, [y,z]), \quad (x,y) =(y,x),
$$ 
i.e. $(*,*)$ is an invariant scalar product on a Lie algebra.

\paragraph{Example.} The space ${\cal C}{\cal L}ie_{3}$ is one dimensional. It is generated by $([v_1, v_2], v_0)$. 
The space ${\cal C}{\cal A}ss_{3}$ is two dimensional, generated  
by $(v_1v_2,v_0)$ and $(v_2v_1, v_0)$. 
\vskip 3mm

Let
${\cal C}{\cal A}ss^*_{m+1}$ (resp. ${\cal C}{\cal L}ie^*_{m+1}$) 
be the $\Q$-vector space dual to ${\cal C}{\cal A}ss_{m+1}$ 
(resp. ${\cal C}{\cal L}ie_{m+1}$).  
We define the subspace of {\it shuffle relations} in ${\cal C}{\cal A}ss^*_{m+1}$ 
as follows. Let $(v_0 v_1 ... v_m)^* \in 
{\cal C}{\cal A}ss^*_{m+1}$ be a functional whose value 
on $(v_1 ... v_m, v_0 )$ is $1$ and on $(v_{i_1} ... v_{i_{m}}, v_0 )$ is 
zero if $\{i_1, ..., i_m\} \not = 
\{1,2, ..., m\}$ as ordered sets. 
The  subspace of shuffle relations is generated by the expressions 
$$
\sum_{\sigma \in \Sigma_{k,m-k}} (v_0 v_{\sigma (1)} ... v_{\sigma (m)})^*; 
\qquad 1 \leq k \leq m-1.
$$

\begin{lemma}\label{MB1}
There is a canonical isomorphism
\begin{equation} \label{MB}
{\cal C}{\cal L}ie^*_{m+1} = \frac{{\cal C}{\cal A}ss^*_{m+1}}
{\mbox{\rm Shuffle relations}}.
\end{equation}
\end{lemma}

\begin{proof} Notice that 
${\cal L}ie_m\{v_1,...,v_m\} \subset {\cal A}ss_m\{v_1,...,v_m\}$ is the  
subspace of primitive elements for the 
coproduct $\Delta$, $\Delta(v_i) = v_i \otimes 1 + 1 \otimes v_i$. Thus  
its dual  ${\cal L}ie^*_m\{v_1,...,v_m\}$ is the quotient of 
${\cal A}ss^*_m\{v_1,...,v_m\}$ by the shuffle relations, see \cite{LV}, Theorem 1.3.9. 
\end{proof}

{\it The tree complex (\cite{K})}. Denote by ${\rm T}_{(m)}^i$ the abelian group generated 
by the isomorphism classes of pairs  $(T, {\rm Or}_T)$ where $T$ is a 
tree (not  a plane tree!) with $2m-i$ edges and with 
$m+1$  ends decorated by the set $\{v_0, ..., v_m\}$. 
Here ${\rm Or}_T$ is an orientation of the tree $T$. The only relation is 
that changing the orientation of a tree $T$ amounts to changing the sign of 
the generator. 
There is a differential 
$
d: {\rm T}_{(m)}^i \lra {\rm T}_{(m)}^{i+1}
$ 
defined by shrinking internal edges of a tree: 
$$
(T, {\rm Or}_{T}) \lms \sum_{\mbox{internal edges $E$ of $T$}} (T/E, {\rm Or}_{T/E}).
$$ 
Here if ${\rm Or}_{T} = E \wedge E_1 \wedge ... $, then  
${\rm Or}_{T/E}:= E_1 \wedge ...$ . 

\vskip 3mm
The group  ${\rm T}_{(m)}^i$ is a free abelian group with a basis; the basis vectors 
are defined up to a sign. Thus there is a perfect pairing 
${\rm T}_{(m)}^i \otimes {\rm T}_{(m)}^i \lra \Z$, and we 
may identify the group ${\rm T}_{(m)}^i$ with its dual $({\rm T}_{(m)}^i)^*:= 
{\rm Hom}({\rm T}_{(m)}^i, \Z)$. 
Consider the map $w$
$$
 \Bigl({\rm T}_{(m)}^1\Bigr)^* \quad \stackrel{w}{\lra} \quad 
{\cal C}{\cal L}ie_{m+1}
\quad \stackrel{i}{\hookrightarrow} \quad {\cal C}{\cal A}ss_{m+1}.
$$
defined as follows. Take an element $(T, {\rm Or}_T)$ represented by a 
 $3$-valent tree $T$ decorated by 
$v_0, v_1, ..., v_m$ and an orientation ${\rm Or}_T$ of $T$. 
Make a decorated tree rooted at $v_0$ out of the tree $T$. 
Define a Lie word 
in $v_1, ..., v_m$ using the oriented tree. We are getting an  element $w(T) \subset {\cal L}{\cal L}ie_{m+1}$. 
Viewing it  as an associative word, we arrive at the element $i \circ w (T)$. 
There are the dual maps
\begin{equation} \label{UYAA}
{\cal C}{\cal A}ss^*_{m+1} \stackrel{i^*}{\lra}  
{\cal C}{\cal L}ie^*_{m+1} \quad \stackrel{w^*}{\lra} \quad 
{\rm T}_{(m)}^1.
\end{equation}

The next lemma offers an explanation of the ``sum over plane trivalent trees'' construction. 

\begin{lemma} \label{7.3.06.1}
The element $w^*\circ i^* (v_0 v_1 \ldots v_m)^*$ is the sum of all plane trivalent trees 
decorated by the cyclic word ${\cal C}(v_0 v_1 \ldots v_m)$, and equipped 
with the canonical orientation induced by the clockwise orientation of the plane. 
\end{lemma}
For any 
$p, q \geq 1$ with $p+q=m$ one has 
 $$
w^*\circ i^* \sum_{\sigma \in \Sigma_{p,q}}(v_0 v_{\sigma(1)} \ldots v_{\sigma(m)})^* = 0, 
$$ 
where the sum is over all $(p,q)$-shuffles. Indeed, we apply $w^*$ to the expression which is 
zero by Lemma \ref{MB1}. 
This plus Lemma \ref{7.3.06.1} immediately imply 

\begin{corollary} \label{7.3.06.2}
The element of ${\rm T}^1_{(m)}$ given by the sum over all shuffles $\sigma \in \Sigma_{p,q}$ 
of the sum of all plane trivalent trees decorated by $(v_0 v_{\sigma(1)} \ldots v_{\sigma(m)})$ is zero. 
\end{corollary}

\paragraph{Proof of Proposition \ref{7.3.06.4}.}
Corollary \ref{7.3.06.2} tells that the integral in (\ref{7.3.06.5}) is  over 
zero cycle.

\paragraph{Remark.} It was proved in \cite{GK} that the following  complex is exact:
\begin{equation} \label {GKap1}
{\cal C}{\cal L}ie^*_{m+1} \quad \stackrel{w^*}{\lra} \quad 
{\rm T}_{(m)}^1 \stackrel{d}{\lra}{\rm T}_{(m)}^2 \stackrel{d}{\lra}
 ... \stackrel{d}{\lra}{\rm T}_{(m)}^{m-2}.
\end{equation}
\section{The twistor transform and the Hodge complex functor} \la{hc3sec}

\subsection{The Dolbeaut complex of a variation of Hodge structures.}

\paragraph{Variations of Hodge structures.} Let $X$ be a complex manifold, given by 
 a $C^\infty$-manifold ${\cal X}$ with an integrable almost complex structure $I$. 
Then  $\overline X:= ({\cal X}, -I)$ is the complex conjugate manifold. 
Denote by ${\cal A}^{p,q}_X$ the space of smooth $(p,q)$-forms on $X$. 
So ${\cal A}^{p,q}_{\overline X} = {\cal A}^{q,p}_{X}$. 

Given a complex $C^\infty$-vector bundle ${\cal E}$ on ${\cal X}$, a structure of a holomorphic vector 
bundle on ${\cal E}$ over the complex manifold $X$ is determined by an integrable 
 $\overline \partial$-connection on ${\cal E}$:
$$
\overline \partial: {\cal E} \lra {\cal E}\otimes {\cal A}^{0,1}_{X}, ~~~~\overline \partial^2=0. 
$$
Denote by $E$ the holomorphic vector bundle on $X$ determined by the $\overline \partial$-connection. Its sheaf of holomorphic sections is given by 
 the sections of ${\cal E}$ 
annihilated by $\overline \partial$. 

A filtration $F^\bullet$ by holomorphic subbundles of $E$ can be described by a decreasing filtration  
$F^{p}{\cal E} \supset F^{p+1}{\cal E}\supset \ldots $ by 
$C^\infty$-subbundles,  preserved by $\overline \partial$:  
$$
\overline \partial: F^p{\cal E} \lra F^p{\cal E} \otimes {\cal A}^{0,1}_{X}. 
$$

Let $({\cal E}, {\bf d})$ be a complex $C^\infty$-vector bundle ${\cal E}$ with a connection ${\bf d}$ 
on a manifold ${\cal X}$. 

If ${\cal X}$ has a structure of  a complex manifold 
 $X$, the connection  has a decomposition
$$
{\bf d} = \partial_X + \overline \partial_X, ~~~~\partial_X: {\cal E} \lra {\cal E} \otimes {\cal A}_X^{1,0}, ~~
\overline \partial_X: {\cal E} \lra {\cal E} \otimes {\cal A}_X^{0,1}.
$$
The same connection, considered 
on the complex conjugate manifold $\overline X$, has a decomposition 
$$
{\bf d} = \partial_{\overline X} + \overline \partial_{\overline X}, ~~~~\partial_{\overline X} = \overline \partial_{X}, ~~
\overline \partial_{\overline X} = \partial_{X}. 
$$

A $C^\infty$-filtration $F^\bullet$ of ${\cal E}$ and a connection ${\bf d}$ on ${\cal E}$ 
satisfy the Griffiths transversality condition on a complex manifold $X = ({\cal X}, I)$ if 
$$
\partial_X: F^p{\cal E} \lra F^{p-1}{\cal E} \otimes {\cal A}^{1,0}_{X}. 
$$

 So a filtered $C^\infty$-vector bundle with connection  $({\cal E}, F^\bullet, {\bf d})$ 
gives rise to a holomorphic vector bundle $E$ on $X$ with a 
filtration  
by holomorphic subbundles satisfying the Griffiths transversality condition   
if and only if
\be \la{GCON}
\overline \partial_{X}^2=0, ~~~~
\overline \partial_{X}: F^p{\cal E} \lra F^{p}{\cal E} \otimes {\cal A}^{0,1}_{X}, ~~~~\partial_{X}: F^p{\cal E} \lra F^{p-1}{\cal E} \otimes {\cal A}^{1,0}_{X}.
\ee
So it makes sense to combine the two conditions on the triple $({\cal E}, F^\bullet, {\bf d})$ 
into a single one, which we call the Griffiths condition.

\bd
A filtered $C^\infty$-vector bundle with connection  $({\cal E}, F^\bullet, {\bf d})$ on a complex manifold $X$ 
satisfies the Griffiths condition if it satisfies (\ref{GCON}).
\ed

Recall that two filtration $F^\bullet$ and $\overline F^\bullet$ on a complex vector space $V$ are $n$-complimentary 
if 
$$
{\rm gr}^p_F{\rm gr}^q_{\overline F}V=0 ~~\mbox{unless $p+q=n$}. 
$$
In this case the vector space $V$ is bigraded of pure weight $n$: 
$$
V = \oplus_{p+q=n} {\rm gr}^p_F{\rm gr}^q_{\overline F}V. 
$$

Given a complex vector space $V$, denote by $\overline V$ the complex conjugate vector space. 
It coincides with $V$ as a real vector space. The new action $\circ$  of $\C$ given by 
$\lambda \circ v:= \overline \lambda v$. 
\bd
A weight $n$ variation of Hodge structures\footnote{What we call a variation of Hodge structures is usually called a 
variation of \underline{complex} Hodge structures.} is a data $(X, {\cal E}, {\bf d}; F^\bullet, \overline F^\bullet)$, 
given by: 

\begin{itemize} 

\item A complex manifold $X$; a complex $C^\infty$-vector bundle ${\cal E}$ with  
a flat connection ${\bf d}$ on 
${ X}$. 

\item 
A filtration $F^\bullet$ by subbundles of ${\cal E}$ satisfying  the Griffiths condition; 
$$
\overline \partial_{X}: F^p{\cal E} \lra F^{p}{\cal E} \otimes {\cal A}^{0,1}_{X}, ~~~~
\partial_{X}: F^p{\cal E} \lra F^{p-1}{\cal E} \otimes {\cal A}^{1,0}_{X}.
$$
\item 
A filtration $\overline F^\bullet$ by subbundles  of $\overline {\cal E}$ satisfying   the Griffiths condition: 
$$
\partial_{X}: \overline F^p\overline {\cal E} \lra \overline F^{p}\overline {\cal E} \otimes {\cal A}^{1,0}_{X}, 
~~~~\overline \partial_{X}: \overline F^p\overline {\cal E} \lra \overline F^{p-1}\overline {\cal E} \otimes {\cal A}^{0,1}_{X}. 
$$
\item Filtrations $F^\bullet$ and $\overline F^\bullet$ are $n$-complimentary. 
\end{itemize}

A variation of  Hodge structures is a direct sum of  variations of different weights. 
\ed

Denote by $\overline {\cal E}$ the complex conjugate vector bundle to ${\cal E}$. 
If $E$ is a holomorphic vector bundle over a complex manifold $X$, then 
$\overline E$ is a holomorphic vector bundle over the  $\overline X$. 
The conditions on $\overline F^\bullet$ just mean that  
it is a filtration of the complex bundle $\overline {\cal E}$ on $\overline X$ 
by 
subbundles satisfying the Griffiths  condition.

The Galois group ${\rm Gal}(\C/\R)$  acts on the category of variations of  complex Hodge structures. 
Its generator acts by  
$$
(X, {\cal E}, {\bf d}; F^\bullet, \overline F^\bullet) \lra 
(\overline X, \overline {\cal E}, {\bf d}; \overline F^\bullet,  F^\bullet).
$$

\bd
The category $\R$-MHS of  real variations of Hodge structures  is the category of 
${\rm Gal}(\C/\R)$-equivariant objects in the category MHS. 
\ed

Equivalently, using the traditional terminology, 
a real weight $n$ variation of Hodge structures is given by a real local system 
${\cal L}$ on a complex manifold $X$, whose complexification is 
equipped with a filtration $F^\bullet$ 
by holomorphic subbundles satisfying Griffiths transversality condition. 
Finally, the filtration $F^\bullet$  and its complex conjugate $\overline F^\bullet$  
are $n$-complementary.

\paragraph{The Dolbeaut complex of a variation of Hodge structures.} 
Let $({\cal L}, {\bf d})$ be a variation 
of Hodge structures 
over  $X$. Here ${\cal L}$ denotes the underlying complex local system. 
Its de Rham complex ${\cal A}^{*}({\cal L}):= 
{\cal A}^{*}\otimes_\C {\cal L}$ has a differential, also denoted by ${\bf d}$. 
Recall the Hodge decomposition
$$
{\cal L} = \oplus_{p,q}{\cal L}^{p,q}.
$$
We define the {\it Dolbeaut bicomplex} of ${\cal L}$ by setting
\be \la{1.20.15.1}
{\cal A}^{*,*}({\cal L}) = \oplus_{s,t}{\cal A}^{s,t}({\cal L}), \qquad 
{\cal A}^{s,t}({\cal L}):= \oplus_{s=a+p, ~t=b+q}
{\cal A}^{a,b}\otimes_{\C} {\cal L}^{p,q}.
 \ee
 It is bigraded
by the {\it Hodge bidegree}  $(s,t)$, and equipped with the   
differential ${\bf d}$.

Now we  use the fact that ${\cal L}$ is a \underline{variation} of 
 Hodge structures. 
The Griffiths  conditions for $F^\bullet$ and $\overline F^\bullet$   
 just mean that  
the components of the differential ${\bf d}$ are of Hodge 
bidegrees $(1,0)$ and $(0,1)$. 

Let ${\bf d} = \partial' + {\partial''}$ be the decomposition 
into  components of Hodge bidegrees $(1,0)$ and $(0,1)$. 
Decomposing  ${\bf d}^2=0$ into the components of Hodge bidegrees 
$(2,0)$, $(1,1)$ and $(0,2)$ we get 
$$
\partial'^2 = \partial''^2= \partial'\partial''+ \partial''\partial'=0.
$$
So we get the Dolbeaut  bicomplex $({\cal A}^{*,*}({\cal L}), \partial', 
\partial'')$ of a variation of Hodge structures. It satisfies the $\partial\overline \partial$-lemma.  
Its total complex is  the Rham complex. 
Set
$$ 
{\bf d}^\C := \partial' - {\partial''}.
$$

Therefore the $({\cal A}^{*,*}({\cal L}), \partial', 
\partial'')$ is a cohomological Dolbeaut complex defined in Section \ref{hc2sec}.

\paragraph{Remark.} Let ${\bf d} = \partial + {\overline \partial}$ be the decomposition 
into  the holomorphic and antiholomorphic components. 
Let $\partial_0$ and $\overline \partial_0$ be the components of 
$\partial$ and $\overline \partial$ preserving the Hodge decomposition. 
Then we have
$$
\partial =\partial_0 +\lambda^{1,0}, \qquad 
\lambda^{1,0}\in {\rm Hom}({\cal L}^{*,*}, {\cal L}^{*-1, *+1})
\otimes {\cal A}^{1,0},
$$
$$
\overline \partial  =\overline \partial_0+\lambda^{0,1}, \qquad \lambda^{0,1}
\in {\rm Hom}({\cal L}^{*, *}, {\cal L}^{*+1, *-1})\otimes {\cal A}^{0, 1}.
$$
$$
\partial' = \partial_0 + \lambda^{0,1}, \quad 
\partial'' = \overline \partial_0 + \lambda^{1,0}; 
$$
However we avoid using  these formulas and operators $\partial$ and ${\overline \partial}$.

\subsection{The Hodge complex functor ${\cal C}_{\cal H}^\bullet(-)$} \la{2.21.15.1}

In Section \ref{2.21.15.1} we define the Hodge complex ${\cal C}_{\cal H}^\bullet({\cal A}^{\ast, \ast})$
of a cohomological Dolbeaut complex
 ${\cal A}^{\ast, \ast}$. Using the natural structure of a Hodge-Deligne complex 
on a cohomological Dolbeaut complex ${\cal A}^{\ast, \ast}$, we define a 
quasiisomorphism 
$$
{\rm RHom}_{{\rm MHS}}(\C(0), {\cal A}^{\ast, \ast}) \stackrel{\rm qis}{\lra} 
{\cal C}_{\cal H}^\bullet({\cal A}^{\ast, \ast}).
$$
Therefore the Hodge complex ${\cal C}_{\cal H}^\bullet({\cal L})$ 
of a variation of Hodge structures ${\cal L}$ on a complex manifold $M$ 
calculates the Hodge cohomology of ${\cal L}$. There are similar results for the $\R$-Hodge complex. 




\paragraph{The Hodge complex ${\cal C}_{\cal H}^{\bullet}({\cal A}^{\ast, \ast})$.}  
Let ${\cal A}^{\ast, \ast}$ be a cohomological Dolbeaut complex,  
bigraded by the Hodge bidegree. Its  cohomological grading is given by the cohomological degree   
${\rm deg}_{\cal A}$. 
The  differential  is a sum ${\bf d}=\partial' +\partial''$ of two anticommuting differentials, where   
$\partial'$ is of Hodge bidegree $(1,0)$, and 
$\partial''$ is of Hodge bidegree $(0,1)$. 

\begin{figure}[ht]
\centerline{\epsfbox{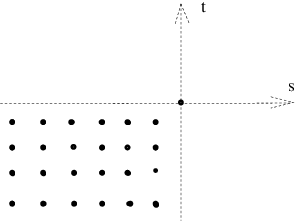}}
\caption{The Hodge bigrading $(s,t)$ of the Hodge complex ${\cal C}_{\cal H}^{\bullet}({\cal A}^{\ast, \ast})$.}
\label{hc7}
\end{figure}

Let ${\cal A}^{0,0}_{\rm cl}$ be the subspace of ${\bf d}$-closed elements in 
${\cal A}^{0,0}$. 
Consider the following subspace of 
${\cal A}^{\ast, \ast}$ 
sitting in non-positive Hodge bidegrees, see Fig \ref{hc7}:
\begin{equation} \label{wpart}
{\cal C}_{\cal H}({\cal A}^{\ast, \ast}):= \bigoplus_{s,t \leq -1}
{\cal A}^{s, t} 
\bigoplus {\cal A}_{\rm cl}^{0,0}.
\end{equation}

Let us define a degree  ${\rm deg}$ 
on $ {\cal C}_{{\cal H}}({\cal A}^{\ast, \ast})$. 
Given an element $\varphi$ of degree ${\rm deg}_{\cal A}(\varphi)$ in ${\cal A}^{s,t}$,  set 
\begin{equation} \label{11.18.ups.2hu}
{\rm deg}(\varphi) := \left \{\begin{array}{ll}
{\rm deg}_{\cal A}(\varphi) & 
\mbox{if $(s,t)=(0,0)$,}   \\
 {\rm deg}_{\cal A}(\varphi) +1 & \mbox{otherwise.} \\
 \end{array} \right.
\ee
In other words, 
\begin{equation} \label{wpart1}
{\cal C}^\bullet_{\cal H}({\cal A}^{\ast, \ast}):= \bigoplus_{s,t \leq -1}
({\cal A}^{s, t})^\bullet[-1] 
\bigoplus ({\cal A}_{\rm cl}^{0,0})^\bullet.
\end{equation}

We define a degree $1$ differential $D$ on ${\cal C}^\bullet_{{\cal H}}({\cal A}^{\ast, \ast})$ 
by setting 
\be \label{11.18.ups.2hus}
D \varphi = \left \{\begin{array}{ll}
-\frac{1}{2}{\bf d}{\bf d}^\C & 
\mbox{ if $(s,t) = (-1, -1)$,}   \\
 {\partial'}\varphi & \mbox{ if $s< - 1, t= -1$,} \\
{\partial''}\varphi & \mbox{ if $s= -1, t < -1$,} \\
{\bf d}\varphi &
 \mbox{ otherwise.}   \\
 \end{array} \right.
\end{equation}
 So unless $(s,t) = (-1, -1)$,  the map $D$ is the 
differential ${\bf d}$   
truncated to fit the rectangle $s, t\leq -1$. 
 
 \paragraph{Conventions.}
Below $\varphi_{s,t}$ denotes a homogeneous element 
of the Hodge bidegree $(-s, -t)$ in the Hodge complex, so $s,t \geq 0$. 
Abusing notation, we still denote by $\partial'$ and $\partial''$ the components 
of $D$ of the Hodge bidegrees $(1,0)$ and $(0,1)$. So 
in the Hodge complex $\partial'\varphi_{1,t}=0$ and $\partial''\varphi_{s,1}=0$, 
while this may not be so  in the original  bicomplex ${\cal A}^{\ast, \ast}$. 
We either omit the cohomological grading, or denote it by $\bullet$. 
The Hodge bidegree is denoted by $(\ast, \ast)$. 
\vskip 3mm

Let $N$ be a linear operator on $ 
{\cal C}^{\bullet}_{{\cal H}}({\cal A}^{\ast, \ast})$ acting
 on a  homogeneous  element $\varphi_{s,t}$ by
$$
N\varphi_{s,t} = \left \{\begin{array}{ll}
0  & \mbox{ if $(s,t) = (0,0)$.}  \\
(s+t-1)\varphi_{s,t} & 
\mbox{otherwise}.
 \end{array} \right. 
$$
We denote by $N_{\varphi}$ the eigenvalue of the operator $N$ on a homogeneous element $\varphi$.

\paragraph{The real Hodge complex functor ${\cal C}^{\bullet}_{{\cal H}, \R}(-)$.}  
Let ${\cal A}^{\ast}_\R$ be a real graded vector space 
with 
a degree $+1$ cohomological differential ${\bf d}$, whose  
complexification  ${\cal A}^{\ast}_\R\otimes_\R\C$ is a 
cohomological Dolbeaut bicomplex equipped with an antiholomorphic involution  $\varphi \lms \overline \varphi$: 
$$
  {\cal A}^{n}_\R\otimes_\R\C = \oplus_{s+t=n}{\cal A}^{s,t}, ~~~~\overline {{\cal A}^{s,t}} = {\cal A}^{t, s}. 
$$
We define the real Hodge complex ${\cal C}^{\bullet}_{{\cal H}; \R}({\cal A}^\ast_\R)$ as a  subcomplex
$$
{\cal C}^{\bullet}_{{\cal H}; \R}({\cal A}^\ast_\R) \subset 
{\cal C}^{\bullet}_{{\cal H}}({\cal A}^{\ast, \ast}). 
$$
Namely, a homogeneous element $\varphi \in  {\cal C}^{\bullet}_{{\cal H}}({\cal A}^{\ast, \ast})$
 is in ${\cal C}^{\bullet}_{{\cal H}; \R}({\cal A}^\ast_\R)$ if and only if 
\be \la{11.19.ups.10}
\overline {\varphi} = \left \{\begin{array}{ll}
\varphi & \mbox{ ~ if $(s,t)=(0,0)$}. \\
-\varphi& \mbox{ ~ otherwise}.   \\
 \end{array} \right.
\ee

\paragraph{The Hodge complex ${\cal C}^{\bullet}_{{\cal H}}({\cal L})$.} 
Applying the functor ${\cal C}^{\bullet}_{{\cal H}}(-)$ to the cohomological 
Dolbeaut bicomplex ${\cal A}^{\ast, \ast}({\cal L})$ provided by a variation of Hodge structures ${\cal L}$, 
we get {\it the Hodge complex ${\cal C}^{\bullet}_{{\cal H}}({\cal L})$ of ${\cal L}$}. 
Given an $(a,b)$-form $\varphi \in {\cal A}^{a,b}({\cal L})$, its cohomological degree   
in $ {\cal C}^{\bullet}_{{\cal H}}({\cal L})$ is given by 
\begin{equation} \label{11.18.ups.2hu}
{\rm deg}(\varphi) := \left \{\begin{array}{ll}
a+b & 
\mbox{if $(s,t)=(0,0)$,}   \\
 a+b+1 & \mbox{otherwise.} \\
 \end{array} \right.
\ee

\paragraph{The real Hodge complex ${\cal C}^{\bullet}_{{\cal H}; \R}({\cal L})$.} 
Let ${\cal L}$ be a variation of \underline{real} Hodge structures. 
Let ${\cal A}^*_\R({\cal L})$  be its real de Rham complex. 
Its complexification is a cohomological Dolbeaut bicomplex.  
Applying the real Hodge complex functor ${\cal C}^{\bullet}_{{\cal H}; \R}(-)$ to ${\cal A}^*_\R({\cal L})$ 
we get the real Hodge complex of the variation of real Hodge structures ${\cal L}$:
$$
{\cal C}^{\bullet}_{{\cal H}; \R}({\cal L}):= {\cal C}^{\bullet}_{{\cal H}; \R}({\cal A}^\ast_\R({\cal L})).
$$

\paragraph{Calculating Hodge cohomology.} A pure weight $w$ Hodge-Deligne complex\footnote{What we call here 
Hodge-Deligne complex traditionally should be called \underline{complex} Hodge-Deligne complex.}  is a 
complex $(K^\bullet, d_K)$ with finite dimensional cohomology and the following structures: 

i) Filtrations $F^\bullet$ and $\overline F^\bullet$ on $K^\bullet$ such that the 
differential $d_K$ is strictly compatible with  them. 

ii) The filtrations $F^\bullet$ and $\overline F^\bullet$ induce a pure weight $w+n$ Hodge structure 
on $H^n(K^\bullet)$.

\paragraph{Example.} Our weight $w$ cohomological Dolbeaut complex ${\cal A}^{\ast, \ast}$ is  an example of a weight $w$ 
Hodge-Deligne complex  for the following standard filtrations:
\be \la{HDCE}
F^p{\cal A}^{\ast, \ast} =  \oplus_{s \geq p} 
{\cal A}^{s,t}, ~~~~
\overline F^q{\cal A}^{\ast, \ast} =  \oplus_{t \geq q} 
{\cal A}^{s,t}.
\ee
The $\partial\overline \partial$-lemma implies that the differential is strictly compatible 
with the filtrations $F^\bullet$ and $\overline F^\bullet$. 

\begin{proposition} \label{4.13.ups.1z} Let ${\cal A}^{\ast, \ast}$ be a cohomological Dolbeaut complex. 
Then the Hodge complex ${\cal C}^{\bullet}_{{\cal H}}({\cal A}^{\ast, \ast})$ is quasiisomorphic to 
 ${\rm RHom}_{{\rm MHS}}(\C(0), {\cal A}^{\ast, \ast})$. 
\end{proposition}

\begin{proof}  Given a Hodge-Deligne complex $K^\bullet$, we have 
\be \la{11.19.ups.1}
{\rm RHom}_{\rm MHS}(\C(0), K^\bullet) = 
{\rm Cone}\Bigl((\overline F^0 \cap W_0)(K^\bullet)\oplus (F^0 \cap W_0)(K^\bullet)
 \lra W_0(K^\bullet)\Bigr).
\ee

So we can apply formula (\ref{11.19.ups.1}) to calculate ${\rm RHom}_{\rm MHS}(\C(0),{\cal A}^{\ast, \ast})$.  

\begin{lemma} \label{11.19.ups.2}
There is a quasiisomorphism
$$
{\rm Cone}\Bigl((\overline F^0 \cap W_0)({\cal A}^{\ast, \ast})\oplus (F^0 \cap W_0)({\cal A}^{\ast, \ast})
 \lra W_0({\cal A}^{\ast, \ast})\Bigr) ~\lra~ {\cal C}^{\bullet}_{{\cal H}}({\cal A}^{\ast, \ast}).
$$
\end{lemma}

\begin{proof} Let $f^*: X^* \to Y^*$ be a morphism of complexes such that 
$f^i$ is injective for $i < k$, isomorphism for $i=k$, and surjective for $i > k$. 
Then there is a complex 
$$
Z^*:= {\rm Coker} f^{<k}[-1] \lra  {\rm Ker}f^{>k}
$$
with a differential $D: {\rm Coker} f^{k-1} \to  {\rm Ker}f^{k+1}$ 
defined by a diagram chase (\cite{G8}, Proposition 2.1). 
Thanks to Lemma 2.2 in {loc. cit.}, the complex $Z^*$ is canonically isomorphic to 
${\rm Cone}(X^* \to Y^*)$.

Let us apply this construction in our case. Using (\ref{HDCE}) we get 
\be 
\begin{split}
&{\rm Ker}\Bigl((\overline F^0 \cap W_0) \oplus (F^0 \cap W_0) \lra W_0\Bigr) = 
\overline F^0 \cap F^0 \cap W_0 = {\cal A}^{0,0}. \\
&{\rm Coker}\Bigl((\overline F^0 \cap W_0) \oplus (F^0 \cap W_0) \lra W_0\Bigr) =\oplus_{s,t\leq -1} {\cal A}^{s,t}. \\
\end{split}
\ee
The differential is the differential $D$. 
The Lemma and hence Proposition \ref{4.13.ups.1z} are proved.
\end{proof} \end{proof}

Similarly if $K^\bullet_\R$ a real  Hodge-Deligne complex, the standard result \cite{B} tells that
\be \la{11.19.ups.1z}
{\rm RHom}_{MHS/\R}(\R(0), K^\bullet_\R) = {\rm Cone}\Bigl(W_0K^\bullet_\R \oplus F^0W_0K^\bullet_\C \lra 
W_0K^\bullet_\C\Bigr).
\ee
Therefore a similar argument shows that there  is a quasiisomorphism
$$
{\rm Cone}\Bigl(W_0{\cal A}_\R^\bullet \oplus F^0W_0{\cal A}_\C^\bullet \lra 
W_0{\cal A}_\R^\bullet\Bigr) ~\lra~ {\cal C}^{\bullet}_{{\cal H}_\R}({\cal A}_\C^\bullet).
$$ 

\begin{figure}[ht]
\centerline{\epsfbox{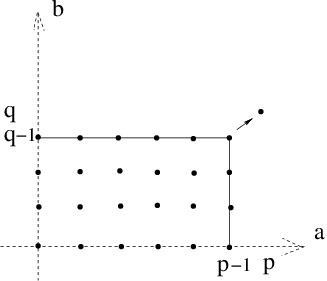}}
\caption{The de Rham bigrading $(a,b)$ 
on the ${\cal L}^{-p,-q}$-part of the Hodge complex ${\cal C}^{\bullet}_{{\cal H}}({\cal L})$.}
\label{feyn6}
\end{figure}

\begin{proposition} \label{4.13.ups.1}
Let ${\cal L}$ be a variation of Hodge structures on a complex manifold $X$. 
Then the complex ${\cal C}^{\bullet}_{{\cal H}; \R}({\cal L})$ is quasiisomorphic to 
 ${\rm RHom}_{{\rm Sh}_{{\rm Hod}(X)}}(\C(0), {\cal L})$. 

Let ${\cal L}_\R$ be a variation of real Hodge structures on a complex manifold $X$. 
Then the complex ${\cal C}^{\bullet}_{{\cal H}; \R}({\cal L}_\R)$ is quasiisomorphic to 
 ${\rm RHom}_{{\rm Sh}_{\R-{\rm Hod}(X)}}(\R(0), {\cal L}_\R)$. 
\end{proposition}

\begin{proof} We calculate 
${\rm RHom}_{{\rm Sh}_{{\rm Hod}(X)}}(\C(0), {\cal L})$ 
in two steps. First, take the   Dolbeaut  
bicomplex ${\cal A}^{\ast, \ast}({\cal L})$ of ${\cal L}$.  
It is a Hodge-Deligne complex of weight given by the weight of ${\cal L}$. 
Then calculate 
${\rm RHom}_{\rm MHS}(\C(0), {\cal A}^{\ast, \ast}({\cal L}))$  using Proposition \ref{4.13.ups.1z}. 
The second part is similar. 
\end{proof}

\subsection{The twistor transform and the DGCom structure on Hodge complexes}

Recall  the twistor plane  $\C^2$ with  the canonical coordinates 
$(z,w)$, and  the twistor line $\C$ given by the equation $z+w=1$, with the  parameter  $u$ such that 
$z=1-u, w=u$. 

Denote by $(\Omega^\bullet_{\C^2}, d_{\C^2})$ the algebraic de Rham complex  on $\C^2$.
Set 
\be \la{invinva}
{\cal A}^{\ast, \ast} \otimes_\C {\Omega}^{\bullet}_{\C^2}.
\ee
There is an action of the group $\C^* \times \C^*$ on (\ref{invinva}) given by the natural action 
on $\C^2$ and the action on the 
 ${\cal A}^{\ast, \ast}$ provided by the Hodge bigrading: an element 
$(u,v)\in \C^*\times \C^*$ acts by 
$$
(z,w) \lms (uz,vw), ~~~~\varphi_{s,t}\lms u^{-s}v^{-t}\varphi_{s,t},
~~~~ \varphi_{s,t}\in {\cal A}^{-s, -t}.
$$

\vskip 3mm
Take the invariants of  the 
${\C^*\times\C^*}$ action:
\be \la{invinvtw}
\Bigl({\cal A}^{\ast, \ast} \otimes_\C \Omega^{\bullet}_{\C^2}\Bigr)^{\C^*\times\C^*}. 
\ee

\bd \la{12.09.ups.15} The \underline{homogeneous twistor transform} is a linear map 
$$
\widehat \gamma: {\cal C}_{\cal H}^{\bullet}({\cal A}^{\ast, \ast})
 \lra \Bigl({\cal A}^{\ast, \ast} \otimes_\C\Omega^{\leq 1}_{\C^2}\Bigr)^{\C^*\times\C^*},
$$
$$
\widehat \gamma: \varphi_{0,0}+ \sum_{s,t\geq 0} \varphi_{s+1,t+1} \lms 
$$
\be \la{twistortr}
\varphi_{0,0}+ \sum_{s,t \geq 0}
z^{s}w^{t}\Bigl(
(zdw-wdz) \wedge (s+t+1)\varphi_{s+1, t+1}  + (w\partial' - z\partial'')\varphi_{s+1,t+1} \Bigr).
\ee

The \underline{twistor transform} $\gamma$ is its restriction to  the twistor line $\C$:
$$
\gamma: {\cal C}_{\cal H}^{\bullet}({\cal A}^{\ast, \ast})
 \lra {\cal A}^{\ast, \ast} \otimes_\C\Omega^{\bullet}_{\C},
$$
\ed

The homogeneous twistor transform  is  injective and functorial.

To formulate the reality condition for the homogeneous twistor transform, recall the complex conjugation $c$ 
and the involution 
$\sigma:(z,w) \lms (\overline w, \overline z)$ of the $\C^2$. So there is an involution
$$
c \otimes \sigma: {\cal A}^{\ast, \ast}\otimes \Omega^\bullet_{\C^2} \lra 
{{\cal A}^{\ast, \ast}}\otimes \Omega^\bullet_{\C^2}, ~~~~ f \otimes \omega\lms c(f \otimes \sigma^*(\omega)).
$$

One easily checks the following. 
\bl
$$
(c\otimes \sigma)^*\widehat \gamma(\varphi) = \widehat \gamma(\varphi) ~~
\mbox{if and only if} ~~ \varphi \in {\cal C}^{\bullet}_{{\cal H}; \R}({\cal A}^{\ast, \ast}).
$$
\el
 
Denote by the subscript $+$  the coinvariants of the involution $c\otimes \sigma$. 
Then we have 
$$
\widehat \gamma: {\cal C}_{{\cal H}; \R}^{\bullet}({\cal A}^{\ast, \ast})
 \lra \Bigl({\cal A}^{\ast, \ast} \otimes\Omega^{\leq 1}_{\C^2}\Bigr)_+^{\C^*\times\C^*},
$$

The following theorem generalizes Theorem \ref{12.8.ups.1s} in the Introduction.  

\bt \la{12.8.ups.1ss} i) The \underline{homogeneous} 
twistor transform $\widehat \gamma$ preserves the degree.

ii) The image of the  {twistor transform} $\gamma$
is closed under the differential. 

iii) Let us assume that ${\cal A}^{\ast, \ast}$ is a cohomological Dolbeaut DGCom. Then the image of 
the \underline{homogeneous} 
twistor transform $\widehat \gamma$ is closed under the product. 

So the image of the  {twistor} transform $\gamma$ 
is a DG subalgebra of the de Rham DGCom of the twistor line:
$$
\gamma\Bigl({\cal C}_{\cal H}^{\bullet}({\cal A}^{\ast, \ast})\Bigr)
 \subset {\cal A}^{\ast, \ast} \otimes_\C\Omega^{\bullet}_{\C}.
$$

iv) The same true for the real twistor transform. 
\et

\paragraph{First steps of the proof of Theorem \ref{12.8.ups.1ss}.} 
 . 

We have to check the following three statements: 

\begin{itemize}

\item (i) {\it The homogeneous twistor transform  $\widehat \gamma$ preserves the degrees}. 

\item (ii) {\it Multiplicativity}. Given cohomological Dolbeaut bicomplexes
 ${\cal A}^{\ast, \ast}$ and ${\cal B}^{\ast, \ast}$, 
there is a  product map
$$
\ast: {\cal C}^{\bullet}_{\cal H}({\cal A}^{\ast, \ast})\otimes 
{\cal C}^{\bullet}_{\cal H}({\cal B}^{\ast, \ast}) \stackrel{}{\lra}
{\cal C}^{\bullet}_{\cal H}({\cal A}^{\ast, \ast}\otimes {\cal B}^{\ast, \ast})
$$
 making 
the following diagram is commutative:
\begin{displaymath} \la{12.09.ups.3}
    \xymatrix{
        {\cal C}^{\bullet}_{\cal H}({\cal A}^{\ast, \ast})\otimes 
{\cal C}^{\bullet}_{\cal H}({\cal B}^{\ast, \ast}) \ar[r]^{\ast} \ar[d]_{\widehat \gamma \otimes \widehat \gamma} & 
{\cal C}^{\bullet}_{\cal H}({\cal A}^{\ast, \ast}\otimes {\cal B}^{\ast, \ast})
\ar[d]^{\widehat \gamma} \\
         {\cal A}^{\ast, \ast}\otimes \Omega^\bullet_{\C^2} \otimes 
{\cal B}^{\ast, \ast}\otimes \Omega^\bullet_{\C^2} \ar[r]^{\wedge}       
& {\cal A}^{\ast, \ast} \otimes {\cal B}^{\ast, \ast}
\otimes \Omega^\bullet_{\C^2} }
\end{displaymath}
The map $\wedge$ at the 
bottom line uses the wedge product map $\Omega^\bullet_{\C^2} \wedge \Omega^\bullet_{\C^2} \to 
\Omega^\bullet_{\C^2}$. 
Thanks to the injectivity of $\widehat \gamma$ this implies that the map 
$\ast$ is commutative and associative.

\item (iii) {\it Differentials.} There is a differential $\delta$ on the 
Hodge complex ${\cal C}^{\bullet}_{\cal H}({\cal A}^{\ast, \ast})$ 
such that the map $\gamma$ transforms $\delta$ 
into the de Rham differential ${\bf d} + d_\C$ 
on the de Rham complex ${\cal A}^{\ast, \ast}\otimes \Omega^\bullet_{\C}$, holomorphic along the twistor line $\C$,  
i.e. we have a commutative diagram:

\begin{displaymath} 
    \xymatrix{
        {\cal C}^{\bullet}_{\cal H}({\cal A}^{\ast, \ast}) \ar[r]^{\delta} \ar[d]_{\gamma } & 
{\cal C}^{\bullet}_{\cal H}({\cal A}^{\ast, \ast})
\ar[d]^{\gamma} \\
         {\cal A}^{\ast, \ast} \otimes \Omega^\bullet_{\C} \ar[r]^{{\bf d}+{\rm d}_\C}       
& {\cal A}^{\ast, \ast} \otimes \Omega^\bullet_{\C}}
\end{displaymath}

In particular, if ${\cal A}^{\ast, \ast}$ is an algebra in the tensor category of 
cohomological Dolbeaut bicomplexes, then ${\cal C}^{\bullet}_{\cal H}({\cal A}^{\ast, \ast})$ 
has a commutative differential graded algebra structure. 
\end{itemize}

The property (i) is clear.  
Indeed, an element $\varphi_{s+1, t+1}$, whose degree in the Hodge complex 
is one plus its cohomological degree in the complex ${\cal A}^{\ast, \ast}$, is either multiplied by a 1-form, 
or acted by a degree 1 differential operator. 
We do not change  the element $\varphi_{0, 0}$
whose degree in the Hodge complex coincides with its cohomological degree in ${\cal A}^{\ast, \ast}$. 
\vskip 3mm

The injectivity of the  twistor transform immediately implies that 
the differential $\delta$ satisfies the Lebniz rule for the $\ast$-product. 

The formulas for the product $\ast$ and the differential $\delta$ 
are not complicated. It is easy to check directly that $\delta^2=0$ 
and that the product is graded commutative. 
However without using the twistor transform   
a proof of the Leibniz rule for the $\delta$ 
requires a messy calculation. 
This is hardly surprising: hiding the role of the twistor line and 
the differential $d_\C$ 
obscures the story.

\paragraph{The $\ast$-product.}
It is handy to introduce a linear map of vector spaces preserving the grading: 
$$
{\bf 1}\ast: {\cal C}^{\bullet}_{{\cal H}}({\cal A}^{\ast, \ast}) \lra 
{\cal A}^{\ast, \ast}, \qquad  
{\bf 1}\ast \varphi:=   \left\{ \begin{array}{lll} 
{\bf d}^\C \varphi& \mbox{ if $N_{\varphi}>0$}, \\ 
\varphi& \mbox{ if $N_{\varphi}=0$}. \\
 \end{array}\right.
$$
Its restriction  to ${\cal C}^{\bullet}_{{\cal H}, \R}({\cal A}^{\ast, \ast})$ satisfies 
the reality condition:
\be \la{11.22.ups.5}
\overline {{\bf 1}\ast\varphi} = {\bf 1}\ast \varphi, ~~~~\varphi \in 
{\cal C}^{\bullet}_{{\cal H}, \R}({\cal A}^{\ast, \ast}). 
\ee

We use a shorthand  $|\varphi|$ for the degree of 
 $\varphi$ in the Hodge complex.  

\begin{definition} \label{3.30.ups.3a} Given $\varphi_i \in {\cal C}_{\cal H}({\cal A}_i^{\ast, \ast})$, we define
$$
\varphi_1 \ast \ldots \ast \varphi_m \in 
{\cal C}_{\cal H}({\cal A}_1^{\ast, \ast} \otimes \ldots \otimes{\cal A}_n^{\ast, \ast}) 
$$
by setting
$$
\varphi_1 \ast \ldots \ast \varphi_m := 
$$
$$
\left\{\begin{array}{ll} 
\varphi_1 \otimes \ldots \otimes \varphi_m &  \mbox{if $\sum N_{\varphi_i}= 0$},\\
N^{-1}\sum_{i=1}^m(-1)^{|\varphi_1|
+ \ldots + |\varphi_{i-1}|}
{\bf 1}\ast\varphi_1 \otimes \ldots \otimes N\varphi_i 
\otimes \ldots 
 \otimes  {\bf 1} \ast\varphi_m&  \mbox{otherwise.}\\
\end{array} \right.
$$
\end{definition}

\begin{proposition} \label{11.14.ups.1} (i) One has 
\begin{equation} \label{3.30.ups.3qa}
{\bf 1} \ast(\varphi_1 \ast \ldots \ast \varphi_m) = {\bf 1} \ast\varphi_1 \otimes 
\ldots \otimes  {\bf 1} \ast\varphi_m. 
\end{equation}

(ii) The product $\ast$ makes the diagram (\ref{12.09.ups.3}) commutative. 

(iii) The product $\ast$  is 
a graded commutative and associative.
 
\end{proposition}

\begin{proof}  (i) Obvious.

(ii) The commutativity of the diagram 
 when $\sum N_{\varphi_i}= 0$ is obvious. 
Let us assume that $\sum N_{\varphi_i}> 0$. 
Take the $\wedge$-product of several expressions (\ref{12.09.ups.1a}).
Look at its $(zdw-wdz)$-component. 
The degree $|\varphi|$ of $\varphi$ in the Hodge complex coincides with the cohomological degree of 
${\bf 1}\ast \varphi$ in the  complex ${\cal A}^{\ast, \ast}$, as well as with the 
 de Rham degree of $(zdw-wdz) \wedge N\varphi_i$, provided $N\varphi_i>0$. 
So moving $(zdw-wdz)$ to the left  we get  
$$
(zdw-wdz) \wedge \sum_{i=1}^k(-1)^{|\varphi_1|  
+ \ldots + |\varphi_{i-1}|}
{\bf 1}\ast\varphi_1 \otimes \ldots \otimes N\varphi_i  
\otimes \ldots 
 \otimes  {\bf 1} \ast\varphi_m.
$$
This is equal to   
$(zdw-wdz)  \wedge N(\varphi_1 \ast \ldots \ast \varphi_m)$. 
So we recovered Definition \ref{3.30.ups.3a} for the $\ast$-product in the case 
$\sum N_{\varphi_i}> 0$. 
The claim that the second term in (\ref{12.09.ups.1a}) 
is multiplicative 
is obvious (and equivalent to (i)).

(iii) Follows from the injectivity of the twistor transform. 
\end{proof}

\paragraph{Explicit formulas.} 
\be \la{12.17.ups.1}
\varphi_1 \ast \varphi_2= 
\left\{\begin{array}{ll} 
\varphi_1 \otimes \varphi_2& \mbox{ if $N_{\varphi_1}= N_{\varphi_2}=0$},\\
    N^{-1}\Bigl(
N\varphi_1 \otimes \varphi_2 + (-1)^{|\varphi_1| }\varphi_1
 \otimes N\varphi_2\Bigr)& 
\mbox{ if $N_{\varphi_1}=0 \mbox{~or}~ N_{\varphi_2}=0$}, \\
& \mbox{ but 
$N_{\varphi_1}+ N_{\varphi_2}>0$}, \\  
N^{-1}\Bigl(N\varphi_1 \otimes {\bf d}^\C \varphi_2  + (-1)^{|\varphi_1|} 
{\bf d}^\C \varphi_1\otimes N\varphi_2 \Bigr) 
&  \mbox{ if $N_{\varphi_1}>0,  N_{\varphi_2}>0$.}\\
\end{array} \right.
\ee

The {twistor transforms} can be written as  
\be \la{12.09.ups.1a}
\widehat \gamma: \varphi = \sum_{s,t} \varphi_{s,t}
\lms \sum_{s,t \geq 0}
z^{s}w^{t}\Bigl(
(zdw-wdz) \wedge (s+t+1)\varphi_{s+1, t+1}  +
({\bf 1}\ast \varphi)_{s,t}\Bigr).
\ee

\paragraph{The differential $\delta$.} Recall the differential  $D$ from 
on $ {\cal C}_{\cal H}({\cal A}^{\ast, \ast})$ defined in (\ref{11.18.ups.2hus}). 

Let 
$\mu$ be a linear map on ${\cal C}^{\bullet}_{\cal H}({\cal A}^{\ast, \ast})$ 
acting on the homogeneous components as 
\be \la{mu}
\mu(\varphi_{0, 0}):= \varphi_{0, 0}; ~~~~
\mu(\varphi_{s+1, t+1}) := {s+t \choose s}^{-1}\varphi_{s+1, t+1} ~~
\mbox{if $ s,t \geq 0$}. 
\ee

\begin{definition} \la{DEL} $\delta:= -\mu^{-1}\circ D \circ \mu$. 
\end{definition}


We are going to give three different proofs of the part ii) of Theorem \ref{12.8.ups.1ss}. 

1) The first is a direct check, which is the most convincing. 
The reader can skip this proof. 

2) The second contains an important characterisation of the 
twistor transform. Yet the corresponding computation is a version of the one in  1). 

3) The third is the most conceptual. It elucidates  
the origin of the map $\mu$  in Definition \ref{DEL}.

\paragraph{1) A direct proof of the part ii) of Theorem \ref{12.8.ups.1ss}.} 
We start with an explicit formula for the differential $\delta$ on ${\cal C}_{\cal H}({{\cal A}^{\ast, \ast}})$ 
from Definition \ref{DEL}. 

\begin{lemma} \la{11.24.ups.10} One has 
\begin{equation} \label{11.18.ups.2inv}
\delta(\varphi_{0,0})= 0, ~~~~~~\delta(\varphi_{s+1,t+1})= \left\{\begin{array}{ll}    
\partial'' \partial' \varphi_{1,1}& 
\mbox{if $(s,t)=(0,0)$}   \\
\frac{-1}{s+t}\Bigl( 
s\partial' + 
t \partial'' \Bigr)\varphi_{s+1, t+1}.&
 \mbox{if $s+t\geq 1$}   \\
 \end{array} \right.
\end{equation}
\end{lemma}

\begin{proof}  The cases when $(s,t)$ is $(0,0)$ or $(1,1)$  are obvious. Otherwise we should have 
\begin{equation} \label{11.18.ups.2}
\delta(\varphi_{s+1, t+1} )\stackrel{?}{=}  
\frac{-1}{s+t}\Bigl( 
s\partial' + 
t \partial'' \Bigr)\varphi_{s+1, t+1}.
\end{equation}
Let us prove the claim for $\partial'$-component.  
Using (\ref{11.18.ups.2}) as a definition, we have  
$$
-\mu\delta'(\varphi_{s+1, t+1})=  
{s+t-1 \choose s-1}^{-1}\frac{s}{s+t}\partial'\varphi_{s+1, t+1} ~~
= 
{s+t \choose s}^{-1}\partial'\varphi_{s+1, t+1} = 
\partial'\mu(\varphi_{s+1, t+1}). 
$$
The case of $\partial''$-component is similar. 
\end{proof}

Let us now prove that 
the twistor transform commutes with the differentials by a direct computation.  
We use the notation ${\rm Res}_{z+w=1}$ for the restriction to the line $z+w=1$. 
One has  
\be \la{11.18.ups.1asti}
({\bf d}+d_{\C}) \gamma(\varphi) = 
{\rm Res}_{z+w=1}~\Bigl(-\sum_{s, t\geq 0}
z^{s}w^{t}(zdw-wdz)\wedge
(s+t+1) {\bf d}  \varphi_{s+1, t+1} 
\ee
\be \la{11.18.ups.1asi}
~~+~~\sum_{s, t\geq 0}{\rm d}_{\C^2}\left[z^{s}w^{t}
(w\partial' - z \partial'')\right]\varphi_{s+1,t+1}
\ee
\be \la{12.12.ups.10}
~~+~~\sum_{s,t\geq 0}
z^{s}w^{t}
{\bf d} (w\partial' - z \partial'')\varphi_{s+1,t+1}\Bigr). 
\ee

i) Let us look first at (\ref{11.18.ups.1asti}) + (\ref{11.18.ups.1asi}), 
which is  the $\Omega^1_{\Bbb C}$-part of $({\bf d}+d_{\C}) \gamma(\varphi)$. 
We examine the $\partial'$-component - 
the other one is similar.  
Since in the Hodge complex $\partial'\varphi_{1, t+1}=0$, we can assume $s>0$, and use
 the summation over $s\geq 1,t\geq 0$:
\be \la{12.1.ups.3}
{\rm Res}_{z+w=1}~\sum_{s\geq 1, t\geq 0}z^{s-1}w^{t}  
\Bigl(-z( s+t+1)(zdw-wdz)  
+swdz + (t+1)zdw)\Bigr)
\partial'\varphi_{s+1, t+1}. 
\ee

 Now we substitute $z=1-u, w=u$. Then $zdw-wdz  =du$, and  we get
$$
{\rm Res}_{{z=1-u, ~w=u}}~\Bigl(-z( s+t+1)(zdw-wdz)  
+swdz + (t+1)zdw)\Bigr) =
$$
$$
 \Bigr(-(1-u)( s+t+1)
-su + (t+1)(1-u)\Bigr)du = -sdu.
$$
Therefore we get
\be \la{12.1.ups.3a}
(\ref{12.1.ups.3}) = 
-\sum_{s\geq 1, ~t \geq 0}(1-u)^{s-1}u^{t}du\wedge  
s\partial'\varphi_{s+1,t+1}.
\ee

On the other hand, we have
$$
\gamma (\delta \varphi) =  {\rm Res}_{z+w=1}\sum_{s,t \geq 0}
z^{s}w^{t}
\Bigl((s+t+1) (zdw-wdz) \wedge (\delta\varphi)_{s+1, t+1}
+ (w\partial' - z\partial'') (\delta\varphi)_{s+1,t+1}\Bigr). 
$$

Using (\ref{11.18.ups.2inv}) we see that the $\Omega^1_{\Bbb C}$-part of its 
$\partial'$-component, reduced to the $z+w=1$ line, is 
\be \la{12.12.ups.4}
-\sum_{s,t \geq 0}
(1-u)^{s}u^{t}
\Bigl((s+t+1) du\wedge \frac{(s+1)}{s+t+1}\partial'\varphi_{s+2, t+1}\Bigr)=
\ee
\be \la{12.12.ups.4*}
-\sum_{s,t \geq 0}
(1-u)^{s}u^{t}
\Bigl(du\wedge (s+1)\partial'\varphi_{s+2, t+1}\Bigr).
\ee
We conclude that  (\ref{12.1.ups.3a})  
coincides with (\ref{12.12.ups.4*}).

\vskip 3mm
ii) Let us  compare the $du$-free part of $({\bf d}+d_{\C})\gamma(\varphi)$, given by (\ref{12.12.ups.10}), with 
 the $du$-free part of $\gamma (\delta\varphi)$. 
If $N_{\varphi}=0$, both  are zero since 
${\bf d}\varphi_{0,0} = 0$. Otherwise this follows from a key formula
\be \la{5.4.14.10}
({\bf d}+d_{\C}) (u\partial' - (1-u)\partial'')\varphi 
 = du {\bf d}\varphi +  \partial''\partial'\varphi.
\ee

Indeed, if $N_{\varphi}>0$, the $({\bf d}+d_{\C})\gamma(\varphi)$ contributes 
${\bf d}(u\partial'-(1-u)\partial'') \varphi ~= ~\partial''\partial'\varphi.$

If $N_{\varphi}=1$, then 
$
\gamma (\delta\varphi_{1,1}) = \delta\varphi_{1,1} = \partial''\partial'  \varphi_{1,1}.
$ 

If $N_{\varphi}>1$, the $\gamma (\delta\varphi)$ contribution, calculated on a homogeneous component 
$\varphi = \varphi_{s,t}$ is:
\be \la{12.23.15.4}
(w\partial' - z \partial'')
\delta\varphi=  
-\frac{(w\partial' - z\partial'')}{s+t-2}\Bigl( 
\frac{(s-1)}{z}\partial' + 
\frac{(t-1)}{w} \partial'' \Bigr)\varphi=  \partial'' \partial'\varphi. 
\ee
Theorem \ref{12.8.ups.1ss} is proved.

\paragraph{2) A proof of Theorem \ref{12.8.ups.1ss}ii) via a characterisation of the twistor transform.} 
The subspace of $\C^*\times \C^*$-invariants 
\be \la{soinv}
\Bigl({\cal A}^{\ast, \ast}\otimes\Omega^{\leq 1}_{\C^2}\Bigr)^{\C^*\times \C^*}
\ee
is not preserved by the differential ${\bf d}+d_{\C^2}$: 
the differential ${\bf d}$ breaks the $\C^*\times \C^*$-invariance. 

Let us consider a  subspace of (\ref{soinv}), denoted by
\be \la{sub}
\Bigl({\cal A}^{\ast, \ast}\otimes\Omega^{\leq 1}_{\C^2}\Bigr)_0^{\C^*\times \C^*},  
\ee
consisting of the elements of the form 
\begin{equation} \label{8.20.05.2asam}
\Phi_{g, \beta}(z,w):=  
\sum_{s,t\geq 0}z^{s}w^{t}\Bigl((s+t+1)(zdw-wdz) g_{s+1,t+1} + {\beta}_{s,t}
  \Bigr).
\end{equation}
The image of the twistor transform is precisely the space of  the 
elements $\Phi_{g, \beta}(z,w)$ such that 
\begin{equation} \label{8.20.05.2asaq}
\beta_{s,t} = ({\bf d}^\C g)_{s,t}, ~~(s,t) \not = (0,0).
\end{equation}

Consider the restriction of the space (\ref{sub}) to the twistor line 
$z+w=1$:
\be \la{soinv1}
{\rm Res}_{z+w=1}~\Bigl({\cal A}^{\ast, \ast}\otimes\Omega^{\leq 1}_{\C^2}\Bigr)_0^{\C^*\times \C^*}.
\ee

Denote by $\Phi_{g, \beta}(u)$ be the restriction to the twistor line of an element  $\Phi_{g, \beta}(z,w)$ 
in (\ref{8.20.05.2asam}). 
\bp \la{12.22.15.1} The $({\bf d}+d_{\C})\Phi_{g, \beta}(u)$ lies in the subspace (\ref{soinv1}) if and only if (\ref{8.20.05.2asaq}) holds.

The image of the twistor transform is the maximal subspace of  (\ref{soinv1}) 
preserved by the  ${\bf d}+d_{\C}$. 
\ep

\begin{proof} The second claim is a reformulation of the first. So let us prove the first one.

The $du$-part of the differential  
$( {\bf d}+d_{\C}) \Phi_{g, \beta}(u)  $ is given by
\be \la{2.7.15.31z}
\sum_{s,t\geq 0}(1-u)^su^tdu \wedge 
\Bigl(-(s+1){\beta }_{s+1, t} + (t+1) {\beta }_{s, t+1} - (s+t+1) {\bf d}{g}_{s+1, t+1} \Bigr).
\ee
We set $g_{s,t}:=0$ if $s< 1$ or $t< 1$, and 
write (\ref{2.7.15.31z}) as a sum of $(s,t)$-components: 
$$
\sum_{s,t\geq 0}(1-u)^{s-1}u^{t-1}du \wedge 
\Bigl((-su+ t (1-u)){\beta}_{s,t}  
-(s+t)  [(1-u)\partial' g_{s+1,t} 
+ u\partial'' g_{s,t+1}] \Bigr).
$$
We stress that 
if $s=0$ then the $\wedge$-factor in brackets is divisible by $(1-u)$, thanks to $g_{0, t+1}=0$. 
Similarly for $t=0$.  
So we have a polynomial in $u$ for every $(s,t)$-component. 

We rewrite the last formula as 
\be \la{2.10.15.1zm}
(\ref{2.7.15.31z}) = \sum_{s,t}(1-u)^{s-1}u^{t-1}du \wedge (uP_{s,t}+Q_{s,t}),
\ee
where we set 
\be \la{2.7.15.20az}
\begin{split}
&P_{s,t}:= (s+t)(-{\beta}_{s,t} + \partial' g_{s+1,t} -  \partial'' g_{s,t+1}), \\
&Q_{s,t}:= t{\beta}_{s,t} 
-(s+t) \partial' g_{s+1,t}.  \\
\end{split}
\ee
Let us compare the 
Hodge bidegree $(s,t)$ parts of the $1$-forms in (\ref{8.20.05.2asam}) and (\ref{2.10.15.1zm}). 
In (\ref{8.20.05.2asam}), it is  $g_{s,t} (1-u)^{s-1}u^{t-1}du $ - a polynomial 
of degree $s+t-2$. In (\ref{2.10.15.1zm}), it is  $u(1-u)^{s-1}u^{t-1}du~P_{s,t}$ - a polynomial of degree $s+t-1$. 
Therefore $P_{s,t}=0$. This just means that, unless $s=t=0$, 
\be \la{12.11.15.1}
{\beta}_{s,t} = \partial' g_{s+1,t} -  \partial'' g_{s,t+1}.
\ee
\end{proof}

\bl
 The twistor transform intertwines the differential $\delta$ with the differential  
${\bf d}+d_{\C}$. 

\el

\begin{proof} 
Substituting (\ref{12.11.15.1}) to the formula for $Q_{s,t}$ in (\ref{2.7.15.20az}), we get 
$$
Q_{s,t}:=  -s\partial' g_{s+1,t} -  t\partial'' g_{s,t+1}. 
$$
Since $P_{s,t}=0$, comparing with (\ref{8.20.05.2asam}), we see that this implies that 
\be \la{12.23.15.1}
(s+t)~\delta(g_{s+1,t+1})= -s~\partial' g_{s+1,t+1} - t~\partial'' g_{s+1,t+1}.
\ee

But this is just the formula for the differential $\delta$ in the Hodge complex 
${\cal C}_{\cal H}({\cal A}^{*,*})$. 

Therefore 
the $du$-part of the  $({\bf d}+d_{\C})\Phi_{g,\beta}$ is given by the $du$-part of 
$\Phi_{\delta g,\beta'}$ for some $\beta'$. 

Let us check  that if $\beta'$ relates 
to $\delta g$ as in the twistor transform, then the $du$-free parts of 
$({\bf d}+d_{\C})\Phi_{g,\beta}$ and $\Phi_{\delta g,\beta'}$ coincide. 
To prove this, notice 
 the following:
\be
\begin{split}
&{\bf d}(w\partial'-z\partial'') g_{s+1,t+1} ~= ~(w+z)\partial''\partial'g_{s+1,t+1},\\
&(w\partial'-z\partial'')\delta g_{s+1,t+1} \stackrel{(\ref{12.23.15.4})}{=} 
\partial''\partial'g_{s+1,t+1},~~~~(s,t) \not = (0,0).\\ 
\end{split}
\ee
 Since 
$\delta g_{1,1} = \partial''\partial'g_{1,1}$, we get also the $(s,t) = (0,0)$ case.
\end{proof}

\paragraph{Remark.} There is  a  natural isomorphism of the \underline{vector spaces}:
\be \la{22.12.15.3}
H^*_{d_{\C^2}}\Bigl({\cal A}^{\ast, \ast}\otimes_\C
(\Omega^0_{\C^2} \stackrel{{d_{\C^2}}}{\lra} \Omega^1_{\C^2}) \Bigr)^{\C^*\times \C^*} \stackrel{{\rm qis}}{\lra} {\cal C}^{\bullet}_{\cal H}({\cal A}^{\ast, \ast}).
\ee
Indeed, the cohomology 
of the complex $\Omega^0_{\C^2} \stackrel{{d_{\C^2}}}{\lra} \Omega^1_{\C^2}$ 
are represented by $1\in \Omega^0_{\C^2}$ and $z^sw^t(zdw-wdz)\in \Omega^1_{\C^2}$. 
So multiplying by ${\cal A}^{\ast, \ast}$ and taking the $(\C^*\times \C^*)$-invariants we get the space 
${\cal C}_{\cal H}({\cal A}^{\ast, \ast})$ as defined in 
(\ref{wpart1}). Precisely, 
$$
1\otimes \varphi_{0,0} \lms \varphi_{0,0}, ~~~~
z^sw^t(zdw-wdz)\otimes \varphi_{s+1,t+1} \lms \varphi_{s+1,t+1}.
$$

  \paragraph{3) A proof of Theorem \ref{12.8.ups.1ss}ii) 
illuminating the map $\mu$, see (\ref{mu}).}  

Consider the operator of integration along the real segment $\Delta$ 
connecting points $(1,0)$ and $(0,1)$ in the twistor plane:
\be \la{5.3.14.1}
{\rm I}:=
 {\cal A}^{\ast, \ast}\otimes\Omega_\C
\lra {\cal A}^{\ast, \ast}.
\ee
The image of the composition ${\rm I}\circ \gamma$ 
lies in the subspace ${\cal C}_{\cal H}^{\bullet}({\cal A}^{\ast, \ast}) \subset {\cal A}^{\ast, \ast}$.
So  
\be \la{5.3.14.1st}
{\rm I}\circ \gamma: {\cal C}_{\cal H}^{\bullet}({\cal A}^{\ast, \ast})
\lra {\cal C}_{\cal H}^{\bullet}({\cal A}^{\ast, \ast}).
\ee

\bl \la{5.5.14.1} i) One has 
${\rm I}\circ \gamma = \mu$ on the $N_\varphi > 0$ subspace.  Precisely, 
\begin{equation} \label{5.3.14.2}
\begin{array}{ll}    
{\rm I}\circ \gamma:~ \varphi_{0, 0} \rightarrow 0, & \\
{\rm I}\circ \gamma:~ \varphi_{s+1, t+1} \rightarrow {s+t \choose s}^{-1} \varphi_{s+1, t+1}&
 \mbox{if $s,t\geq 0$}.   \\
 \end{array} 
\end{equation}
  
ii) The restriction of the operator ${\rm I}$ to the image of 
  $\gamma$ 
intertwines  differentials: 
$$
{\rm I}\circ ({\bf d}+{\rm d}_\C) \circ  \gamma = \delta.
$$
Equivalently, there is  a commutative diagram 

\begin{displaymath} \la{1.08.08.2az}
    \xymatrix{
        {\cal A}^{\ast, \ast}\otimes\Omega_\C   ~\ar[r]^{{\bf d}+{\rm d}_{\C}} & 
~~{\cal A}^{\ast, \ast}
\otimes\Omega_\C
\ar[d]^{{\rm I} } \\
{\cal C}_{\cal H}({{\cal A}^{\ast, \ast}}) \ar[u]^{\gamma}  \ar[r]^{\delta~~~~~~} & {\cal C}_{\cal H}({{\cal A}^{\ast, \ast}})\subset {\cal A}^{\ast, \ast}
}
\end{displaymath}
\el

\begin{proof} i) This follows from (\ref{2.19.15.1}), which tells that
\be \la{2.19.15.1a}
\int^1_{0}(1-u)^{s}u^{t}
du  = \frac{1}{s+t+1}{s+t \choose s}^{-1}.
\ee 

ii) Given a fibration  whose fibers are compact manifolds without boundary, 
the integration along the fibres 
commutes with the differential. 
However in our case the fibers do have boundary. 

We claim that for any $s,t \geq 0$,  we have 
$$
{\rm I} \circ {\rm Res}_{z+w=1} {\rm d}_{\C^2}\Bigl(z^sw^t (w \partial' - z\partial'')\varphi_{s+1, t+1}\Bigr)=0.
$$
Here $\varphi_{s+1, t+1}$ is a homogeneous element of the Hodge complex 
${\cal C}_{\cal H}({{\cal A}^{\ast, \ast}})$, and  
$\partial', \partial''$ are  differentials in the \underline{Hodge} complex. 
The boundary contributions are zero 
for the following reasons. 

If $s>0, t>0$ then 
$z^sw^t(w\partial' - z\partial'')\varphi_{s+1, t+1}$ restricts to zero 
when $z=0$ or $w=0$.

If $s=0, t>0$ then 
$w^t(w\partial' - z\partial'')\varphi_{1, t+1}$ restricts to zero 
when $w=0$. It also restricts to zero when $z=0$ since 
$\partial' \varphi_{1, t+1} =0$ in the Hodge complex. 
The case $s>0, t=0$ is similar. 

Finally, in the Hodge complex  $(w\partial' - z \partial'')\varphi_{1, 1}=0$.
\end{proof} 

Lemma \ref{5.5.14.1} implies that the map ${\rm I} $ commutes with the 
differentials on the $\Omega_{\C}^1$-component. 
We  use Proposition \ref{12.22.15.1}  to complete the proof.

\paragraph{Conclusion.} 
 
We have defined a tensor functor 
$$
{\cal C}_{\cal H}: \{\mbox{\rm cohomological Dolbeaut bicomplexes}\} \lra \{\mbox{\rm Complexes}\}. 
$$
In particular, it restricts to a functor 
$$
{\cal C}_{\cal H}: \{\mbox{\rm Dolbeaut DGCom's}\} \lra \{\mbox{\rm DGCom's}\}. 
$$

\subsection{Concluding remarks}

\paragraph{A) The differential $\tau$.}  The constant function $1$ on $X$ provides
 a canonical element
\be \la{11.23.ups.1}
{\bf 1}_{\R(1)}:= 1 \otimes (2\pi i)
\in {\cal C}^{1}_{{\cal H}, \R}(\R(1)).
\ee
Given a variation of $\R$-Hodge structures 
${\cal L}$, 
there is a map  provided by the $\ast$-product with 
${\bf 1}_{\R(1)}$:
 \be \la{11.23.ups.2}
\tau :  
 {\cal C}^{\bullet}_{{\cal H}, \R}({\cal L}) \lra 
{\cal C}^{\bullet+1}_{{\cal H}, \R}({\cal L}(1)), 
\qquad 
\tau  {\varphi}:= \left\{ \begin{array}{lll} 
{\bf d}^\C {\varphi} \otimes 2\pi i& \mbox{ if 
${\varphi}$ homogeneous, $N_{\varphi}>0$}, \\ 
{\varphi}\otimes 2\pi i & \mbox{ if ${\varphi}$ homogeneous, $N_{\varphi}=0$}. \\
 \end{array}\right.
\ee

The linear map ${\bf 1}\ast$ is 
the $\ast$-product with ${\bf 1}_{\R(1)}$ followed by 
multiplication by $(2\pi i)^{-1}$ and the tautological linear embedding 
${\cal C}^{\bullet+1}_{{\cal H}, \R}({\cal L}) 
\hra {\cal A}^{\bullet}({\cal L})$.

\bl \la{12.09.ups.8} One has $\tau^2=0$. 
The map $\tau$ commutes with the differential $\delta$. 
So the Hodge cohomology groups form a complex: 
$$
\ldots \stackrel{\tau }{\lra} {\rm H}^{\bullet}_{\cal H}({\cal L}) \stackrel{\tau }{\lra} 
{\rm H}^{\bullet+1}_{\cal H}({\cal L}(1)) \stackrel{\tau  }{\lra} 
{\rm H}^{\bullet+2}_{\cal H}({\cal L}(2)) \stackrel{\tau  }{\lra} \ldots 
$$
\el

\begin{proof}  Clearly 
${\bf 1}_{\R(1)}\ast {\bf 1}_{\R(1)} = {\bf d}^\C 1 \otimes (2\pi i)^2 = 0$. 
This plus associativity of the $\ast$-product implies $\tau^2=0$. 
Since $\delta {\bf 1}_{\R(1)}=0$, the second statement follows from the Leibniz rule. 
\end{proof}

A similar, yet different, complex of
 Weil-\'etale cohomologies is 
used by S. Lichtenbaum \cite{Li}.

\paragraph{B)  Algebra structures on ${\cal A}_X^*[-1]$.}
Here is a family of graded commutative 
algebra structures on 
the shifted $C^\infty$ de Rham complex ${\cal A}_X^*[-1]$ of 
a complex manifold 
$X$.  
Let $N$ be a linear operator  on ${\cal A}_X^*$ acting by a scalar 
on the subspace of forms of given degree: 
for a homogeneous form $\alpha$ we have 
 $N(\alpha) = N_{\alpha}\cdot \alpha$. Let us assume that  
for any  homogeneous forms $\alpha_1, \alpha_2$ we have 
\begin{equation} \label{3.30.ups.1}
N_{\alpha_1 \circ \alpha_2} = N_{\alpha_1 } + N_{\alpha_2}. 
\end{equation} 
Let  
$|\alpha|$ be the degree of 
$\alpha$ in ${\cal A}_X^*[-1]$. 
We define  a product $\circ$ on  ${\cal A}_X^*[-1]$ by
\begin{equation} \label{3.30.ups.3}
N(\alpha_1 \circ \ldots \circ \alpha_m) := \sum_{i=1}^k(-1)^{|\alpha_i| (
|\alpha_1| + \ldots + |\alpha_{i-1}|)}N(\alpha_i) 
 \wedge d^{\C}\alpha_1 \wedge \ldots \wedge 
\widehat {N(\alpha_{i})} \wedge  \ldots \wedge d^{\C}\alpha_m.
\end{equation}

\begin{lemma} \label{3.30.ups.2} (i) One has 
$
d^{\C}(\alpha_1 \circ \ldots \circ \alpha_m) = d^{\C}\alpha_1 \wedge 
\ldots \wedge  d^{\C}\alpha_m. 
$

\noindent
(ii) The product $\circ$ provides ${\cal A}_X^*[-1]$ with a structure of a 
 graded commutative associative algebra. 
\end{lemma}

\begin{proof}  (i) Obvious thanks to (\ref{3.30.ups.1}). 

\vskip 3mm
(ii) The product is supercommutative by the very definition. 
One has, using (i), 
$$ 
N ((\alpha_1 \circ \alpha_2)\circ \alpha_3) = 
N (\alpha_1 \circ \alpha_2) \wedge d^{\C}\alpha_3 +(-1)^{ |{\alpha_3}| 
(|{\alpha_1}| + |{\alpha_2}|)} 
N \alpha_3 \wedge  d^{\C}(\alpha_1 \circ \alpha_2)= 
N (\alpha_1 \circ \alpha_2\circ \alpha_3). 
$$
Similarly $
N (\alpha_1 \circ \alpha_2\circ \alpha_3) = 
N (\alpha_1 \circ (\alpha_2\circ \alpha_3)). 
$ 
\end{proof}

\section{Twistor connections and variations of mixed Hodge structures} \la{hc4sec} 

Below we assume that $X$ is a compact complex manifold. 

All results of Section \ref{hc4sec} admit a rather straightforward 
 generalization when $X$ is a smooth complex variety, 
not necessarily projective. 
We use a regular compactification  of $X$  
with a normal crossing divisor at infinity 
and the standard technique to extend our results.

\bd A variation of mixed Hodge structures\footnote{What we call a variation of mixed Hodge structures is usually called a variation of \underline{complex} mixed Hodge structures.} over a 
complex manifold $X$ is  a complex local 
system $({\cal L}, \nabla)$ on $X$ 
equipped with three filtrations: $W_{\bullet}$, $F^{\bullet}$, $\overline F^{\bullet}$
 such that 
$$
{\rm gr}_{F}^p{\rm gr}_{\overline F}^q{\rm gr}^{W}_n=0~~\mbox{unless $p+q=n$}. 
$$
\be \la{2.7.15.100}
\nabla (F^p) \subset F^{p-1}\otimes \Omega^1_{X} \oplus F^{p}\otimes \overline \Omega^1_{X}. 
\ee
\be \la{2.7.15.101}
\nabla (\overline F^q) \subset \overline F^{q-1}\otimes 
\overline \Omega^1_{X} \oplus F^{q}\otimes \Omega^1_{X}. 
\ee
\ed
The condition  (\ref{2.7.15.100}) is the {\it Griffiths condition for $F^{\bullet}$}. 
It is equivalent to the Griffiths transversality 
 for $F^{\bullet}$ plus holomorphicity of $F^{\bullet}$. 
Similarly   (\ref{2.7.15.101}) is the {\it Griffiths condition for $\overline F^{\bullet}$}. 

\vskip 3mm
The  Galois group ${\rm Gal}(\C/\R)$ acts on the 
category ${\rm MHod_X}$ of variations of  mixed Hodge structures. Its generator acts by a functor 
$$
(X, {\cal L}, \nabla; W_{\bullet}, F^{\bullet}, \overline F^{\bullet}) \lra
 (\overline X, {\cal L}, \nabla; W_{\bullet}, \overline F^{\bullet}, F^{\bullet}). 
$$

The category ${\rm MHod^\R_X}$ is the category of ${\rm Gal}(\C/\R)$-equivariant objects in the category 
${\rm MHod_X}$.

\vskip 3mm
The 
category ${\rm Hod}_X$ of variations of  
Hodge structures on $X$  is a  semi-simple abelian  
tensor category. 
The 
category ${\rm MHod}_X$ of variations of mixed  Hodge structures on $X$ is an abelian  
tensor category. It is a mixed category, equipped with  a canonical tensor functor 
$$
{\rm gr}^W_{\bullet}: {\rm MHod}_X \lra {\rm Hod}_X; \quad 
 {\cal L} \lms {\rm gr}^W_{\bullet}{\cal L}.
$$ 
There is a canonical bigrading 
\begin{equation} \label{8.20.05.1}
{\rm gr}^W_{\bullet}{\cal L} = \oplus_{p,q}
({\rm gr}^W_{\bullet}{\cal L})^{p,q}. 
\end{equation}
Our goal is to equip the  object (\ref{8.20.05.1}) with 
an additional  {\it Green datum}, 
 which allows to recover the original 
variation of  mixed Hodge structures in a functorial way, compatible with the  ${\rm Gal}(\C/\R)$-action. This  
implies a similar description of variations of mixed $\R$-Hodge structures.

\subsection{Green data and twistor connections} 

\paragraph{The DG Lie algebra ${\cal C}_{\cal H}^{\bullet}({\rm End}{\cal V})$.}
 Let ${\cal L}$ be 
variation of Hodge structures  on $X$ with a Lie algebra structure. 
Then the Hodge complex ${\cal C}_{\cal H}^{\bullet}({\cal L})$ has a DG Lie 
algebra structure: the commutator is given by the composition
$$
{\cal C}^\bullet_{{\cal H}}({\cal L}) \otimes {\cal C}^\bullet_{{\cal H}}({\cal L}) \stackrel{\mbox{$\ast$-product}}{\lra}
{\cal C}^\bullet_{{\cal H}}({\cal L}\otimes {\cal L})
\stackrel{[\ast, \ast]}{\lra} {\cal C}^\bullet_{{\cal H}}({\cal L}). 
$$
To prove the Jacobi identity, observe that a triple commutator is induced by the composition
$$
{\cal C}^\bullet_{{\cal H}}({\cal L}) \otimes {\cal C}^\bullet_{{\cal H}}({\cal L}) \otimes 
{\cal C}^\bullet_{{\cal H}}({\cal L}) \stackrel{\ast}{\lra}
{\cal C}^\bullet_{{\cal H}}({\cal L}\otimes {\cal L}\otimes {\cal L})\stackrel{[[\ast,\ast], \ast]}{\lra} 
{\cal C}^\bullet_{{\cal H}}({\cal L}). 
$$
The Jacoby identity follows from three facts:

i) The $\ast$-product is associative.

ii) The Hodge complex ${\cal C}^\bullet_{{\cal H}}({\cal L})$ is functorial. 

iii) The commutator $[, ]$ is a morphism of Hodge structures. It satisfies 
the Jacoby identity.  
\vskip 3mm

Similarly, if ${\cal L}$ is a 
variation of Hodge structures  on $X$ with an algebra structure, then 
the Hodge complex ${\cal C}_{\cal H}^{\bullet}({\cal L})$ has a DG 
algebra structure. 
\vskip 3mm

Let ${\cal V}$ be a 
variation of Hodge structures  on $X$. 
 Let us introduce a DG  algebra 
of derived Hodge endomorphisms of ${\cal V}$. 
Let ${\rm End}{\cal V}$ be 
the local system of endomorphisms of the local system 
${\cal V}$. 
Since it is an algebra, 
the Hodge complex ${\cal C}_{\cal H}^{\bullet}({\rm End}\cal V)$ is a DG 
algebra.

\bd \la{11.17.ups.4} 
A {\rm Green datum} on a 
variation of Hodge structures ${\cal V}$  is a degree $1$
 element 
\be \la{11.17.ups.2}
{\cal G} \in {\cal C}_{\cal H}^1({\rm End}{\cal V}) 
\ee  
of the  DG Lie algebra ${\cal C}_{\cal H}^{\bullet}({\rm End}{\cal V})$ 
which satisfies a Maurer-Cartan equation 
\begin{equation} \label{10.14.05.1MC}
\delta {\cal G} + {\cal G}\wedge {\cal G} =0. 
\end{equation}
\ed

\paragraph{Explicit formulas.} 
Let us elaborate Definition \ref{11.17.ups.4}. 
The bigrading on ${\cal V}$ induces a bigrading 
on the local system ${\rm End}{\cal V}$: 
its bidegree $(-p, -q)$  component ${\rm End}^{-p,-q}{\cal V}$ consists 
of linear maps $A$ such that $A {\cal V}^{a,b} \subset {\cal V}^{a-p, b-q}$. 

\bl \la{11.10.ups.15} An element ${\cal G} \in {\cal C}_{\cal H}^1({\rm End}{\cal V})$  is given by 
the following datum:

\begin{itemize}
\item 
A collection of linear operators 
\begin{equation} \label{8.20.05.6a}
G_{s,t}\in {\rm End}^{-s,-t}{\cal V}_{\infty}, \quad s,t \geq 1,
\end{equation}
\item 
A closed $1$-form 
with values in ${\rm End}{\cal V}$ 
\begin{equation} \label{1fa}
g_{0,0} 
\in ~{\cal A}^{1, 0}\otimes {\rm End}^{-1, 0}{\cal V}_{\infty} ~\oplus ~
 {\cal A}^{0,1}\otimes {\rm End}^{0, -1}{\cal V}_{\infty}, 
\quad 
 {\bf d}g_{0,0} = 0. 
\end{equation}
\end{itemize}
\el
We call the operator $G = \sum_{s,t\geq 1} G_{s,t}$ the {\it Green operator}. 

\begin{proof} An element of 
${\cal C}_{\cal H}^{1}({\rm End}{\cal V})$ is given by 
a pair $\{G, g_{0,0}\}$, where $G$ is a 
function, $g_{0,0}$ is a 
$1$-form with the values ${\rm End}{\cal V}$. 
 Since $G$ is a function,  
its components have the Hodge bidegree $(-s,-t)$, where  $(s,t) \geq (1,1)$.  
So we recover conditions (\ref{8.20.05.6a}).

Since  $g_{0,0}$ is a $1$-form 
of degree $1$ in the complex ${\cal C}_{\cal H}^{\bullet}({\rm End}{\cal V})$, its 
 Hodge bidegree is $(0,0)$. So it is closed.   If 
a component  of the $1$-form $g_{0,0}$ lies in ${\cal A}^{a,b}\otimes {\rm End}^{-p, -q}{\cal V}$ then 
$(0,0) = (a-p, b-q)$. This imples   conditions (\ref{1fa}). 
Clearly these conditions are not only necessary, but also sufficient.  
\end{proof}

\paragraph{Twistor connections.} 
 Recall  the projection $\pi: X \times \C^2 \to X$. 
An element ${\cal G}\in {\cal C}_{\cal H}^1({\rm End}{\cal V})$  
gives rise, via the twistor transform $\gamma$, to a {connection} $\nabla_{\cal G}$ on 
$\pi^*{\cal V}_{\infty}$:
\begin{equation} \label{8.20.05.2zzz}
\begin{split}
&\nabla_{\cal G}:= {\bf d}+ d_{\C^2}+ \gamma({\cal G}) \\
&= {\bf d}+  d_{\C^2} + g_{0,0}+ 
\sum_{s,t\geq 0}z^{s}w^{t}\Bigl((s+t+1)G_{s+1,t+1}(zdw-wdz) + (w\partial' - z\partial'')G_{s+1,t+1} \Bigr).\\ 
\end{split}
\end{equation}
Recall the twistor line $\C \subset \C^2$, given by $z+w=1$. The crucial fact about it is the following:

\begin{proposition} \label{8.20.05.3}  
An  element ${\cal G}\in {\cal C}_{\cal H}^1({\rm End}{\cal V})$  on  a 
variation of Hodge structures ${\cal V}$ satisfies the Maurer-Cartan 
equation (\ref{10.14.05.1MC}) if and only if 
the  connection $\nabla_{\cal G}$ on $X \times \C$ is flat.
\end{proposition}

\begin{proof} Follows immediately from Theorem \ref{12.8.ups.1s}. \end{proof}

\bd \la{12.26.15.10} The connections $\nabla_{\cal G}$ on $X \times \C^2$ as in (\ref{8.20.05.2zzz}) 
such that the restriction 
${\rm Res}_{z+w=1}\nabla_{\cal G}$ is flat arew called twistor connection. 
\ed

The Green data and twistor connections are two ways to look at the same object. 

Tensor product of two flat twistor connections is again 
a flat twistor connection. Therefore flat 
twistor connections, or, equivalently, Green data,  form an abelian 
tensor category ${\cal G}_{\cal H}(X)$.    

Explicitly, the  objects of the category ${\cal G}_{\cal H}(X)$ 
are pairs  
$({\cal V}, {\cal G})$, 
where ${\cal V}$ is a 
variation  of Hodge structures on $X$ and ${\cal G}$ 
is a     Green datum on it. 
The morphisms are morphisms of variations of 
Hodge structures ${\cal V}_1\to {\cal V}_2$ 
commuting with the Green data. 
The tensor product is given by 
$
({\cal V}_1, {\cal G}_1)\otimes ({\cal V}_2, {\cal G}_2) = 
({\cal V}_1\otimes {\cal V}_2,  {\cal G}_1 \otimes {\rm Id}_{{\cal V}_2} +
{\rm Id}_{{\cal V}_1} \otimes {\cal G}_2).
$

\begin{lemma} \la{12.17.ups.10}
The   
category ${\cal G}_{\cal H}(X)$ 
is canonically equivalent to the category 
of comodules over the 
Lie coalgebra ${\cal L}_{{\cal H}; X}$ in the category ${\rm Hod}_X$. 
\el

We prove this lemma in Section \ref{Sec6.5} below. We will not use it before.

\begin{theorem} \label{8.20.05.7} One has a canonical  equivalence 
between the tensor category 
of variations of  mixed Hodge structures on $X$ 
and the tensor category ${\cal G}_{\cal H}(X)$ of twistor connections on $X$. 
\end{theorem}

\subsection{Deligne's description of triples of opposite filtrations} 

For the convenience of the reader we reproduce Deligne's approach to mixed (complex) Hodge structures, 
presenting it, following \cite{D2},  in a bit more general set-up. 

Let ${\cal A}$ be an abelian category. Let $A$ be an object of ${\cal A}$ with three finite {\it complimentary} filtrations 
$W_\bullet$, $F^\bullet$, $\overline F^\bullet$,  which means that 
$$
{\rm Gr}^p_F{\rm Gr}^q_{\overline F}{\rm Gr}_n^WA =0 ~~\mbox{if $n \not = p+q$}. 
$$
Here the $W_\bullet$ is an increaing filtration, and the other two are decreasing. 
This just means that  filtrations $F^\bullet$, $\overline F^\bullet$ induce $n$-opposite filtrations 
on the 
 ${\rm Gr}_n^WA$, providing a decomposition
$$
{\rm Gr}_n^WA = \bigoplus_{p+q=n}A_W^{p,q}
$$
such that the filtrations on ${\rm Gr}_n^WA$ induced by the ones $F^\bullet$ and $\overline F^\bullet$ 
are given by
$$
F^i{\rm Gr}_n^WA = \bigoplus_{p+q=n, ~p\geq i}A_W^{p,q}, ~~~~
\overline F^i{\rm Gr}_n^WA = \bigoplus_{p+q=n, ~q\geq i }A_W^{p,q}.
$$
Equivalently, one has $A_W^{p,q} = F^p{\rm Gr}_n^WA \cap \overline F^q{\rm Gr}_n^WA$. 

Any bigraded object $B = \bigoplus_{p+q}B^{p,q}$ has three canonical 
filtrations $W$, $F_W$, and $\overline F_W$: 
$$
W_n:= \bigoplus_{p+q\leq n}B^{p,q}, ~~~~
F^i_W:= \bigoplus_{p \geq i}B^{p,q}, ~~~~ \overline F^j_W:= \bigoplus_{q \geq j}B^{p,q}, ~~~~ 
$$
The filtrations $F, \overline F$ induce on the bigraded object 
${\rm Gr}^WA = \bigoplus_{p, q}A_W^{p,q}$ the filtrations $F_W, \overline F_W$. 

Deligne considers the following subobjects of the object $A$:
\be
\begin{split}
&A_F^{p,q}:= (W_{p+q}\cap F^p) \cap \Bigl((W_{p+q}\cap \overline F^q) + 
\sum_{i \geq 0}(W_{p+q-i} \cap  \overline F^{q-i+1})\Bigr),\\
&A_{\overline F}^{p,q}:= (W_{p+q}\cap  \overline F^q) \cap \Bigl((W_{p+q}\cap F^p) + \sum_{i \geq 0}(W_{p+q-i} \cap  F^{p-i+1})\Bigr).\\
\end{split}
\ee
It is proved in \cite[1.2.8]{D1},  that if  $n=p+q$ then the canonical projection $W_nA \to {\rm Gr}^W_nA$ induces 
an isomorphism of $A^{pq}_F \stackrel{\sim}{\to} A^{pq}_W$.\footnote{Let us scetch the argument when
 ${\cal A}$ is a category of finite dimensional vector spaces. One has 
$$
W_{p+q-1} = W_{p+q-1}\cap F^p \oplus ( 
\sum_{i \geq 1}(W_{p+q-i} \cap  \overline F^{q-i+1}).
$$
This implies that ${\rm dim}A_F^{p,q} = {\rm dim}A_W^{p,q}$, and 
the map  $A^{pq}_F \stackrel{}{\to} A^{pq}_W$ 
is injective.  Therefore it is an isomorphism. } 
This implies that the subobjects $\{A^{pq}_F\}$ provide a bigrading of $A$. 
Therefore the sum of isomorphisms $a_F: A_F^{p,q} \stackrel{\sim}{\to}A_W^{p,q}$ gives rise to a bifiltered 
isomorphism
$$
a_F: (A; W, F) \lra ({\rm Gr}^WA; W, F_W).  
$$
Similarly  the subobjects $\{A^{pq}_{\overline F}\}$ give rise to a bifiltered 
isomorphism
$$
a_{\overline F}: (A; W, \overline F) \lra ({\rm Gr}^WA; W, \overline F_W).  
$$
Then the automophism ${\rm d}=  a_{\overline F}a^{-1}_{F}$ of ${\rm Gr}^WA$ satisfies 
\be \la{autd}
({\rm d}-1)(A^{pq}_W) \subset \bigoplus_{r<p, s<q}A^{rs}_W.
\ee
\bt\la{autd1}
The functor $A \lra {\rm Gr}^WA$ provides an equivalence of the category of 
objects of ${\cal A}$ equipped with three 
complimentary filtrations $(W_\bullet, F^\bullet, \overline F^\bullet)$ with the category of 
bigraded objects of ${\cal A}$ 
equipped with an automorphism ${\rm d}$ satisfying (\ref{autd}). 
\et

\subsection{Proof  of Theorem \ref{8.20.05.7}} 
\paragraph{1. Constructing a tensor functor ${\cal F}$.} Let us define a functor  
$$
{\cal F}:  {\cal G}_{{\cal H}}(X) \lra 
{\rm MHod}_X.
$$ 
Recall a twistor connection
 $\nabla_{\cal G}$ on $\pi^*{\cal V}_{\infty}$: 
\begin{equation} \label{8.20.05.2}
\nabla_{\cal G}= {\bf d}+ {\rm d}_{\C^2}+ g_{0,0}+ 
\sum_{s,t\geq 0}z^{s}w^{t}\Bigl((s+t+1)G_{s+1,t+1}(zdw-wdz) + (w\partial'-z\partial'')G_{s+1,t+1}\Bigr). 
\end{equation}
Let us define  
a variation of mixed Hodge structures on $X$ starting from the 
 connection  $\nabla_{\cal G}$.

\vskip 3mm 
(i) {\it The flat connection}. Restricting the 
 connection $\nabla_{\cal G}$  
to $X \times \{\frac{1}{2}, \frac{1}{2}\}$ we get 
a flat connection $\nabla^{\{\frac{1}{2}\}}_{\cal G}$ on the vector bundle ${\cal V}_{\infty}$ 
on  $X$. 
 
\vskip 3mm 
(ii) {\it The weight filtration $W_\bullet$}. It is the standard weight filtration on a bigraded object: 
$$
W_n{\cal V}_\C:= \oplus_{p+q\leq n}{\cal V}^{p,q}.
$$

The connection 
$\nabla^{\{\frac{1}{2}\}}_{\cal G}$ 
preserves the weight filtration. 
\vskip 3mm 
(iii) {\it The Hodge filtrations $F^\bullet$ and $\overline F^\bullet$}. 
Consider the standard Hodge filtration 
on the restriction of $\pi^*{\cal V}_{\infty}$ to $X \times \{0,1\}$: 
\be \la{st}
F^p_{\rm st}:= \oplus_{i\geq p}{\cal V}^{i, *}.
\ee
Consider the standard conjugate Hodge filtration 
on the restriction of $\pi^*{\cal V}_{\infty}$ to $X \times \{1, 0\}$:
\be \la{st1}
\overline F^q_{\rm st}:= \oplus_{j\geq q}{\cal V}^{\ast, j}.
\ee
Since the restriction of the twistor connection to 
the coordinate cross $zw=0$ is trivial, the weight filtration and the standard Hodge filtrations 
are extended canonically to the fiber over the coordinate cross. 

Let $P$ be the operator of parallel transport 
for the connection $\nabla_{\cal G}$ along the twistor 
line $z+w=1$ from $X \times \{0,1\}$ to $X \times \{\frac{1}{2}, \frac{1}{2}\}$, and $\overline P$ 
a similar  one  
from $X \times \{1, 0\}$ to $X \times \{\frac{1}{2}, \frac{1}{2}\}$. 

\bd The Hodge filtrations $F^{\bullet}$ and $\overline F^{\bullet}$ on 
the ${\cal V}$ are 
the images  of the corresponding standard Hodge filtrations by the  operators $P$ and $\overline P$: 
\be \la{stho}
F^{p} := P(F^p_{\rm st}), ~~~~
\overline F^{q} := \overline P(\overline F^q_{\rm st}).
\ee
\ed

\bl
The data $({\cal V}, \nabla^{\{\frac{1}{2}\}}_{\cal G}; W_{\bullet}, 
F^{\bullet}, \overline F^{\bullet})$ is a variation of mixed Hodge structures on $X$. 
\el

\begin{proof} The only thing we have to prove is that the Hodge filtrations satisfy 
the Griffiths   condition. 
We employ the following general observation. 

\bl \la{12.20.ups.1} 
Let $(E, \nabla)$ be a smooth flat bundle 
on $X \times \R$, and ${\cal F}^{\bullet}$ a  
filtration on 
$E$ invariant under the parallel transport 
along the $\R$-lines. Then if 
the restriction of $(E, {\cal F}^{\bullet}, \nabla)$ 
to $X \times \{s\}$ 
satisfies the Griffiths   
condition for a single $s \in \R$,  
the same is true 
for any $s$.  
\el

We apply Lemma \ref{12.20.ups.1} to the bundle $E= \pi^*{\cal V}_{\infty}$ 
with the connection $\nabla_{\cal G}$ restricted to the  line $z+w=1$ 
and  filtrations ${\cal F}^{\bullet}$ and $\overline {\cal F}^{\bullet}$ 
on $X$ 
given by  Hodge filtrations 
$F^{\bullet}$ and $\overline F^{\bullet}$ in (\ref{stho}). 

The standard Hodge filtration (\ref{st}) on the restriction $\nabla^{\{0,1\}}_{\cal G}$ of the connection 
$\nabla_{\cal G}$  to $X \times \{0,1\}$ satisfies the Griffiths   
condition. This just means 
that  the Hodge bidegree $(-s,-t)$ component of $\nabla^{\{0,1\}}_{\cal G}$ is 
zero if   $s>0$. To check this notice that in the restriction  
to $X \times \{0,1\}$ monomials $z^sw^t$ with $s>0$ vanish.
 
Similarly the Griffiths   condition for $\overline F^\bullet$ is 
clear for the restriction to $X \times \{1,0\}$.  
 
Since the connection $\nabla_{\cal G}$ is flat over $X \times \{z+w=1\}$, by Lemma \ref{12.20.ups.1} 
the filtrations $F^\bullet$ and  $\overline F^\bullet$ in (\ref{stho}) at $\{\frac{1}{2}, \frac{1}{2}\}$ 
 satisfy the Griffiths  
condition. 
\end{proof}

\paragraph{2. The functor ${\cal F}$ is an equivalence.}
The functor ${\cal F}$ is clearly a tensor functor. 
To show that the functor ${\cal F}$ 
is an equivalence of categories we have to show that it
is an isomorphism on ${\rm Hom}$'s, and 
every variation  of mixed Hodge structures is  
obtained this way. 
The former claim is clear since it is
 known when $X$ is a point. So it remains to prove the latter. 

\begin{figure}[ht]
\centerline{\epsfbox{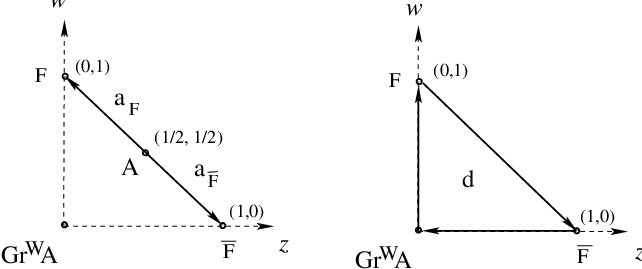}}
\caption{Deligne's operators $a_F$, $a_{\overline F}$, ${\rm d}$ via the homogeneous twistor connection}
\label{hc50}
\end{figure}

Fig \ref{hc50} tells how Deligne's operators $a_F$, $a_{\overline F}$, ${\rm d}$ 
are expressed via the twistor connection: 
\newline 
The object $A$ is the fiber over $\{\frac{1}{2}, \frac{1}{2}\}$. 
\newline 
The associate graded ${\rm Gr}^WA$ is the fiber over $(0,0)$. 
Since the restriction of the connection to the 
coordinate cross $zw=0$  is zero, the ${\rm Gr}^WA$  is  identified 
the fiber over the coordinate cross. 
\newline
 The operator $a_F: A \to {\rm Gr}^WA$ is the parallel transport from 
$\{\frac{1}{2}, \frac{1}{2}\} \to \{0, 1\}$.  
\newline
The operator $a_{\overline F}: A \to {\rm Gr}^WA$ is the parallel transport from 
$\{\frac{1}{2}, \frac{1}{2}\} \to \{1,0\}$.
\newline
Deligne's operator ${\rm d}:=a_{\overline F}a^{-1}_F: {\rm Gr}^WA \to {\rm Gr}^WA$ is 
the clockwise monodromy around the triangle 
$\{0,0\} \to\{0,1\}\to \{1,0\}\to \{0,0\}$. 
\newline Identifying  ${\rm Gr}^WA$  with 
the fiber over  $zw=0$, the  ${\rm d}$
 is the parallel transport  $\{0,1\} \to \{1,0\}$.

\bp \la{12.25.15.5} For any variation ${\cal L}$ of  MHS  
there  is a unique smooth family of Green operators $G_x$ on 
${\rm gr}^W{\cal L}_x$ 
describing the MHS's fiberwise. 
\ep

\begin{proof}
Over a point 
Deligne's operator ${\rm d}$ is recovered from $G$, 
and vice versa. Indeed, let ${\Bbb P}$ be the operator of parallel transport from $\{1,0\}$ to $\{0,1\}$. Then 
$$
{\rm d} = {\Bbb P}.
$$ 
The operator ${\Bbb P}$ is given explicitly 
by 
\be \la{P}
{\Bbb P} = 1 + \sum_{n>0}\Bigl(\int_{0}^1\varphi_{p_1, q_1} \circ \ldots 
\circ \varphi_{p_n, q_n} \Bigr) 
G_{p_1+1, q_1+1}\ldots 
 G_{p_n+1, q_n+1}.
\ee
Here $\varphi_{p,q}:= (p+q+1)(1-u)^{p}u^{q}du$, and 
$\int_{0}^1\varphi_{p_1, q_1} \circ \ldots \circ \varphi_{p_n, q_n}$ is a positive rational number given by 
the iterated integral of the 
$1$-forms $\varphi_{p_1, q_1},  \ldots , \varphi_{p_n, q_n}$. 
Therefore we have 
$$
g_{p,q} := \log({\rm d}) = \log({\Bbb P})_{p,q} = \int_{0}^1\varphi_{p-1, q-1} G_{p,q} + \ldots 
$$
where $\ldots$ denotes a polynomial  with rational 
coefficients of 
degrees $\geq 2$ in $G_{*,*}$, of the Hodge bidegree $(-p,-q)$. 
By the induction $G_{p,q}$ is recovered  from $g_{*,*}$. 
\end{proof}

Thus there is   a unique isomorphisms of fibrations 
$$
\eta: {\rm gr}^W{\cal L} \stackrel{\sim}{\to} {\cal L}
$$
 such that ${\rm gr}^W\eta: 
{\rm gr}^W{\cal L} \stackrel{}{\to} {\rm gr}^W{\cal L}$ is the identity map, and 
which identifies pointwise   
the MHS provided by the pair $({\rm gr}^W{\cal L}_x, G_x)$ with the  one 
 ${\cal L}_x$. 

Let ${\bf d} = \partial' + \partial''$ be the connection 
on ${\rm gr}^W{\cal L}$ induced by the connection on ${\cal L}$. 
The isomorphism $\eta^{-1}$ 
transforms the connection on ${\cal L}$ to another connection $\nabla$ on 
${\rm gr}^W{\cal L}$. 
We have to show that:

\vskip 2mm
{\it The connection on ${\cal L}$ 
determines uniquely a $1$-form $g_{0,0}$ as in (\ref{1fa}). 

The Griffiths   conditions for
  $F^\bullet$ and $\overline F^\bullet$ 
are equivalent to the following two conditions:

1. One has \be \la{11.5.ups.1a}
\nabla = \nabla^{\frac{1}{2}}_{\cal G}.      
\ee

2. The operators $G_x$ and the $1$-form $g_{0,0}$, organised 
into a Green data ${\cal G} = (G, g_{0,0})$,  
satisfy 
\be \la{11.5.ups.1}
\delta{\cal G} + {\cal G}\wedge {\cal G} =0.
\ee}

Write the connection $\nabla$  as follows, where 
the ${\cal B}_{s,t}$ is a   $1$-form 
of the  Hodge bidegree $(-s,-t)$:
$$
\nabla = {\bf d} +{\cal B}; ~~~~ {\cal B}= \sum_{s,t}2^{-(s+t)}{\cal B}_{s,t}.
$$ 

Consider  the  third quadrant 
in  the $(s,t)$-plane,  shown on 
Fig. \ref{hc11}:
$$
B= \{(s,t)| s,t \leq 0\}.
$$

\begin{figure}[ht]
\centerline{\epsfbox{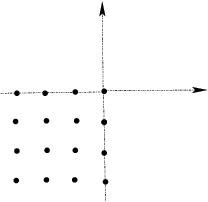}}
\caption{The Hodge bidegrees of the components of the $1$-form ${\cal B}$ must be in the domain $B$.}
\label{hc11}
\end{figure}

\begin{lemma} \la{12.22.ups.4}
${\cal B}_{s,t}=0$  unless $(s,t)\in -B$. 
\el

\begin{proof}
Since the connection $\nabla$ preserves the weight filtration, 
its components of positive weights are zero. 
The weight zero component 
of $\nabla$  is the 
connection ${\bf d}$ on ${\rm gr}^W{\cal L}$.

The Griffiths   condition for the filtrations $F^\bullet$ and
 $\overline F^\bullet$ 
on the 
bundle with connection ${\bf d} +{\cal B}$ 
just mean that 
\be \la{12.22.ups.3}
(P^{-1}{\bf d} P  + P ^{-1}{\cal B} P )_{s,t} =0 \quad \mbox{if $s>0$}. 
\ee
\be \la{12.22.ups.3h}
(\overline P^{-1}{\bf d} \overline P  + \overline P^{-1}{\cal B} \overline P )_{s,t} =0 \quad \mbox{if $t>0$}. 
\ee

It is clear from (\ref{P}) that: 

a) Operators $P - {\rm Id}$ and $P ^{-1}- {\rm Id}$  are sums  
of the terms with  
$s, t\geq 1$. 

b) The components of $P ^{-1} {\bf d}P$ are inside of the domain $-B$.

This plus similar arguments for the filtration $\overline F^\bullet$ using (\ref{12.22.ups.3h}) imply the   Lemma. 
\end{proof}

Let us prove   (\ref{11.5.ups.1a}).    Set  $g_{0,0} := {\cal B}_{0,0}$. The
 Hodge bidegree $(0,0)$ case  follows from
   Lemma \ref{12.22.ups.4}. 

Consider the following cousin of the twistor  connection on $X\times \C$: 
\begin{equation} \label{8.20.05.2asa}
\nabla_{{\cal B}, {\cal G}} = {\bf d}+{ d}_{\C}+ 
\sum_{s,t\geq 0}(1-u)^{s}u^{t}\Bigl((s+t+1)G_{s+1,t+1}du + {\cal B}_{s,t}
  \Bigr). 
\end{equation}
Its curvature  has 
two components: the horizontal one, and the vertical one, which contains $du$.  
\begin{proposition} \la{12.22.ups.1}
The $du$-part of the curvature of the connection 
$\nabla_{{\cal B}, {\cal G}}$ is zero if and only if 
the Griffiths   condition for filtrations $F^\bullet$  and $\overline F^\bullet$  
for the connection ${\bf d} +{\cal B}$ 
hold. 
\ep

\begin{proof} Since the vertical components of  connections 
$\nabla_{{\cal B}, {\cal G}}$ and $\nabla_{{\cal G}}$ 
are the same,  the operator $P$ is also the operator of  
parallel transform from $\{0,1\}$ to $\{\frac{1}{2}, \frac{1}{2}\}$ 
for  the connection 
 $\nabla_{{\cal B}, {\cal G}}$.   
The connection ${\bf d} + {\cal B}$ is the restriction of 
$\nabla_{{\cal B}, {\cal G}}$ to $X \times \{\frac{1}{2}, \frac{1}{2}\}$. 
The value of $P^{-1} ({\bf d} + {\cal B})P$ on a tangent vector $v$ at 
$(x; 0,1)$ is 
a linear endomorphism of the fiber of $\pi^*{\cal V}$ at $(x; 0,1)$. It 
is  
 the composition of the 
parallel transform from $\{0,1\}$ to $\{\frac{1}{2}, \frac{1}{2}\}$, followed by the infinitesimal 
parallel transform along $v$,  and then by the 
parallel transform back from $\{\frac{1}{2}, \frac{1}{2}\}$ to $\{0,1\}$ for the connection 
$\nabla_{{\cal B}, {\cal G}}$. 
So it is the parallel transport 
for the connection $\nabla_{{\cal B}, {\cal G}}$ 
along the three sides of an infinitesimal rectangle $R$
 with the sides $v$ and $x \times [0,1/2]$. Thus it differs 
from the value of the connection $\nabla_{{\cal B}, {\cal G}}$ on $v$  
by the integral 
of the curvature over $R$. 

The restriction of the connection $\nabla_{{\cal B}, {\cal G}}$ 
to $X \times \{0,1\}$ (respectively $X \times \{1, 0\}$) 
satisfies the Griffiths condition for $F^\bullet$  (respectively $\overline F^\bullet$) 
since monomials $z^sw^t$ with $s>0$ 
(respectively $t>0$) restrict to zero. 
Therefore using (\ref{12.22.ups.3}) (respectively  (\ref{12.22.ups.3h})) we see that  the 
Griffiths condition for $F^\bullet$  (respectively $\overline F^\bullet$) is  equivalent 
to the condition that the integral of the $du$-part of the curvature 
$(\nabla_{{\cal B}, {\cal G}})^2$ over the interval $ \frac{1}{2} \leq u \leq 1$ (respectively $ 0 \leq u \leq \frac{1}{2} $)
has zero $(s,t)$-components  if $s>0$ (respectively $t>0$). 

The argument below is reminiscent of the proof 2) of Theorem \ref{12.8.ups.1ss}ii). 

The $du$-part of the curvature 
$(\nabla_{{\cal B}, {\cal G}})^2$ is given by
\be \la{2.7.15.31}
\sum_{s,t\geq 0}(1-u)^su^tdu \wedge 
\Bigl(-(s+1){\cal B}_{s+1, t} + (t+1) {\cal B}_{s, t+1} - (s+t+1) {\bf d}{G}_{s+1, t+1} 
\ee
\be \la{2.7.15.21dfd}
+ \sum_{s'+s''=s, ~t'+t''=t}(s'+t'+1) G_{s'+1, t'+1}  {\cal B}_{s'', t''} \Bigr). 
\ee
We set $G_{s,t}=0$ if $s< 1$ or $t< 1$,  
and rewrite this as a sum of $(s,t)$-components: 
\be \la{2.10.15.1}
(\ref{2.7.15.31}) + (\ref{2.7.15.21dfd}) = \sum_{s,t}(1-u)^{s-1}u^{t-1}du \wedge (uP_{s,t}+Q_{s,t}),
\ee
where we set 
\be \la{2.7.15.20ab}
\begin{split}
&P_{s,t}:= (s+t)(-{\cal B}_{s,t} + \partial' G_{s+1,t} -  \partial'' G_{s,t+1}), \\
&Q_{s,t}:= t{\cal B}_{s,t} 
-(s+t) \partial' G_{s+1,t} 
+\sum_{s'+s''=s, ~t'+t''=t}(s'+t'-1) G_{s', t'}  {\cal B}_{s'', t''}. \\
\end{split}
\ee

Now we have two conditions: 

$F$: For each $(s,t)$ such that $s>0$, the integral of  (\ref{2.10.15.1}) over  
$\frac{1}{2} \leq u \leq 1$ is zero. 

$\overline F$: For each $(s,t)$ such that $t>0$, 
the integral of  (\ref{2.10.15.1}) over $0\leq u \leq \frac{1}{2} $  is  zero. 

\vskip 3mm
Let us assume first that $s,t>0$. Set 
$
a_{s,t}:= \int_{1/2}^1(1-u)^{s-1}u^{t-1}du>0.
$ 
Then the  conditions $F$ and $\overline F$ are equivalent to 
 the following equations on  $P_{s,t}, Q_{s,t}$: 
\be \la{12.15.25.11}
\begin{split}
&a_{s,t+1}P_{s,t} + a_{s,t}Q_{s,t} =0, \\
&a_{t+1,s}P_{s,t} + a_{t,s}Q_{s,t} =0.\\
\end{split}
\ee
Observe that  $a_{s,t+1}/a_{s,t}>1/2$. Indeed, the integrand $(1-u)^{s-1}u^{t}$ of 
$a_{s,t+1}$ is obtained by multiplying 
the integrand $(1-u)^{s-1}u^{t-1}$ of $a_{s,t}$ by $u$, and $u>1/2$ on the integration interval.  
Similarly, $a_{t+1,s}/a_{t,s}<1/2$. Therefore system (\ref{12.15.25.11}) has a unique solution: 
\be \la{12.25.15.1}
P_{s,t} = Q_{s,t}=0~~\mbox{ if $s,t>0$}.
\ee
If $t=0$, then $Q_{s,0}=0$. If $s=0$, then $P_{0,t}=-Q_{0,t}$. So 
the claim is  evident in both cases. 
 \end{proof}

Thanks to the Griffiths   conditions  
for the connection ${\bf d} +{\cal B}$ and Proposition \ref{12.22.ups.1}, 
the $du$ part of the curvature of 
$\nabla_{{\cal B}, {\cal G}}$ is zero. 
In particular we have  $P_{s,t} = 0$, which just means that
\be \la{2.7.15.20a}
{\cal B}_{s,t} = 
\partial' G_{s+1,t} 
-\partial'' G_{s,t+1} ~~~~\mbox{where}~s>0 ~\mbox{or}~ t>0. 
\ee
This plus $g_{0,0} = {\cal B}_{0,0}$ proves (\ref{11.5.ups.1a}). 
The horisontal part of the connection is zero by the very definition.  
Since the $du$ part of the curvature of 
$\nabla_{{\cal B}, {\cal G}}$ is zero, 
the connection $\nabla_{{\cal B}, {\cal G}}$ is flat. 
Then (\ref{11.5.ups.1}) follows from Proposition \ref{8.20.05.3}. 
Theorem \ref{8.20.05.7} is proved. 

\paragraph{Explicit differential equations.}
Let $\psi_{s,t} = 2^{-s-t}({\bf d}^\C G)_{s,t}$ if $s+t >0$, and  $\psi_{0,0}= g_{0,0}$. 
We decompose $\psi_{s,t}$ into 
the $(1,0)$ and $(0,1)$-components: $\psi_{s,t} = \psi_{s,t}'+\psi_{s,t}''$. 
Then we have: 

\begin{equation} \label{10.14.05.1}
\begin{split}
&s~ {\partial}'G_{s+1, t} = 
\sum_{s'+s''=s, ~~t'+t''=t}(s'+t'-1)\left[G_{s', t'}, \psi'_{s'', t''}\right ], \\
&t~{\partial}''G_{s,t+1} =  
\sum_{s'+s''=s, ~~t'+t''=t}(s'+t'-1)\left[G_{s', t'}, 
\psi''_{s'', t''}\right],\\
&{\partial}\overline {\partial}G_{1,1} =  g_{0,0}\wedge g_{0,0} \ .\\
\end{split}
\end{equation}
Indeed, substituting (\ref{2.7.15.20a}) and ${\cal B}_{0,0} = g_{0,0}$ to (\ref{2.7.15.20ab}) and using 
$Q_{s,t}=0$, see (\ref{12.25.15.1}), we get  the first two differential equations. 
The last one follows from the flatness of the horisontal part of the connection. 
Notice that equations (\ref{10.14.05.1}) impose no 
conditions on 
$\partial' G_{1,t}$ and $\partial'' G_{s,1}$ for $s,t >1$.

\bl \la{11.18.ups.1}
Maurer-Cartan equation (\ref{10.14.05.1MC}) 
is equivalent to differential equations 
(\ref{10.14.05.1}).
\el

\begin{proof} Differential equations (\ref{10.14.05.1}) just mean that the twistor connection 
$\nabla_{\cal G}$ is flat. By Proposition \ref{8.20.05.3}, this is  
equivalent to Maurer-Cartan equation $\delta {\cal G}+ {\cal G}\wedge {\cal G}=0$, see
 (\ref{10.14.05.1MC}). \end{proof}

\subsection{Conclusions} 

\bd \la{TCDEF}
Given a complex manifold $X$, a twistor connection is a 

\begin{itemize}

\item
$\C^*\times \C^*$-equivariant connections 
on $X\times \C^2$ whose restriction to $X\times \{z+w=1\}$ is flat. 
\end{itemize}
\ed

Definition \ref{TCDEF} is equivalent to Definition \ref{12.26.15.10} thanks to the following Lemma.

\bl \la{12.31.15.5}
A twistor connection has a unique representative in its gauge class, given by 
\begin{equation} \label{8.20.05.2as}
\nabla_{\cal G}= {\bf d}+  g_{0,0}+ 
\sum_{s,t\geq 0}z^{s}w^{t}\Bigl((s+t+1)G_{s+1,t+1}(zdw-wdz) + (w\partial' - z\partial'')G_{s+1,t+1} \Bigr). 
\end{equation}
\el
\begin{proof} A $\C^*\times \C^*$-equivariant connection looks as follows: 
$$
\nabla = {\bf d} + \sum_{s,t\geq 0}\Bigl(z^{s-1}w^t A'_{s,t}dz + 
z^{s}w^{t-1} A''_{s,t}dw\Bigr) + \sum_{s,t\geq 0}z^sw^t{\cal B}_{s,t}.  
$$
Using a gauge transformation in the fibers $\C^2$, it can be uniquely transformed 
into 
$$
\nabla = {\bf d} +  \sum_{s,t\geq 0}z^{s-1}w^{t-1} (s+t-1)G_{s,t}(zdw - wdz) + 
\sum_{s,t\geq 0}z^{s}w^{t}{\cal B}_{s,t}.  
$$
The $du$-part of the curvature of its restriction to $X \times \C$ 
is given by 
\be \la{12.31.15.1}
\sum_{s,t}(1-u)^{s-1}u^{t-1}du \wedge (uP_{s,t}+Q_{s,t}),
\ee
where 
the $P_{s,t}$ and $Q_{s,t}$ are given by (\ref{2.7.15.20ab}).
Just as in the proof of Lemma \ref{12.22.15.1}, this implies that $P_{s,t}=0$. 
Indeed, the $P_{s,t}$ is the leading coefficient of the polynomial in $u$ giving the 
$(-s, -t)$-part. This immediately implies that $Q_{s,t}=0$. The $P_{s,t}=0$ just means that 
$$
\sum_{s,t\geq 0}z^sw^t{\cal B}_{s,t} = {\cal B}_{0,0}+ \sum_{s,t\geq 0}z^sw^t
 (w\partial' - z\partial'')G_{s+1,t+1}. 
$$
 It remains to set $g_{0,0}= {\cal B}_{0,0}$. 
\end{proof} 

Using Lemma \ref{12.31.15.5} one can restate Theorem 
\ref{8.20.05.7} as follows. 

\bt
The tensor category of variations of mixed Hodge structures on a 
complex manifold $X$ is equivalent to the tensor category 
of twistor connections on $X$. 
\et

Deligne's operator ${\rm d}$ is the monodromy of the twistor connection  
over the triangle with vertices $\{0,0\}, \{0,1\}, \{1,0\}$ in the twistor plane $\C^2$, 
while the Green operator $(s+t+1)G_{s+1, t+1}$ is the 
$(s+1,t+1)$-homogeneous component of the  curvature form in the twistor plane.

\subsection{Variations of real mixed Hodge structures.} 

The category of variations of real mixed Hodge structures as 
the category of ${\rm Gal}(\C/\R)$-equivariant objects 
in the category of variations of mixed Hodge structures. 
Here is an elaborate description.

\begin{definition} 
A variation of mixed $\R$-Hodge structures over a 
complex manifold $X$ is given by a real local 
system $({\cal L}, \nabla)$ on $X$ 
equipped with the following data: 

\begin{itemize}

\item
An  increasing weight filtration $W_{\bullet}{\cal L}$ 
on the local system  ${\cal L}$,

\item A decreasing  Hodge filtration 
$F^{\bullet}$ on the holomorphic vector bundle  
${\cal L}_{\cal O}:= {\cal L}\otimes_{\R}{\cal O}_{X}$, 
\end{itemize}

satisfying the following conditions: 

\begin{itemize}

\item The two filtrations induce at each fiber a mixed $\R$-Hodge structure, 

\item The Griffiths transversality condition: 
$
\nabla^{1,0} (F^p) \subset F^{p-1}\otimes \Omega^1_{M}.
$
\end{itemize}
\end{definition}

A {\rm Green datum} on a 
variation of real Hodge structures ${\cal V}$  is a Maurer-Cartan 
 element:  
\be \la{11.17.ups.2gR}
{\cal G} \in {\cal C}_{{\cal H}, \R}^1({\rm End}{\cal V}), ~~~~ 
\delta {\cal G} + {\cal G}\wedge {\cal G} =0. 
\end{equation}
Explicitly, ${\cal G}$ is given by the data $(G_{p,q}, g_{0,0})$ as in 
Lemma \ref{11.10.ups.15}, satisfying the reality condition: 
$$
\overline G_{p,q} = - G_{p,q}, ~~\mbox{i.e.} ~~\overline G = -G, ~~~~~~\mbox{and} 
~~\overline {g_{0,0}} = g_{0,0}.
$$

Equivalently, a {\it real twistor connection} is a twistor connection which is 
invariant under the involution $c\circ \sigma$, where 
$\sigma: (z,w) \lms (\overline w, \overline z)$, and $c$ is the complex conjugation. 

The category of real twistor connections on $X$ 
is equivalent to the category variations of real mixed Hodge structures. 

 The complex conjugation interchanges the first two differential 
equations in (\ref{10.14.05.1}). 

\subsection{DG generalizations} \la{Sec6.5}

Recall that  
a DG-module over a DG Lie algebra  is 
a module over a Lie 
algebra in the category of complexes rather then vector spaces. 
So a DG module over a 
DG Lie algebra  $(L_{\bullet}, \partial)$ 
is a complex $M_{\bullet}$ plus an action 
of the 
graded Lie algebra $L_{\bullet}$ on the graded object 
$M_{\bullet}$ given by a map of complexes  
$\mu: L_{\bullet}\otimes M_{\bullet} \to M_{\bullet}$. 

Given two DG-modules $M_{\bullet}$ and $N_{\bullet}$, 
${\rm Hom}^{\bullet}_{L_{\bullet}}(M_{\bullet}, N_{\bullet})$ is a complex, consisting 
of elements of the complex ${\rm Hom}^{\bullet}(M_{\bullet}, N_{\bullet})$ 
commuting with the action of 
$L_{\bullet}$. Taking $H^0$'s of these ${\rm Hom}$-complexes   
we arrive at the homotopy category of $DG$-modules. Here 
the  morphisms 
are morphisms of $L_{\bullet}$-modules $M_{\bullet} \to N_{\bullet}$ 
commuting with the differentials in $M_{\bullet}$ and $N_{\bullet}$, 
up to homotopies given 
by morphisms of $L_{\bullet}$-modules $M_{\bullet} \to N_{\bullet}[-1]$. 
Localizing by the quasiisomorphisms, we 
get the derived category ${\cal D}(L_{\bullet})$. 
One easily translates this to the world of DG Lie coalgebras.

\vskip 3mm
Recall the DG Lie coalgebra ${\cal L}_{{\cal H}; X}^*$ in the pure category 
${\rm Hod}_X$ defined in Section \ref{hc1.6}. 
\bd
The DG Hodge category ${\cal G}^*_{\cal H}(X)$ of a regular 
complex projective variety $X$ is the derived category 
${\cal D}({\cal L}^*_{{\cal H}; X})$ 
of DG comodules over the DG coalgebra 
${\cal L}_{{\cal H}; X}^*$ in the pure category 
${\rm Hod}_X$. 
\ed
The precise form of Conjecture \ref{12.1.ups.7} is 
\bcon \la{12.1.ups.7pr}
The category of smooth 
complexes of 
real Hodge sheaves on $X$ 
is equivalent to the derived category ${\cal D}({\cal L}^*_{{\cal H}; X})$.
\econ
\vskip 3mm

Let us  describe  DG comodules over the 
DG Lie coalgebra ${\cal L}_{{\cal H}; X}^*$. 
Consider the following data: 

\begin{itemize} 
\item A bounded complex ${\cal V} =  \{ 
 \ldots \lra {\cal V}^1 \stackrel{\partial}{\lra} {\cal V}^2 
\stackrel{\partial}{\lra} 
 {\cal V}^3 \stackrel{\partial}{\lra}  {\cal V}^4 
\stackrel{\partial}{\lra} \ldots\} 
$   in ${\rm Hod}_X$.

\item A Maurer-Cartan element:
\begin{equation} \label{3/29/07/1}
{\cal G} \in {\cal C}^1_{\cal H}({\rm End}{\cal V}) 
\quad \mbox{such that} \quad 
\delta{\cal G} + {\cal G}\wedge {\cal G} =0,  
\end{equation}
defined modulo elements of type 
$\delta B + [B,{\cal G}]$, where 
$B \in  {\cal C}^0_{\cal H}({\rm End}{\cal V})$. 
\end{itemize} 
Here 
${\rm End}{\cal V} = 
{\cal V}^{\vee} \otimes {\cal V}$  
is a complex in ${\rm Hod}_X$. The complex   
${\cal C}^*_{\cal H}({\rm End}{\cal V})$ is the total complex of the 
bicomplex ${\cal C}^{\bullet}_{\cal H}({\rm End}{\cal V})$. 
So 
an element (\ref{3/29/07/1}) includes a function on $X$ with values in 
${\rm End}^0{\cal V}$, a $1$-form with values in 
${\rm End}^{-1}{\cal V}$, etc. 
Thanks to Theorem \ref{8.20.05.7} 
it provides the cohomology $H^*({\cal V})$ with a variation of 
$\R$-MHS's data.

\bl \la{12.02.ups.1}
(i) DG comodules over the DG Lie coalgebra 
${\cal C}^*_{\cal H}({\rm End}{\cal V})$ are the same thing 
as Maurer-Cartan elements (\ref{3/29/07/1}). 

(ii) Homotopy equivalent DG comodules correspond to the equivalent 
Maurer-Cartan elements. 
\el

\begin{proof}
(i) 
Neglecting the differential, the graded Lie coalgebra
${\cal L}^*_{{\cal H}; X}$ is a free 
Lie coalgebra cogenerated by the graded object ${\cal D}_X[1]$ in the category 
${\rm Hod}_X$. 
Therefore a comodule over ${\cal L}^*_{{\cal H}; X}$ is 
determined by the coaction of the cogenerators. The latter 
is given 
by a graded object 
${\cal V}$ plus a linear map 
$ {\cal V} \lra {\cal V} \otimes_\R {\cal D}_X[1]. 
$ 
It can be interpreted as a degree zero element  in 
$  
{\cal D}_X[1]\otimes_\R  {\rm End}{\cal V},  
$ 
which is the same as 
an element ${\cal G}\in {\cal C}^1_{\cal H}({\rm End}{\cal V})$. 

ii) Similarly a homotopy is 
given by an element $ {\cal V} \lra {\cal V} \otimes_\R {\cal D}_X$, 
which can be interpreted as an element  
$B \in  {\cal C}^0_{\cal H}({\rm End}{\cal V})$. 
\end{proof}

Assume that $L_k=0$ for $k>0$. Then 
subcategory of DG-modules concentrated in the degree zero 
is canonically identified with the category of modules over the 
Lie algebra $H_0(L_{\bullet})$.  Indeed, given such 
a module $M_0$, $L_k$ where $k \not = 0$ acts on 
it trivially since $M_0$ a graded 
module over a graded Lie algebra. Further, since  
the action map $\mu$ is 
a map of complexes, $\partial(L_{-1})$ acts on $M_0$ trivially. 

There is a similar picture for DG comodules over a DG Lie 
coalgebra concentrated in degrees $\geq 0$, which we use below.

\paragraph{Proof of Lemma \ref{12.17.ups.10}.}  
Follows immediately from Lemma \ref{12.02.ups.1}.

\vskip 3mm
A  {\it DG-connection} on $X$ 
is a first order linear differential operator 
${\bf d} + {\cal A}$ where ${\cal A}$ is an  ${\rm End}{\cal V}$-valued 
differential form on $X$  of total degree $1$. 
A DG-connection is {\it flat} 
if $({\bf d} + {\cal A})^2=0$.

\bd \la{11.22.ups.1}
Let ${\cal V}$ be a 
variation of real Hodge structures  on $X$, and 
${\cal G} \in {\cal C}^{1}_{\cal H}({\rm End}{\cal V})$. 
 The {\it twistor DG-connection of ${\cal G}$} is a 
DG-connection on $X$ given by
$\nabla_{\cal G}:= {\bf d} +
\gamma({\cal G})$. 
\ed

\bp \la{12.09.ups.5} \la{12.12.ups.1}
A form ${\cal G}\in {\cal C}^{1}_{\cal H}({\rm End}{\cal V})$ 
satisfies the Maurer-Cartan equation
$
\delta{{\cal G}} + {\cal G}\wedge {\cal G} =0
$ 
if and only if the twistor 
 DG-connection $\nabla_{{\cal G}}$ is flat. 
\ep

\begin{proof} Follows immediately from Theorem \ref{12.8.ups.1s}. 
\end{proof}


\section{DG Lie algebra ${\cal C}_{H, S}$ and $L_\infty$-algebra of plane decorated trees} \label{hc5sec}

Let  $H^{\vee}$ be a  finite dimensional symplectic vector 
space, and $S$ a finite set. 
Set
\be \la{VHS}
{\rm V}^{\vee}_{H,S}:= H^{\vee} \oplus \C[S]. 
\ee

Let 
${\cal C}^{\vee}_{H, S}:= {\cal C}{\rm T}({\rm V}^{\vee}_{H,S})$ 
be the cyclic envelope 
of the tensor algebra of ${\rm V}^{\vee}_{H,S}$. 
In Section \ref{SSSec7.1} 
we define a Lie coalgebra structure on ${\cal C}^{\vee}_{H, S}$.  
The coproduct  generalizes the coproduct for the dihedral Lie coalgebras from 
\cite{G4} -- the latter corresponds to the $H^{\vee}=0$ case. 
The 
Lie coalgebra structure on ${\cal C}^{\vee}_{H, S}$ provides the dual vector space 
${\cal C}_{H, S}$   with a Lie algebra structure. 
The Lie algebra  ${\cal C}_{H, S}$ 
plays a key role in our story.  We give a different interpretation of the  
Lie algebra ${\cal C}_{H, S}$ in Section \ref{hc7sec}.

In Section \ref{SSSec7.2} we define an $L_\infty$-coalgebra 
$({\cal T}^{\vee, \bullet}_{H, S}[1], \partial)$ of the plane ${\rm V}^{\vee}_{H,S}$-decorated trees. 
 
Its dual is an $L_\infty$-algebra ${\cal T}^{\bullet}_{H, S}[-1]$ of  
plane ${\rm V}_{H,S}$-decorated trees. 

The standard cochain complex of the $L_\infty$-algebra 
${\cal T}^{\bullet}_{H, S}[-1]$ is identified with a commutative differential graded  algebra 
${\cal F}^{\vee, \bullet}_{H, S}$ of 
plane ${\rm V}^{\vee}_{H,S}$-decorated forests:
$$
S^\bullet({\cal T}^{\vee, \bullet}_{H, S}) = {\cal F}^{\vee, \bullet}_{H, S}. 
$$

The $L_\infty$-coalgebra structure on ${\cal T}^{\vee, \bullet}_{H, S}[1]$  generalises both 
the Kontsevich-Boardman differential $\partial_\Delta$ in the graph complex, and the differential 
defined in \cite{G5}.  The differential $\partial_\Delta$ 
is indentified with  
the differential on the ${\cal T}^{\vee, \bullet}_{H, S}[1]$. 
Our definition was suggested by computation 
of the differential of the Hodge correlator  in  Theorem \ref{5.16.06.12g}.

The sum over all plane trivalent  trees with a given decoration  
gives rise to an injective map
$$
F: {\cal C}^{\vee}_{H, S} \hra {\cal T}^{\vee, 1}_{H, S}[1]. 
$$
The $L_\infty$-algebra ${\cal T}^{\bullet}_{H, S}[-1]$ with the differential $\partial_\Delta$ 
is a resolution of the Lie algebra ${\cal C}_{H, S}$. So
$$
{\cal C}_{H, S} = H_0({\cal T}^{\bullet}_{H, S}[-1]).
$$

\subsection{A Lie algebra structure on ${\cal C}_{H, S}$}  \la{SSSec7.1}
Let us define a cobracket 
$\delta: {\cal C}^{\vee}_{H, S} \lra \Lambda^2{\cal C}^{\vee}_{H, S}$. It is  
a sum of two maps: 
$$
\delta = \delta_{\rm Cas}+\delta_S.
$$ Let us define these maps.
We picture a cyclic word $W$ on an oriented  circle. 

\begin{figure}[ht]
\centerline{\epsfbox{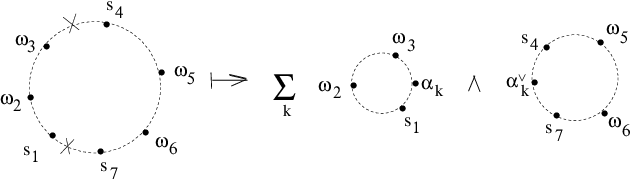}}
\caption{The component $\delta_{\rm Cas}$ of the differential.}
\label{feyn23}
\end{figure}

(i) {\it The map $\delta_{\rm Cas}$}. Cut two different arcs of the circle. We get two semicircles. 
Make an oriented  circle out of each of them. Each of them comes with a special point on it,
 obtained by gluing the ends of the  semicircle. Decorate the special point of the first circle by 
$\alpha_k$, and the special point of the second circle by  $\alpha_k^{\vee}$, 
take the sum over $k$, and put the sign. This is illustrated on Fig \ref{feyn23}. 
The map  $\delta_{\rm Cas}$ is obtained by taking the sum over 
all possible cuts. 

Observe that the wedge product comes with the plus sign, and changing the order of the cuts (and hence the order of the terms in the wedge product)
we do not change it, since if $\{\alpha_k\}$, $\{\alpha_k^{\vee}\}$  
is a  basis and its symplectic dual, then the same is true for 
 $\{-\alpha_k^{\vee}\}$, $\{\alpha_k\}$. 

\begin{figure}[ht]
\centerline{\epsfbox{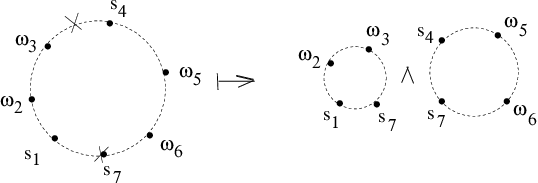}}
\caption{The component $\delta_S$ of the differential.}
\label{feyn24}
\end{figure}

(ii) {\it The map $\delta_S$}. Cut the circle at an arc and at 
an $S$-decorated point, which is not at the end of the arc, 
and make two circles as before. Their special points inherit 
an $S$-decoration from the 
$S$-decorated point, see Fig \ref{feyn24}. 
Make the sum over 
all possible cuts. 

This time the order of the terms in the wedge product matters. 
We put the semicircle which, going according to the orientation of the circle, 
has the $S$-cut at the end  as the first term. 

\begin{figure}[ht]
\centerline{\epsfbox{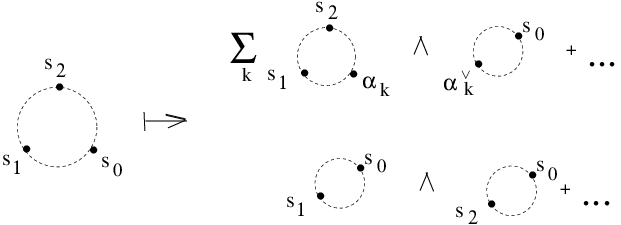}}
\caption{The coproduct of ${\cal C}(\{s_0\}\otimes \{s_1\}\otimes \{s_2\})$.  }
\label{feyn25}
\end{figure}

\paragraph{Examples.} One has 
\begin{equation} \label{5.16.06.h2}
\delta {\cal C}(\{s_0\}\otimes \{s_1\}) = {\cal C}(\{s_0\}\otimes \alpha_k) \wedge 
{\cal C}(\alpha_k^{\vee}\otimes \{s_1\}). 
\end{equation}
$$
\delta {\cal C}(\{s_0\}\otimes \{s_1\}\otimes \{s_2\}) = 
$$
\begin{equation} \label{5.16.06.h1}
{\rm Cycle}_{0,1,2}\left ( 
{\cal C}(\{s_0\}\otimes \{s_1\}\otimes \alpha_k) \wedge 
{\cal C}(\alpha_k^{\vee}\otimes \{s_2\}) + {\cal C}(\{s_0\}\otimes \{s_1\}) \wedge 
{\cal C}(\{s_1\}\otimes \{s_2\}\right). 
\end{equation}
where ${\rm Cycle}_{0,1,2}$ means the cyclic sum. See Fig \ref{feyn25} 
illustrating the last formula. 

\begin{proposition} \label{6.18.05.1} One has $\delta^2 =0$. 
\end{proposition}

Proposition \ref{6.18.05.1} is deduced from the construction 
of the $L_\infty$-coalgebra of plane decorated trees, see 
 the end of Section \ref{hc5sec}.  

\subsection{A commutative differential graded algebra of plane decorated forests} \la{SSSec7.2}

We start from a background on DG Lie coalgebras, 
$L_\infty$-coalgebras, decorated trees and forests.  

\paragraph{DG Lie coalgebras.} 
Let ${\cal L}^\bullet$ be a DG Lie coalgebra with a differential $d$. Its cobracet 
$\delta: {\cal L}^\bullet \lra \Lambda^2{\cal L}^\bullet$ is a map of complexes: $\delta d = d \delta$. 
 The map $\delta$  is understood as   
another differential 
$$
\delta: {\cal L}^\bullet[-1] \lra S^2({\cal L}^\bullet[-1]).
$$
We alter the signs of $\delta$, so that $\delta d + d \delta=0$. Then 
the symmetric algebra $S^\bullet({\cal L}^\bullet[-1])$ is a commutative DGA  
with a differential $\partial$, given by 
the map $d+\delta$ on the space of cogenerators ${\cal L}^\bullet[-1]$. 
Conversely, consider  a commutative DGA $S^\bullet({\cal L}^\bullet[-1])$ with a differential $\partial$, such that 
\be \la{mapD}
\partial: {\cal L}^\bullet[-1] \lra {\cal L}^\bullet[-1] \oplus S^2({\cal L}^\bullet[-1]). 
\ee
Then there is  DG Lie coalgebra 
structure on ${\cal L}^\bullet$ defined as follows. Write the map (\ref{mapD}) as $\partial = d+\delta$, where $d$ and $\delta$ 
are the two natural components, where
$d$ is the one preserving ${\cal L}^\bullet[-1]$. Then $d$ is a differential on ${\cal L}^\bullet$, and  
$\delta$ is a Lie coalgebra coproduct on ${\cal L}^\bullet$. 

\paragraph{$L_\infty$-coalgebras.} Generaliseing the last construction, a  graded $L_\infty$-coalgebra structure 
on a graded vector space ${\cal L}^\bullet$ 
is given by a differential $\partial$ of the graded commutative algebra 
$S^\bullet({\cal L}^\bullet[-1])$, $\partial^2=0$. It is determined by its action on the cogenerators: 
$$
\partial: {\cal L}^\bullet[-1] \lra S^\bullet({\cal L}^\bullet[-1]). 
$$
Write it as a sum of the $S^\bullet$-components: $\partial = \sum_{i\geq 1}\partial_i$. The
 component $\partial_1: {\cal L}^\bullet[-1] \lra {\cal L}^\bullet[-1]$ is a differential. 
If $\partial_m =0$ for $m>2$ then it is a DG Lie coalgebra with the coproduct $\partial_2$.

\paragraph{Trees, forests, decorations.} 
Below a plane tree is a tree with internal vertices of valency $\geq
3$. A tree may have no internal 
vertices, i.e. just two external vertices. 
A {\it plane forest} is a disjoint union of plane trees. 
An $R$-decorated plane forest 
is a plane forest decorated by a set $R$. If ${F}  = {T}_1 \cup \ldots 
\cup {T}_k$ is
a forest presented as a union of trees ${T}_i$, its orientation torsor
is the product 
of the tree orientation torsors: ${\rm or}_{F} 
= {\rm or}_{{T}_1}\otimes \ldots \otimes {\rm or}_{{T}_k}$. 

{\it Decorations by cyclic words}. 
A decoration of a tree $T$ by vectors of a vector space
$V$ gives rise to a decoration of $T$ by an element ${\cal C}{\rm T}(V)$: 
if the external edges are decorated clockwise 
 by the vectors
$v_1, ..., v_n$, the 
resulting decoration is a cyclic word  $W:= {\cal C}(v_1 \otimes \ldots
\otimes v_n)$. We say then that {\it $T$ is decorated by the cyclic word
$W$}. Different decorations of external edges of $T$ by vectors $v_i$ 
can lead to the same decoration by a cyclic word $W$. 

\bd 
A decoration of a forest is the symmetric 
product $W_{T_1} \circ \ldots \circ W_{T_n}$ of cyclic words $W_{T_i}$ decorating its connected components $T_i$. 
 \ed

\paragraph{The graded space of decorated  trees.} 
Let ${\cal T}^{\vee, \bullet}_{H, S}$ be the graded $\C$-vector space generated by triples $(T, W; {\rm Or}_T)$ 
where $T$ is a plane tree, ${\rm Or}_T$ is an element of the orientation torsor of $T$, 
and $W$ is a ${\rm V}^{\vee}_{H,S}$-decoration of $T$. The relations between the generators 
are the following: 
$$
(T, W; {\rm Or}_T) = - (T, W; -{\rm Or}_T), ~~~~(T, \lambda W_1+ \mu W_2; {\rm Or}_T) = \lambda (T, W_1; {\rm Or}_T) + 
\mu (T, W_2; {\rm Or}_T). 
$$
The degree  is defined as follows,  where $v$ runs 
through the set of all  internal vertices of $T$: 
\be \la{12.18.15.1}
{\rm deg}(T, W; {\rm Or}_T):= 1+ \sum_v ({\rm val} (v)-3).
\ee 
The extra $+1$ in the formula for the degree is necessary to have a differetial $\partial$, defined below. 

\paragraph{The graded commutative algebra of decorated  forests.} 
Consider the symmetric algebra of 
the graded space of ${\cal T}^{\vee, \bullet}_{H, S}$ of ${\rm V}^{\vee}_{H,S}$-decorated plane trees:
\be \la{11.30.13.1}
{\cal F}^{\vee, \bullet}_{H, S}:= S^\bullet\Bigl({\cal T}^{\vee, \bullet}_{H, S}\Bigr). 
\ee
It is a graded commutative algebra. 
The vector space ${\cal F}^{\vee, \bullet}_{H, S}$  is spanned  
by oriented forests decorated by symmetric products
of cyclic words in ${\rm V}^{\vee}_{H, S}$. 
The grading is given by the degree 
\be \la{degF}
{\rm deg}(F,  W; {\rm Or}_F):= \sum_{\mbox{vertices $v$ of $F$}}({\rm val}(v)-3) + \pi_0(F). 
\ee
The product 
in ${\cal F}^{\vee, \bullet}_{H, S}$ is described on the generators by 
$$
(F_1, W_1; {\rm Or}_{T_1}) \ast (F_2, W_2; {\rm Or}_{T_2}) =
(F_1\cup F_2, W_1\circ W_2; {\rm Or}_{T_1}\wedge {\rm Or}_{T_2}). 
$$

\paragraph{Commutative DGA ${\cal F}^{\vee, \bullet}_{H, S}$ of plane decorated forests.}  Let us define a differential
 $$
\partial: {\cal F}^{\vee, \bullet}_{H, S} \to 
{\cal F}^{\vee, \bullet+1}_{H, S}.
$$
 When $H^\vee  = 0$ it is
the differential defined in \cite{G5}. We define it on trees, 
and then extend by the Leibniz rule. It  
has three components: 
$$
\partial= \partial_\Delta+\partial_{\rm Cas}+\partial_{\rm S}.
$$ 

\begin{figure}[ht]
\centerline{\epsfbox{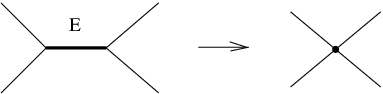}}
\caption{$\partial_\Delta$: The $\delta_X$-term for an internal edge  $E$.}
\label{feyn17}
\end{figure}

(i) {\it The map $\partial_\Delta$}. Let $(T, W_T; {\rm Or}_T)$ be a generator
 and $E$ an {\it internal} edge of $T$, see Fig \ref{feyn17}.
Let $T/E$ be the tree
 obtained by contraction of the edge $E$. 
So it has one less edge, and one less vertex. The orientation ${\rm Or}_T$ 
of the tree $T$ induces an orientation ${\rm Or}_{T/E}$ 
of the tree $T/E$. Namely, if ${\rm Or}_T = E \wedge E_1 \wedge E_{2
  }\wedge  ...$ then ${\rm Or}_{T/E} = E_1 \wedge E_{2 }\wedge
... $. 
Constructing an internal 
edge we do not touch the decoration. Set 
$$
\partial_\Delta(T,W_T; {\rm Or}_T):= \sum_{\mbox{$E$: internal edges of $T$}}(T/E, W_T; {\rm Or}_{T/E}).
$$
It 
is (a decorated version of) the graph complex differential 
of  Boardman and Kontsevich. 

\begin{figure}[ht]
\centerline{\epsfbox{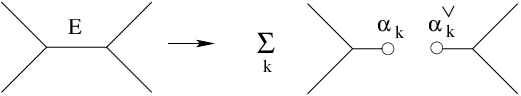}}
\caption{The Casimir map $\partial_{\rm Cas}$ for an internal edge $E$.}
\label{feyn18}
\end{figure}

\begin{figure}[ht]
\centerline{\epsfbox{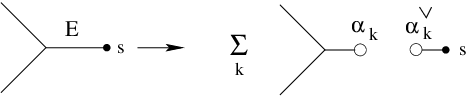}}
\caption{$\partial_{\rm Cas}$: 
Contribution of the Casimir term for an external edge $E$.}
\label{feyn19}
\end{figure}

(ii) {\it The map $\partial_{\rm Cas}$}. Let $E$ be an edge of $T$. 
Cut the tree $T$ along the edge $E$, getting two trees, 
$T_1$ and $T_2$, see Fig \ref{feyn18} and Fig \ref{feyn19}. 
Choose their orientations ${\rm Or}_{T_1}$ and 
 ${\rm Or}_{T_2}$ so that 
$
{\rm Or}_{T} = E \wedge {\rm Or}_{T_1}\wedge {\rm Or}_{T_2}. 
$ 
The trees $T_i$ inherit partial decorations by 
the (non-cyclic!) words $W_i'$, so that 
$W_T = {\cal C}(W_1' \otimes W_2')$. 
Denote by $E_1$ and $E_2$ the new external edges of the trees $T_1$ and
$T_2$ obtained by cutting $E$. 

Below we use the identity (Casimir) element ${\rm Id}: H \to H$. 
We write it as follows. Choose a basis $\{\alpha_k\}$ of $H$. 
Denote by $\{\alpha^{\vee}_k\}$ the dual basis:
$(\alpha_k, \alpha^{\vee}_l) = \delta_{kl}$. Then 
 ${\rm Id} = \sum_k \alpha^{\vee}_k \otimes \alpha_k$. 

We decorate $E_1$ by $\alpha_k$, 
$E_2$ by $\alpha^{\vee}_k$, getting a decoration ${\cal C}(W_1 \otimes
\alpha_k)$ 
of the tree $T_1$, and a decoration ${\cal C}(\alpha^{\vee}_k \otimes W_2)$ 
of the tree $T_2$. Then 
\begin{equation} \label{7.29.06.1}
\partial_{\rm Cas}(T, W_T; {\rm Or}_T):= \sum_{E} \sum_k (T_1, {\cal C}(W'_1 \otimes
\alpha_k); {\rm Or}_{T_1})\wedge (T_2, {\cal C}(\alpha_k^{\vee} 
\otimes W_2'); {\rm Or}_{T_2})
\end{equation}
where the first sum is over all internal edges $E$ of $T$.

\begin{figure}[ht]
\centerline{\epsfbox{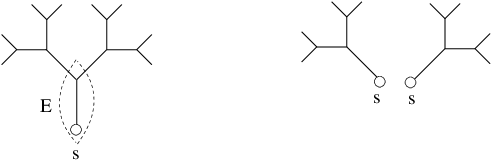}}
\caption{The map $\partial_S$: contribution of  an external edge $E$.}
\label{feyn20a}
\end{figure}

\vskip 3mm
(iii) {\it The map $\partial_S$}. Let $E$ be an external 
$S$-decorated edge of $T$, 
see Fig \ref{feyn20a}. 
Remove it together with a little neighborhood of its vertices. 
If  
$E\not =T$,  
one of the vertices has the valency $v \geq 3$. 
The tree $T$ is replaced by 
$v-1$ trees $T_1, \ldots , T_{v-1}$:  
each of the  internal edges sharing the vertex with $E$
produces a new tree. Choose their orientations ${\rm Or}_{T_i}$
so that 
$ 
{\rm Or}_{T} = E \wedge {\rm Or}_{T_1} \wedge \ldots \wedge 
{\rm Or}_{T_{v-1}}. 
$ 
The  decoration of the tree $T$ provides a 
 decorations $W_{T_i}$ of the  new trees $T_i$, 
so that the new external
 vertex of each of the trees $T_i$ inherits the $S$-decoration of $E$. 
Set
\be
\begin{split}
&\frac{\partial}{\partial E}(T, W_T; {\rm Or}_{T}) = 
(T_1, W_{T_1}; {\rm Or}_{T_1})\wedge \ldots \wedge 
(T_{v-1}, W_{T_{v-1}}; {\rm Or}_{T_{v-1}}), \\
&\partial_S:= 
\sum_{\mbox{$E$: $S$-decorated}}\frac{\partial}{\partial E}.\\
\end{split}
\ee

\bt \label{6.20.06.1}
The map $\partial = \partial_\Delta + \partial_{\rm Cas} +\partial_{\rm S}$ is a differential 
on the graded commutative algebra  ${\cal F}^{\vee, \bullet}_{H, S}$: 
$\partial^2=0$. More specifically, $\partial_{\rm Cas}$ anticommutes with 
$\partial_\Delta+\partial_{\rm S}$, and one has 
$$
\partial_\Delta^2=0, ~~\partial_{\rm Cas}^2=0, ~~
  ~~(\partial_\Delta  + \partial_{\rm S})^2  =0.
$$ 
\et

\begin{figure}[ht]
\centerline{\epsfbox{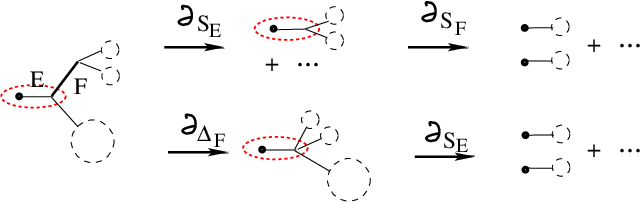}}
\caption{Illustrating $(\partial_{\rm S}+ \partial_\Delta)^2=0$.  
The ${\rm S}$-decorated external vertex of the edge $E$ and its descendants are the boldface vertices. 
We also trace contributions 
of the  internal edge $F$  to $\partial_\Delta$. }
\label{feyn55}
\end{figure}

\begin{proof} Thanks to the term $\pi_0(F)$ in  formula (\ref{12.18.15.1}) for the degree, the map 
$\partial$ has degree $+1$. Notice that without the   $\pi_0(F)$ the map 
$\partial_{\rm Cas}$ has degree $0$, and 
the map $\partial_{\rm S}$ is not homogeneous. 

It is well known that $\partial_\Delta^2=0$: the $\partial_\Delta$ is the Kontsevich-Boardman 
differential. 

To check that 
$\partial_{\rm Cas}^2=0$ 
observe that $\partial_{\rm Cas}^2(T, W; {\rm Or}_{T})$ is a sum of two terms, one of whom is 
\begin{equation} \label{7.29.06.2}
\sum_{E_1, E_2} \sum_{k, l} (T_1, {\cal C}(W_1 \otimes
\alpha_k); {\rm Or}_{T_1})\wedge (T_2, {\cal C}(\alpha_k^{\vee} 
\otimes W_2 \otimes \beta_l); {\rm Or}_{T_2}) \wedge (T_3, {\cal C}(\beta_l^{\vee} 
\otimes W_3); {\rm Or}_{T_3}) 
\end{equation}
and the other has the opposite sign, see 
Fig \ref{feyn37}. 
\begin{figure}[ht]
\centerline{\epsfbox{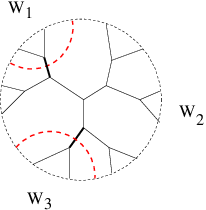}}
\caption{The Casimir map is a differential: $\partial_{\rm Cas}^2=0$.}
\label{feyn37}
\end{figure}
Here we cut the tree $T$ along the edges $E_1$ and $E_2$, shown by  the thick segments. 
The dotted arcs show how cutting along these edges we get three new trees. 
Computing the signs we see that the differentials 
$\partial_{\rm Cas},\partial_\Delta$ anticommute.  

Although $\partial_{\rm S}^2 \not = 0$, we have 
$(\partial_{\rm S} + \partial_\Delta)^2  =0$. 
To show this it is sufficient to trace the contribution of a single 
${\rm S}$-decorated external edge $E$. For any internal edge $F$ not incident to $E$, 
the naive operations of shrinking of the edge $F$ and removing the edge $E$ evidently commute. 
So taking into account the orientation torsors, we see that  
the $\partial_{\rm S}$-component at $E$ and the $\partial_\Delta$-component at $F$ anticommute. 
The contributions of the incident edges $E$ and $F$ are depicted on  Fig \ref{feyn55}. 
Naively, $\partial_{S_F} \circ\partial_{S_E}$ gives the same result as 
$\partial_{\Delta_E} \circ\partial_{S_F}$. So taking into account the orientation torsor, 
we get $\partial_{S_F} \partial_{S_E} + \partial_{\Delta_E} \partial_{S_F}=0$.

Similarly, the maps 
$\partial_{\rm Cas},\partial_{\rm S}$ anticommute: notice that  $\partial_{\rm Cas}$ does not 
create new ${\rm S}$-decorated vertices. 
\end{proof}

\paragraph{$L_\infty$-coalgebra  ${\cal T}^{\vee, \bullet}_{H, S}[1]$ of plane decorated trees.} 

The map $\partial_{\Delta}$ is the component of the differential $\partial$ on the 
DGCom ${\cal F}^{\vee, \bullet}_{H, S}$ 
preserving
 the subspace ${\cal T}^{\vee, \bullet}_{H, S}$. So we get a 
differential
$$
\partial_{\Delta}: {\cal T}^{\vee, \bullet}_{H, S} \lra {\cal T}^{\vee, \bullet}_{H, S}. 
$$
 The differential $\partial$ defines a $L_\infty$-coproduct 
$$
\partial: {\cal T}^{\vee, \bullet}_{H, S} \lra S^\bullet\Bigl({\cal T}^{\vee, \bullet}_{H, S}\Bigr).
$$
Therefore 
${\cal T}^{\vee, \bullet}_{H, S}[1]$ is an $L_\infty$-coalgebra
with a differential $\partial_{\Delta}$ and a $L_\infty$-cobracket $\partial$. 
The commutative DGA ${\cal F}^{\vee, \bullet}_{H, S}$ is the standard cochain complex of the $L_\infty$-coalgebra 
${\cal T}^{\vee, \bullet}_{H, S}[1]$.

\begin{figure}[ht]
\centerline{\epsfbox{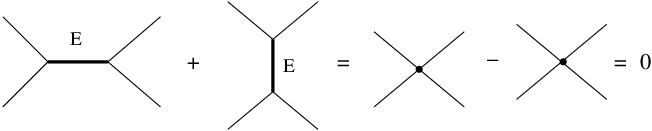}}
\caption{The two terms for the internal edges $E$ cancel each other.}
\label{feyn17a}
\end{figure}
\begin{figure}[ht]
\centerline{\epsfbox{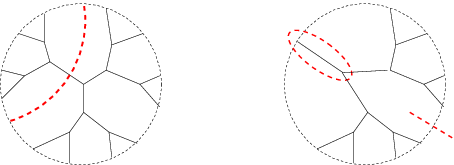}}
\caption{Compatibility of the differential $\partial_{\rm Cas}$ (on the left) and $\partial_{\rm S}$ 
(on the right) on cyclic words and trees.}
\label{feyn38}
\end{figure}

\paragraph{The sum over  plane trivalent  trees map.} 
Let us define an injective 
linear map 
$$
F: {\cal C}^{\vee}_{H, S} \hra {\cal T}^{\vee, 1}_{H, S}[1], \qquad 
F(W):= \sum_T(T, W; {\rm Or}_T).  
$$
The sum is over all plane trivalent trees decorated by the
cyclic word $W$, and ${\rm Or}_T$ is the canonical orientation of the
trivalent tree $T$ corresponding to the clockwise orientation of the
plane. 

\begin{lemma} \label{6.19.06.1}
i) The differential $\partial_\Delta$ 
is zero on the image of the map $F$, i.e. $\partial_{\Delta} F = 0$. 

ii) The map $F$ commutes with the differentials:  $\partial F = 
F \delta$. 

\end{lemma}

\begin{proof}  i) Let $E$ be an internal edge. There are exactly two 
contributions to $\partial_\Delta$ corresponding to the trees
 shown on  Fig \ref{feyn17a} 
(the parts of these graphs which are not shown are the same). 
They cancel each other since the corresponding trees
 with one $4$-valent vertex   
obtained by shrinking the edge $E$   
inherit different orientations. Thus $\partial_\Delta$ 
is zero on the image of the map $F$. 

ii) To check that the differential $\partial_{\rm Cas}$ 
on ${\cal T}^{\vee, \bullet}_{H, S}$
match the differential $\delta_{\rm Cas}$ for ${\cal C}^{\vee}_{H, S}$,  
look at the left picture on 
Fig \ref{feyn38}. 
Take a cut  
shown by a punctured arc, and  vary the trees 
in the two obtained domains keeping their external vertices untouched.   

To check the similar statement for the map $\partial_{\rm S}$ 
on ${\cal T}^{\vee, \bullet}_{H, S}$ and the differential $\delta_{\rm S}$ on ${\cal C}^{\vee}_{H, S}$,  
look at the right picture on 
Fig \ref{feyn38}. 
Cut out an edge $E$, cut the circle inside the arc separating the two obtained trees, and vary the trees, keeping 
the external vertices  untouched.   

Now $\partial F = F \delta$ follows from $\partial_{\rm Cas} F = F \delta_{\rm Cas}$, 
$\partial_{\rm S} F = F \delta_{\rm S}$, and $\partial_{\Delta} F = 0$. 
\end{proof}

\paragraph{Proof of Proposition \ref{6.18.05.1}.}  
Follows   from 
the injectivity of the sum over plane trivalent trees 
map $F$, Lemma \ref{6.19.06.1}ii) and Theorem \ref{6.20.06.1}.

\bt \la{1.15.08.2}
The 
graded $L_\infty$-coalgebra ${\cal T}_{H, S}^{\vee, \bullet}[1]$  
is a resolution of the Lie coalgebra ${\cal C}_{H, S}^{\vee}$. So
$$
{\cal C}_{H, S}^{\vee} = H^0_{\partial_{\Delta}}({\cal T}_{H, S}^{\vee, \bullet}[1]).  
$$
\et

\begin{proof} Follows from the exactness of (\ref{GKap1}) and Lemma \ref{6.19.06.1}. 
\end{proof}
\section{The Hodge correlator twistor connection} \la{hc6sec}

The main result of this section is Theorem \ref{MDET}. 
 
\subsection{Preliminary constructions} \la{hc6.1sec}

\paragraph{The Hodge theoretic set-up over a base.}

Let $p:X \to B$ 
be a smooth family of complex projective curves over $B$, 
and $S \subset X$  a smooth divisor in $X$ over $B$. 
Let $s_0$ be a single component of $S$, and $S^* = S - s_0$. 
The divisor $S$ intersects each fiber at a  finite collection of points, containing a special point $s_0$. 
Below we usually suppress the subscript $B$ from the notation like $X_{/B}, S_{/B}$. There are 
variations of $\R$-Hodge structures over $B$: 
\be \la{VHS*}
{\rm V}^{\vee}_{X,S^*} := {\rm gr}^WR^1p_*(X-S, \R) = H^1(X_{/B}, \R) \oplus \R[S_{/B}^*](-1), 
\ee
\be \la{VHS**}
{\rm V}_{X,S^*}:= {\rm gr}^WR^1p_*(X-S, \R)^\vee = H_1(X_{/B}, \R) \oplus \R[S_{/B}^*](1).
\ee
All  objects considered  in  Section \ref{hc5sec} 
can be upgraded to variations of  $\R$-Hodge structures over $B$. Indeed, they are obtained by 
linear algebra constructions from the variation 
(\ref{VHS*}). 
In particular:

1. Consider the object 
$$
{\cal C}^{\vee}_{X, S^*}(1):= {\cal C}{\rm T}({\rm V}^\vee_{X,S^*})(1).
$$ 
It 
is a Lie coalgebra in the category of variations of $\R$-Hodge structures over $B$. 
The twist  by the $\R(1)$ is necessary to make the coproduct a 
morphism in this category. 

2. The dual object is a Lie algebra in the same category:
\be \la{vvaarrrrhhss}
{\cal C}_{X, S^*}(-1):= {\cal C}{\rm T}({\rm V}_{X,S^*})(-1).
\ee

3. The $L_\infty$-coalgebra ${\cal T}^\bullet_{X, S^*}(-1)$ of ${\rm V}_{X,S^*}$-decorated trees
 is a resolution of the Lie algebra ${\cal C}_{X, S^*}(-1)$ in the same category, provided by    
the ``sum over plane trivalent trees''  map: 
$$
F: {\cal C}_{X, S^*}(-1) \stackrel{\rm qis}{\hra} {\cal T}^\bullet_{X, S^*}[-1](-1).
$$

There is a motivic set up, which works exactly the same way, starting  
with pure motives 
$$
{\rm V}^{\vee}_{X,S^*}:= {\rm gr}^WR^1p_*(X-S) ,~~~~
{\rm V}_{X,S^*}:= {\rm gr}^WR^1p_*(X-S)^\vee.
$$

\paragraph{De Rham description of the 
Gauss-Manin connection.} Given a smooth map $p: X \to B$, 
the Gauss-Manin connection $\nabla_{GM}$ in $R^kp_*\C$ 
can be described as follows.  Take a section $s$ of $R^kp_*\C$. 
Choose a smooth family of differential forms $\alpha_t$ 
on the fibers $X_t$ whose cohomology 
classes form the section $s$.  Choose any differential form 
$\beta$ on $X$ whose restriction to every fiber $X_t$ 
is the form $\alpha_t$. Then 
$d\beta$ gives rise to $\nabla_{GM}(s)$ as follows. 
The de Rham differential of $\beta$ along the fibers of $p$ is zero since 
$d\alpha_t=0$ for every $t$. Therefore 
$d\beta$ provides a  $1$-form on the base with values 
in ${\cal A}_{\rm cl}^{k-1}({X_{/B}})$. 
The cohomology classes of the closed $(k-1)$-forms on the fibers  do not depend on the choice of 
the  form $\beta$ representing the section $s$. 
The obtained $1$-form on $B$ with values in $R^{k-1}p_*\C$ 
equals $\nabla_{GM}(s)$. 

\paragraph{A $1$-form $\nu$.} 
The ${\cal H}:= R^1p_*(\R)^\vee$ is a local system on $B$ with the Gauss-Manin connection $\nabla_{GM}$ and a fiber $H_1(X_t, \R)$ over $t \in B$. 
By  definition, $H_1(X_t, \Z)$ is flat for the connection. 

The Albanese variety $A_{X_t}$ of 
${X_t}$ is  identified as a manifold with $H_1(X_t, \R)/H_1(X_t, \Z)$. 
The fibers of the local system 
${\cal H}$ are the tangent spaces 
to Albanese varieties $A_{X_t}$.

Let $\pi: X \times_B X\to B$ be the fibered product of $X$ over $B$. 
 Then there is a $C^{\infty}$ $1$-form 
\be \la{nu}
\nu\in {\cal A}_{X\times_B X}^1\otimes \pi^*{\cal H}
\ee 
on $X\times_B X$ with values in the local system $\pi^*{\cal H}$. 
Its value at a tangent vector $(v_1, v_2)$ at  $(x_1, x_2)$ 
is obtained as follows. 
Let $(x_1(t), x_2(t))$ be a smooth path 
from $(x_1, x_2)$ in the direction $(v_1, v_2)$. 
Assuming that 
the points $(x_1(t), x_2(t))$ are in the same fiber ${X_t}$, 
we get a section $(x_1(t) - x_2(t)) 
\in A_{X_t}$. 
The Gauss-Manin connection provides another 
section of the Albanese fibration passing through 
$(x_1(0) - x_2(0))$. We define $\nu(v_1, v_2)$ to be 
the difference of the tangents to these sections with the velocity
$(v_1, v_2)$. 

\paragraph{Green functions for a family of curves.} 
Let  
${\cal H}_{\C}:= {\cal H}\otimes {\C}= 
{\cal H}^{-1,0}_{\C}\oplus {\cal H}^{0,-1}_{\C}$ be the Hodge  decomposition. 
The intersection pairing 
$
H_1(X_t, \Z) \wedge  H_1(X_t, \Z) \lra  \Z(1)
$ 
provides a pairing 
$$
\langle *, *\rangle: 
{\cal H}^{-1,0}_{\C}\otimes {\cal H}^{0,-1}_{\C}\lra \C. 
$$
The wedge product $\nu \wedge \nu$, 
followed by the 
pairing $\langle *, *\rangle$, provides a canonical $(1,1)$-form 
$$
\langle\nu\wedge \nu\rangle \in {\cal A}^{1,1}_{X\times_B X}.
$$
Denote by $\langle\nu\wedge \nu\rangle^{[1\times 1]}$ its $1 \times 1$-component 
 on $X\times_B X$. 
Let $\delta_{\Delta}$ be the $\delta$-function of the relative diagonal 
$\Delta \subset X \times_B X$. 
A section $a: B \to X$ provides a divisor 
$S_{a}\subset X \times_B X$, whose fiber at 
$t \in {\cal B}$ is $a_t \times X_t \cup X_t\times a_t$.  
It gives rise to a $(1,1)$-current $\delta_{S_{a}}$ on 
$X$. 
\begin{definition} \label{5.17.07.2} Let $p: X\to B$ be 
a proper map and $a: B \to X$ is its  section.  
A Green current $G_{a}(x,y)$ is a $0$-current  on 
$X \times_B X$ satisfying the equation 
\be \la{10.22.07.3}
\overline \partial \partial G_{a}(x,y)= \delta_{\Delta} - 
\delta_{S_{a}} - 
\langle\nu\wedge \nu\rangle^{[1\times 1]}.
\ee
\end{definition}
The  Green function exists fiberwise by the 
$\overline \partial  \partial$-lemma. It determines a $0$-current 
$G_{a}(x,y)$ uniquely up to a function lifted from the base $B(\C)$. 
Let us choose a non vanishing $v \in T_aX$. It provides a uniquely defined
 normalized 
Green function  $G_{v}(x,y)$. 

Below we choose a section $v_0$ of the fibration  $T_{s_0}X_{/B}$, providing us 
a  Green function $G_{v_0}(x,y)$.

\subsection{The Hodge correlator twistor connection}

Let us define the Hodge correlator map in a bit more general set-up. 

Let $T$ be a plane tree, not necessarily trivalent.   
We assume that a tree does not have two-valent vertices. We proceed similarly to  Section 2.3. 
We start with a decorated tree $(T, W)$. 

i) Suppose that the  $T$ has more then one edge. 
Euivalently, it has an internal vertex. 
Consider the fibered product over 
the base $B$ of the copies of $X$ parametrised by the vertices of  $T$.
\be \la{12.12.13.1}
{\cal X}_T:= X_{/B}^{\{\mbox{vertices of $T$}\}}. 
\ee 
Let  $\{E_0, \ldots, E_{r}, F_1, \ldots , F_k\}$ be the   edges of $T$
numbered so that the edges $E_0, \ldots, E_r$ are internal or $S^*$-decorated, 
and the edges $F_1, \ldots , F_k$  are decorated by sections   
$\omega_{F_1}, \ldots ,  \omega_{F_k}$ of 
the local system ${\rm gr}_WH^1(X-S, \R)$. 
We assume that the order $F_1, \ldots , F_k$ is compatible with the cyclic order provided 
by the orientation of the plane containing the tree $T$. 
Choose 
a $1$-form $\alpha_{F_j}$ on $X$ representing a section $\omega_{F_j}$: its restriction to 
a fiber $X_t$ is a  holomorphic/antiholomorphic 
form.

Consider an element 
of the orientation torsor of the tree $T$:
\be \la{5.21.15.1}
(E_0 \wedge \ldots \wedge E_r) \wedge (F_{1} \wedge \ldots \wedge F_{k}) 
\in {\rm or}_T.
\ee
The difference between the element (\ref{5.21.15.1}) and 
 a given orientation  ${\rm Or}_T \in {\rm or}_T$ is denoted by 
$$
{\rm sgn}_{{\rm Or}_T}(E_0 \wedge \ldots \wedge E_r \wedge F_{1} \wedge \ldots \wedge F_{k}):=  
(\ref{5.21.15.1})/{\rm Or}_T \in \{\pm 1\}.
$$

We assign to each edge $E_i$ the Green current $G_{E_i}$ on $X\times_BX$, where the copies of $X$ match the vertices of the edge $E_i$. Since the Green current is symmetric, it 
does not depend on the choice of an orientation of the edge $E_i$. 

Recall the twistor plane $\C^2$ with coordinates $(z,w)$. 
Using the forms  $\alpha_{F_j}$, the Green functions $G_{E_i}$ at the edges 
of $T$, and the 
orientation ${\rm Or}_T$, let us cook up a current 
$\widehat \kappa(T, W; {\rm Or}_T)$ on 
\be \la{12.12.13.1a}
{\cal X}_T \times \C^2.
\ee 
 Put the form  $\alpha_{F_j}$ to the 
copy of $X/B$ assigned to the unique internal vertex of the edge $F_j$. 
Abusing notation, we identify it with its pull back to ${\cal X}_T \times \C^2$. 
Recall 
$$
{\rm D}^\C:= (w\partial - z \overline \partial) + (zdw-wdz). 
$$

\bd \la{MDHCCL} We define a differential form  on (\ref{12.12.13.1a}) by setting
$$
\widehat \kappa(T, W; {\rm Or}_T):=
$$
\begin{equation} \label{3.8.05.15}
{\rm sgn}_{{\rm Or}_T}(E_0\wedge \ldots \wedge E_{r} \wedge F_1 \wedge \ldots \wedge F_k)~{\rm D}^\C G_{E_0}\wedge \ldots \wedge 
{\rm D}^\C G_{E_r} 
\wedge \alpha_{F_{1}}\wedge \ldots \wedge \alpha_{F_{k}}.
\end{equation} 
\ed

\bt  The form $\widehat \kappa(T, W; {\rm Or}_T)$ a current on ${\cal X}_T \times \C^2$, 
 algebraic along the twistor plane:
$$
\widehat \kappa(T, W; {\rm Or}_T)\in {\cal D}^{*,*}({\cal X}_T) \otimes_\C\Omega^{\leq 1}_{\C^2}.
$$
\et

\begin{proof}
Just as 
in Lemma 2.2.
\end{proof}

There is a canonical projection 
$$
p_T: {\cal X}_T^{\rm int} \times \C^2 := X_{/B}^{\{\mbox{internal vertices of $T$}\}}\times \C^2  \lra  B\times \C^2. 
$$ 

\bd Suppressing $v_0$ from the notation, set
\be \la{12.15.13.1}
{\rm Cor}_{\cal H}(T; {\rm Or}_T)(W) := 
{p_T}_*\Bigl(\widehat \kappa(T, W; {\rm Or}_T) \Bigr) \in {\cal D}^{*,*}(B) \otimes_\C\Omega^{\leq 1}_{\C^2}.
\ee
\ed

When $B$ is a point, we recover the Hodge correlator map from Definition 2.6.
In this case 
the form (\ref{12.15.13.1}) does not depend on the choice of  forms $\alpha_{F_j}$. Indeed, 
if $\widetilde \alpha_{F_j}$ is a different choice, then 
the restriction of 
$\widetilde \alpha_{F_j} - \alpha_{F_j}$ to 
the fiber curve is zero. 
Indeed, the fibers are 
smooth complex projective curves, and 
so there are canonical harmonic representatives for 
$H^1$ given by holomorphic/antiholomorphic forms. 
The same is true for the components of the Hodge correlator which are functions on the base. 
In particular, this is the case in the crucial for us case 
when $T$ is a plane trivalent tree, as we will see in Proposition \ref{12.16.15.1}.

\begin{figure}[ht]
\centerline{\epsfbox{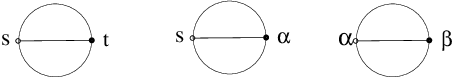}}
\caption{The three possible decorations of a single edge tree.}
\label{hc8a}
\end{figure}

ii) Let $T$ be a tree with a single edge $E$. Since there are no internal vertices,  the 
 definition of the Hodge correlator involves no integration. 
There are three possible decorations, 
see Fig \ref{hc8a}. 

1. $T=T_1$ is decorated by sections $s, t\in S^*$. 
The sections provide  a map $(s, t): B \to X\times_B X$. 
We assign to $T_1$ the pull back of the Green current 
$G(s, t):= (s, t)^*G$ to the base $B$: 
\be \la{2.12.13.10asas}
{\rm Cor}_{\cal H}(T_1)(\{s\}\otimes \{t\}) :=   G(s, t)\in {\cal D}_{B}^{0}\otimes \R(-1).
\ee 
Dualising,  and since the dual to $\{s\}\otimes \{t\}$ is a section of the variation 
 $\R(2)$ on $B$,  we get
\be \la{8.12.13.1}
{\rm Cor}_{\cal H}(T_1) :=   G(s,t)\otimes (\{s\} \otimes \{t\})^\vee\in 
{\cal D}_{B }^{0}\otimes \R(1)\subset {\cal C}_{{\cal H}_\R}^1({\cal C}_{X, S^*}). 
\ee 
To check the last inclusion, notice that  $\R(1)$-valued $0$-forms  have degree $1$ in the Hodge complex 
according to (\ref{11.18.ups.2hu}), and that  
the Green function $G(s, t)$ is imaginary,  in agreement with 
the reality condition (\ref{11.19.ups.10}) on the Hodge complex. The Hodge bidegree is $(-1,-1)$.

2. $T=T_2$ is a single edge, decorated by a section $s: B \to X$ and a 1-form $\alpha$ on $X$.  
We set 
\be \la{2.12.13.10}
{\rm Cor}_{\cal H}(T_2)(\{s\}\otimes\alpha) := s^*\alpha\in {\cal A}_{B }^{1}.
\ee 
Let us dualise this element. Let $s^*\alpha | h$ be the canonical Casimir element, defined as follows. 
Choose a basis $\alpha_1, ..., \alpha_{2g}$ in ${\cal H}^\vee = R^1p_*(\R)$,  the dual basis 
$h_{1}, ..., h_{{2g}}$ in ${\cal H} = R^1p_*(\R)^\vee$, and set 
$$
s^*\alpha | h := \sum_{i=1}^{2g}s^*\alpha_i \otimes h_{i}\in {\cal H}^\vee \otimes {\cal H}. 
$$
Dualising  map (\ref{2.12.13.10}),  using the fact that the dual to $\{s\}\otimes \alpha$ is a section of the variation 
of real Hodge structures 
 ${\cal H}(1)$ on $B$, and twisting by $\R(-1)$, 
we get an element
\be \la{8.12.13.2}
{\rm Cor}_{\cal H}(T_2):= s^*\alpha | h \otimes \R(-1) \in {\cal A}^{1,0}_{B }\otimes {\cal H}^{-1,0} \oplus 
{\cal A}^{0, 1}_{B }\otimes {\cal H}^{0, -1} \subset {\cal C}_{{\cal H}_\R}^1({\cal C}_{X, S^*}). 
\ee 
Each summand has degree $1$ in the Hodge complex, see (\ref{11.18.ups.2hu}). The 
Hodge bidegree is $(0,0)$.  
It is imaginary since  $s^*\alpha | h$ is real,  in agreement with 
 reality condition (\ref{11.19.ups.10}) on the Hodge complex.

3. $T$ is a single edge decorated by forms $\alpha, \beta$. We assign to it zero.

\paragraph{The Hodge correlator class.} 
Let  $W= (x_{1} \otimes \ldots 
\otimes x_{m})_{\cal C}$. Let $|{\rm Aut}(W)|$ be the order of its automorphism group. 
Let ${\cal C}$ be the operator of weighted projection 
on the cyclic tensor algebra: 
\be \la{WCP}
{\cal C}(x_1 \otimes \ldots \otimes x_m):= \frac{1}{|{\rm Aut}(W)|}
 (x_{1} \otimes \ldots 
\otimes x_{m})_{\cal C}.
\ee

Given a collection of cohomology classes $\gamma_i\in V^\vee_{X, S^*}$ and homology classes $h_i\in V_{X, S^*}$,
 we define the 
Hodge correlator map of the 
cyclic tensor product of the expressions $\gamma_i | h_i$ as a cyclic 
product of their homology factors 
taken with the coefficient 
given by the Hodge correlator map applied to the cyclic product 
of the corresponding  differential form factors. 

Here is an example. 
Let $\alpha_1, \alpha_2, \alpha_3$ be forms on $X$. Then 
\be \la{HCC}
{\rm Cor}_{\cal H}(\alpha_1 | h_1 \otimes \alpha_2 | h_2 \otimes 
\alpha_3 | h_3 ) :=
\ee
$$
{\rm Cor}_{\cal H}(\alpha_1 \otimes 
\alpha_2\otimes \alpha_3)
\cdot ( h_{1}\otimes h_{2} \otimes h_{3})_{\cal C} = 
\int_{X/B} (\alpha_1 \wedge 
\alpha_2\wedge \alpha_3)
\cdot ( h_{1}\otimes h_{2} \otimes h_{3})_{\cal C}.
$$

Choose a basis $\{\gamma_k\}$  in $V^\vee_{X, S^*}$. 
Let $h_{\gamma_k}$ be the dual
 basis in $V_{X, S^*}$.  
Set
\be \la{CASx1}
\gamma|h:= \sum_k \gamma_k | h_{\gamma_k}.
\ee
Take a plane oriented tree $(T, {\rm Or}_T)$, not necessarily trivalent, 
with $m$ external vertices. 
Decorate the tree $T$ by the cyclic word ${\cal C}(\gamma|h\otimes \ldots \otimes \gamma|h)$, 
the $m$-fold cyclic tensor product of the $\gamma|h$'s.  
By this we mean decorating the $T$ by the $m$-fold tensor product 
$\gamma|h\otimes \ldots \otimes \gamma|h$ and all its cyclic shifts. We apply 
the Hodge correlator map to this tree, 
decorated by the ${\cal C}(\gamma|h\otimes \ldots \otimes \gamma|h)$:
\be \la{12.15.15.1}
{\rm Cor}_{\cal H}(T, {\rm Or}_T){\cal C}(\gamma|h\otimes \ldots \otimes \gamma|h).
\ee

A plane trivalent tree $T$ has a canonical orientation ${\rm Or}_T$, 
induced by the plane orientation. 

\bd \la{12.15.15.100} The Hodge correlator class ${\bf G}$ is obtained by taking
 the sum of  Hodge correlators (\ref{12.15.15.1}) over all plane trivalent trees $T$, 
  with the canonical orientation ${\rm Or}_T$, twisted by $\R(-1)$:
\be \la{12.15.15.2}
{\bf G} := ~\sum_{m \geq 2}\sum_T\sum_{W}\frac{1}{|{\rm Aut}(W)|}
{\rm Cor}_{\cal H}(T, {\rm Or}_T)(\gamma_1|h_{\gamma_1} 
\otimes \ldots \otimes \gamma_m| 
h_{\gamma_m})(-1).
\ee
Here $\{W\} = \{\gamma_1 
\otimes \ldots \otimes \gamma_m\}$   is a basis  in ${\cal C}({\rm T}(V^\vee_{X, S^*}))$. 
\ed
One can rephraise Definition \ref{12.15.15.100} as follows. Set
\be \la{CASx}
{\rm I}:=  1 + \gamma|h_{\gamma} + \gamma|h\otimes \gamma|h
+ \gamma|h\otimes \gamma|h \otimes\gamma|h + \ldots.
\ee
Applying the operator ${\cal C}$, see (\ref{WCP}), we get a cyclic Casimir element:
\be \la{CYCCAS}
{\cal C}({\rm I}):= {\cal C}\Bigl( 1 + \gamma|h + \gamma|h\otimes \gamma|h 
+ \gamma|h\otimes \gamma|h \otimes\gamma|h + \ldots \Bigr).
\ee

\bd The  Hodge correlator class is obtained by 
applying to (\ref{CYCCAS})  the  Hodge correlator map, twisted by $\R(-1)$, and taking the sum over all plane 
trivalent trees $T$:
\be
\begin{split}
&{\bf G} := \sum_{T}{\rm Cor}_{\cal H}(T, {\rm Or}_T)\Bigl( {\cal C}({\rm I}) (-1)\Bigr). \\
\end{split} 
\ee
\ed
Explicitly, the Hodge correlator class  is given by: 
\be \la{12.16.15.2}
\begin{split}
&{\bf G} = {\bf G }_{0,0} +\sum_{s,t \geq 0}
z^{s}w^{t}\Bigl(
(zdw-wdz)\wedge (s+t+1){\bf G}_{s+1, t+1} +
(w\partial'-z\partial''){\bf G}_{s+1,t+1}\Bigr),\\
&{\bf G }_{0,0}: = 
{\rm Cor}_{\cal H}(s^*\alpha | h_{\alpha}) +  
{\rm Cor}_{\cal H}(\alpha_1|h_{\alpha_1} \otimes \alpha_2|h_{\alpha_2} \otimes  \alpha_3|h_{\alpha_3})_{\cal C}.\\
\end{split}
\ee
For example,  ${\bf G}_{1,1} = G(a,b)$ is the Green function. 

Recall the differential ${\rm d}_{\C^2}$ on the twistor plane $\C^2$. Recall  the projection 
$$
\pi: B\times \C^2 \to B.
$$
 The  Hodge correlator class ${\bf G}$ is a current on $B \times \C^2$ with values in 
the local system of Lie algebras $\pi^*{\cal C}_{X,S^*}(-1)$ 
with a flat connection  ${\bf d}+{\rm d}_{\C^2}$. This local system
 is  a Lie algebra 
in the category of variations of 
$\R$-Hodge structures on $B \times \C^2$.

Recall the Lie algebra $\pi^*{\cal C}{\cal L}ie_{X, S^*}(-1)$, see (\ref{mvlca}).  
It is obtained by linear algebra 
constructions (\ref{9.25.13.1}), (\ref{9.25.13.2}) from variations of  $\R$-Hodge structures on $B\times \C^2$. So it 
is a Lie subalgebra of $\pi^*{\cal C}_{X,S^*}(-1)$ 
in the same category. Its underlying   
 flat connection is ${\bf d}+{\rm d}_{\C^2}$.

\bl \la{12.16.15.1}
The Hodge correlator  ${\bf G}$ is a $\pi^*({\cal C}{\cal L}ie_{X, S^*})(-1)$-valued 
$1$-form on $B \times \C^2$. 
\el

\begin{proof} 
The Hodge correlator  ${\bf G}$ satisfies
the shuffle relations, see Section \ref{sec2.3x}. This just means that  ${\bf G}$ 
is a $\pi^*({\cal C}{\cal L}ie_{X, S^*})(-1)$-valued current on $B \times \C^2$. 
A plane trivalent tree with $n$ vertices has $2n+1$ edges. 
In the integrand (\ref{3.8.05.15}) we assign a $1$-form to each edge. Since $(2n+1) - 2n =1$, 
after the integration we get a $1$-form on $B \times \C^2$. 
\end{proof}

\paragraph{The Hodge correlator twistor connection $\nabla_{\bf G}$.} 
Lemma \ref{12.16.15.1} implies that the Hodge correlator  ${\bf G}$ gives 
rise to a connection  on $B \times \C^2$ 
with values in the Lie algebra $\pi^*({\cal C}{\cal L}ie_{X, S^*}(-1))$ 
in the category of variations of real Hodge structures: 
\be\la{12.15.14.22}
\nabla_{\bf G}:= {\bf d}+ {\rm d}_{\C^2}+ {\bf G}.
\ee
We call it {\it the Hodge correlator twistor connection $\nabla_{\bf G}$}. The name is justified by 
Theorem \ref{MDET}.

\begin{theorem} \la{MDET} The connection $\nabla_{\bf G}$ is a twistor connection 
on $B \times \C^2$. This just means that:

\begin{itemize}

\item The connection $\nabla_{\bf G}$  is $\C^* \times \C^*$-equivariant real connection 
on $B \times \C^2$.

\item  The restriction of the connection $\nabla_{\bf G}$ 
to $B \times \C$ is flat. Equivalently, the restriction ${\bf G}'$ of the Hodge correlator   
 ${\bf G} $ to the twistor line $z+w=1$ satisfies the Maurer-Cartan 
equation
\be \la{12.15.15.300}
({\bf d}+ {\rm d}_{\C}){\bf  G}'  + {\bf G}' \wedge {\bf G}'  =0.
\ee
\end{itemize}
\end{theorem}

The Hodge correlator class ${\bf G}$ in (\ref{12.16.15.2}) is nothing else but the 
 twistor transform of the data 
$$
{\bf G}^* := ({\bf G}_{0,0}, \{{\bf G}_{s,t}\}_{s,t\geq 1}).
$$ 
Theorem \ref{MDET} just means that  ${\bf G}^*$  is a Green data. 
Since $\pi^*({\cal C}{\cal L}ie_{X, S^*}(-1))$ is a Lie algebra 
in the category of variations of $\R$-Hodge structures, its  Hodge complex 
${\cal C}^{\bullet}_{{{\cal H}_\R}}({\cal C}{\cal L}ie_{X,S^*}(-1))$  is defined, and it has  
a DG Lie algebra structure, see Section 4.1. So we can restate Theorem \ref{MDET}  as follows:

\bt \la{2.12.13.11}
The data ${\bf G}^*$ defines a 1-cocycle   in the Hodge complex:
\be \la{12.31.ups.1aas}
{\bf G}^*  
\in {\cal C}_{{\cal H}_\R}^1({\cal C}{\cal L}ie_{X,S^*}(-1)).
\ee 
\et

One of the benefits of  Definition \ref{MDHCCL} 
is that it delivers the twistor connection $\nabla_{\bf G }$  directly.  
So the twistor transform is 
built naturally in the construction of Hodge correlators.

The construction from Section 2 delivers the Green data ${\bf G}^*$, and we need to apply to it the twistor transform, 
which seems to come out of the blue, 
  to get the twistor connection $\nabla_{\bf G }$. 

 \begin{figure}[ht]
\centerline{\epsfbox{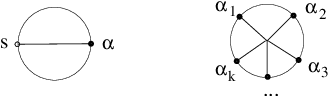}}
\caption{The decorated trees providing Hodge correlators of 
the Hodge bidegree $(0,0)$.}
\label{hc8}
\end{figure}

\paragraph{Remark.} Denote by $e_{T, W}$ the number of edges of a plane decorated tree $(T, W)$, not necessarily 
trivalent,  
which are not external edges with a $1$-form in the decoration. It
 is the number of the edges of $T$ to which we assign Green currents.  
Then $e_{T, W}$ is the N-degree of the Hodge correlator ${\rm Cor}_{\cal H}(T, W)$.  
Indeed, the $e_{T, W}$ is the total degree of the ${\rm Cor}_{\cal H}(T, W)$ as a 
polynomial in $(z,w)$, e.g. the degree of $zdw-wdz$ is $1$, which coinsides with the 
N-degree of the Hodge correlator ${\rm Cor}_{\cal H}(T, W)$. 
In particular, the  ${\bf G }_{0,0}$
 is the contribution of the decorated plane trivalent trees with $e_{T, W}=0$. Here are couple of examples. 

If $e_{T, W}=0$, then $(T, W)$ is one of the two decorated trees on Fig \ref{hc8}.  

\begin{figure}[ht]
\centerline{\epsfbox{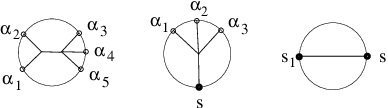}}
\caption{The decorated trees providing Hodge correlators of 
the Hodge bidegree $(-1,-1)$.}
\label{hc10}
\end{figure}

If $e_{T, W}=1$, then $(T, W)$ is one of the three decorated trees on Fig \ref{hc10}.  

\subsection{Proof of Theorem \ref{MDET}}  \label{6.3}

We prove Theorem \ref{MDET}  in three steps. 

1) Establish  Proposition \ref{12.16.15.1}. 
Show that the connection $\nabla_{\bf G}$ is $\C^*\times \C^*$-equivariant and real. 

2) Show that the Hodge correlator  ${\rm Cor}_{\cal H}$  
is multiplicative. 

3) Prove the Maurer-Cartan formula (\ref{12.15.15.300}) for the differential 
of the 
Hodge correlator ${\bf G}'$.

\paragraph{1) The $\nabla_{\bf G}$ is a $\C^*\times \C^*$-equivariant  real connection.} 
First, let us consider the first component of  ${\bf G }_{0,0}$. It is the restriction 
of the real Casimir element $\gamma | h$ to a section $s$, see (\ref{8.12.13.2}). 
So it is real, and its
 Hodge bidegree is $(0,0)$, see (\ref{8.12.13.2}). The rest is done as follows.

The Casimir element $\gamma | h$ is of the Hodge bidegree $(0,0)$. The form 
${\rm D}^\C G$ is of the Hodge bidegree $(1,1)$. For any plane trivalent tree $T$ the number of the edges 
$E$ on which we put the form ${\rm D}^\C G$ is one less then the number of the internal vertices. 
Indeed, this is true if all external edges are decorated by 1-forms; adding an 
$s$-decorated external edge amounts to  adding one internal edge, and one internal vertex. 
Since integration over the relative curve $X/B$ has the Hodge bidegree $(-1,-1)$, and since 
we twisted the result by $\R(-1)$, the total bidegree is $(0,0)$. This just means that the 
connection $\nabla_{\bf G}$ is $\C^*\times \C^*$-equivariant.

Recall the antiholomorphic involution of the twistor plane 
$\sigma: (z,w) \lra (\overline w, \overline z)$, 
and the complex conjugatoon $c$. 
For the 
involution $c\circ \sigma$, the both the operator ${\rm D}^\C$ and the Green current are imaginary, and 
therefore ${\rm D}^\C G$ is real:
$$
(c\circ \sigma)^*{\rm D}^\C = - {\rm D}^\C, ~~~~(c\circ \sigma)^*G = -G, ~~~~
 (c\circ \sigma)^*{\rm D}^\C G = {\rm D}^\C G. 
$$
The Casimir element $\gamma | h$ is evidently real. 
Therefore the integrand (\ref{3.8.05.15}) is 
real. This implies that the Hodge correlator class ${\bf G}$ is real. 
We conclude that the connection $\nabla_{\bf G}$ is  real. 

So if $s,t \geq 1$, the   ${\bf G}_{s,t}$  is a 
${\cal C}{\cal L}ie_{X, S^*}(-1)$-valued function  on $B$ of the Hodge bidegree $(-s,-t)$.   
The ${\bf G}_{0,0}$ is a ${\cal C}{\cal L}ie_{X, S^*}(-1)$-valued $1$-form on $B$.

\paragraph{2. Multiplicativity of the Hodge correlators.} 

The Hodge correlator construction can be applied to an oriented forest $(F; {\rm Or}_{F})$, 
resulting an element in the Hodge complex:  
\be \la{1.12.08.1}
{\rm Cor} _{\cal H}(F; {\rm Or}_F)\in 
 {\cal C}^{\bullet}_{\cal H}(S^\bullet({\cal C}_{X,S^*}(-1))).
\ee
Since $S^\bullet({\cal C}_{X,S^*}(-1))$ is a graded commutative 
algebra in the category of $\R$-Hodge structures, its Hodge complex 
${\cal C}^{\bullet}_{{\cal H}_\R}(S^\bullet({\cal C}_{X,S^*}(-1)))$ is a 
graded commutative  algebra,  with the product $\ast$. 

\bp \la{FII}
The  Hodge correlator map (\ref{1.12.08.1}) 
is multiplicative: 
$$
{\rm Cor}_{\cal H}((F'; {\rm Or}_{F'}) \circ (F''; {\rm Or}_{F''})) = {\rm Cor}_{\cal H}(F'; {\rm Or}_{F'}) \ast 
{\rm Cor}_{\cal H}(F''; {\rm Or}_{F''}). 
$$
\ep

\begin{proof} One has 
$
{\cal X}_{F'\times F''} = {\cal X}_{F'}  \times {\cal X}_{F''}.
$ 
There is a canonical direct product map of distributions
$$
{\cal D}^{*,*}({\cal X}_{F'}) \times {\cal D}^{*,*}({\cal X}_{F''}) \lra {\cal D}^{*,*}({\cal X}_{F'\times F''}). 
$$
The current $\widehat \kappa(F, W; {\rm Or}_F)$ assigned to a decorated forest $F$, see (\ref{3.8.05.15}), 
is evidently multiplicative: 
$$
\widehat \kappa(F', W'; {\rm Or}_{F'}) \times \widehat \kappa(F'', W''; {\rm Or}_{F''}) = 
\widehat \kappa(F'\times F'', W'\times W''; {\rm Or}_{F'\times F''}). 
$$
Recall the projection 
$
p_{F}: {\cal X}^{\rm int}_{F} \times \C^2 \lra B \times \C^2.
$ 
Provided that the $S^*$-decorations of $F'$ and $F''$ are disjoint, 
the multiplicativity implies that 
$$
p_{F'}\left(\widehat \kappa(F', W'; {\rm Or}_{F'})\right) \wedge  
p_{F''}\left(\widehat \kappa(F'', W''; {\rm Or}_{F''})\right)  = 
p_{F' \times F''}\left(\widehat \kappa(F'\times F'', W'\times W''; {\rm Or}_{F'\times F''})\right) . 
$$
Since the Hodge complex functor ${\cal C}_{\cal H}$ is multiplicative, this implies the claim. \end{proof}

\paragraph{3. Maurer-Cartan differential equations for the Hodge correlators.}  
Recall that the Hodge correlator  ${\rm Cor}_{\cal H}(T; {\rm Or}_T)$  
 is the push forward to $B \times \C^2$ 
of the current 
\be \la{1.08.08.4}
{\rm D}^\C G_{E_0}\wedge ... \wedge {\rm D}^\C G_{E_r} \wedge \gamma|h \wedge 
\ldots \wedge \gamma|h. 
\ee

\begin{theorem} \label{5.16.06.12g} The following diagram  is commutative: 
\begin{displaymath} \la{1.08.08.2a}
    \xymatrix{
        {\cal F}^{\vee, {\bullet}}_{X, S^*}(1)  \ar[r]^{\partial} \ar[d]_{{\rm Cor}_{\cal H}} & {\cal F}^{\vee, {\bullet}}_{X, S^*}(1)
\ar[d]^{{\rm Cor}_{\cal H}} \\
         {\cal D}^{*,*}_{B }\times \Omega^{\bullet}_{\C}  \ar[r]^{{\rm d}_B+{\rm d}_\C}       
& {\cal D}^{*,*}_{B }\times \Omega^{\bullet}_{\C}}
\end{displaymath}
This means that for any decorated oriented forest $(F, W; {\rm Or}_F)$ we have 
\be \la{1.08.08.10}
({\rm d}_B+{\rm d}_\C) {\rm Cor}_{\cal H}(F, W;  {\rm Or}_F)
= {\rm Cor}_{\cal H} \circ \partial(F, W;  {\rm Or}_F).
\ee 
\end{theorem}

\begin{proof}  It is sufficient to prove the result 
for an oriented
 decorated tree $(T, {\rm Or}_T)$. To calculate the differential ${\rm d}_B + {\rm d}_\C$ on $B \times \C^2$, 
we apply the de Rham differential 
${\rm d}+{\rm d}_\C$ to the integrand, and use the fact that the de 
Rham differential commutes with the proper push forward. 
The Casimir element $\gamma|h$ represents the identity map, 
and thus is annihilated by the differential. Therefore the calculation of the value of 
$({\rm d}+{\rm d}_\C)  {\rm Cor}_{\cal H}(T, W; {\rm Or}_T)$ 
boils down  to  the calculation
of    the differential 
\be \la{22.5.15.1}
({\rm d}+{\rm d}_\C) ~{\rm D}^\C G_{E_0}\wedge ... \wedge {\rm D}^\C G_{E_r}. 
\ee
Recall the Laplacian $\Delta = \partial''\partial'$. 
We employ formula (\ref{5.22.15.2}) to calculate (\ref{22.5.15.1}). Namely:
$$
({\rm d}+{\rm d}_\C) ~{\rm D}^\C G_{E_0}\wedge ... \wedge {\rm D}^\C G_{E_r} = 
\Delta G_{E_0}\wedge {\rm D}^\C G_{E_1}... \wedge {\rm D}^\C G_{E_r} + \ldots .
$$
We use formula 
(\ref{7.3.00.1}) for the  Laplacian:
$$
\Delta G_{s_0}(x,y) = \delta_{\Delta_X} - 
\Bigl(p_1^*\delta_{s_0} + p_2^*\delta_{s_0} - \frac{i}{2}\sum_{k=1}^g (p_1^*\alpha_k \wedge p_2^*\overline \alpha_k +
p_2^*\alpha_k \wedge p_1^*\overline \alpha_k) \Bigr).
$$ 
There are three terms in the formula for the Laplacian:

(a) The $\delta_{\Delta_X}$-term.

(b) The volume term $p_1^*\delta_{s_0} + p_2^*\delta_{s_0}$.

(c) The reduced Casimir term -- the remaining term. 

We call the sum of the terms (b)+(c) the Casimir term. 

Let us calculate 
the contribution of the Laplacian corresponding to an edge $E$. 
We will show that each term of the contribution matches 
certain term in the form ${\rm  Cor}_{\cal H} \circ \partial$.

\vskip 3mm
1. {\it  $E$ is an internal edge}. 
Then there are three cases, matching the terms 
in   $\Delta G$.

(1a) The $\delta_{\Delta_X}$-term contributes
the diagram with a special vertex,  obtained by shrinking the edge $E$, see the right of Fig \ref{feyn17}. 
This matches the graph complex differential $\partial_\Delta$.

(1b) Let us show that the contribution of the volume term 
 is zero.  Let $T_1$ and $T_2$ 
be the two trees obtained by cutting $T$ along the edge $E$. 
We assign the current $\delta_{s_0}$ to the end of one of them, say $T_1$. 
Denote by $W_1$ the induced decoration of $T_1$. 
Let us show that 
${\rm Cor}_{{\cal H}, v_0}(T_1, W_1; {\rm Or}_{T_1}) =0$. 

Recall the regularisation procedure called the specialization at a tangent vector, Section \ref{2.2}.  
One defines  specializations of 
forms by taking specializations of their coefficients. 
\bl \la{1.11.08.1}
Let $W = {\cal C}(\{s\} \otimes W_0)$. Then 
the specialization ${\rm Sp}^{s\to s_0}_{v_0}{\rm Cor}_{{\cal H}, v_0}(W)=0$. 
\el

\begin{proof} If $W$ is different from 
${\cal C}(\{s\} \otimes \alpha)$ this follows from 
the property (\ref{1.11.08.2}) of the Green function 
assigned to the edge decorated by $\{s\}$. The restriction of  
${\rm Cor}_{{\cal H}} {\cal C}(\{s\} \otimes \alpha)$ 
to the divisor $s=s_0$ is zero by the very definition. 
\end{proof}  

The decoration 
$W_1 = \{a\} \otimes W'_1$  delivers an {\it a priori} divergent 
correlator ${\rm Cor}_{{\cal H}, a}(\{a\} \otimes W'_1)$. 
We understand it via the specialization 
of a correlator ${\rm Cor}_{{\cal H}, s_0}(\{a\} \otimes W'_1)$ 
when $a \to s_0$ at the tangent vector $v_0$ at $s_0$. 
This specialization is zero by Lemma \ref{1.11.08.1}. 
This settles the (1b) case.

(1c) The contribution of the reduced Casimir term matches the form 
assigned by the map ${\rm Cor}_{\cal H}$ to the minus of the 
sum of  forests on 
Fig \ref{feyn18}. 
The latter sum is the contribution of the edge $E$ to the map $\partial_{\rm Cas}$.

\begin{figure}[ht]
\centerline{\epsfbox{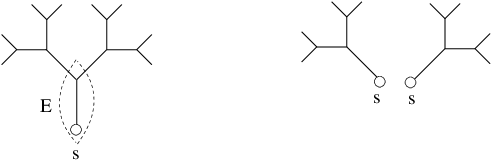}}
\caption{Contribution of the $\delta_{\Delta_X}$-term for
an  external ${\rm S}$-decorated edge.}
\label{feyn20}
\end{figure}

2. {\it $E$ is an external edge decorated by a section $s$}.
Then there are the following cases. 

(2a) The contribution of the $\delta_{\Delta_X}$-term 
in $\Delta G_E$ for a tree $T$ is calculated as follows. 
Let $s$ be the decoration at the external vertex $v_E$ of $E$. Let us cut out a small neighborhood of $E$ 
from $T$, getting  connected trees $T_1, \ldots , T_k$, see Fig \ref{feyn20}. 
Each of them has a new external 
vertex, 
and the other external vertices are the ones of $T$ minus $v_E$. 
The tree $T_i$ is decorated:  
the new vertex is decorated by $s$, 
and the others inherit their decorations from $T$. 
This matches the component $\partial_S$ of the differential $\partial$. 

(2b) 
The contribution of the volume term 
in $\Delta G_E$ is zero. It is similar 
to  (1b), except that one needs to add an argument for 
the contribution of $\delta_{s_0}$ related to an external vertex $v$ of $E$, decorated by 
$s$. The contribution is $\delta(s_0-s)$.  
It contributes zero since the divisor $S$ is smooth over $B$, 
so, moving over  $B$, the point $s$ can
 not collide with $s_0$.

(2c) It is similar to  (1c). 
The contribution of the reduced Casimir term 
in $\Delta G_E$ 
is  negative of the form on $X ^{\{\mbox{$S$-vertices}\}}$ 
given by the following recepee, 
illustrated on Fig \ref{feyn19}. 
Cutting the edge $E$, we 
get a  connected tree $T'$, and a single edge graph. 
The latter is decorated by a form 
$\alpha_k^{\vee}$ and $s$.  The form $\alpha_k^{\vee}$ is viewed as 
 a form 
 on the $s$-copy of $X $. 
The tree $T'$ is decorated: the inherited from $T$ external vertices 
 inherit their decoration, and 
the new external vertex 
is decorated by $\alpha_k$. 
It 
 gives rise to a Hodge correlator, which is a form on  
$$
X ^{\{\mbox{$S$-decorated vertices of $T$}\} - \{s\}}.
$$ 
The sum  over $k$ of 
the product of this correlator and the form above delivers the desired form.

\end{proof}

Denote by $\delta_{\rm Colie}$ the cocommutator in the Lie coalgebra 
${\cal C}_{X, S^*}(-1)$.

\bc \la{12.15.21.1}
The following diagram is commutative:   
\begin{displaymath} \la{1.08.08.2a}
    \xymatrix{
        {\cal C}^{\vee}_{X, S^*}[-1](1)  ~\ar[r]^{\delta_{\rm Colie}} \ar[d]_{F} & 
~~S^2({\cal C}^{\vee}_{X, S^*}[-1](1))
\ar[d]^{S^2F} \\
{\cal T}^{\vee, {\bullet}}_{X, S^*}(1)  \ar[r]^{\partial} \ar[d]_{{\rm Cor}_{\cal H}} & {\cal F}^{\vee, {\bullet}}_{X, S^*}(1)
\ar[d]^{{\rm Cor}_{\cal H}} \\
         {\cal D}^{*,*}_{B }\times \Omega^{\bullet}_{\C}  \ar[r]^{{\rm d}+{\rm d}_\C}       
& {\cal D}^{*,*}_{B }\times \Omega^{\bullet}_{\C}}
\end{displaymath}
\ec

\begin{proof}
The bottom square is commutative by Theorem \ref{5.16.06.12g}, the top  by 
Lemma \ref{6.19.06.1}. 
\end{proof}
 
By the very definition, the element ${\bf G}\in {\cal C}^1_{{\cal H}_\R}({\cal C}_{X,S^*}(-1))$ 
is obtained by dualising the composition 
of the left vertical maps. 

Recall 
the Hodge complex differential  $\delta_{\rm Hod}$. Similarly, the element $-\delta_{\rm Hod}{\bf G}$ 
is 
obtained by dualising  the  triple composition: 
$$ 
 {\cal C}^{\vee}_{X, S^*}[-1](1)\stackrel{F}{\lra} {\cal T}^{\vee, {\bullet}}_{X, S^*}(1)  
\stackrel{{\rm Cor}_{\cal H}}{\lra} {\cal D}^{*,*}_{B }\otimes\Omega_{\C^2}^\bullet 
\stackrel{{\rm d}_B+{\rm d}_\C}{\lra} {\cal D}^{*,*}_{B }\otimes\Omega_{\C^2}^\bullet. 
$$ 
By Corollary \ref{12.15.21.1}, the triple composition coincides with the other triple composition:
$$ 
 {\cal C}^{\vee}_{X, S^*}[-1](1)\stackrel{\delta_{\rm Colie}}{\lra} S^2({\cal C}^{\vee}_{X, S^*}[-1](1))
\stackrel{S^2F}{\lra} {\cal F}^{\vee, {\bullet}}_{X, S^*}(1) 
\stackrel{{\rm Cor}_{\cal H}}{\lra} {\cal D}^{*,*}_{B }\otimes\Omega_{\C^2}^\bullet. 
$$ 
 By Proposition \ref{FII}  the product 
 ${\bf G} \wedge {\bf G} \in {\cal C}^2_{{\cal H}_\R}({\cal C}_{X,S^*}(-1))$ 
is obtained by dualising the latter composition.   This implies that
$\delta_{\rm Hod}{\bf G} + {\bf G} \wedge {\bf G}=0$. 
Theorem \ref{MDET} is proved.

\section{The Lie algebra of special derivations} \label{hc7sec}

The pronilpotent completion  $\pi^{\rm nil}_1(X-S, v_0)$ was described in Section 1.2. 
It is a pro-Lie algebra over $\Q$. 
It carries a weight filtration. The corresponding associate 
graded ${\rm L}_{X, S}$ 
has a simple description in terms of the cohomology of  $X$ and  $S$. 
We describe the 
Lie algebra of all special derivations of ${\rm L}_{X, S}$, and 
identify it with the  Lie algebra 
${\cal C}_{X,S}$ from Section \ref{hc5sec}.

\subsection{The linear algebra set-up} 
Let $H$ be a 
finite dimensional 
vector space with a symplectic structure 
$\omega^{\vee} \in \Lambda^2H^{\vee}$. 
Let $\{p_i, q_i\}$ be a symplectic basis in $H$: $\omega^{\vee} (p_i, q_i)=1, \omega^{\vee} (p_i, p_j)= \omega^{\vee} (q_i, q_j)= 0$. 
There is the dual  of $\omega^{\vee}$: 
$$
w = \sum_{i}p_i \wedge q_i \in \Lambda^2 H.
$$

Let $S$ be a finite set.  We need the following 
associative/Lie algebras:

\begin{enumerate}

\item 
$A_{H,S}$: the tensor algebra of  $H \oplus \Q[S]$. 
Then  $\omega$ provides an element 
$
[p,q]:= \sum_i [p_i, q_i] \in A_{H,S}.
$ 
Denote by  $X_s$ the generator corresponding to $s \in S$.

\item $L_{H,S}$: the free Lie algebra generated by $H\oplus \Q[S]$. 
The algebra $A_{H,S}$ is identified with its universal enveloping 
algebra. 

\item $\overline A_{H,S}$: the quotient of the algebra $A_{H,S}$ 
by the two-sided ideal generated by the element 
\begin{equation} \label{relk}
[p,q] + \sum_{s \in S} X_s. 
\end{equation}

\item  $\overline L_{H,S}$: the  quotient of $L_{H, S}$ by 
the ideal generated by (\ref{relk}). 
\end{enumerate}
When $S$ is empty, set    $A_H:= A_{H,\emptyset}, 
L_H:= L_{H, \emptyset}$, etc. 

Let $S$ assume from now on that the set $S$ is non-empty.  
Let $0$ be an element of $S$. Set $S^* := S - \{0\}$. 
Then there are canonical isomorphisms 
\begin{equation} \label{6.26.00.10}
L_{H, S^*}\quad \stackrel{\sim}{\lra} \quad \overline L_{H,S}, \qquad  A_{H, S^*}\quad \stackrel{\sim}{\lra} \quad \overline A_{H,S}. 
 \end{equation} 
Indeed, one can express $X_0$ via the other generators:
\begin{equation} \label{6.26.00.10*}
X_0 = -\sum_{s \in S^*}X_s - [p,q].
\end{equation}

\paragraph{Example.} Let $X$ be a smooth 
complex compact curve, $H:= H_1(X; \Q)$. The symplectic form $\omega^{\vee}$ 
is the intersection form 
on $H_1(X; \Q)$. Let $S$ be a non-empty subset of points of $X$. 
Then $\overline L_{H,S}$ is 
isomorphic, non-canonically, to 
$\pi_1^{\rm nil}(X-S, v_0)$. 
The $\overline L_{H,S}$  is canonically 
identified with the associate graded for the weight filtration ${\rm L}_{X, S}$ 
of the $\pi_1^{\rm nil}(X-S, v_0)$.

\subsection{Special derivations of the algebra $\overline A_{H,S}$}  \la{sec7.2}
Constructions of Section \ref{sec7.2} generalize the ones of Drinfeld  [Dr] 
and Kontsevich [K].  
A derivation ${\cal D}$ of the algebra 
$\overline A_{H,S} \stackrel{(\ref{6.26.00.10})}{=} A_{H,S^*}$ 
is {\it special} if there exist elements $B_s$, 
$s \in S^*$ such that 
$$
 {\cal D}(X_0) =0, \quad {\cal D}(X_s) = [B_s, X_s], \quad s \in S^*. 
$$
Special derivations form a Lie algebra denoted 
${\rm Der}^SA_{H,S^*}$. 
Let us emphasize that we define it 
using the isomorphism (\ref{6.26.00.10}).

\paragraph{Remark.} If  $H \not =0$,  a 
collection of  elements $ \{B_{s}\}$  
does not determine 
a special derivation. 
For instance if $S^*$ is empty this collection is also empty  
while,  as we show below, the space of special derivations is not.  If $H =0$ it does, see 
Section 4 in [G4].

\vskip 3mm
 Let $A$ be an associative algebra. Recall  the projection
${\cal C}: A \lms {\cal C}(A) = A/[A,A]$. 

If $B$ is a free associative algebra generated by a finite 
set ${\cal S}$, then 
there are linear maps
$$
{\partial }/{\partial X_s}: {\cal C}(B) \lra B; \qquad 
{\cal C}(X_{s_1}X_{s_2} \cdot \cdot \cdot  X_{s_k}) \lms \sum_{s_i = s} X_{s_{i+1}}X_{s_{i+2}} \cdot \cdot \cdot X_{s_{i-1}}.
$$
These maps are components of the non-commutative differential map
\begin{equation} \label{6.27.00.10ds}
{\rm D}: {\cal C}(B) \lra B\otimes \Q[S], \qquad {\rm D}(F) = 
\sum_{s\in S}\frac{\partial F}{\partial X_s}\otimes X_s.
\end{equation}

According to   formula (96) in Section 4 of [G4] (an easy exercise) one has 
\begin{equation} \label{6.27.00.10}
\sum_{s\in S}[\frac {\partial F}{\partial X_s}, X_s] =0.
\end{equation}
 
In particular, there are linear maps
$
{\partial }/{\partial p_i}, {\partial }/{\partial q_i}, {\partial }/{\partial X_s}: 
{\cal C}(A_{H, S^*} ) \to A_{H, S^*}$,  $s \in S^*$.

\paragraph{Remark.}
 We can not define these maps for 
$s \in S$ acting on $\overline A_{H,S}$ since relation (\ref{relk}) will not be killed. 
\vskip 3mm

The vector space ${\cal C}(A_{H,S^*})$ is  decomposed into a direct sum  of 
$\Q\cdot 1$ and of ${\cal C}^+(A_{H,S^*})$. 

\begin{lemma} \label{10.27.04.1} There is a map 
$$
\kappa: {\cal C}^+(A_{H,S^*}) \lra {\rm Der}^S(A_{H,S^*}).
$$
such that given an element $F \in {\cal C}^+(A_{H,S^*})$, 
 the derivation $\kappa_F$ acts on the generators 
by
\begin{equation} \label{12.2.04.1}
X_0 \lms 0; \quad p_i \lms -\frac{\partial F}{\partial q_i}; \quad q_i \lms \frac{\partial F}{\partial p_i}; \quad 
X_s \lms [ X_s, \frac{\partial F}{\partial X_s}], \qquad s \in S^*.
\end{equation}\end{lemma}

\begin{proof}
It follows from (\ref{6.27.00.10}) that 
$
\sum_i[\frac{\partial F}{\partial q_i}, q_i] + \sum_i[\frac{\partial F}{\partial p_i}, p_i] + \sum_{s \in S^*}
[\frac{\partial F}{\partial X_s}, X_s] =0.
$ 
This is equivalent to 
$
\kappa_F\Bigl([p,q]+ \sum_{s \in S^*}X_s\Bigr) =0
$,  
 which just means that $\kappa_F(X_0) =0$, and thus  $\kappa_F$ 
is a special derivation. 
\end{proof}

The space ${\cal C}(A_{H,S^*})$ is identified, as a vector space, 
with the Lie algebra $  {\cal C}_{H,S^*}$ from Section 
\ref{hc5sec}. So it  
inherits a Lie algebra structure. 

\bp \la{7.4.08.1}
a) The Lie bracket on ${\cal C}(A_{H,S^*})$ 
is given by the following formula: 
\begin{equation} \label{6.20.00.11}
\{F, G\}:= \quad {\cal C}\left( 2\sum_{s \in S^*}[\frac{\partial F}{\partial X_s}, 
\frac{\partial G}{\partial X_s}]\cdot X_s  \quad + \quad \sum_i\Bigl(\frac{\partial F}{\partial p_i}\cdot 
\frac{\partial G}{\partial q_i} - \frac{\partial F}{\partial q_i}\cdot 
\frac{\partial G}{\partial p_i}\Bigr)\right).
\end{equation}

b)  The map $\kappa$  is a Lie algebra morphism. 
\ep

Notice that (\ref{6.20.00.11}) does not depend on the choice of the symplectic basis $\{p_i, q_i\}$ in $H$. 

\begin{proof} a) Follows immediately from the definitions. 
The first term in formula (\ref{6.20.00.11}) corresponds to the 
map $\delta_{S}$, and the second $\delta_{\rm Cas}$. In particular this shows that 
the bracket  (\ref{6.20.00.11}) satisfies the Jacobi identity. 

b) We will check it in the more 
general motivic set-up in Section \ref{hc8sec}. 
\end{proof}

Let ${\cal C}'(A_{H,S^*})$ be the quotient 
of the vector space ${\cal C}(A_{H,S^*})$ by $\oplus_{s\in S^*}\Q[X_s]$. 
The Lie bracket  (\ref{6.20.00.11}) descends to a Lie bracket on 
${\cal C}'(A_{H,S^*})$,

\begin{proposition} \label{ssw}
The map $\kappa': {\cal C}'({A_{H,S^*}}) \lra 
{\rm Der}^S{A_{H,S^*}}$ is a Lie algebra  isomorphism. 
\end{proposition}

\begin{proof}  Since the algebra ${\rm A}_{H, S^*}$ is free, the centralizer of $X_s$  
is the algebra $\Q[X_s]$ of polynomials in $X_s$. 
So ${\rm Ker} \kappa = \oplus_{s} \Q[ X_s]$. 
We get  an injective Lie algebra map
$ 
 \kappa': {\cal C}'({\rm A}_{H, S^*})
\hookrightarrow {\rm Der}^S{\rm A}_{H, S^*}.
$ 
So it remains to check that $\kappa'$ is surjective. 
We reduce this to the known case when ${\rm dim}H=0$. 
Observe that a system $B = (B_{q_i}, B_{p_i}, B_s)$ of elements of 
$A_{S^*}$, where 
$i = 1, ..., {\rm dim} H$ and $s \in S$, satisfying  
\begin{equation} \label{10.31.04.2}
\sum_{i=1}^{{\rm dim}H}([B_{q_i}, q_i] + 
[B_{p_i}, p_i]) +  \sum_{s \in S}[B_s,X_s]) =0
\end{equation}
defines a special derivation $D_B$ of the algebra 
$A_{S^*}$, acting on the generators as follows: 
$$
D_B: p_i \lms - B_{q_i}, \quad q_i \lms  B_{p_i}, \quad X_s \lms  B_s. 
$$
Indeed, $D_B([p,q] + \sum_{s \in S}X_s)$ is given by the left hand side of  
(\ref{10.31.04.2}). So by Proposition 5.1 of \cite{G4}, 
there exists an $F \in {\cal C}(A_{S^*})$ such that 
$B_{q_i} = -\partial F/\partial q_i$, 
$B_{p_i} = \partial F/\partial p_i$, and 
$X_{s} = \partial F/\partial X_s$.
\end{proof}

Just like in (\ref{mvla}), the Lie algebra ${\cal C}{\cal L}ie_{H,S^*}$ 
is a Lie subalgebra of ${\cal C}({A_{H,S^*}})$, defined 
as follows:
\be \la{mvlca1}
{\cal C}{{\cal L}ie}^{\vee}_{X, S^*}:= 
\frac{ {\cal C}({\rm A}^{\vee}_{H, S^*})}{\mbox{Shuffle relations}};~~~~
{\cal C}{{\cal L}ie}_{H, S^*}:= \mbox{the dual of 
${\cal C}{{\cal L}ie}^{\vee}_{H, S^*}$}.
\ee

\begin{proposition} \label{sswas}
The map $\kappa: {\cal C}{\cal L}ie_{H,S^*} \lra 
{\rm Der}^S{L_{H,S^*}}$ is a Lie algebra  isomorphism. 
\end{proposition}

\begin{proof}   Let us consider a Casimir element, given by a formal infinite sum
$$
\xi:= \sum_W \frac{1}{|{\rm Aut}(W)|}W^{\vee} \otimes W. 
$$
Here the sum is over a basis $\{W\}$ in  ${\cal C}_{H, S^*}$, and 
$\{W^{\vee}\}$ is the dual basis. Here ${\rm Aut}(W)$ is the automorphism 
group of  the cyclic word $W$. 
There is a coproduct $\Delta$ on the algebra ${\rm A}_{H, S^*}$, making it into a Hopf algebra, 
 determined by 
the property that the generators of the algebra are primitive. 
It dualizes the shuffle product $\circ_{\rm Sh}$. 
Let $\overline \Delta(Z):= \Delta(Z) - (1\otimes Z + Z \otimes 1)$ 
be the reduced coproduct. 
Then ${\rm L}_{H, S^*} = {\rm Ker}\overline \Delta$. 
Choose a basis $\{Y_i\}$ in ${\rm V}_{H, S^*}$.
 Given a basis vector $Y_i$, consider the expression 
\be \la{1.16.08.1}
{\rm id} \otimes (\overline \Delta \circ {\cal D}_{Y_i}(\xi)) \in 
{\cal C}^{\vee}_{H, S^*} \otimes {\rm A}_{H, S^*}.
\ee
Choose two elements $A,B$ of the induced basis in ${\rm A}_{H, S^*}$. Then 
$A\otimes B \in {\rm A}_{H, S^*}\otimes {\rm A}_{H, S^*}$ appears in 
(\ref{1.16.08.1}) with coefficient ${\cal C}(Y(A\circ_{\rm Sh} B))$. 
Notice that 
the factor $1/|{\rm Aut}(W)|$ is necessary. 
\end{proof}

\subsection{A symmetric description 
of the Lie algebra  of special derivations} \la{hc7.3sec}

To define a special derivation we have to choose a specific element, $0$, of the set $S$. 
Let us show that 
 the Lie algebra of all special derivations can be 
defined using 
only the symplectic vector space $H$ and the set $S$. However the way it acts 
by special derivations of $A_{H,S^*}$ depends on the choice of 
$0 \in S$. 

Let us  define a Lie algebra structure on the vector space 
${\cal C}(A_{H,S})$ via formula (\ref{6.20.00.11}) where 
$S^*$ is replaced by  $S$. 
Then  we prove the following result.

\begin{proposition} \label{6.27.00.30}
Take the two-sided ideal in the algebra $A_{H,S}$,  
generated by $[p,q] + \sum_{s\in S} X_s$.  Its  image in the Lie algebra 
${\cal C}(A_{H,S})$ is an ideal. 
So the quotient  ${\cal C}(\overline A_{H,S})$ of ${\cal C}(A_{H,S})$ 
by this ideal is a Lie algebra. 
The canonical map $
i_0: {\cal C}(\overline A_{H,S}) \lra {\cal C}(A_{H,S^*})$ 
is a Lie algebra isomorphism.  
\end{proposition}

\begin{proof}  It is an easy exercise left to the reader.\footnote{A natural proof of 
is given by deducing it from a general result in the DG version of the story.} 
\end{proof}

So for any element $s \in S$ there is a Lie 
algebra isomorphism 
\be \la{1.26.08.1}
\kappa_s: {\cal C}(\overline A_{H,S}) \stackrel{\sim}{\lra }
{\rm Der}^{S}A_{H, S - \{s\}}.
\ee
The isomorphism 
$
i_s: A_{H, S -  \{s\}} \to \overline A_{H,S}
$ 
provides a Lie algebra isomorphism
$
{\cal C}(A_{H, S -  \{s\}}) \to {\cal C}(\overline A_{H,S}).
$ 
Furthermore, for any $s \in S$ isomorphism (\ref{1.26.08.1}) 
restricts to a Lie 
algebra isomorphism 
\be \la{1.26.08.2}
\kappa_s: {\cal C}{\cal L}ie(\overline L_{H,S}) \stackrel{\sim}{\lra }
{\rm Der}^{S}L_{H, S - \{s\}}.
\ee
There is a Lie algebra isomorphism
$
{\cal C}{\cal L}ie(L_{H, S -  \{s\}}) \to {\cal C}(\overline L_{H,S}).
$ 

\paragraph{The Lie algebra of outer semi-special derivations.}   
We say that a derivation ${\cal D}$ of the associative algebra 
$\overline A_{H,S}$ is {\it semi-special} 
if it preserves the conjugacy classes 
of each of the generators $X_s$, $s \in S$, i.e. if there exist 
elements $B_s$ of the algebra such that 
$$
{\cal D}(X_s) = [B_s, X_s], \quad s\in S.
$$
The Lie algebra of semi-special derivations is denoted by 
${\rm Der}^{SS}\overline A_{H,S}$. An inner derivation is a particular example 
of a semi-special derivation. We define the Lie algebra of 
{\it outer semi-special} derivations of the algebra 
$\overline A_{H,S}$ by taking the quotient modulo the Lie subalgebra of inner derivations: 
$$
{\rm ODer}^{SS}\overline A_{H,S}:= \frac{{\rm Der}^{SS}\overline A_{H,S}}{{\rm In}{\rm Der}^{SS}\overline A_{H,S}}.
$$

For an $s \in S$, let 
$
{\rm Der}_s^{S}\overline A_{H,S}\subset {\rm Der}^{SS}\overline A_{H,S}
$  be the Lie subalgebra 
of all semi-special derivations which are special for the generator $X_s$, i.e. kill $X_s$. 
Thus there is a Lie algebra map
$$
p_s: {\rm Der}_s^{S}\overline A_{H,S} \lra {\rm ODer}^{SS}\overline A_{H,S}.
$$

\begin{lemma} \label{12.15.04.1}
The map $p_s$ is surjective. Its kernel is spanned 
by the inner derivation $[X^n_s, *]$. 
\end{lemma}

\begin{proof}  Let ${\cal D}$ be a semi-special derivation,  such that 
${\cal D}(X_s) = [B_s, X_s]$. Subtracting the inner derivation 
$[B_s, *]$ from ${\cal D}$ we get a derivation special with respect to 
the generator $X_s$. So $p_s$ is surjective. Since 
$\overline A_{H, S-\{s\}}$ is a free associative algebra, we get the second claim. 
\end{proof} 

\section{Variations of mixed 
$\R$-Hodge structures via Hodge correlators} \label{hc9sec}

 \paragraph{A variation of $\R$-MHS on $\pi^{\rm nil}(X-S, v_0)$ via Hodge correlators.} 
Recall that $S = S^* \cup \{s_0\}$. Recall (Section \ref{sec1.9}) 
the enhanced moduli space 
${\cal M}'_{g, n}$ 
parametrising collections $S$ of 
$n$ distinct points on a genus $g$ complex curve $X$ plus a tangent vector at $s_0$. 
There is a local system ${\cal L}$ over ${\cal M}'_{g, n}$ 
with the fiber 
\be \la{20:33}
{\rm L}_{X,S^*; v_0}\stackrel{\sim}{=} {\rm gr}^W{\pi^{\rm nil}_1}(X-S, v_0).
\ee 
It is equipped with the Gauss-Manin connection.

A tangent vector $v_0$ at the  point $s_0\in X$ 
provides the  
normalized Green function $G_{v_0}(x,y)$ (Section \ref{2.2}). 
It gives rise to  the  Hodge correlator  map, and hence to an element 
$$
\widetilde {\bf G}_{v_0} \in {\cal C}^1_{\cal H}({\cal C}{\cal L}ie_{X,S}). 
$$
The Lie algebra ${\cal C}_{X,S}$ acts 
by 
special derivations of the Lie algebra 
${\rm L}_{X,S^*; v_0}$ (Section \ref{hc7sec}). So we arrive at an endomorphism 
${\bf G}_{v_0}$ of the smooth bundle ${\cal L}_{\infty}$. 
We proved in Section \ref{hc6sec}   that it satisfies the Maurer-Cartan 
differential equation 
$\delta {\bf G}_{v_0} + [{\bf G}_{v_0}, {\bf G}_{v_0}] =0$. Therefore 
the  construction of Section \ref{hc4sec} provides 
a variation of mixed $\R$-Hodge structures. 
Let us recall Theorem \ref{9.16.05.2}:
 
\begin{theorem} \label{9.16.05.2aa} When the data 
$(X, S, v_0)$ varies, 
the 
variation of mixed $\R$-Hodge structures given by the Green datum 
${\bf G}_{v_0}$ 
is isomorphic to the classical 
variation of the  
mixed $\R$-Hodge structures on $\pi_1^{\rm nil}(X -S, v_0)$. 
\end{theorem}

\paragraph{A simple proof of Theorem \ref{9.16.05.2aa} for rational curves.} 
Denote by  $G^{\rm Cor}_{p,p}$ the 
Green operators for the Hodge correlators on ${\Bbb P}^1-S$. 
The standard variation of  mixed $\R$-Hodge structures 
$\pi_1^{\cal H}({\Bbb P}^1-S, v_0)$ is described 
by the Green operators denoted by $G^{\rm st}_{p,p}$. 
Let us prove by induction on $p$ that 
$G^{\rm Cor}_{p,p} = G^{\rm st}_{p,p}$. 
The base of the induction is given by the following straightforward Lemma: 
\begin{lemma}
One has $G^{\rm Cor}_{1,1} = G^{\rm st}_{1,1}$.  
\end{lemma}

Suppose that the claim is proved for $p<n$. The systems of differential equations 
for Green operators describing variations of $\R$-MHS implies 
that the endomorphisms 
$G^{\rm Cor}_{n,n}$ and $G^{\rm st}_{n,n}$ satisfy the same differential 
equations, and  that  
$G^{\rm Cor}_{n,n} - G^{\rm st}_{n,n}$  
is annihilated by both $\partial$ and $\overline \partial$, 
and thus is a constant. 
To check that this constant is zero we use the shuffle relations 
from Proposition \ref{7.3.06.4}, 
which imply that a non-zero 
multiple of this constant is zero. 
This proves the theorem. 

\vskip 3mm

This argument  does not work 
in the non-Tate case since we do not impose $\partial$-differential equations on $G_{1,q}$, as well as 
$\overline \partial$-differential equations on $G_{p, 1}$. 
Let us address now the general case of Theorem \ref{9.16.05.2}.

\paragraph{$\R$-MHS on $\pi_1^{\rm nil}(X-S, {a})$ with a regular base point $a$ via Hodge correlators.} 
Let $a\in X-S$. 
Choose a non-zero tangent vector $v_a$ at $a$. There is a canonical isomorphism 
$$
\pi_1^{\rm nil}(X-S, {a}) = \pi_1^{\rm nil}(X-S\cup \{a\}, v_{a})/(v_a) 
$$
where $(v_a)$ is the ideal generated by $v_a$. 
The Green operator  on ${\rm gr}^W\pi_1^{\rm nil}(X-S\cup \{a\}, v_{a})$ is
a special derivation, and thus preserves the ideal 
${\rm gr}^W(v_a)$. So it induces an operator ${\bf G}_{v_a}$ 
on the quotient, which gives rise to an $\R$-MHS on $\pi_1^{\rm nil}(X-S, {a})$.

\paragraph{Good unipotent variations.} R. Hain and S. Zucker \cite{HZ} introduced a notion of a 
{\it good unipotent variation of mixed Hodge structures} on an open algebraic variety 
$U$, and proved that 
the category of  such variations  is equivalent to the category  
of representations  $\pi_1^{\rm nil}(U, u) \otimes V \lra V$ 
which are morphism of mixed Hodge structures.  Equivalently, 
they are representations  in the category mixed Hodge structures of 
the universal enveloping algebra ${\rm A}^{\cal H}(U, u)$ of 
$\pi_1^{\rm nil}(U, u)$ with the induced mixed Hodge structures. 

Recall (\cite{H1}, Section 7) that a unipotent variation of mixed Hodge structures 
${\cal V}$ over a smooth complex curve $X = \overline X -D$, where 
 $\overline X$ is a smooth compactification, and $D$ a finite number of points, 
is a good unipotent one if the following conditions at infinity hold:   

\begin{itemize} 

\item (i) Let $\overline {\cal V} \to \overline X$ be Deligne's canonical extension 
of the local system ${\cal V}$. The Hodge bundles ${F}^p{\cal V}$ 
are required to extend to holomorphic subbundles ${F}^p\overline {\cal V}$
of $\overline {\cal V}$. 

\item (ii) Let $N_P:= \frac{1}{2\pi i}\log T_P$ be the logarithm of 
the local monodromy $T_P$ around  $P \in D$. Then 
$$
N_P(W_lV_x) \subset W_{l-2}V_x.
$$
\end{itemize}
The second condition is so simple because we assumed that the 
global monodromy is unipotent.

\bt \la{2.27.08.3}
The variation ${\cal L}$ 
of $\R$-MHS on $\pi_1^{\rm nil}(X-S, {a})$  defined by the 
operator ${\bf G}_{v_a}$ acting on (\ref{20:33}), 
 is a good unipotent variation of mixed Hodge structures over $a\in X-S$.
\et

We introduce a DGA of forms with logarithmic singularities, and use 
it to prove Theorem \ref{2.27.08.3}. 
Then we show how Theorem \ref{2.27.08.3} implies Theorem \ref{9.16.05.2}.

\paragraph{Differential forms on $X$ with tame logarithmic 
singularities.} Let $D$ be a normal crossing divisor in a smooth complex projective variety 
$\overline X$, and $X:= \overline X - D$. Let us define a DG subalgebra 
${\cal A}_{\log}(X)$ of differential forms on $X$ with {\it tame logarithmic 
singularities} on $X$. It is a slight modification of the space 
used in 1.3.10 in \cite{L}. It is a DG subalgebra of the de Rham DG algebra 
of smooth forms on $X$. Choose local equations $z_i$ of the 
irreducible components $D_i$ of the divisor $D$. 
The space ${\cal A}_{\log}(X)$ consists of 
forms $\omega$ which can be represented as polynomials in 
$\log|z_i|$, $\partial \log|z_j|$, and $\overline \partial \log|z_k|$,  
$$
\omega = \sum_{a_i, \varepsilon'_j, \varepsilon''_k} \omega_{a_i, \varepsilon'_j, \varepsilon''_k}
\log^{a_i}|z_i| \wedge \bigwedge_j(\partial \log|z_j|)^{\varepsilon'_j} \wedge 
\bigwedge_k(\overline \partial \log|z_k|)^{\varepsilon''_k}, 
$$
whose coefficients are smooth functions on $\overline X$ with the following property: 
a coefficient of a monomial containing (in positive degree) any of the three 
expressions $\log|z_i|$, $\partial \log|z_i|$, and $\overline \partial \log|z_i|$ 
related to the stratum $D_i$ vanishes at the stratum $D_i$. 
For example, $z_1 \log|z_1|$ belongs ${\cal A}_{\log}(X)$, while 
$z_1 \log|z_2|$ does not. The space ${\cal A}_{\log}(X)$ clearly does not depend on the choice of local equations $z_i$. The space ${\cal A}_{\log}(X)$ is closed under 
$\partial$ and $\overline \partial$. 

\bp \la{3.8.08.1}
Choose a local parameter $t$ at $s$ such that $dt$ is dual to 
the tangent vector $v_s$ at $s$. Then ${\bf G}_{v_a} - \log t \cdot {\rm ad}(X_s)$  
is a smooth function 
 with tame logarithmic singularities near the point $s$. 
\ep

\bc \la{popopop}
The logarithm of the monodromy of the local system 
${\cal L}$ around  $s$ is ${\rm ad}(X_s)$.
\ec

\begin{proof} Indeed, $d^\C\log|t| = i d\arg(t)$, and $d^\C$ of the other components 
of ${\bf G}_{v_a}$ vanish at $t=0$ since these components have tame logarithmic singularities
at $t=0$. \end{proof}

We prove Proposition \ref{3.8.08.1} later on. Let us show first how it implies 
 Theorems \ref{2.27.08.3} and  \ref{9.16.05.2}. 

\paragraph{Proof of Theorem \ref{2.27.08.3}.}
The associate graded ${\rm gr}^W{\cal L}$ 
is a constant variation. So the  
local system ${\cal L}$ is unipotent. Thus we must check only  
conditions (i) and (ii).  

Condition (ii) follows immediately from Corollary \ref{popopop}. 

Condition (i) is just a condition on the 
Hodge filtration on the local system on $X$ -- 
it does not involve the weight filtration. 
So, using 
flatness of the twistor connection,  
it is sufficient to prove it for the twistor connection $\nabla^{(1)}_{\cal G}$ restricted to 
$u=1$ instead of $u=0$. 
The Hodge filtration for the connection $\nabla^{(1)}_{\cal G}$ is constant.
This implies the claim. Theorem \ref{2.27.08.3} is proved.

\paragraph{Proof of Theorem \ref{9.16.05.2aa}.}
By  Theorem \ref{2.27.08.3} there is 
 a good  unipotent variation of 
mixed $\R$-Hodge structures on $X-S$ given  by 
the universal enveloping algebras ${\rm A}(X-S, {a})$ of 
(\ref{20:33}), 
when $a$ varies.  
 By the Hain-Zucker theorem \cite{HZ}, 
it is described by a morphism of $\R$-MHS
$$
\theta: {\rm A}^{\cal H}(X-S, a) \otimes {\rm A}(X-S, {a}) \lra {\rm A}(X-S, {a}).  
$$
Observe that ${\rm A}(X-S, {a})$ is a free rank one 
${\rm A}^{\cal H}(X-S, a)$-module with a canonical generator given by  
the unit $e \in {\rm A}(X-S, {a})$. 

\bl \la{2.13.08.1}
The action $\theta$ coincides with the standard one. 
\el

\begin{proof} The action $\theta$ is a morphism of mixed $\R$-Hodge structures. 
We understand them as bigraded objects equipped 
with the action of the Hodge Galois Lie algebra. In particular $\theta$ 
preserves the bigrading and hence the weight. 

The action $\theta$ is determined 
by the images of the generators. The latter are operators of weight $-1$ 
or $-2$. The weight $-1$ generators are determined by the weight $-1$ component
of the connection $\nabla^0_{\cal G}$ on ${\cal L}_{\infty}$. Indeed, the weight $-1$ 
component of the parallel transform is determined by  the weight $-1$ component
of the connection, and the same is true for the weight $-1$ 
component of the 
exponential of the action of the generators. 
The weight $-1$ 
component of our connection is given by the canonical $1$-form 
$\nu$ (Section 2.2), and the same is true for the standard one. Thus the weight $-1$ 
generators act in the standard way. 

Corollary \ref{popopop} implies the claim for the weight $-2$ generators $X_s$.
\end{proof}

The Hodge Galois group acts by automorphisms, and hence acts trivially on the unit $e$. 
So there is an isomorphism of 
mixed $\R$-Hodge structures 
${\rm A}^{\cal H}(X-S, a) \otimes e \stackrel{\sim}{\lra} A(X-S, {a})$. 
Theorem \ref{9.16.05.2} is proved.

\vskip 3mm
The proof of Proposition \ref{3.8.08.1} is a bit technical. 
Here is a crucial step.

\paragraph{Dependence of ${\bf G}_{v_a}$ on $a$.} Below the point $s$,
 a puncture on the curve, is fixed. 
Let us calculate 
the part of ${\bf G}_{v_a}$ depending on the point $a$, and investigate what 
happens when $a$ approaches to 
a puncture $s$. 
Recall (\ref{7.3.00.1ds}) that
$$
G_{a}(x_1, x_2) =  
G_{\rm Ar}(x_1, x_2) - G_{\rm Ar}(a, x_2) - G_{\rm Ar}(x_1, a). 
$$

\begin{figure}[ht]
\centerline{\epsfbox{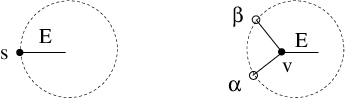}}
\caption{Changing the Green function assigned to the edge  
$E$ may alter the Hodge correlator.}
\label{hc12}
\end{figure}

\bl \la{22.1.08.1}
Assume that none of the vertices of an edge $E$ of a decorated tree $T$ 
is as on Fig \ref{hc12}, i.e. either $S$-decorated 
or decorated by a pair of $1$-forms. Then replacing
 the Green 
function $G_{v_a}$ assigned to $E$ by  $G_{\rm Ar}(x_1, x_2)$ 
we do not change the  
Hodge correlator integral assigned to $T$. 
\el

\begin{proof}  Let $v_1, v_2$ be the vertices of $E$. 
Each of them is a $3$-valent vertex. 
Denote by $x_1, x_2$ 
the points on the curves assigned to them. Let us assign 
the function $G_{\rm Ar}(x_1, a)$ to the edge $E$. 
Then the integral over $x_2$ is 
proportional to one of the 
following: 
$$
\int_X d^\C G_{a}(y,x_2) \wedge d^\C G_{a}(z, x_2), \quad \mbox{or} \quad 
\int_X d^\C G_{a}(y,x_2) \wedge \alpha.
$$ 
Each of them is zero: we integrate  $d^\C$-exact forms. Similarly for 
$G_{\rm Ar}( a, x_2)$. 
\end{proof}

Let $\{\alpha_i, \alpha_j^{\vee}\}$ be a Darboux
 basis in the symplectic space $\Omega^1_X \oplus \overline \Omega^1_X$, so that  
$\int_X \alpha_i\wedge \alpha^{\vee}_j = \delta_{ij}$. Let 
$\{h_i, h^{\vee}_j\}$ be the dual basis in $H_1(X, \C)$. 
Then the Hodge correlator is given by 
\be \la{5.21.08.1}
{\bf G}_{a} = 
\sum_{W}\frac{1}{|{\rm Aut}(W)|}{\rm Cor}_{{\cal H}, a}(\beta_1 
\otimes \ldots \otimes \beta_m)\cdot 
{\cal C}(h_{\beta_1} \otimes \ldots \otimes h_{\beta_m}), \quad \beta_i \in 
\{\alpha_i, \alpha_i^{\vee}\}
\ee
The sum is over a basis $W = {\cal C}(\beta_1 \otimes \ldots \otimes \beta_m)$ in the cyclic envelope of the tensor algebra 
of $H^1(X, \C)$, and $\{h_\beta\}$  is the dual basis to the basis $\{\beta\}$.

\paragraph{Contribution of an edge $E$ as on Fig \ref{hc12}.} Let us calculate the depending 
on the point $a$ contribution to the integral 
of an edge $E$ as on Fig \ref{hc12}. Take a tree $T$ with such an edge $E$. 
Cut out the ``small'' tree growing from 
$E$, getting a new tree. Let $W$ be its decoration 
inherited from $T$. For example, 
the decoration of the  left tree on 
Fig \ref{hc12} is ${\cal C}(s\otimes W)$. 
We assert that the depending on $a$ part of the contribution of the edge $E$ 
to the sum (\ref{5.21.08.1})   is 
$$
{\rm Cor}_{{\cal H}, {\rm Ar}}
(\{a\}\otimes  W)\cdot {\cal C}(X_{a}\otimes  W^{\vee}),
$$
where $W^{\vee}$ is the cyclic word in $H_1$ dual to the cyclic 
words $W$ in $H^1$. 
Indeed, take first 
the tree on the right of Fig \ref{hc12}, and assign to the edge $E$ the Green function $G_{\rm Ar}(a,y)$. Then its contribution to the  integral is 
$$
{\rm Cor}_{{\cal H}, {\rm Ar}}(\{a\} \otimes W)\int_X \alpha \otimes \beta. 
$$
Substituting basis elements 
$\alpha_i, \alpha^{\vee}_j$ for $\alpha$ as well as for $\beta$, and taking 
the sums over $i,j$, we get
$$
{\rm Cor}_{{\cal H}, {\rm Ar}}(\{a\}\otimes  W)
 \cdot {\cal C}(\sum_i[h_i, h_i^{\vee}] \otimes W^{\vee}). 
$$
Similarly, the tree on the left on Fig \ref{hc12} contributes 
$$
{\rm Cor}_{{\cal H}, {\rm Ar}}(\{a\}\otimes  W)
\cdot {\cal C}(\sum_{s\in S}X_s \otimes W^{\vee}). 
$$ 
Since $X_a = \sum_{s\in S}X_s + \sum_i [h_i. h_i^{\vee}]$, the total contribution is 
${\rm Cor}_{{\cal H}, {\rm Ar}}(\{a\}\otimes  W)
\cdot {\cal C}(X_a \otimes W^{\vee})$. 

\vskip 3mm
{\bf Remark}. A similar argument shows that 
the Green operator induced on ${\rm L}(X-\{a\}, v_{a})/(v_a)$
does not depend on the choice of a constant in 
the Green function $G_a(x,y)$. Indeed, 
we may assume that an  edge $E$  as on Fig \ref{hc12} 
contributes to the integral of the Green function but not its derivative. 
Then adding a constant to $G_a(x,y)$ amounts to changing 
 a factor in (\ref{5.21.08.1})  by
$$
\sum_{s\in S} X_s + \sum_{i=1}^g [h_i,  h^{\vee}_i] =0.
$$

{\it Conclusion}. Since we may have several edges like on Fig \ref{hc12}, 
the depending on $a$  part of ${\bf G}_{v_a}$  is  a 
sum of expressions 
$$
{\rm Cor}_{{\cal H}, {\rm Ar}}(\{a\}\otimes  W_1 \otimes \ldots \otimes \{a\}\otimes  W_m) 
\cdot {\cal C}(X_a\otimes  W^{\vee}_1 \otimes 
\ldots \otimes X_a \otimes  W^{\vee}_m). 
$$

\paragraph{Proof of Proposition \ref{3.8.08.1}.} 
Denote by ${\bf G}_{p,q}$ the $(p,q)$-component of ${\bf G}_{v_a}$.

Let us investigate ${\bf G}_{1,1}$. There are only two cyclic words
contributing to the depending on $a$ part of  
${\bf G}_{1,1}$: $W' = {\cal C}(\{a\}\otimes \{s\})$ 
and $W'' = {\cal C}(\{a\}\otimes \alpha \otimes \overline \beta)$. The 
first appears with the coefficient given by the Green function
 $G_{\rm Ar, v_a}(a,s)$, and has the singularity $\log t$. It delivers the 
singular term $\log t \cdot {\rm ad}(X_s)$. 
The second does not depend on $s$. So the claim about ${\bf G}_{1,1}$ follows. 

Let us show that ${\bf G}_{p,q} \in {\cal A}_{\log}(X)$ for $(p,q) = (1,q), (p, 1), (2,2)$ provided that 
$(p,q) \not = (1,1)$.  
After this the Lemma follows immediately by 
the induction using the differential equations on ${\bf G}_{p,q}$, see (\ref{1.08.08.10}). 
Notice that since $\overline {\bf G}_{1, q} = - {\bf G}_{q, 1}$, the claim about 
${\bf G}_{1,q}$ implies the one about ${\bf G}_{q,1}$. 

Let us prove that ${\bf G}_{1,q} \in {\cal A}_{\log}(X)$ for $q>1$.
The only cyclic words contributing 
to the depending on $a$ part of ${\bf G}_{1,q}$ are 
${\cal C}(\{a\} \otimes \{s\} \otimes 
\overline \alpha_1 \otimes \ldots \otimes \overline \alpha_{q-1})$, where 
$\overline \alpha_i$ are 
antiholomorphic $1$-forms on $\overline X$. The crucial case is $q=2$. 
Let us show that 
${\rm Cor}_{\cal H}(\{a\} \otimes \{s\} \otimes \overline \alpha) \in {\cal A}_{\log}(X)$. 
We have to investigate the behavior of the following integral at $a \to s$:
$$
\int_{\overline X} G(a,x)\partial G(s,x) \wedge \overline \alpha. 
$$ 
Choose a small neighborhood $U_s$ of $s$. We may assume that it contains $a$. 
Then the integral over its complement is clearly a smooth function in $a$. 
So we have to show that the integral over $U_s$, as a function of $a \in U_s$, 
has a tame logarithmic singularity at $a=s$. The problem is local. 
So it is sufficient to prove it for a single differential form 
$\overline \alpha$ in $U_s$ 
which does not vanish at $s$.   Since $G(0,x)$ has singularity of type $\log|x|$ near $x=0$,  
and since $\overline \alpha \sim \overline \partial \log |x-1|$ near $x=0$, 
we may assume $U_s$ is a small neighborhood of $0$ in $\C$, and 
consider an integral 
$$
\int_{\C {\Bbb P}^1} \log |x| \partial \log |x-a| \wedge \overline \partial \log |x-1|.
$$
This integral equals to a non-zero multiple of 
of the Bloch-Wigner dilogarithm function ${\cal L}_2(a)$, see Section 10.1.2, 
which is independent of this Section. Since  
$$
d{\cal L}_2(a) = \log|1-a| d\arg a - \log|a| d\arg (1-a), 
$$
it 
 has a tame logarithmic singularity at $a=0$. 
This implies the claim. The claim for $q>2$ follows easily from the $q=2$ case. 

Finally, ${\bf G}_{2,2}$ is given by either 
${\rm Cor}_{\cal H}(\{a\} \otimes \{s\} \otimes 
\overline \alpha \otimes \beta)$, which reduces to the previous case, 
or by ${\rm Cor}_{\cal H}(\{a\} \otimes \{s\} \otimes \otimes \{x\})$, 
which is completely similar to the ${\bf G}_{1,2}$ case, and done by the reduction to the dilogarithm. 
Proposition \ref{3.8.08.1} is proved.

\section{Motivic correlators on curves} \label{hc8sec}

In Section \ref{hc7sec} we described  
the Lie algebra of special derivations 
in the Betti realization. 
In Section \ref{hc8sec} we present it in the motivic set-up. Using this we define  
certain elements of the motivic Lie coalgebra, 
{\it motivic correlators on curves}. 
Section \ref{hc9sec} implies 
that the Hodge correlators on curves are their periods.

\vskip 3mm
Let  $X$ be an irreducible 
 regular projective curve over a field $F$, and $S$ a finite 
subset  of $X(F)$. We explore  
the motivic  fundamental Lie algebra 
${\rm L}(X-S, v_0)$ of the curve $X-S$ 
with a tangential base point $v_0$ over $F$.

We start by describing 
several different frameworks: Hodge, $l$-adic, motivic, etc. 
In each of them we supposed to have a Lie algebra 
in the corresponding category. It is equipped with the 
weight filtration $W_{\bullet}$, and    
${\rm gr}^W{\rm L}(X-S, v_0)$ is expressed via 
${\rm gr}^WH^1(X-S)$. 

The $l$-adic realization ${\rm L}(X-S, v_0)$ is an $l$-adic 
pro-Lie algebra ${\rm L}^{(l)}(X-S, v_0)$ -- 
the Lie algebra of the 
pro-$l$ completion  $\pi^{(l)}_1(X-S, v_0)$ of the 
fundamental group. 
The Galois group of $F$ acts by automorphisms 
of $\pi^{(l)}_1(X-S, v_0)$, and hence by derivations of the 
Lie algebra ${\rm L}^{(l)}(X-S, v_0)$,  preserving some data, mostly  
related to the set $S$. We call such derivations {\it special}.

\subsection{The motivic framework} \la{hc8.2sec}

Let $F$ be a field. Below  we  
work in one of the following five categories  ${\cal C}$:

\vskip 3mm

\begin{enumerate}

\item  {\it Motivic}. The hypothetical abelian category of mixed 
motives over a field $F$. 

\item  {\it Hodge}. $F=\C$, and  ${\cal C}$ is the category of mixed 
$\Q$- or $\R$-Hodge structures. 

\item  {\it Mixed $l$-adic}. 
$F$ is a field such that $\mu_{l^{\infty}} \not \in F$, and 
 ${\cal C}$ 
is 
 the mixed category of 
$l$-adic ${\rm Gal}(\overline F/F)$--modules with a filtration 
$W_{\bullet}$ indexed by integers,  
such that ${\rm gr}^W_n$ is a pure of weight $n$.  

\item 
 {\it Motivic Tate}. $F$ is a number field, 
${\cal C}$ is the abelian category of mixed 
Tate motives over $F$, equipped with the Hodge and $l$-adic 
realization functors, c.f. \cite{DG}. 
 
\item  {\it Variations of mixed Hodge structures}. 
$F=\C$, and  ${\cal C}$ is the category of mixed 
$\Q$- or $\R$-Hodge structures  over a base $B$. 

\end{enumerate}

The setup 1) is conjectural. The other four are well defined.  
A   category ${\cal C}$ from the list  is an $K$--category, 
where $K = \Q$ in 1), 4), 5); $K=\Q$ or $\R$ in 2);  and $K=\Q_l$ in 3). 

Each of the categories has an invertible object, the {\it Tate object}, 
which is denoted, abusing notation, 
$\Q(1)$ in all cases, although ${\rm Hom}_{\cal C}(\Q(1), \Q(1)) = K$. 
 So, for instance, in the setup iii) $\Q(1)$ 
denotes the Tate module $\Q_l(1)$. We set $\Q(n):= \Q(1)^{\otimes n}$.

Each object in ${\cal C}$ 
 carries a canonical weight filtration 
$W_{\bullet}$,  morphisms in ${\cal C}$ are strictly 
compatible with this filtration. The weight of $\Q(1)$ is $-2$. 

An object $M$ of ${\cal C}$ is  {\it pure} if ${\rm gr}^W_nM$ 
is zero for all $n$ but possibly one. 
Tensor powers and direct summands of pure objects are pure. 
Let $M$ be a pure object in ${\cal C}$. The  pure subcategory  
of ${\cal C}$ generated by $M$ is the smallest subcategory containing $M$, 
and closed under operations of taking tensor products, 
duals, finite sums and  direct summands. The mixed subcategory of ${\cal C}$ generated by $M$ 
is the smallest subcategory of ${\cal C}$ closed under extensions 
and containing the pure subcategory generated by $M$.

\vskip 3mm 
Let us briefly recall some features of the Tannakian formalism in 
the setup of mixed categories which is used below, see \cite{G3} 
for the proofs. See also \cite{G7} for the  mixed Tate case. 
 
There is a canonical fiber functor 
to the pure tensor category ${\cal P}^{\cal C}$ of pure objects in ${\cal C}$: 
$$
\omega: {\cal C} 
\lra {\cal P}^{\cal C}, \quad X \lms \oplus_n 
{\rm gr}^W_{n}X.
$$ 
The ind-object ${\rm End}^\otimes( \omega)$ of  endomorphisms of this functor respecting the tensor product 
 has a natural structure of a  Hopf algebra in 
the tensor category ${\cal P}^{\cal C}$. 
Let ${\cal A}_{\bullet}(\cal C)$ be its graded dual. It is  
a commutative graded Hopf algebra. 
The functor $\omega$ provides a 
canonical equivalence between the category ${\cal C}$ and the 
category of graded ${\cal A}_{\bullet}(\cal C)$--comodules in the category 
${\cal P}^{\cal C}$. 
There is 
the a Lie coalgebra 
$$
{\cal L}_{\bullet}({\cal C}) = \frac{{\cal A}_{\bullet}({\cal C})}{
{\cal A}_{>0}({\cal C})^2}
$$
with the cobracket induced by the coproduct in 
${\cal A}_{\bullet}({\cal C})$. Let ${\rm L}_{\bullet}({\cal C})$ 
be the dual Lie algebra.

\begin{definition} \label{11.28.04.1}
Let $X$ be a regular projective curve over $F$. Then 
${\cal P}_X$ (respectively  ${\cal C}_X$) is the 
subcategory of ${\cal P}^{\cal C}$  (respectively  ${\cal C}$)   generated by 
$H^1(X)$.
\end{definition}
Restricting the fiber functor $\omega$ to ${\cal C}_X$ we get an equivalence 
between the category  ${\cal C}_X$ and the category of graded 
${\cal A}_{\bullet}({\cal C}_X)$-comodules in ${\cal P}_X$. 
There is 
the corresponding Lie coalgebra 
$$
{\cal L}_{\bullet}({\cal C}_X) = \frac{{\cal A}_{\bullet}({\cal C}_X)}{
{\cal A}_{>0}({\cal C}_X)^2}\subset {\cal L}_{\bullet}({\cal C}) 
$$
 and the dual Lie algebra ${\rm L}_{\bullet}({\cal C}_X)$, a quotient of 
${\rm L}_{\bullet}({\cal C})$. 

The  Hopf algebra 
${\cal A}_{\bullet}(\cal C)$ 
can be defined using the framed objects in ${\cal C}$, see 
\cite{G3}.

The Lie coalgebra ${\cal L}_{\bullet}({\cal C})$ is 
the direct sum 
of its isotypical components ${\cal L}_{M}({\cal C})$ 
over the set of all  isomorphism classes of simple objects $M$ in ${\cal C}$:
\begin{equation} \label{5.22.05.1}
{\cal L}_{\bullet}({\cal C}) = \oplus_M{\cal L}_{M}({\cal C})\bigotimes M^{\vee}. 
\end{equation}

\subsection{${\cal C}$-motivic fundamental algebras of curves} \label{hc8.2sec}
Below we recall main features of the 
fundamental Lie/Hopf algebras of an irreducible  smooth  open curve 
in the motivic framework and in the realizations. 
The classical point of view is presented in Section 1.2 and
Section \ref{hc7sec}. The $l$-adic one 
is recalled in the end. 

\vskip 3mm
Let $S= \{s_0, ..., s_n\}\subset X(F)$, and $v_0$ is 
a tangent vector at 
$s_0 \in S$ defined over  $F$. 
One should have a Hopf algebra ${\rm A}^{\cal C}(X-S, v_0)$ in the 
category ${\cal C}$, the {\it ${\cal C}$-motivic 
fundamental Hopf algebra of $X-S$ with the tangential base point 
$v_0$}. 
It  was defined in the setups ii) and iii) in \cite{D1}, and 
in the setup iv) in \cite{DG}. For the setup v) see \cite{HZ}.
Set 
\be \la{9.29.1.13}
{\Bbb A}^{\cal C}_{X,S^*}:= {\rm gr}_{\bullet}^W{\rm A}^{\cal C}(X-S, v_0).
\ee
We 
use notation ${\rm A}$ and ${\Bbb A}$ for 
${\rm A}^{\cal C}(X-S, v_0)$ and ${\Bbb A}^{\cal C}_{X,S^*}$. 
The tangential base point $v_0$ gives rise to 
a canonical morphism, the ``canonical 
loop around $s_0$'': 
$$
X_{s_0}: \Q(1) \lra  {\rm A}. 
$$
The punctures $s \in S^*$ provide  
conjugacy classes  in ${\rm A}$, 
corresponding to ``loops around $s$''. 
By passing to the associated graded 
we arrive at canonical morphisms
$$
X_{s}: \Q(1) \lra  {\Bbb A},  \quad s\in S^*.
$$ 
A derivation ${\cal D}$ of the algebra 
${\Bbb A}$ is called {\it special} if there exist objects  
$B_s$ in ${\Bbb A}$
  such that 
$$
{\cal D}(X_s) = [B_s, X_s], \quad s \in S^*;  \qquad 
\mbox{and} \quad {\cal D}(X_{s_0}) =0.
$$
The formula on the left means that ${\cal D}(X_s(\Q(1)))$ is isomorphic to 
the image of 
the commutator map applied to the object 
$B_s \otimes X_s(\Q(1))$. 
Special derivations form a Lie algebra, denoted by ${\rm Der}^S{\Bbb A}$. 

Set  
\be \la{notatioonA}
{\Bbb S}^*= \Q[S^*] \otimes H_2(X), \quad {\Bbb H}= H_1(X); \qquad 
{\Bbb V}:=  {\Bbb S}^*\oplus {\Bbb H} \stackrel{\sim}{=} {\rm gr}^W H_1(X-S).
\ee
Notice that here 
$H_2(X)=\Q(1)$ tells the weight grading. In the linear algebra considerations 
before the additional data given by 
the weight grading 
was not taking into account.  

Let ${\rm T}_{\Bbb V}$ be the tensor algebra of ${\Bbb V}$. 

\begin{lemma} \label{12.1.04.1}
The Hopf algebra ${\Bbb A}$  is isomorphic to the 
tensor 
algebra ${\rm T}_{{\Bbb V}}$.  
\end{lemma}

\begin{proof}
The Hopf algebra (\ref{10/25/04/1}) is identified with the 
Betti realization of the motivic fundamental Hopf algebra of $X-S$. 
This implies Lemma \ref{12.1.04.1}. 
\end{proof}

\paragraph{An example: $l$-adic fundamental algebras.} 
Let $\pi_1^{(l)}:= \pi^{(l)}_1(X-S, v_0)$ be the pro-$l$ completion of the fundamental group of $X-S$ with the 
tangential base point at $v_0$. It gives rise to a pronilpotent Lie algebra ${\rm L}^{(l)}(X-S, v_0)$ over $\Q_l$ 
as follows. Let  $\pi_1^{(l)}(k)$ be the lower central series of $\pi_1^{(l)}$. Then 
$\pi_1^{(l)}/\pi_1^{(l)}(k)$ is an $l$-adic Lie group, and 
$$
{\rm L}^{(l)}(X-S, v_0):= \lim_{{\longleftarrow}}{\rm Lie}\Bigl(\pi_1^{(l)}/\pi_1^{(l)}(k)\Bigr).
$$ 
It is the $l$-adic realization of the fundamental Lie algebra 
${\rm L}(X-S, v_0)$. 

If $X$ is defined over $\C$, 
the comparison theorem between the Betti and $l$-adic realizations 
gives rise to a Lie algebra isomorphism
$
{\rm L}^{(l)}(X-S, v_0) = \pi_1^{\rm nil}(X(\C)-S, v_0)\otimes \Q_l.
$ 

If $X$ is defined over a field $F$, the Galois group 
${\rm Gal}(\overline F/F)$ acts by automorphisms of 
${\rm L}^{(l)}(X-S, v_0)$, giving rise to a canonical homomorphism
$$
{\rm Gal}(\overline F/F) \lra {\rm Aut}{\rm L}^{(l)}(X-S, v_0).
$$
The tangent vector $v_0$ provides a morphism of  
Galois modules, the 
``canonical loop around $s_0$'': $$
X_{s_0}: \Q_l(1) \lra {\rm L}^{(l)}(X-S, v_0).
$$
Thus every puncture $s\in S^*$  gives rise to a 
conjugacy class in ${\rm L}^{(l)}(X-S, v_0)$ 
preserved by the Galois group. 
An automorphism of ${\rm L}^{(l)}(X-S, v_0)$ is {\it special} if 
it preserves $X_{s_0}(\Q_l(1))$ and the conjugacy classes around 
all other punctures. Denote by ${\rm Aut}^S{\rm L}^{(l)}(X-S, v_0)$ the 
group of all special automorphisms of ${\rm L}^{(l)}(X-S, v_0)$. We get a map
$$
{\rm Gal}(\overline F/F(\zeta_{l^{\infty}})) \lra {\rm Aut}^S{\rm L}^{(l)}(X-S, v_0).
$$

\subsection{The Lie algebra of special derivations in the motivic set-up}
\label{hc8.3sec} 
We use a shorthand 
$$
{\cal C}_{{\Bbb V}}:= \frac{{\rm T}_{{\Bbb V}}}{
[{\rm T}_{{\Bbb V}}, 
{\rm T}_{{\Bbb V}}]}, \qquad {\cal H}_X = {\cal H}:= H^2(X).
$$
We are going to define an action of ${\cal C}_{{\Bbb V}}\otimes {\cal H}$  
by special derivations 
of the algebra ${\rm T}_{{\Bbb V}}$, and  introduce 
a  Lie algebra structure on 
${\cal C}_{{\Bbb V}}\otimes {\cal H} $, making this action into an action of 
a Lie algebra.

Since ${\rm T}_{{\Bbb V}}$ is the free associative algebra in a semi-simple  
category generated by 
${\Bbb V}$, there is  a non-commutative differential map 
$$
{\Bbb D}: {\cal C}_{{\Bbb V}} \lra {\rm T}_{{\Bbb V}}\otimes {\Bbb V}, 
\qquad  {\Bbb D}{\cal C}(a_0 \otimes a_1 \otimes ... \otimes a_m):= 
{\rm Cyc}_{m+1}(a_0 \otimes ... \otimes a_{m-1})\otimes a_m. 
$$
Here ${\rm Cyc}_{m+1}$  means the cyclic sum. So on the the right there is  a sum of $m+1$ 
summands obtained by the cyclic shifts of $(a_0 \otimes ... \otimes a_{m-1})\otimes a_m$. 

Let us define a  canonical map 
\be \la{19:29}
\langle \ast \cap {\cal H}\cap \ast \rangle: 
{\Bbb V} \otimes {\cal H} \otimes {\Bbb V} \lra \Q \oplus {\Bbb S}^*. 
\ee
Its only non-zero components are the following:
\be \la{19:26}
\langle \ast \cap {\cal H}\cap \ast \rangle_{{\Bbb H}}: 
{\Bbb H} \otimes {\cal H} \otimes {\Bbb H} \lra \Q, \qquad 
\langle \ast \cap {\cal H}\cap \ast \rangle_{{\Bbb S}^*}: 
{\Bbb S}^* \otimes {\cal H} \otimes {\Bbb S}^* \lra {\Bbb S}^*. 
\ee
The first one is the intersection pairing on $H_1(X)$. The second 
is given by $\langle X_s \cap {\cal H}\cap X_t \rangle_{{\Bbb S}^*} = \delta_{st}X_s$. 
 \vskip 3mm

\paragraph{An action of 
${\cal C}_{{\Bbb V}}\otimes {\cal H}$ by special derivations of  
${\rm T}_{{\Bbb V}}$.}  
We are going to define a map
\be \la{19:55}
\kappa: {\cal C}_{{\Bbb V}}\otimes {\cal H}\lra {\rm Der}^S({\rm T}_{{\Bbb V}}). 
\ee
Since the algebra ${\rm T}_{{\Bbb V}}$ is free, to define 
its derivation we just define it on the generators, i.e.
 produce an element in
${\rm T}_{{\Bbb V}}\otimes {\Bbb V}^{\vee} $. 
We get it as a composition 
$$
{\cal C}_{{\Bbb V}}\otimes {\cal H} \stackrel{{\Bbb D}\otimes {\rm Id}}{\lra} 
{\rm T}_{{\Bbb V}} \otimes ({\Bbb H}\oplus {\Bbb S^*})\otimes {\cal H} 
\stackrel{(\ref{19:26})^*}{\lra} 
{\rm T}_{{\Bbb V}} \otimes ({\Bbb H}^{\vee} \oplus  {\Bbb S^*} 
\otimes {\Bbb S^*}^{\vee}) \stackrel{\sim}{=}
{\rm T}_{{\Bbb V}}\otimes {\Bbb H}^{\vee}  \bigoplus   
{\Bbb S^*}\otimes {\rm T}_{{\Bbb V}} \otimes {\Bbb S^*}^{\vee}   
$$
followed by the commutator map 
${\Bbb S^*}\otimes {\rm T}_{{\Bbb V}} \lra {\rm T}_{{\Bbb V}}$. 
The map ${\rm T}_{{\Bbb V}} \otimes {\Bbb S^*}\otimes {\cal H} \lra 
{\rm T}_{{\Bbb V}} \otimes  {\Bbb S^*} 
\otimes {\Bbb S^*}^{\vee}$ in the second arrow 
is obtained using the map ${\Bbb S}^* \otimes {\cal H} \otimes {\Bbb S}^* \lra {\Bbb S}^*$ from (\ref{19:26}).

\paragraph{A Lie bracket on ${\cal C}_{{\Bbb V}}\otimes {\cal H}$.} 
There is  a canonical map 
$
\theta: {\cal C}_{{\Bbb V}}\otimes {\cal H}\otimes  {\cal C}_{{\Bbb V}} \lra 
{\rm T}_{{\Bbb V}} 
$ 
 given as 
$$
{\cal C}_{{\Bbb V}}\otimes {\cal H}\otimes  {\cal C}_{{\Bbb V}} 
\stackrel{{\Bbb D} \otimes {\rm Id} \otimes {\Bbb D}}{\lra} 
{\rm T}_{{\Bbb V}}\otimes {\Bbb V} \otimes {\cal H}
\otimes {\rm T}_{{\Bbb V}}\otimes {\Bbb V}{\lra}
{\rm T}_{{\Bbb V}}\otimes {\rm T}_{{\Bbb V}}\otimes {\Bbb V} \otimes {\cal H}
\otimes {\Bbb V}
\stackrel{(\ref{19:29})}{\lra}
{\rm T}_{{\Bbb V}}\otimes {\rm T}_{{\Bbb V}}\otimes (\Q \oplus {\Bbb S}^*)
\lra 
{\rm T}_{{\Bbb V}}.
$$
The last map is the product for the $\Q$-component, and  
the commutator map $[,]: {\rm T}_{{\Bbb V}} \otimes {\rm T}_{{\Bbb V}} \to 
{\rm T}_{{\Bbb V}}$ followed by the product with ${\Bbb S^*}$
 for the second component. 
We define a Lie bracket on ${\cal C}_{{\Bbb V}}\otimes {\cal H}$ as follows: 
\begin{equation} \label{6.20.00.11c}
\{F\otimes {\cal H}, G\otimes {\cal H}\}:= {\cal C}\circ 
\theta\Bigl(F\otimes {\cal H}\otimes G 
- G\otimes {\cal H}\otimes F\Bigr)
\otimes 
 {\cal H} 
\end{equation}
 
\paragraph{Explicit formulas.}
 Let ${\Bbb H} = \oplus_i p_i$ 
be a decomposition into simple objects. 
Then, given an $F \in {\cal C}_{\Bbb V}$, we have 
$$
{\Bbb D}F = \sum_i \frac{\partial F}{\partial p_i} \otimes p_i + 
\sum_{s\in S^*} \frac{\partial F}{\partial X_s} \otimes X_s.
$$
Now the action of the derivation $\kappa_{F\otimes {\cal H}}$ obtained by applying the map $\kappa$ to 
${F\otimes {\cal H}}$ is given by 
$$
\kappa_{F\otimes {\cal H}}: p \lms 
\frac{\partial F}{\partial q}\langle q \cap {\cal H} \cap p \rangle, 
\quad X_{s_0} \lms 0, \quad  
X_s \lms \left[X_s, \frac{\partial F}{\partial X_s}\right], \quad s\in S^*.
$$
$$
\theta: F\otimes {\cal H}\otimes G \lms 
\sum_{i,j}\frac{\partial F}{\partial p_i} \otimes 
\frac{\partial G}{\partial p_j}\otimes \langle p_i \cap {\cal H}\cap p_j\rangle 
+ \sum_{s \in S^*}
 \left[\frac{\partial F}{\partial X_s}, 
\frac{\partial G}{\partial X_s}\right]\otimes\langle 
X_s\cap {\cal H} \cap X_s\rangle .
$$
$$
\{F\otimes {\cal H}, G\otimes {\cal H}\} = 2\sum_{s\in S^*}
 {\cal C}\Bigl([\frac{\partial F}{\partial X_s},
\frac{\partial G}{\partial {X_s}}] \otimes X_s \Bigr) \otimes {\cal H}~+  
$$
\begin{equation} \label{6.20.00.11cer}
\sum_{p, q}\langle p \cap {\cal H}\cap q\rangle \otimes 
 {\cal C}\Bigl(\frac{\partial F}{\partial p} \otimes
\frac{\partial G}{\partial {q}} ~~ - ~~\frac{\partial G}{\partial p} 
\otimes
\frac{\partial F}{\partial {q}}\Bigr) \otimes {\cal H}. 
\end{equation}

Just like in Section \ref{hc7sec}, we define 
${\cal C}'_{{\Bbb V}}$ as the quotient of ${\cal C}_{{\Bbb V}}$ by 
$\oplus_{s\in S^*}\Q[X_s]$. 

\begin{theorem} \label{3/11/07/10} 
The map $\kappa$ is a morphism of Lie algebras. 
It induces an isomorphism 
\begin{equation} \label{3/11/07/11} 
\kappa': {\cal C}'_{{\Bbb V}}\otimes {\cal H}
 \stackrel{\sim}{\lra} {\rm Der}^{S}({\rm T}_{{\Bbb V}}).
\end{equation} 
The map $\kappa$ provides an isomorphism of Lie algebras
\begin{equation} \label{3/11/07/11xs} 
\kappa: {\cal C}{\cal L}ie_{{\Bbb V}}\otimes {\cal H}
 \stackrel{\sim}{\lra} {\rm Der}^{S}({\rm L}_{{\Bbb V}}).
\end{equation}
\end{theorem}

\begin{figure}[ht]
\centerline{\epsfbox{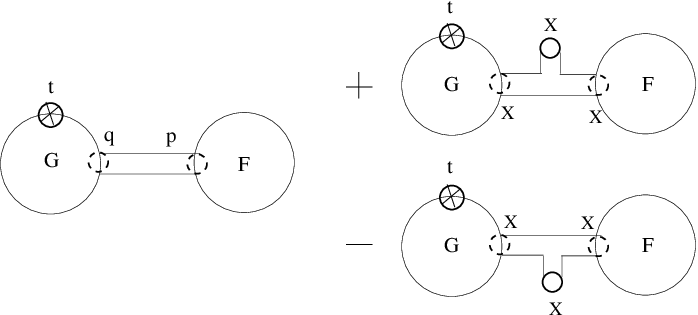}}
\caption{Calculation of $\kappa_{F\otimes {\cal H}} \circ \kappa_{G\otimes {\cal H}}(h)$.}
\label{fey100}
\end{figure}

\begin{proof}
Let us show that the map (\ref{3/11/07/11}) respects the brackets. 
Let us calculate how  the commutator 
$[\kappa_{F\otimes{\cal H}}, \kappa_{G\otimes{\cal H}}]$ acts on 
$h \subset {\Bbb H}$. To calculate the ${\Bbb H}$-part of the action, 
we have to insert 
${\partial F}/{\partial p}$ instead of each factor $q$ in
 $\partial G /\partial {t}$ -- this operation is denoted by
 $\stackrel{q}{\to}$. We get:
$$
\kappa_{F\otimes{\cal H}} \circ \kappa_{G\otimes{\cal H}}: 
h \lms \sum_{t}
\kappa_{F\otimes{\cal H}}(\frac{\partial G}{\partial t}) \otimes 
\langle  t\cap {\cal H}\cap h\rangle= 
$$
$$
\sum_{ t}
 \Bigl(\sum_{ p, q}
\frac{\partial F}{\partial p} \stackrel{q}{\to}
\frac{\partial G}{\partial {t}}\otimes\langle p \cap {\cal H}\cap q\rangle ~+~ 
\sum_{ s\in S^*}\left [\frac{\partial F}{\partial X_s}, 
\frac{\partial G}{\partial X_s}\right ]X_s\Bigr)
\otimes \langle t \cap {\cal H} \cap h\rangle.
$$
The (left) sums are over bases $\{t\}, \{p\}, \{q\}$ in ${\Bbb H}$. 
The result is visualized on 
Fig \ref{fey100}. Big circles show the cyclic 
words $F$ and $G$; little circles show the variables 
added, and the punctured 
little circles show the variables removed  in the process of calculation.  
The added / removed  factors 
shown nearby.  
 The bridge between big circles 
shows how we make  the product. The special role of 
 $t$ is emphasized by a cross. 
It is clear from this that $[\kappa_{F\otimes{\cal H}}, \kappa_{G\otimes{\cal H}}]$ 
acts on $h$ just like the commutator
$\kappa_{\{F\otimes{\cal H}, G\otimes{\cal H} \}}$. 
The action on $X_s$ is treated similarly -- the only difference is that we have to 
replace the operation of removing of $t$ by the commutator with $X_s$. 

The rest follows from this and Proposition \ref{sswas} by going 
to the Betti realization, or repeating the argument in the motivic set-up. 
\end{proof}

\subsection{Motivic correlators on curves} \label{hc8.4sec} 

We define three versions of the motivic correlators on curves, 
which differ in our treatment of the base points: 
base point motivic correlators, symmetric motivic correlators, and 
averaged  base point motivic correlators.

\paragraph{1. Base point motivic correlators on curves.} Similar to (\ref{9.29.1.13}), set 
\be \la{9.29.1.13.1}
{\Bbb L}_{X, S^*}:= {\rm gr}_{\bullet}^W{\rm L}^{\cal C}(X-S, v_0).
\ee 
Since the fundamental Lie algebra ${\rm L}^{\cal C}(X- S, v_0)$ of the curve 
is a pro-object in the category ${\cal C}_X$, the motivic Galois Lie algebra ${\rm L}_{\bullet}({\cal C}_X)$ 
acts by the special derivations on its associate graded for the weight filtration (\ref{9.29.1.13.1}), and so we get a map 
\begin{equation} \label{12.3.04.a2}
{\rm L}_{\bullet}({\cal C}_X) \lra {\rm Der}^S{\Bbb L}_{X, S^*}. 
\end{equation} 
Identifying Lie algebras ${\rm Der}^S{\Bbb L}_{X, S^*}$ and 
${\cal C}{\cal L}ie_{{\Bbb V}}\otimes H^2(X)$ via isomorphism 
(\ref{3/11/07/11xs}),  and dualizing   
map (\ref{12.3.04.a2}), we arrive at a Lie coalgebra
map 
\begin{equation} \label{12.1.04.3}
\Psi_{v_0}: {\cal C}{\cal L}ie^{\vee}_{X, S^*}\otimes H_2(X) \lra 
{\cal L}_{\bullet}({\cal C}_X).
\end{equation} 
Recall notations (\ref{notatioonA}). Abusing notation, write $X^{\vee}_s$ for $X_s(\Q(1))^{\vee}$. 
There is an isomorphism
\begin{equation} \label{12.3.04.2}
\Phi: \bigoplus_{m=1}^{\infty}\bigoplus_{a_i \in S^*}
({\Bbb A}_{{\Bbb H}^{\vee}} \otimes X_{a_0}^{\vee}) \otimes ({\Bbb A}_{{\Bbb H}^{\vee}} \otimes X_{a_1}^{\vee}) 
\otimes ... \otimes ({\Bbb A}_{{\Bbb H}^{\vee}} \otimes X_{a_m}^{\vee}) 
 \otimes {\Bbb A}_{{\Bbb H}^{\vee}}   \stackrel{\sim}{\longrightarrow}  
{\Bbb A}_{{\Bbb V}^{\vee}}.
\end{equation} 
Recall that ${\cal C}{\cal L}ie^{\vee}_{X, S^*}$ is the quotient of the cyclic envelope 
of ${\rm A}_{{\Bbb V}^{\vee}}$ by the shuffle relations: 
\be \la{mvlca1aq}
{\cal C}{{\cal L}ie}^{\vee}_{X, S^*}:= 
\frac{ {\cal C}({\rm A}_{{\Bbb V}^{\vee}})}{\mbox{Shuffle relations}}.
\ee

We denote by $\Phi_{\rm cyc}$ the map $\Phi$ 
followed by the projection to the cyclic envelope of 
${\Bbb A}_{{\Bbb V}^{\vee}}$. Projecting then to the quotient (\ref{mvlca1aq}) 
by the shuffle relations, and applying the map $\Psi_{v_0}$, see (\ref{12.1.04.3}), 
we arrive at an element of the 
motivic Lie coalgebra ${\cal L}_{\bullet}({\cal C}_X)$. This is  our first definition:
\begin{definition} \label{12.1.04.7} Given points $a_i \in S^*$  and 
simple summands 
$M_0, ..., M_m$ of ${\Bbb A}_{{\Bbb H}^{\vee}}$, we set 
 $w = 2m + \sum {\rm wt}(M_i)$, and define 
the {\rm base point motivic correlator}
\begin{equation} \label{12.1.04.15a}
{\rm Cor}_{X;v_0} \Bigl(\{a_0\} \otimes M_0 
\otimes \{a_1\} \otimes M_1 \otimes \ldots 
\otimes \{a_m\} \otimes M_m \Bigr)(1) \subset {{\cal L}_{w}({\cal C}_X)}
\end{equation}
as the image of $
{\cal C}\Bigl(X_{a_0}^{\vee} \otimes 
M_0 \otimes ... \otimes 
X_{a_m}^{\vee} \otimes  M_m\Bigr)(1)$ under the map $\Psi_{v_0}\circ\Phi_{\rm cyc}$.\footnote{The 
twist by $\Q(1)$ comes from the factor $H_2(X)$ in (\ref{12.1.04.3}).} 
\end{definition} 
\bl
The object 
(\ref{12.1.04.15a}) lies in the 
$M_0^{\vee} \otimes ...  \otimes M^{\vee}_m(m)$-isotypical component 
of ${\cal L}_{\bullet}({\cal C}_X)$.
\el

\begin{proof} Recall decomposition (\ref{5.22.05.1}). 
The object on the left in (\ref{12.1.04.15a}) is 
$M_0 \otimes ...  \otimes M_m(-m)$. 
The isotypical component of ${\cal L}_{\bullet}({\cal C}_X)$ is labeled by its dual. 
\end{proof}

Using description of the Hopf algebra ${\cal A}_{\bullet}({\cal C}_X)$ 
via minimal framed 
objects in ${\cal C}_X$ (\cite{G3}), 
and a similar description of the Lie coalgebra 
${\cal L}_{\bullet}({\cal C}_X)$, we conclude that the object 
(\ref{12.1.04.15a}) 
provides a minimal 
 framed object in ${\cal C}_X$, well-defined  modulo products 
of non-trivial objects in ${\cal C}_X$.

\paragraph{2. Symmetric motivic correlators on curves.} \la{9.4.2ref}
Let us employ the symmetric description of the Lie algebra of 
special derivations (Section \ref{hc7.3sec}) in the motivic set-up. 
There is a Lie algebra $\overline {\cal C}_{X, S}$ and its Lie subalgebra 
$\overline {{\cal C}{\cal L}ie}_{X, S}$ in the ${\cal C}$-motivic category,  
whose Betti realizations are the Lie algebras  
${\cal C}{\cal L}ie(\overline {\rm A}_{X,S})$ and 
${\cal C}{\cal L}ie(\overline {\rm L}_{X,S})$. Denote by 
$\overline {\cal C}^{\vee}_{X, S}$ 
and 
$\overline {{\cal C}{\cal L}ie}^{\vee}_{X, S}$ the dual Lie coalgebras. 

There is a Lie coalgebra
map 
\begin{equation} \label{12.1.04.3a}
\Psi^{\rm sym}: \overline {{\cal C}{\cal L}ie}^{\vee}_{X, S}\otimes H^2(X) \lra 
{\cal L}_{\bullet}({\cal C}_X).
\end{equation}

\bd \la{1.26.07.4}
A {\rm symmetric motivic correlator} is the image of 
$W \in \overline {{\cal C}{\cal L}ie}^{\vee}_{X, S}(1)$
 under map (\ref{12.1.04.3a}): 
$$
{\rm Cor}^{\rm sym}_X(W):= \Psi^{\rm sym}(W) \in {\cal L}({\cal C}_X).
$$ 
\ed

The Lie coalgebras  in question are related as follows:
 $$
\overline {{\cal C}{\cal L}ie}^{\vee}_{X, S} 
\longleftarrow \overline {\cal C}^{\vee}_{X, S} \hra {\cal C}^{\vee}_{X, S}. 
 $$
The left map is an epimorphism. 
In general 
there seems to be no simple description of the subspace 
$\overline {\cal C}^{\vee}_{X, S}$ in 
${\cal C}^{\vee}_{X, S}$. However when 
$X$ is either a rational or an elliptic curve, 
such a description exists. 

\bl \la{1.26.08.7} Let $X$ be the projective line ${\Bbb P}^1$. Then an element 
\be \la{1.26.08.9}
\sum_{a_0, ..., a_m\in S} c_{a_0, \ldots, a_m}{\cal C}
\Bigl(X^{\vee}_{a_0} \otimes \ldots \otimes X^{\vee}_{a_m}\Bigr) 
\in {\cal C}^{\vee}_{{\Bbb P}^1, S}
\ee
lies in the subspace $\overline {\cal C}^{\vee}_{{\Bbb P}^1, S}$ if and only if 
$\sum_{s \in S}c_{a_0, \ldots, a_{k-1}, s, a_{k+1}, \ldots , a_m} =0$ for any $k=0, ..., m$. 
\el

\begin{proof}  Elements  (\ref{1.26.08.9}) for arbitrary 
coefficients $c_{a_0, \ldots, a_m}\in \Q$
give a basis in  ${\cal C}^{\vee}_{{\Bbb P}^1, S}$. 
The orthogonality to the two-sided ideal generated 
by $\sum_{s\in S} X_s$ is precisely the condition 
on the coefficients. 
\end{proof} 

\bl \la{1.26.08.11} Let $E$ be an elliptic  curve, ${\Bbb H}^{\vee} = H^1(E)$. Then 
the subspace $\overline {\cal C}^{\vee}_{E, S}$ consists of 
elements 
\be \la{1.26.08.13}
\sum_{a_i\in S} c_{a_0, n_0;  \ldots, a_m, n_m}{\cal C}
\Bigl(X^{\vee}_{a_0} \otimes S^{n_0}{\Bbb H}^{\vee} \otimes \ldots 
\otimes X^{\vee}_{a_m}\otimes S^{n_m}{\Bbb H}^{\vee}\Bigr) 
\in {\cal C}^{\vee}_{E, S}, \qquad 
\ee
such that for any $k=0, ..., m$ one has 
$
\sum_{a_k \in S}c_{a_0,n_0;  \ldots, a_m, n_m} =0.
$ 
\el

\begin{proof} We spell the argument in the Betti realization. 
So ${\Bbb H}^{\vee}$ has a symplectic basis $p,q$. 
A basis of ${\cal C}^{\vee}_{E, S}$ is given by 
cyclic words 
$
{\cal C}(a_0 \otimes F_0(p,q)  \otimes \ldots \otimes a_m \otimes 
F_m(p,q)), 
$ 
where $a_i \in S$ and $F_i(p,q)$ are non-commutative polynomials in $p,q$. 
The defining ideal is generated 
by the element $[p, q] + \sum_{a\in S}X_a$. So expressing 
$[p, q]$ as $-\sum_{a\in S}X_a$ we can replace the polynomials 
$F_i(p,q)$ by commutative polynomials $G_i(p,q)$. We may assume 
that $G_i$ are homogeneous 
of degree $n_i$. Since $G_i(p,q)$ are commutative polynomials, 
they are orthogonal to $[p,q]$. So the orthogonality 
to the element  $[p, q] + \sum_{a\in S}X_a$ boils down to the orthogonality 
to $\sum_{a\in S}X_a$. \end{proof} 

{\bf Remark}. The  symmetric motivic correlator map 
does not depend on the choice of the base point. 
It has the same image  in  the motivic 
Lie coalgebra 
as the one defined by using a base point. 
This apparent ``contradiction'' is resolved as follows. 
A choice of the base point provides a way to 
choose elements in the image of the motivic correlator map, without changing 
the latter. Different base points 
provide  different 
representations, built using 
different supplies of elements in the same space. For example, 
as we see in Section \ref{9.5ref}.1, the elements of 
the Jacobian $J_X\otimes \Q$ of $X$ provided by the symmetric motivic correlators 
are $s_i-s_j$. The ones provided by the motivic correlators corresponding 
to a tangential base point at $s_0$ are $s_i - s_0$. 
These are two different collection of elements of 
$J_X \otimes \Q$; however they span the same  subspace. 

Similarly the symmetric Hodge correlators are complex numbers defined 
by the curve $X-S$. The $\Q$-vector space spanned by them  
 does not depend on 
the choice of the base point, or a Green function. 
However choosing different base points we get 
different real MHS's --  the latter depend not only on the 
$\Q$-vector subspace spanned by the periods, but also on the way 
we view its elements as the matrix coefficients.

\paragraph{3. Averaged base point motivic correlators.} 
Let us  choose for every point $s\in S$ a tangent vector $v_s$ at $s$. 
There are natural choices for the elliptic and modular curves. 

Recall the canonical morphism of Lie coalgebras 
assigned to a tangential base point $v_s$ at $s$: 
\be \la{der}
{\rm Cor}_{X, S-\{s\}; v_s}: {\cal C}{\rm T}({\Bbb H}^{\vee}_{X,S-\{s\}}) 
\otimes H_2(X)
\stackrel{}\lra  {\cal L}_{\rm Mot}. 
\ee
The collection of maps (\ref{der}), for different choices of  $s\in S$, 
is organized
into a Lie coalgebra map  
\be \la{der2}
\oplus_{s \in S}{\cal C}{\rm T}({\Bbb H}^{\vee}_{X,S-\{s\}})\otimes H_2(X) 
\lra {\cal L}_{\rm Mot}. 
\ee
Here on the left stands a direct sum of the Lie coalgebras. 
It is handy extend $S-\{s\}$ to $S$, and 
enlarge the Lie coalgebra on the left to a Lie coalgebra 
\be \la{der3}
{\cal C}{\rm T}({\Bbb H}^{\vee}_{X,S})\otimes H_2(X) \otimes \Q[S], 
\ee
where the added elements have the coproduct zero. So 
the Lie coalgebra on the left of (\ref{der2}) is 
the  quotient of (\ref{der3}). 
The composition (\ref{der3}) $\to$ (\ref{der2}) 
provides a Lie coalgebra map 
\be \la{der1}
{\rm Cor}^{\rm av}_{X-S}: {\cal C}{\rm T}({\Bbb H}^{\vee}_{X,S})\otimes H_2(X) \otimes \Q[S]
\stackrel{}\lra  {\cal L}_{\rm Mot},  
\ee
so that map (\ref{der}) is the component of map (\ref{der1}) 
assigned  to the factor $\{s\} \in \Q[S]$.

The averaged base point motivic correlator map arises when we take 
$\frac{1}{|S|}\sum_{s\in S}\{s\} \in \Q[S]$. 
It allows to work in towers of curves, see Section \ref{sec11.2}

\subsection{Examples of motivic correlators} 
\begin{figure}[ht]
\centerline{\epsfbox{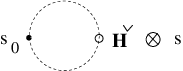}}
\caption{The motivic correlator for ${\cal C}(\{s_0\} \otimes {\Bbb H}^{\vee})$ is the point $s-s_0$ of the Jacobian of $X$.}
\label{hc2}
\end{figure}

\la{9.5ref}
\paragraph{1. The Jacobian of $X$.} Let $J_X(F)$ be the Jacobian of the curve $X$. 
It is well known that in the motivic set-up (i) one should have 
$$
{\rm Ext}_{\rm Mot}^1(\Q(0), {\Bbb H}) = {\cal L}_1({\cal C}_X) = J_X(F)\otimes \Q.
$$
The motivic correlator 
${\rm Cor}_{X, s_0}(\{a\} \otimes {\Bbb H}^\vee)$, see Fig. \ref{hc2},  
is given by the point $\{a\}-\{s_0\}$ in the Jacobian $J_X(F)$. 
Notice the importance of the base point $s_0$. 
In the symmetric description the correlators 
are 
$$
{\rm Cor}^{\rm sym}_{X}{\cal C}((\{a\}-\{b\})  
\otimes {\Bbb H}^\vee) = a-b \in J_X\otimes \Q.
$$
They are defined only on degree zero divisors. 
In any case we get the same supply of elements of the motivic Lie coalgebra.  
In the $l$-adic set-up (iii) the target group is $J_X(F)\otimes \Q_l$. 

\paragraph{ 2. The rational curve case.} Let  $X={\Bbb  P}^1$, the set $S$ includes 
$\{0\}\cup \{\infty\}$, 
and $v=\partial /\partial t$  
is the canonical vector at $0$. Denote by $dt/t$ the cohomology class 
given by the form $dt/t$ in the De Rham realization. The motivic 
correlators in this case are 
\be \la{14:47sas}
{\Bbb L}^{\cal M}_{n_0, \ldots , n_m}(a_0; ... ; a_m):= 
{\rm Cor}_{{\Bbb P}^1; v}\Bigl(
\{a_0\} \otimes \{0\}^{\otimes n_0}\otimes \ldots \otimes  \{a_m\} 
\otimes \{0\}^{\otimes n_m} \Bigr)
\subset {{\cal L}_{w}({\cal C}_T)}.
\ee 
Here $a_i \in S^*$, $w= 2(n_0+ ... +n_m+m)$ is  
the weight, and $m$ is the depth (Section 9.5.9). 
We call them the {\it cyclic motivic multiple polylogarithms}. 
Their periods coincide with the 
corresponding multiple polylogarithms  
up to the lower depth terms.  In Section  \ref{hc10sec}.1.4 we 
derive an integral formula expressing their real periods 
as Hodge correlator type integrals of the classical polylogarithms. 

\vskip 3mm
{\bf Remark}. In \cite{G4}, given an abelian group $G$, 
we defined dihedral Lie coalgebras, which are  
generated, as vector spaces, by symbols 
${\rm L}_{n_0, \ldots , n_m}(g_0; ... ; g_m)$, $g_i \in G$. 
When $G = \C^*$,  
we related them to multiple polylogarithms, although to achieve  this 
we had to 
work modulo the depth filtration.  
It is now clear that 
those generators match the motivic elements (\ref{14:47sas}), 
and so their periods are the corresponding Hodge correlators 
rather then the single valued versions of the multiple polylogarithms.  
See Section \ref{hc10sec} for more information.
\vskip 3mm

An especially interesting case is when $S=\mu_N \cup \{0\}\cup \{\infty\}$.
In this case we call the elements (\ref{14:47sas})  
the depth $m$ motivic multiple Dirichlet $L$-values. 
For instance when $S=\{0\}\cup \{1\}\cup \{\infty\}$, we get the depth 
$m$ motivic multiple $\zeta$-values.

\paragraph{3. The elliptic curves case.} Let  $X=E$ be an 
elliptic curve 
over $F$. 
If $E$ is not CM, the objects ${\rm S}^k{\Bbb H}$ are simple. 
The  {\it symmetric motivic multiple elliptic polylogarithms}  
are defined as 
\begin{equation} \label{12.1.04.15sd}
{\rm Cor}^{\rm sym}_{E; n_0, \ldots , n_m}(D_0; ... ; D_m):= {\rm Cor}_{E}\Bigl(
D_0 \otimes {\rm S}^{n_0}{\Bbb H}^{\vee} \otimes \ldots \otimes  D_m 
\otimes {\rm S}^{n_m}{\Bbb H}^{\vee}  \Bigr)(1)
\subset {{\cal L}_{w}({\cal C}_E)}.
\end{equation} 
Here $D_i$ are degree zero divisors supported on 
$S$ and  $w= n_0+ ... +n_m+2m$. 
When $E$ is a complex elliptic curve, their periods  
 are given by the generalized 
Eisenstein-Kronecker series (Section \ref{hc10sec}.2).

\paragraph{4. Motivic correlators of torsion points on elliptic curves.} 
Let $X=E$ be an elliptic curve, and $S=E[N]$ is 
the set of its $N$-torsion points. There is an almost canonical choice $v_\Delta$ of 
the tangent vector at $0$: its dual is given by $12$-th root of the 
section of $\Omega_{{\cal E}/{\cal M}}^{12}$, 
where ${\cal E}$ is the universal 
elliptic curve over the modular curve  ${\cal M}$, 
 given by the modular form $\Delta= q \prod(1-q^n)^{24}$.  
We employ the corresponding $E$-invariant vector field $v_\Delta$ on $E$, 
and assign to every missing point  on $E-E[N]$ 
the tangent vector $\frac{1}{N}v_\Delta$. We 
 arrive at the averaged base point motivic correlator map. It 
is obviously $E[N]$-invariant, and thus descends to a map 
of  $E[N]$-coinvariants: 
$$
{\rm Cor}^{\rm av}_{E-E[N]}: \Bigl({\cal C}{\rm T}({\Bbb H}^{\vee}_{E, E[N]})\otimes 
{\rm Meas}({E[N]})(1)\Bigr)_{{E[N]}} 
\lra {\cal L}_{\rm Mot}.
$$
There is a unique ${E[N]}$-invariant volume $1$ measure $\mu_{E}^0$ 
on ${E[N]}$. 
So we get a canonical map 
$$
{\rm Cor}^0_{{E-E[N]}}: {\cal C}{\rm T}({\Bbb H}^{\vee}_{E, E[N]})(1)_{{E[N]}}\otimes \mu_{E}^0
\lra {\cal L}_{\rm Mot}.
$$
These are the motivic correlators with 
the ``averaged tangent vector'' at the $N$-torsion points.

\paragraph{5. Motivic double elliptic logarithms at torsion points.} 
They were defined and studied in \cite{G10}. 
Here is an alternative definition, coming from 
the action of the Galois group on the motivic fundamental group 
of an elliptic curve minus the torsion points. 
We use a shorthand 
$$
\theta_E(a_0: \ldots : a_n):= 
{\rm Cor}^0_{{E}}(\{a_0\} \otimes 
\ldots \otimes \{a_n\})(1), \quad a_i \in E_{\rm tors}.
$$ 
It agrees with the notation in \cite{G10}. 
It is translation invariant. 

Here is an $l$-adic / Hodge version of formula 
(\ref{7.3.00.1ds}) expressing 
the Green function $G_a(x,y)$ via the Arakelov Green function. 
We spell it in the $l$-adic set-up. 

\bl \la{theta1} For any $x \in E(F)$, and any 
$E$-invariant vector field $v$ on $E$  one has 
$$
\theta_x(a_0:a_1)= \theta_E(a_0:a_1) - \theta_E(x:a_1)  
- \theta_E(a_0:x) + C, \qquad C \in {\rm Ext}^1_{\cal C}(\Q_l(0), \Q_l(2)).
$$
\el

\begin{proof} 
Indeed,  
$$
\delta \theta_x(a_0:a_1)= {\rm Cor}_{{E}, x}(\{a_0\} \otimes {\Bbb H}^{\vee})(1) 
\wedge {\rm Cor}_{{E}, x}(\{a_1\} \otimes {\Bbb H}^{\vee})(1) = 
(a_0- x) \wedge (a_1-x) \in \Lambda^2 E(F). 
$$
The coproduct of the right hand side modulo $N$-torsion is the same. 
Therefore the difference lies in 
${\rm Ext}^1_{\cal C}(\Q_l(0), \Q_l(2))$. 
The latter  $l$-adic  ${\rm Ext}$-group is rigid. Thus,  
as a variation in $x, a_0, a_1$, it is constant variation.  
\end{proof} 

 The constant $C$ for the canonical vector $v_\Delta$ is zero. 
Indeed, $v_\Delta$, viewed as a section over the modular curve, 
 is defined over $\Q$, and 
we can take $E, x,a_i$'s defined over $\Q$, and use the fact that 
$C$ is of geometric origin, and 
$K_3^{\rm ind}(\Q) =0$. 

\bl \la{theta} Assume that $a_i \in E[N]$. Then, in the $l$-adic or Hodge set-ups, 
one has
\be \la{theta2}
\delta \theta_E(a_0: a_1 : a_2) = {\rm Cycle}_{\{0,1,2\}} 
\Bigl(\theta_E(a_0: a_1) \wedge \theta_E(a_1 : a_2)\Bigr) \quad 
\mbox{modulo $N$-torsion}. 
\ee
\el

\begin{proof} We have $N\delta_{\rm Cas}\theta_E(a_0: a_1 : a_2)=0$. 
Indeed, 
$$
N\delta_{\rm Cas}\theta_E(a_0: a_1 : a_2) = N{\rm Cycle}_{\{0,1,2\}}
\Bigl({\rm Cor}^0_{{E}}(\{a_0\} \otimes {\Bbb H}^{\vee})(1) \wedge 
{\rm Cor}^0_{{E}}({\Bbb H}^{\vee} \otimes \{a_1\} \otimes \{a_2\})(1)\Bigr),  
$$
and $N{\rm Cor}^0_{{E}}(\{a_0\} \otimes {\Bbb H}^{\vee})=Na_0 =0$. 
Using Lemma \ref{theta1}, $\delta_{S}\theta_E(a_0: a_1 : a_2)$  
 equals 
$$
{\rm Cycle}_{\{0,1,2\}}
\frac{1}{N^2}\sum_{x \in E[N]}
\theta_x(a_0: a_1) \wedge \theta_x(a_0: a_2) = 
{\rm Cycle}_{\{0,1,2\}}
\Bigl(\theta_E(a_0: a_1) \wedge \theta_E(a_0: a_2)\Bigr). 
$$
\end{proof}

The elements $\theta_E(a_0: a_1 : a_2)$ from \cite{G10} 
have the same coproduct, and thus coincide with the ones above up to an element 
of ${\rm Ext}^1_{\cal C}(\Q_l(0), \Q_l(2))$. 
Each of them is skew-symmetric in $a_0, a_1, a_2$. Since the difference between them 
is constant on the modular curve, this implies that it is zero.

\paragraph{6. Motivic correlators of  torsion points on a CM elliptic curve.} 
Suppose that  $E$ is a CM curve. 
Then it has a complex multiplication by an order in the ring of integers in an imaginary quadratic field $K$. Extending the scalars 
from $L$ to $L\otimes K$, we have a decomposition 
 ${\Bbb H}^{\vee} = \psi\oplus \overline \psi$ into a sum 
of pure motives corresponding 
to the two Hecke Gr\"ossencharacters.  Set $\psi^n:= \psi^{\otimes n}$ 
and $ \overline \psi^n:= \overline \psi^{\otimes n}$. 
Then 
$
{\rm S}^n{\Bbb H}^{\vee} = \oplus_{n'+n''=n}\psi^{n'}\overline \psi^{n''}.
$ 
So we decompose the  motivic correlators accordingly.

Let $a_i$ be torsion points on $E$, 
and $(n'_i, n''_i)$ non-negative integers.   
The corresponding
 {\it motivic Hecke  Gr\"ossencharacters multiple $L$-values} 
$L_{E; n'_0,n''_0; \ldots ; n'_m,n''_m}(a_0: \ldots : a_m)$  
are motivic correlators 
$$
{\rm Cor}^0_E\Bigl( \{a_0\}\otimes \psi^{n'_0}\overline \psi^{n''_0}\otimes 
\{a_1\}\otimes \psi^{n'_1}\overline \psi^{n''_1}\otimes \ldots 
\otimes \{a_m\}\otimes \psi^{n'_m}\overline \psi^{n''_m}\Bigr)(1). 
$$

\paragraph{ 7. The Fermat curves case.} Let ${\Bbb F}_N$ be 
the projective Fermat curve $x^N+y^N=z^N$. 
Let $S$ be the intersection of ${\Bbb F}_N$ with the coordinate triangle. 
The group 
$\mu_N^2 = \mu_N^3/\mu_N$ acts on the curve ${\Bbb F}_N$. 
The motive $H^1({\Bbb F}_N)\otimes \Z[\zeta_N]$ is 
decomposed into a direct sum of rank one motives $\psi_{\chi}$ 
parametrized by the characters $\chi$ of $\mu_N^2$. 
These are Weil's  Jacobi sums Gr\"ossencharacter motives. 

We use the symmetric 
motivic correlators. 
Just like in  the CM 
elliptic curve case, by using the motives $\psi_{\chi}$
 and their powers, we get symmetric motivic correlators 
parametrised by cyclic tensor products 
of $\psi^n_{\chi}$ and the Tate motives assigned to the points of $S$. 

\paragraph{ 8. The modular curves case.} This is one of the most 
interesting cases, see Section \ref{hc11sec}.

\paragraph{9. The depth filtration.} 
Suppose that $s_0 \in S' \subset S \subset X$.  
Then the inclusion $X-S \subset X-S'$ 
provides a projection  of algebras 
$
 {\rm A}(X-S, v_0) \lra {\rm A}^{\cal C}(X-S', v_0).  
$ 
 A {\it depth filtration} on the Hopf algebra ${\rm A}(X-S, v_0)$ 
is a filtration by powers of  its kernel. 
One defines similarly a depth filtration on the Lie algebra 
${\rm L}(X-S, v_0)$. It is induced by the embedding 
${\rm L}(X-S, v_0)\hookrightarrow {\rm A}(X-S, v_0)$. 
Here are two examples. 

(i) Assume that $X\not = {\Bbb P}^1$. 
The embedding $X-S \hookrightarrow X$  provides a canonical projection 
\begin{equation} \label{6.16.00.1g}
{\rm A}(X-S, v_0) \lra {\rm A}(X, s_0). 
\end{equation} 
Let ${\rm I}_{X,S}$ be its kernel. 
A {depth filtration} $D$ on the Hopf algebra ${\rm A}(X-S, v_0)$ 
is the filtration  
$  
{\rm A}(X-S, v_0) \supset {\rm I}_{X,S}\supset 
{\rm I}_{X,S}^2\supset  \ldots
$ 
by powers of  the ideal ${\rm I}_{X,S}$, 
indexed by $m=0, 1, 2, ... $. 
One has 
$$
{\rm Gr}^D_{m}{\rm A}(X-S, v_0) = \oplus_{s_i \in S^*}{\rm A}_X\otimes X_{s_0}
\otimes {\rm A}_X \otimes ... \otimes X_{s_m} \otimes {\rm A}_X. 
$$

(ii) The depth filtration for  ${\Bbb G}_m-S$ is defined 
using the embedding ${\Bbb G}_m-S \hra {\Bbb G}_m$. 
\section{Examples of Hodge correlators} \la{hc10sec}

\subsection{Hodge correlators on $\C{\Bbb P}^1-S$ and polylogarithms}

\paragraph{1. Multiple Green functions.}  
Let $a_0, \ldots, a_{m}$ be points of a complex curve $X$. 
Let $\mu$ be a volume one measure on $X$. 
The Hodge correlator of a cyclic word 
${\cal C}(\{a_0\} \otimes \ldots \otimes \{a_m\})$
 is the {\it depth $m$
multiple Green function} (Section 9 of \cite{G1}):
\begin{equation} \label{7.3.00.1a}
\begin{split}
&G_\mu(a_0, \ldots, a_{m}):= 
\sum_{\mbox{T}} {\rm sgn}(E_1 \wedge \ldots \wedge E_{2m-1}) \int_{X^{m-1}}
\omega_{2m-2}(G_{E_1} \wedge \ldots \wedge G_{E_{2m-1}}). \\
\end{split}
\end{equation}
The sum is over all plane trivalent trees $T$ decorated by $a_0, \ldots, a_{m}$. The integral 
is  over the product of copies of $X$ 
parametrised by the internal vertices of $T$. 
The function $G_\mu(a_0, \ldots, a_{m})$ enjoys

\begin{itemize}
\item  The dihedral symmetry: 
$
G_\mu(a_0, \ldots, a_{m}) = G_\mu(a_1, \ldots, a_{m}, a_0) = (-1)^{m+1}G_\mu(a_{m}, \ldots, a_0).
$ 
\item  The shuffle relations: 
$
\sum_{\sigma \in \Sigma_{p,q}}G_\mu(a_0, a_{\sigma(1)}, \ldots, 
a_{\sigma(m)}) = 0.
$ 
\end{itemize}

The first property is clear from the definition. The second 
follows from Proposition \ref{7.3.06.4}.

\paragraph{2. The Bloch-Wigner function.} 
 The simplest multiple Green function appears when 
$X = {\Bbb P}^1$, $\mu = \delta_a$ 
 and $m=3$. It is described by a single Feynman diagram,  
shown on the left of Fig \ref{feyn10}. 

Recall the dilogarithm function ${\rm Li}_2(z)$, its single valued version, the Bloch-Wigner function  
${\cal L}_2(z)$, and the 
the cross-ratio $r(z_1, \ldots, z_4)$ of four points $z_1, \ldots, z_4$ 
on ${\Bbb P}^1$ normalized by $r(\infty, 0, 1, z) = z$, see Section 1.3. 
\begin{lemma} \label{3.7.05.1}
One has 
$$
G_{a}(a_0, a_1, a_2)= -\frac{1}{(2\pi i)^2}{\cal L}_2(r(a, a_0, a_1, a_2)).
$$
\end{lemma}

\begin{proof}  We may assume without loss of generality that 
$a = \infty$. Then, up to an additive constant,  
$
G_{\infty}(x,y) = \log|x-y|.
$ 
Therefore 
$$
G_{\infty}(a_0, a_1, a_2)= \frac{1}{(2\pi i)^3}\int_{\C{\Bbb P}^1}\omega_2\Bigl(\log|x-a_0|\wedge 
\log|x-a_1|\wedge \log|x-a_2|\Bigr).
$$
It remains to use the differential equations for 
the Hodge correlator from Section \ref{hc6sec} and ${\cal L}_2(z)$.
\end{proof}

\begin{figure}[ht]
\centerline{\epsfbox{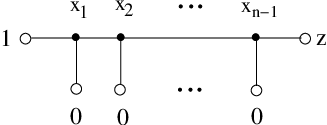}}
\caption{The Feynman diagram for the classical $n$-logarithm.}
\label{feyn15}
\end{figure}

\paragraph{3. The classical polylogarithms.} Let $X = \C {\Bbb P}^1$ and 
$s_0=\infty$. Consider the following cyclic word 
of length $n+1$: 
\begin{equation} \label{3.21.05.1x}
W_n= {\cal C}\Bigl(\{1\} \otimes \{z\} \otimes 
\{0\} \otimes \ldots \otimes \{0\} \Bigr).
\end{equation}
A $W_n$-decorated Feynman diagram  with  
an internal vertex incident to  two $\{0\}$-decorated external vertices   
evidently contributes the zero integral. 
There is a unique $W_n$-decorated Feynman diagram 
with no internal vertices incident to two $\{0\}$'s, see  
Fig \ref{feyn15}. 
Denote by 
${\bf L}_n(z)$ the Hodge correlator for $W_n$, which uses the 
clockwise orientation  of the tree. We can package them into the generating series
$$
{\bf L}_n(z|u):= \sum_{n=0}^\infty{\bf L}_n(z)\frac{u^n}{n!}.
$$

Recall that the classical $n$-logarithm is a multivalued function 
on $\C{\Bbb P}^1-\{0, 1, \infty\}$,  
defined by an iterated integral: 
$$
Li_n(z) := 
\int_{0 \leq x_1 \leq \ldots \leq x_n\leq z}\frac{dx_1}{1-x_1}\wedge 
\frac{dx_2}{x_2} \wedge \ldots \wedge \frac{dx_n}{x_n}.
$$ 
Here we integrate along a  simplex consisting 
of points $(x_1, \ldots, x_n)$ on a path $\gamma$ connecting $0$ and $z$, following 
each other on the path. It 
has a  single-valued
version  (\cite{Z1}, \cite{BD}): 
\begin{eqnarray*} 
{\cal L}_{n}(z) &:=& \begin{array}{ll} 
{\rm Re} & (n:\ {\rm odd}) \\ 
{\rm Im} & (n: \ {\rm even}) \end{array} 
\left( \sum^{n-1}_{k=0} \beta_k
\log^{k}\vert z\vert \cdot Li_{n-k}(z)\right)\; , \quad n\geq 
2. \\ 
\end{eqnarray*}           
Here   $\frac{2x}{e^{2x} -1}  =
\sum_{k=0}^{\infty}\beta_k x^k $, so $\beta_k = \frac{2^kB_k}{k!}$ 
where the $B_k$ are the Bernoulli numbers. For example ${\cal L}_2(z)$
is the Bloch - Wigner function.

Consider a modification ${\Bbb L}^*_{n}(z)$ of the 
 function ${\cal L}_n(z)$ studied  by A.Levin in \cite{L}:
$$
{\Bbb L}^*_{n}(z):=  4^{-(n-1)}
\sum_{\mbox{$k$ even;  $0 \leq k \leq n-2$}}
{2n-k-3\choose n-1}\frac{2^{k+1}}{(k+1)!}
 {\cal L}_{n-k}(z)\log^k|z|.
$$
It is also handy to have a slight modification 
of this function:
$$
{\Bbb L}_{n}(z) = 4^{(n-1)}{2n-2\choose n-1}^{-1}{\Bbb L}^*_{n}(z).
$$
So ${\Bbb L}_{n}(z) = {\cal L}_{n}(z)$  
if and only if $n \leq 3$. 
\begin {lemma}  \label{clpoly} One has 
$$
-{\bf L}_n(z) = (2\pi i)^{-n}{\Bbb L}_{n}(z).
$$ 
\end {lemma}

\begin{proof}
A direct integration carried out in Proposition 4.4.1 of \cite{L} tells
$$
-(2\pi i)^{-n}
{\Bbb L}^*_{n}(z) = 
$$
\be \la{1.08.04.1}
(2\pi i)^{-(2n-1)}\int_{(\C{\Bbb P}^1)^{n-1}}\log|1-x_1| \bigwedge_{i=1}^{n-2}\Bigl(d\log |x_i| \wedge 
d\log |x_i - x_{i+1}|\Bigr) \wedge d\log |x_{n-1}|\wedge d\log |x_{n-1} - z|.
\ee
By  \cite[Proposition 6.2]{GZ},\footnote{notice that we  use the form $-\omega_m$ in \cite{GZ} instead of $\omega_m$} 
 given functions 
 $\varphi_i$ on a complex $(n-1)$-dimensional manifold $M$, 
$$
\int_M\omega_{2n-2}(\varphi_0 \wedge ... \wedge \varphi_{2n-2}) = 
(-4)^{n-1}{2n-2\choose n-1}^{-1}
\int_M\varphi_0  d\varphi_1 \wedge ... \wedge d\varphi_{2n-2}. 
$$  
Using  $\varphi_i= \log|f_i|$ we see that 
integral (\ref{1.08.04.1}) coincides up to a factor with 
the Hodge correlator integral 
for the counterclockwise orientation of the tree on Fig \ref{feyn15}. 
The factor $(-1)^{n-1}$ tells 
the difference between the counterclockwise  and 
clockwise orientations of the tree. 
\end{proof}

Alternatively, one can 
 check using the differential equations for the function ${\cal L}_n(z)$ that 
the function ${\Bbb L}^*_{n}(z)$ satisfies the differential equation
$$
d {\Bbb L}^*_{n}(z) = \frac{2n-3}{2n-2} {\Bbb L}^*_{n-1}(z)d^\C\log |z| - 
\frac{1}{2n-2}\log |z| d^\C{\Bbb L}^*_{n-1}(z).
$$
The function ${\bf  L}_{n}(z)$ is a Hodge correlator, and thus
 satisfies the same 
differential equation (Section \ref{hc6sec}), and is zero if $z=0$. 
This implies Proposition \ref{clpoly} by induction, since ${\Bbb  L}^*_{n}(0)=0$. 

So the Hodge correlator ${\rm Cor}_{\cal H}$ delivers a single-valued 
version of the classical polylogarithm, while multiplying it 
by $4^{-(n-1)}{2n-2\choose n-1}$, i.e. using the 
normalized Hodge correlator ${\rm Cor}^*_{\cal H}$ (Section \ref{hc6sec}),  
we make the differential equation nicer. 
\vskip 3mm

\paragraph{4. Cyclic multiple polylogarithms and Hodge correlators.} 
Let us set $s_0=\infty$, and 
\be \la{2.5.08.3}
{\Bbb L}_{k_0; \ldots ; k_m}(a_0:\ldots  :a_m):= 
{\rm Cor}_{\cal H, \infty}\Bigl( \{a_0\} \otimes 
\{0\}^{\otimes k_0} \otimes \ldots 
\otimes \{a_m\} \otimes \{0\}^{\otimes k_m}\Bigr), \quad a_i \in \C^*. 
\ee
We package these numbers into a generating series
$$
{\Bbb L}(a_0:\ldots  :a_m|u_0: \ldots : u_m):= 
\sum^{\infty}_{k_i = 0}{\Bbb L}_{k_1; \ldots ; k_m}(a_0:\ldots :a_m)
\frac{u_0^{k_0}}{k_0!}
\ldots 
\frac{u_m^{k_m}}{k_m!}.
$$

\begin{figure}[ht]
\centerline{\epsfbox{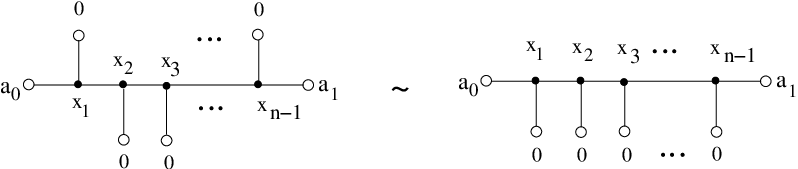}}
\caption{Feynman diagrams for the  depth one Hodge correlators.}
\label{hc16}
\end{figure}

The 
Hodge correlator integral assigned to the right diagram on 
Fig \ref{hc16} equals ${\bf L}_n(a_1/a_0)$.

\bl \la{deps1c}
One has 
\be \la{2.08.03.10}
{\Bbb L}(a_0: a_1|u_0: u_1) =  {\bf L}(\frac{a_1}{a_0}| u_1-u_0).
\ee
\el

\begin{proof}
The left hand side is a sum 
of the integrals attached to the diagrams 
on the left of Fig \ref{hc16}. 
The integral attached to the left Feynman diagram  
on Fig \ref{hc16}  coincides
up to a sign, by its very definition,  with the integral 
for the right Feynman diagram on Fig \ref{hc16}. 
The sign is 
$$
(-1)^{\mbox{the number of legs on the left diagram looking up}}.
$$ 
Equivalently, the exponent is the number of legs which we need to flip 
in order to transform one diagram to the other. 
Indeed, there is  a natural bijection 
between the edges of these two graphs. It does not respect 
their canonical orientations. The sign measures the difference: 
 flipping an edge we change the sign of the canonical orientation of the tree. 
 It follows that
$$
{\Bbb L}(a_0: a_1|u_0: u_1) \stackrel{\rm def}{=} 
\sum^{\infty}_{k_0, k_1 = 0}{\Bbb L}_{k_0;  k_1}(a_0 :a_1)
\frac{u_0^{k_0}}{k_0!} 
\frac{u_1^{k_1}}{k_1!} = 
$$
$$
\sum^{\infty}_{k_0, k_1 = 0}(-1)^{k_0}{\bf L}_{k_0+  k_1}(a_0 :a_1)\frac{u_0^{k_0}}{k_0!}
\frac{u_1^{k_1}}{k_1!} =
\sum^{\infty}_{k = 0}{\bf L}_{k}(a_0: a_1)
\frac{(u_1-u_0)^{k}}{k!} = {\bf L}(\frac{a_1}{a_0}| u_1-u_0).
$$
The second equality was explained above. 
\end{proof}

\begin{figure}[ht]
\centerline{\epsfbox{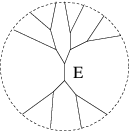}}
\caption{An edge $E$ of a tree provides two trees rooted at $E$.}
\label{mp1bbarc}
\end{figure}

\paragraph{Cyclic multiple polylogarithms.} We are going 
to define them, following Section 8 of \cite{G1}, as a yet another 
 generating series 
$$
{\bf L}(a_0: \ldots :a_{m}| t_0, \ldots, t_{m}), \qquad 
a_i\in \C^*, \quad t_i \in H_1(\C^*, \R), 
\quad 
t_0 + \ldots + t_{m}=0.
$$
Consider a plane trivalent tree $T$ 
decorated by 
$m+1$ pairs 
$(a_0, t_0), \ldots , (a_{m}, t_{m})$.  
We picture $a_i$'s outside, and $t_i$'s inside of the 
circle. 
\begin{figure}[ht]
\centerline{\epsfbox{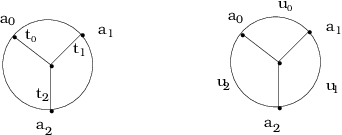}}
\caption{Decorations  for the cyclic multiple polylogarithms.}
\label{feyn}
\end{figure}
Each oriented edge $\stackrel{\to}{E}$ of the tree $T$ 
provides an element $t(\stackrel{\to}{E})$ in $H_1(\C^*, \R)$ 
defined as follows. The edge $E$ determines two trees rooted at $E$, see 
Fig \ref{mp1bbarc}. 
An orientation of the edge $E$ allows to choose one of these trees: 
the one obtained by going in the direction shown by 
the orientation. 
The union of the incoming (i.e. different from  $E$) legs 
of these rooted trees coincides with the set of all legs of the initial tree. 
Let $t(\stackrel{\to}{E})$ be the sum of all $t_i$'s 
corresponding to the incoming legs  of one of these trees. 
 The opposite orientation of the edge $E$ produces 
the element $-t(\stackrel{\to}{E})$. 
We define the generating series for the cyclic multiple polylogarithms  
just like  a Hodge correlator, with the 
generating series (\ref{2.08.03.10}) for the classical polylogarithms 
serving as Green functions. 
Namely, given a plane trivalent tree $T$ decorated as above, and an edge $E$ of  $T$, 
let us define a generating series assigned to the pair $(E, T)$:
\begin{equation} \la{9.15.13.1}
{\bf L}(x^{E}_1: x^{E}_2|t_1^{E}, t_2^{E}).
\ee
We define first the arguments $x^{E}_1, x^{E}_2, t_1^{E}, t_2^{E}$. 
Choose an orientation of the edge $E$. It provides 
an ordered pair $(x^{E}_1, x^{E}_2)$ of the vertices of $E$, as well as 
a pair of vectors  $(t^{E}_1, t^{E}_2):= 
(t(\stackrel{\to}{E}), -t(\stackrel{\to}{E}))$. Then 
${\bf L}(x^{E}_1: x^{E}_2|t_1^{E}, t_2^{E})$
 is given by the generating series  (\ref{2.08.03.10}) evaluated at these arguments: 
$$
{\bf L}(x^{E}_1: x^{E}_2|t_1^{E}, t_2^{E}):= {\bf L}(\frac{x^{E}_1}{x^{E}_2}|t_1^{E}).
$$
Thanks to  (\ref{2.08.03.10})
and the symmetry ${\bf L}(a|t) = {\bf L}(a^{-1}|-t)$, 
the function ${\bf L}(x^{E}_1: x^{E}_2|t_1^{E}, t_2^{E})$  
does not change when we change the orientation of the edge $E$:
$$
{\bf L}(x^{E}_1: x^{E}_2|t_1^{E}, t_2^{E}) = {\bf L}(x^{E}_2: x^{E}_1|t_2^{E}, t_1^{E}).
$$ 

Now we are going to apply the map $\omega_m$ to the wedge product 
of the generating series (\ref{9.15.13.1}) over the edges $E$ of the tree $T$. 
The map $\omega_m$ is usually applied, in the case of curves, to a wedge product  of functions. 
In our case (\ref{9.15.13.1}) is a generating series whose coefficients are functions, and  
we apply $\omega_m$  to the coefficients. Here is an example. Let $a(u):= \sum a_iu^i$ and $b(v):= \sum b_jv^j$ 
where $a_i, b_j$ are functions on a curve, and $u,v$ are formal variables. Then 
$$
\omega_2(a(u), b(v)) = \omega_2(\sum_{i}a_i u^i, \sum_{j}b_j v^j):= \sum_{i,j}\omega(a_i, b_j)u^iv^j. 
$$
Then we integrate over the product of $\C{\Bbb P}^1$'s 
labeled by the internal vertices 
of $T$. 
Finally, we take the sum over all plane $3$-valent trees 
$T$ with the given decoration. We arrive at the following definition. 

\bd \la{CMP}
\begin{equation} \label{7.1.00.2qw}
{\bf L}(a_0: \ldots: a_{m}|t_0, \ldots, t_{m}):= (2\pi i)^{-(2m-1)}\sum_{T}
 {\rm sgn}(E_1 \wedge \ldots \wedge E_{2m-1})\cdot
\end{equation}
$$
 \int_{(\C{\Bbb P}^1)^{m-1}}
\omega_{2m-2}\left({\bf L}(x^{E_1}_1: x^{E_1}_2|t_1^{E_1}, t_2^{E_1}) \wedge \ldots \wedge 
{\bf L}(x^{E_{2m-1}}_1: x^{E_{2m-1}}_2|t_1^{E_{2m-1}}, t_2^{E_{2m-1}}) \right). 
$$
\ed

\bt \la{1.08.03.2}
$
{\Bbb L}(a_0: \ldots :a_m|u_0: \ldots : u_m) = 
{\bf L}(a_0:\ldots  :a_m|u_0- u_m, u_1- u_0, \ldots , u_m- u_{m-1}).   
$ 
\et

\begin{proof} {\it The depth $1$ case}. This is Lemma \ref{deps1c}.  
To proceed further, we need temporarily 
another  generating series, whose combinatorics is illustrated on Fig \ref{hc17}:
\be \la{1.08.03.1}
\begin{split}&{\Bbb L}'(a_0: \ldots :a_m|s_0, t_0; \ldots ; s_m, t_m):= \\
&\sum^{\infty}_{n'_i, n''_i=0}{\Bbb L}_{n'_0, n''_0; \ldots ; n'_m, n''_m}(a_0:\ldots :a_m)
\frac{s_0^{n'_0}}{n'_0!}\frac{t_0^{n''_0}}{n''_0!}
\ldots
\frac{s_m^{n_m}}{n_m!}\frac{t_m^{n''_m}}{n''_m!}.\\
\end{split}
\ee
\begin{figure}[ht]
\centerline{\epsfbox{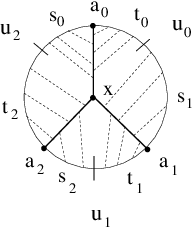}}
\caption{A Feynman diagram for the double polylogarithm.}
\label{hc17}
\end{figure}

{\it The depth $2$ case}. 
The Hodge correlator integral is given by 
a Feynman diagram as on Fig \ref{hc17}. One has 
$$
{\Bbb L}'(a_0: a_1 :a_2|u_2, u_0; u_0, u_1;  u_{1}, u_2) = 
{\Bbb L}(a_0:a_1 :a_2|u_0, u_1 , u_2). 
$$
Indeed, 
to get the right hand side we have to substitute 
$t_0 = s_1 = u_0$, 
$t_1 = s_2 = u_1$, $t_2 = s_0 = u_0$ into the sum (\ref{1.08.03.1}),
 getting the left hand side.  Furthermore, the same arguments 
about flipping the edges as in 
the proof of Lemma \ref{deps1c}  show that  
$$
{\Bbb L}'(a_0: a_1 :a_2|s_0, t_0; s_1, t_1;  s_{2}, t_2)=
{\bf L}(a_0:a_1 :a_2|t_0- s_0, t_1- s_1, t_2- s_2).   
$$
Here we integrate first over the internal vertices incident 
to the external $\{0\}$-decorated vertices. Thanks to Lemma \ref{deps1c} 
we get the ${\bf L}$-propagators assigned to the three ``thick'' edges 
of the tree. The remaining  we get precisely the definition 
of the depth two function ${\bf L}$. 
Combining, we get 
$$
{\Bbb L}(a_0:a_1 :a_2|u_0, u_1, u_2) = 
{\bf L}(a_0:a_1 :a_2|u_0- u_2, u_1- u_0, u_2- u_1).   
$$

{\it The general case}. It proceeds just like the depth $2$ case. 
The case $m=3$ is illustrated on Fig \ref{hc19}. For example, 
we assign $u_2-u_0 = (u_2-u_1) + (u_1-u_0)$ to the internal edge $E$ 
of the tree on Fig \ref{hc19},   
in agreement with the recepee for $t(\stackrel{\to}E)$. 
\begin{figure}[ht]
\centerline{\epsfbox{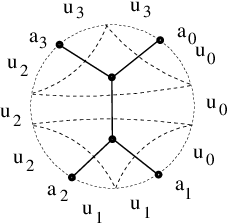}}
\caption{A Feynman diagram for the triple polylogarithm.}
\label{hc19}
\end{figure}
\end{proof}

\paragraph{5. Double shuffle relations.} 
We have defined 
the motivic cyclic multiple polylogarithms,  
denoted below  ${\rm L}_{n_0, \ldots , n_m}^{\cal M}(a_0: \ldots : a_m)$, 
in (\ref{14:47sas}). Recall that they are elements of the motivic Lie algebra of the category 
of mixed Tate motives. They are linear combinations of the traditional 
motivic multiple polylogarithms functions, but satisfy cleaner identities, as Lemma \ref{7.1.00.3} shows.

When the ground field is $\C$, and elements $a_i \in \C$ vary, 
we can apply the Hodge realisation functor, followed by the canonical real period map. 
Then we get well defined real valued functions depending on the complex variables $a_i$.
We prove that they are nothing else but the 
 {\it single valued cyclic multiple polylogarithm functions} defined above. They, of course, satisfy all 
the identities which their motivic counterparts have.

The traditional weight $m$ 
motivic multiple polylogarithms ${\rm L}_{n_1, \ldots , n_m}^{\cal M}(a_1, \ldots , a_m)$ live in the motivic Hopf algebra. 
Therefore, when the elements $a_i \in \C$ vary, 
applying the Hodge realisation functor to them, we 
get a variation of framed mixed $\Q$-Hodge structures, which gives rise to a multivalued analytic function 
on the parameter space. Precisely, to get the functions we have to choose a splitting of the weight filtration, 
which exists only locally, and then use the framings.  

 Since motivic cyclic multiple polylogarithms live in the motivic Lie coalgebra rather then 
in the motivic Hopf algebra, they gives rise to multivalued analytic functions well defined only modulo products 
of similar lower weight functions. 
We would like to stress that, in contrast with this, 
 the single valued cyclic multiple polylogarithms functions 
are well defined functions.

\begin{lemma} \label{7.1.00.3}
i) The cyclic multiple polylogarithms  enjoy 
the dihedral symmetry relations
\be \la{1.03.08.11}
{\Bbb L}(a_0: \ldots : a_{m}| t_0: \ldots : t_{m}) \quad = \quad 
{\Bbb L}(a_1: \ldots : a_{m}: a_0|t_1: \ldots : t_{m}: t_0) \quad = 
\ee
$$
 (-1)^{m+1}{\Bbb L}(a_{m}: \ldots : a_0|t_{m}: \ldots : t_0).
$$
and the shuffle relations 
\be \la{1.03.08.12}
\sum_{\sigma \in \Sigma_{p,q}}
{\Bbb L}(a_{0}: a_{\sigma(1)}\ldots: a_{\sigma(m)}| 
t_{0}: t_{\sigma(1)}: \ldots: t_{\sigma(m)}) = 0.
\ee

ii) The motivic 
cyclic multiple polylogs 
${\rm L}^{\cal M}(a_0: \ldots : a_{m}| t_0: \ldots : t_{m})$ satisfy the same relations. 
\end{lemma}

\begin{proof}  i) The first is clear from the definition. 
The second follows from Proposition \ref{7.3.06.4}, which shows 
the shuffle relations for Hodge correlators. 
ii) is similar. \end{proof}

{\bf Remark}. The shuffle relations (\ref{1.03.08.12}) do not contain any products since 
the cyclic multiple polylogarithms by the very definition live in the motivic Lie coalgebra, and the products 
are killed by passing to the motivic Lie coalgebra. 
\vskip 2mm

The definitions below were motivated by the results of \cite{G4}. Set
$$
{\Bbb L}(a_0: \ldots : a_m):= {\Bbb L}_{0, \ldots , 0}(a_0: \ldots : a_m),   
$$
$$
{\Bbb L}(a_0, \ldots , a_m):= 
{\Bbb L}(1: a_0, a_0a_1, 
\ldots , a_0a_1 ... a_{m-1}), \qquad a_0 ... a_m =1. 
$$
We also need their motivic versions 
${\rm L}^{\cal M}(a_0: \ldots : a_m)$ and ${\rm L}^{\cal M}(a_0, \ldots , a_m)$.  
\bp \la{2.05.08.1} i) For any $p+q=m, p,q \geq 1$ one has the second 
shuffle relations
\be \la{2.05.08.2}
\sum_{\sigma \in \Sigma_{p,q}}
{\Bbb L}(a_0, a_{\sigma(1)}, \ldots , a_{\sigma(m)}) = 0.
\ee

ii) The same is true for the ${\cal C}$-motivic multiple logarithms 
${\rm L}^{\cal M}(a_0, \ldots , a_m)$. 
\ep

\begin{proof} ii) 
The elements 
${\rm L}^{\cal M}(a_0: \ldots : a_m)$
satisfy the same coproduct formula 
as the generators of the dihedral Lie coalgebra, see (81) in \cite{G4}: 
$$
\delta {\rm L}^{\cal M}(a_0: \ldots : a_m) = 
{\rm Cycle}_{m+1}\Bigl(\sum_{k=1}^{m-1} {\rm L}^{\cal M}(a_0: \ldots : a_k)
\wedge {\rm L}^{\cal M}(a_k: \ldots : a_{m})\Bigl). 
$$
Indeed, this formula is equivalent to the one for the differential $\delta_S$ 
in Section \ref{hc5sec}.2, which was made motivic in Section \ref{hc8sec}.

Let us prove the shuffle relations by the induction on $m$. 
Denote by $s_2^{p,q}(a_0, \ldots , a_m)$ the left hand side in 
(\ref{2.05.08.2}). Theorem 4.3 in \cite{G4} 
tells that these elements span a coideal. This implies 
by the induction that 
$\delta s_2^{p,q}(a_0, \ldots , a_m)=0$. 
Therefore  
$s_2^{p,q}(a_0, \ldots , a_m) \in {\rm Ext}^1(\Q(0), \Q(m))$ 
in the category of variations 
over the base parametrising the $a_i$'s. Since $m\geq 2$, 
${\rm Ext}^1(\Q(0), \Q(m))$ is rigid, and thus is a constant. 
Setting $a_i=1$ in (\ref{1.03.08.12}) we get  
${\rm L}^{\cal M}(1: \ldots : 1)=0$. Thus
$s_2^{p,q}(1, \ldots , 1) = 0$, and hence the constant above is zero. 
The statement is proved. 

i) Follows from ii) by applying the period map, or by 
using the differential equations.
\end{proof}

So the cyclic multiple logarithms as well as their motivic 
avatars  
satisfy the double shuffle relations from Section 4 of \cite{G4} on the nose, 
without lower depth corrections or products. 
The absence of the products is a general feature of the 
Lie-type elements. It is also a general feature of the Lie-type periods, since the latter 
are obtained by applied the period map to the Hodge counterpart of the motivic Lie coalgebra. 

The absence of the lower depth corrections 
is a remarkable fact. 
It implies that there is a homomorphism from the diagonal part of 
the dihedral Lie coalgebra of $\C^*$ or $F^*$ ({\it loc. cit.})  
to the motivic Lie coalgebra which sends the standard generators of the former 
to the motivic cyclic multiple logarithms. 
Notice that the construction of a similar homomorphism in 
\cite{G4} required the associate graded for the depth filtration 
of the motivic Lie coalgebra.  
Specifying to the roots of unity, we conclude that the 
mysterious 
connection between the geometry of modular varieties and 
motivic cyclic multiple logarithms at roots of unity is valid without 
going to the associated graded for the depth filtration.

So the generators of the dihedral Lie coalgebras 
are related to the cyclic multiple polylogarithms, 
rather then the usual multiple 
polylogarithms, which differ by the lower depth terms. 

\vskip 3mm 
{\bf Remark}. The double shuffle relations for the depth 2 motivic 
cyclic multiple polylogarithms follow from (\ref{1.03.08.11}), since 
in this case they are equivalent 
to the dihedral symmetry relations. Their proof in general 
is more involved since the coproduct of the motivic 
cyclic multiple polylogarithms has the lower depth terms.

\vskip 3mm
\subsection{Hodge correlators on elliptic curves are 
multiple Eisenstein-Kronecker series} 

\paragraph{1. The classical Eisenstein-Kronecker series.} 
Let $E$ be a complex elliptic 
curve with a lattice of periods $\Gamma$, so that $E = \C/\Gamma$. 
The intersection form $\Lambda^2\Gamma \to 2 \pi i \Z$ 
leads to a  pairing $\chi: E \times \Gamma \to S^1$. 
So a point   $a \in E$ provides a character $\chi_a: \Gamma \lra S^1$. 
We denote by $dz$ the differential on $E$ provided by the coordinate $z$ on $\C$.
We normalize the lattice so that $\Gamma = \Z \oplus \Z\tau$, ${\rm Im}\tau >0$. 

Consider the generating series for the classical Eisenstein-Kronecker series
\be \la{KLEK}
K(a| t):= \frac{\tau - \overline \tau}{2\pi i}\sum_{\gamma \in \Gamma - \{0\}}
 \frac{\chi_{a}(\gamma)}{|\gamma-t|^2}.
\ee
It depends on a point $a$ of the elliptic curve 
and an element $t$ in  a formal neighborhood of zero in $H_1(E, \R)$. It is invariant
 under the involution $a \lms -a, t \lms -t$. 
Expanding it into the 
power series in $t$ and $\overline t$ we get 
the classical Eisenstein-Kronecker series as the coefficients: 
\be \la{CEKS}
K(a| t) = \sum_{p,q \geq 1}\Bigl(\frac{\tau - \overline \tau}{2\pi i}
\sum_{\gamma \in \Gamma - \{0\}}
 \frac{\chi_{a}(\gamma)}{\gamma^p\overline \gamma^q}\Bigr)t^{p-1}\overline t^{q-1}.
\ee

Set ${\rm Sym}_{p+q}F(x_1, ..., x_{p+q}):= 
\sum_{\sigma} F(x_{\sigma(1)}, ..., x_{\sigma(p+q)})$. 
Consider the following cyclic word:
\begin{equation} \label{3.21.05.1}
W_{p,q}:= {\cal C}\Bigl(\{0\} \otimes \{a\}\otimes 
\frac{d\overline z^{p}}{p!} \cdot  \frac{d{z}^{q}}{q!} \Bigr) := 
{\rm Sym}_{p+q}{\cal C}\Bigl(\{0\} \otimes \{a\}\otimes 
\frac{d\overline z^{\otimes p}}{p!} \otimes  \frac{d{z}^{\otimes q}}{q!} \Bigr).
\end{equation}
The corresponding Hodge correlator is given by 
the Feynman diagram  on 
Fig \ref{feyn16}. 
\begin{figure}[ht]
\centerline{\epsfbox{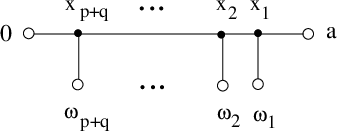}}
\caption{A Feynman diagram for the 
classical Eisenstein-Kronecker series; $\omega_i = dz$ or 
$d\overline {z}.$}
\label{feyn16}
\end{figure}
\bp \la{2.06.08.1} i) The Hodge correlator of $W_{p,q}$ is 
the classical Eisenstein-Kronecker series:
$$
{\rm Cor}_{\cal H}(W_{p,q}) = 
\frac{(-1)^p}{2\pi i}\Bigl(\frac{\tau - \overline \tau}{2\pi i}\Bigr)^{p+q+1}
\sum_{\gamma \not = 0}
 \frac{\chi_{a}(\gamma)}{\gamma^{p+1}\overline \gamma^{q+1}}.
$$
ii) The function $K(a|t)$ coincides with 
the generating function of the Hodge correlators of (\ref{3.21.05.1}). 
\ep

\begin{proof}
i) The normalized Green function $g(x,y)$ for the invariant volume form 
${\rm vol}_E$ on $E$ is a distribution   given by  $g(x,y):= g(x-y)$ where 
\be \la{disgz}
\chi_{z}(\gamma) := {\rm exp}
(\frac{2\pi i(z\overline \gamma 
- \overline z \gamma)}{\tau - \overline \tau }), 
\qquad g(z) := \frac{\tau - \overline \tau}{2\pi i}\sum_{\gamma \not = 0}
 \frac{\chi_{z}(\gamma) }{|\gamma|^2}, \qquad 
{(2\pi i)}^{-1}\overline \partial \partial g(z) = 
 \delta_0 - {\rm vol}_E.
\ee
Therefore the Hodge correlator  of (\ref{3.21.05.1}) 
is obtained by integrating the following form over $E^{p+q}$:
$$
\Omega_{p,q}:= \frac{(-1)^p}{(2\pi i)^{(p+q+1)}}
 g(a-x_{1}) 
\bigwedge_{s=1}^p
d \overline x_s\wedge \partial g(x_{s}-x_{s+1})  \wedge 
\bigwedge_{t=1}^qd x_{p+t}\wedge \overline \partial g(x_{p+t}-x_{p+t+1}). 
$$
Here $x_{p+q+1}:=0$. We used the clockwise orientation of the tree 
to make the product. Notice that the component of the form $\omega_{p+q}$ 
responsible for this form comes with the  coefficient  $(-1)^pp!q!/(p+q)!$.  
The integral does not depend on the 
normalization of the Green function. 
The term by term differentiation of the sum defining the distribution $g(z)$ in (\ref{disgz}) gives 
$$
\partial \log g(z) = 
\sum_{\gamma \not = 0}
 \frac{\chi_{z}(\gamma) }{\gamma}d z, \qquad 
\overline \partial g(z) = -  
\sum_{\gamma\not = 0}
 \frac{\chi_{z}(\gamma) }{\overline \gamma} d \overline z.
$$
Thus we get
$$
\int_{E^{p+q}}\Omega_{p,q} = \frac{(-1)^p}{2\pi i}\Bigl(\frac{\tau - \overline \tau}{2\pi i}\Bigr)^{p+q+1}
\sum_{\gamma \not = 0}
 \frac{\chi_{a}(\gamma)}{\gamma^{p+1}\overline \gamma^{q+1}}. 
$$
Indeed, we have a convolution of the Green function with its 
first derivatives, which amounts to a product after the Fourier transform. 

\vskip 3mm
ii) Let $l, \overline l\in H_1(E, \C)$ be the homology classes dual to 
$dz$ and  $d\overline z$. The Hodge correlators of the elements 
(\ref{3.21.05.1})  are naturally  
packaged 
into a generating function
\be \la{KLEK1}
\frac{(-1)^p}{2\pi i} \frac{\tau - \overline \tau}{2\pi i}
\sum_{p,q=0}^{\infty}\Bigl(\frac{\tau - \overline \tau}{2\pi i}\Bigr)^{p+q}
\sum_{\gamma \not = 0}
 \frac{\chi_{a}(\gamma)}{\gamma^{p+1}\overline \gamma^{q+1}}\overline l^p l^q.
\ee 
This is naturally a function on 
a formal neighborhood of the origin in $H^1_{DR}(E, \C)(1)^+$, where $+$ means 
the complex conjugation acting on the forms on $E$. 
On the other hand $K(a|t)$ is naturally a function 
on $t \in \C = H_1(\Gamma) \otimes \R = H_1(E, \R)$. 
Indeed, $\gamma$ and $t$ in formula (\ref{KLEK})  lie in the same space! 

The intersection pairing 
$\langle \ast, \ast\rangle$ on $H_1$ provides 
an isomorphism  $i: H^1_{DR}(E, \C)(1)^+ \to H_1(E, \R)$ 
such that  $\langle i(\omega), h\rangle  := (\omega, h)$. It follows that
\be \la{vart}
\langle l, \overline l \rangle = 
\frac{2\pi i}{\tau - \overline \tau}, \qquad 
i: dz \lms  t:= -\frac{\tau - \overline \tau}{2\pi i}\overline l, \qquad 
i: d\overline z \lms \overline t:= \frac{\tau - \overline \tau}{2\pi i}l.
\ee
Using this we get 
$$ 
 (\ref{KLEK1}) = \frac{1}{2\pi i} \frac{\tau - \overline \tau}{2\pi i}
\sum_{\gamma \not = 0}
 \frac{\chi_{a}(\gamma)}{|\gamma- t|^2} = \frac{1}{2\pi i} K(a|t).
$$
\end{proof}

\paragraph{Example.}  It is easy to see that 
when $a\in E$ is a torsion point, the coproduct of the 
averaged base point motivic correlator of 
${\cal C}(\{0\} \otimes \{a\} \otimes {\rm Sym}^{n}{\Bbb H})$ (Section 
\ref{9.5ref}.4) is zero.
 So we get 
an element  
$$
{\rm Cor}^0_{E-E_{\rm tors}}({\cal C}(\{0\} \otimes \{a\} 
\otimes {\rm Sym}^{n}{\Bbb H}))  \in {\rm Ext}^1_{\cal C}(\Q(0), 
{\rm Sym}^{n}{\Bbb H}(1)). 
$$
Its real period is given by the Eisenstein-Kronecker series. 
When $E$ is a CM curve it gives 
the special values of the Hecke L-series with  
Gr\"ossencharacters. Summarizing, we get the motivic 
cohomology class corresponding to this special value. 
See \cite{B}, \cite{Den1}, \cite{BL} for different approaches. 

\paragraph{2. Hodge correlators on an elliptic curve.} Let 
$\mu$  be the unique invariant volume form on $E$ of volume $1$. 
The Green function corresponding to $\mu$ 
is normalized by 
\begin{equation}  \label{clpolysd}
\int_{E(\C)} G_\mu(s,z)dz\wedge d\overline z  =0.
\end{equation} 
Below we consider the  Hodge correlators decorated by the elements of 
\be \la{12.29.08.1}
{\cal C}\Bigl(\{a_0\} \otimes 
S^{k_0}(\Omega^1_E \oplus \overline \Omega^1_E) \otimes \ldots \otimes \{a_m\} \otimes 
S^{k_m}(\Omega^1_E \oplus \overline \Omega^1_E)\Bigr). 
\ee
The symmetric Hodge correlators on an elliptic curve 
are linear combinations of those which appear when 
$\{a_i\}$ are replaced by degree zero divisors  on $E$. 
Abusing language, we call any Hodge 
correlator on an elliptic curve decorated by (\ref{12.29.08.1}) and defined 
by using the normalized Green function $G_\mu(x,y)$   
a symmetric Hodge correlator. 
The {\it depth} of the Hodge correlator of $W$ is the number of 
the $S$-factors  of $W$ minus one. 

We package depth $m$ symmetric Hodge correlators 
on $E$ into non-holomorphic 
 generating series:
\be \la{1.29.08.2q}
\begin{split}
&K_{p_0, q_0; \ldots ; p_{m}, q_{m}}(a_0: \ldots : a_{m}):= 
{\rm Cor}_{\cal H}{\cal C}\Bigl(\{a_0\} \otimes 
\frac{dz^{p_0}d\overline z^{q_0}}{p_0!q_0!} 
\otimes \ldots \otimes \{a_m\} \otimes 
\frac{dz^{p_m}d\overline z^{q_m}}{p_m!q_m!} \Bigr). \\
&K(a_0: \ldots : a_{m}| l_0: \ldots: l_{m}) = 
\sum_{p_i, q_i\geq 0}K_{p_0, q_0; \ldots ; p_{m}, q_{m}}(a_0: \ldots : a_{m})
\overline l_0^{p_0} l_0^{q_0}\ldots 
\overline l_m^{p_m} l_m^{q_m}.\\
\end{split}
\ee
Following (\ref{vart}), we introduce the variables 
$$
t_i: =-\frac{\tau - \overline \tau}{2\pi i}\overline l_i, \qquad 
\overline t_i: =\frac{\tau - \overline \tau}{2\pi i}l_i. \qquad 
t_i, \overline t_i \in H_1(E, \R). 
$$

\paragraph{3. The multiple Eisenstein-Kronecker series.} 
Our next goal is to present the symmetric Hodge 
correlators on elliptic curves as 
integrals of the generating series $K(a| t)$. 
We rewrite  the function 
(\ref{CEKS}) as 
\be \la{plGreen}
{\bf K}(a_1:a_2| t_1, t_2):= K(a_1-a_2| t_1); \qquad t_1+t_2 =0.
\ee
Let us define, following Section 8 of \cite{G1}, 
the 
{\it depth $m$ multiple Eisenstein-Kronecker series} 
$$
{\bf K}(a_0: \ldots :a_{m}| t_0, \ldots, t_{m}); 
\qquad a_i\in E, \quad t_i \in H_1(E, \R), 
\quad 
t_0 + \ldots + t_{m}=0, 
$$
Consider a plane trivalent tree $T$, see Fig \ref{feyn},  
decorated by 
$m+1$ pairs 
$(a_0, t_0), \ldots , (a_{m}, t_{m})$. 
Each oriented edge $\stackrel{\to}{E}$ of the tree $T$ 
provides an element $t(\stackrel{\to}{E}) \in H_1(E, \R)$ 
just like in Section 10.1. 
We define the multiple Eisenstein-Kronecker series 
just like  the 
multiple Green functions, with the 
generating series 
(\ref{plGreen}) for the classical Eisenstein-Kronecker
 series serving as Green functions: 

\begin{equation} \label{7.1.00.2}
\begin{split}
&{\bf K}(a_1: \ldots: a_{m+1}|t_1, \ldots, t_{m+1}):= (2\pi i)^{-(2m-1)}\sum_{\mbox{T}}
 {\rm sgn}(E_1 \wedge \ldots \wedge E_{2m-1})~\cdot\\
&\int_{E^{m-1}}
\omega_{2m-2}\left({\bf K}(x^{E_1}_1: x^{E_1}_2|t_1^{E_1}, t_2^{E_1}) \wedge \ldots \wedge 
{\bf K}(x^{E_{2m-1}}_1: x^{E_{2m-1}}_2|t_1^{E_{2m-1}}, t_2^{E_{2m-1}}) \right). \\
\end{split}
\end{equation}

Here ${\bf K}(x^{E}_1: x^{E}_2|t_1^{E}, t_2^{E})$ 
is the function (\ref{plGreen}) assigned to the edge $E$, defined 
just like in the rational case. 
The sum in (\ref{7.1.00.2}) is over all plane $3$-valent trees $T$ cyclically 
labeled by the pairs  $(a_i, t_i)$. 
The integral is over the product of copies of $E$ labeled by the internal vertices 
of $T$.

Clearly 
$$
{\bf K}(a+a_0: \ldots :a+a_{m}| t_0, \ldots, t_{m}) = 
{\bf K}(a_0: \ldots :a_{m}| t_0, \ldots, t_{m}).
$$
The multiple Eisenstein-Kronecker series  enjoy 
the dihedral symmetry relations (\ref{1.03.08.11}) 
and the shuffle relations (\ref{1.03.08.12}) 
It is proved just the same way as Lemma \ref{7.1.00.3}.

\bt \la{1.29.08.2}
The generalized Eisenstein-Kronecker series 
are the generating series for the symmetric Hodge correlators 
on the elliptic curve $E$:
\be \la{1.29.08.2q1}
K(a_0: \ldots : a_{m}| u_0: \ldots: u_{m}):= 
{\bf K}(a_0: \ldots : a_{m}| u_m-u_0, u_0-u_1,  \ldots: u_{m-1}-u_m). 
\ee
\et
We visualize this as follows. The points $a_i$ sit on a circle, and 
cut it into $m+1$ arcs. We assign the variables $u_i$ 
 to the arcs so that the vertex $a_i$ shares 
the arcs labeled by $u_{i-1}$ and $u_{i}$, see Fig \ref{hc21}.

\begin{figure}[ht]
\centerline{\epsfbox{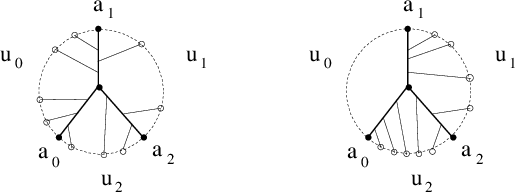}}
\caption{Feynman diagrams for symmetric Hodge correlators on an elliptic curve.}
\label{hc21}
\end{figure}

 \begin{proof}
A  Feynman diagram  
decorated by a vector in (\ref{12.29.08.1}), and  which 
has  
an internal vertex incident to  two $\omega$-decorated external vertices,  
contributes zero. Indeed, the correlator integral involves 
$\int_E (\omega_1 \wedge \omega_2 +\omega_2 \wedge \omega_1) =0$.   
It follows that the only Feynman diagrams 
contributing to the Hodge correlators we consider are 
as on Fig \ref{hc21}.

Symmetric Hodge correlators 
of depths $-1$ and $0$ are zero. Indeed, 
every Feynman diagram with less then two $S$-decorated vertices 
has an internal vertex sharing two $\omega$-decorated 
edges. 

The depth $1$ case follows from Proposition \ref{2.06.08.1}. 
There are other depth $1$ symmetric cyclic words, 
similar to the one on the left of Fig \ref{hc16}, where $\{0\}$'s 
are replaced by $\omega_i$'s. The integral for each of them 
coincides, up to a sign,  
with the one for the standard diagram on Fig \ref{feyn16}. 
Namely,  just flip the $\omega$-decorated edges looking up, 
thus transforming the diagram 
to the one on Fig \ref{feyn16}. 


The rest of the proof repeats the proof of Theorem \ref{1.08.03.2} 
with a slight modification stemming from the fact that 
$H_1(E, \R)$ is two-dimensional. \end{proof}

\paragraph{Remark.} Let us show directly that the generating 
series (\ref{1.29.08.2q}) are invariant under the shift 
$u_i \lms u_i - u_0$. We  use the same argument as in the proof of 
Lemma \ref{deps1c}.  
Each external edge touching the  arc labeled by $u_0$ 
is flipped by moving its external vertex to the other arcs.  
For example, the right picture on Fig \ref{hc21} represents 
the result of flipping of all external edges 
touching the $u_0$-arc on the left one. 
The correlator integral changes the sign under a flip of a single edge.  
Indeed, a flip changes only the  orientation  of the tree. 
This proves the claim.

\paragraph{4. Polylogarithms on curves.} 
It is naturally to define 
{\it polylogarithms on curves} as the motivic/Hodge correlators of depth one.  
For the rational and elliptic curves 
we get this way the classical and elliptic polylogarithms. 
For the higher genus curves there is a yet smaller 
class of functions, defined in \cite{G1} and discussed below. 
Their motivic avatars, however,  
do not span a Lie coalgebra. 

Let $X$ be a smooth compact complex 
curve  of genus $g \geq 1$. 
 and $\mu$ is a volume one measure on $X$. Set 
$$
{\Bbb H}_\R:= H_1(X, \R), \qquad {\Bbb H}_{\C}:= 
{\Bbb H}_\R\otimes \C = {\Bbb H}_{-1,0} \oplus {\Bbb H}_{0,-1}.
$$
For each integer 
$n \geq 1$ we define a $0$-current $G_n(x,y; \mu)$ on $X \times X$ 
with values in 
\begin{equation} \label{7.6.00.1}
{\rm Sym}^{n-1}{\Bbb H}_{\C} (1) = \oplus_{s+t=n-1}
S^{p}{\Bbb H}_{-1,0} \otimes S^{q}{\Bbb H}_{0,-1}\otimes \R(1). 
\end{equation}
The current $G_{1}(x,y; \mu)$ is the  Green function $G_\mu(x,y)$, 
The current $G_{n}(x,y; \mu)$ for $n>1$ is given by a function on $X \times X$. 
To get complex valued functions out of the vector 
function $G_n(x,y; \mu)$ we proceed as follows. 
Let 
$$
\overline \Omega_p \in S^{p}\overline\Omega^1, \quad 
\Omega_q \in S^{q}\Omega^1.
$$
Then $\overline \Omega_p \otimes \Omega_q$ is an element of the dual to 
(\ref{7.6.00.1}). Let us define the pairing 
$\langle G_n(x,y; \mu), \overline \Omega_p\cdot \Omega_q \rangle $.  
Let 
$$
\Omega_p = \omega_{\alpha_1}\cdot \ldots \cdot \omega_{\alpha_{p}}, \qquad 
\Omega_q = \omega_{\beta_1}\cdot \ldots \cdot \omega_{\beta_{q}};\qquad  
\omega_{*} \in \Omega^1.
$$
Denote by $p_i: X^{n-1} \lra X$ the projection on $i$-th factor. 
\begin{definition} \label{7.7.00.1} The $n$-logarithm 
on the curve $X$ is given by
$$
\langle G_n(x,y;a), \overline \Omega_p\cdot \Omega_q \rangle :=\qquad 
$$
$$\frac{1}{(2\pi i)^{n}}{\rm Alt}_{\{t_1, \ldots, t_{n-1}\}}\Bigl(\int_{X^{n-1}}\omega_{n-1}
\Bigl(G(x, t_1) \wedge G(t_1, t_2) \wedge \ldots \wedge G(t_{n-1}, y)\Bigr)  \wedge 
\bigwedge_{i=1}^{p}p_i^*\overline\omega_{\alpha_i} 
\wedge\bigwedge_{j=1}^{t}p^*_{p+j} \omega_{\beta_j}\Bigr).
$$
\end{definition}
Here we skewsymmetrize the integrand with respect to 
$t_1, \ldots, t_{n-1}$. Before the skewsymmetrization 
the integrand depends on the element 
$\overline \omega_{\alpha_1}\otimes \ldots \otimes
 \overline \omega_{\alpha_{p}}$ from $\otimes^{p}\overline \Omega^1$;  
after 
it depends only on its image in $S^{p}\overline \Omega^1$; 
similarly for the second factor. 

It 
is the Hodge correlator for the Feynman diagram 
on Fig \ref{feyn16}. For example, 
the dilogarithm on a curve $X$ is described by the second 
Feynman diagram on Fig \ref{feyn10}. It is given by the integral
$$
\langle G_2(x,y;\mu), \overline \omega_{\alpha}\rangle := \frac{1}{(2\pi i)^{2}}
\quad \int_{X}\omega_1\Bigl(G_\mu(x, t) 
\wedge G_\mu(t, y)\Bigr) \wedge \overline \omega_{\alpha}(t).
$$

\paragraph{Remark.} Unlike in the genus $\leq 1$ case, there are non-trivial Hodge correlators of depths $-1$ and $0$. The simplest depth $-1$ Hodge correlator 
is illustrated on Fig \ref{feyn32}. Furthermore, there are 
other depth one Hodge correlators then the polylogarithms defined above.

\begin{figure}[ht]
\centerline{\epsfbox{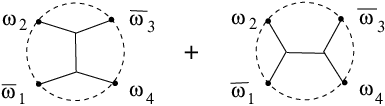}}
\caption{The simplest depth $-1$ Hodge correlator.}
\label{feyn32}
\end{figure}

Just like in the rational and elliptic cases, one can introduce 
cyclic multiple polylogarithms by taking the generating function 
of the polylogarithms on curves as the propagator assigned to an edge 
of a trivalent tree, and repeating the construction of 
Definition \ref{CMP}. 

One proves just the same way as in  Theorem \ref{1.08.03.2}, that these way we get 
all Hodge correlators corresponding 
the elements of 
${\cal C}\Bigl(\{a_0\} \otimes S^{k_0}(\Omega^1_X \oplus \overline \Omega^1_X) 
\otimes \ldots \otimes \{a_m\} \otimes 
S^{k_m}(\Omega^1_X \oplus \overline \Omega^1_X)\Bigr)$.

However, unlike in the cases of rational or  elliptic curves, 
we do not get all Hodge correlators on curves this way, 
and the corresponding motivic correlators are not closed under the coproduct. 
The reason is that when $E$ is an elliptic curve, the direct summands of 
the twisted symmetric powers ${\rm Sym}^mH^1(E)(n)$ of the motive $H^1(E)$ 
give all simple object of the abelian tensor category of pure motives generated by $H^1(E)$. 
However when $X$ is regular projective curve of genus $\geq 2$, the direct summands 
of ${\rm Sym}^mH^1(X)(n)$ do not deliver all simple objects of the tensor category generated by 
$H^1(X)$. 
\section{Motivic correlators on modular curves} \la{hc11sec}

In the Section $X$ is a modular curve, and $S$ is the set of cusps. 

\subsection{Hodge correlators on modular curves and generalized 
Rankin-Selberg integrals}\la{sec11.1ref}

\paragraph{Rankin-Selberg integrals as Hodge correlators on modular curves.} \la{sec11.1}
For an arbitrary curve $X$, 
there are three different types of the Hodge correlators for 
the length three cyclic words. 
They correspond to the 
three Feynman diagrams on Fig. \ref{feyn10}. 
\begin{figure}[ht]
\centerline{\epsfbox{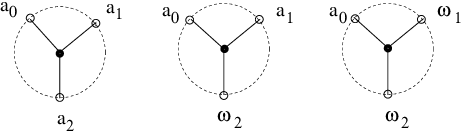}}
\caption{The three possible Feynman diagrams for the $m=2$ correlators.}
\label{feyn10}
\end{figure}
The first is the double Green function  from Section 10.1. Let us 
interpret the other two in the case when $X$ is a modular curve, and $S$ is the set of its cusps. 
Let $\omega$ and $\omega_1, \omega_2$ be cuspidal weight two Hecke eigenforms.  They are  holomorphic $1$-forms  on $X$. 
Let $s \in S$ be a cuspidal divisor on the modular curve $X$. Then 
$G(s, t)$ is an Eisenstein series, see, for example \cite{GZ}. 
By the Manin-Drinfeld theorem \cite{Ma}-\cite{Dr3} if ${\rm deg}(s)=0$, there exists 
a rational function $g_s$ on $X$ such that ${\rm div}g_s = N\cdot s$. So $G(s, t) = N \log|g_s|^2$.  
Setting 
$a \sim_{\Q^*} b$ if $a/b \in \Q^*$, we get 
\be
\begin{split}
&{\rm Cor}_{\cal H}\Bigl((\{a_0\}\otimes \{a_1\} \otimes \omega)\Bigr) = 
\int_{X(\C)} G(a_0, t) d^\C G(a_1, t) \wedge \omega ~\sim_{\Q^*}~
\int_{X(\C)} \log|g_{a_0}| d^\C \log|g_{a_1}| \wedge \omega, \\
&{\rm Cor}_{\cal H}
\Bigl({\cal C}(\{a_0\}\otimes \omega_1 \otimes \overline \omega_2)\Bigr) = 
\int_{X(\C)} G(a_0, t)\omega_1 \wedge \overline \omega_2  ~\sim_{\Q^*} ~\int_{X(\C)}
 \log|g_a| \omega \wedge \overline \omega.  \\
\end{split}
\ee
Let  $L(\omega, s)$ be the $L$-function of the rank two weight $-1$ 
motive $M_\omega$ corresponding to the Hecke eigenform $\omega$, 
and $L(\omega_1\otimes \omega_2, s)$ the $L$-function of the product of 
motives corresponding to different Hecke eigenforms $\omega_1$ and  $\omega_2$. 
According to the classical Rankin-Selberg method  
the first integral is proportional, up to some standard factors, to $L(\omega, 2)$, 
and the second one to 
$L(\omega_1\otimes \omega_2, 2)$. Moreover, there exist degree zero divisors $a_0$, $a_1$ such that the proportionally 
coefficient is not zero \cite{B}.

\paragraph{Generalized Rankin-Selberg integrals.} 
Let $X$ be a modular curve, given together with  
its uniformization $p_X: {\cal H} \to X$  by the hyperbolic plane ${\cal H}$. 
Choose a hyperbolic metric of curvature $-1$ on $X$. Denote by 
$\mu_X$ the corresponding volume $1$ volume form on $X$. So 
$p_X^*\mu_X$ is the standard volume form $dxdy/y^2$ on the hyperbolic plane. 
Therefore if $\pi: Y \to X$ is a natural map of modular curves, i.e. a map 
 commuting with 
the uniformization maps, 
then the volume forms $\mu_X$ and $\mu_Y$ are compatible: 
 $\pi_*\mu_Y = \mu_X$. 

\begin{lemma} \label{7.9.06.3ed} Let $\pi: Y\to X$ be a map of curves 
and $\pi_*(\mu_{Y}) = \mu_X$ for certain measures $\mu_Y, \mu_X$ on them. 
Let $\pi_1:Y\times Y \to X \times Y$ and $\pi_2: X\times Y \to X \times X$ 
be two natural maps. Then 
$$
\pi_{1*}G_{\mu_Y}(y_1, y_2) = \pi^*_{2}G_{\mu_X}(x_1, x_2) +C. 
$$
\end{lemma}

\begin{proof}  The Green function is defined uniquely 
up to adding a constant by the differential equation 
(\ref{7.3.00.1}). The latter 
involves the volume form and the identity maps on the cohomology 
realized by the 
delta-currents $\delta_\Delta$ and the Casimir elements 
for the space of holomorphic/antiholomorphic  $1$-forms. Clearly  
$\pi_{1*}\delta_{\Delta_Y} = \pi^*_{2}\delta_{\Delta_X}$.
The latter fact  implies a  similar identity for the Casimirs. Finally, 
\be
\begin{split}
&\pi_{1*}(1_Y\otimes \mu_Y ) = {\rm deg}\pi~  (1_X\otimes \mu_Y)  =  
\pi^*_{2}(1_X\otimes \mu_X ); \\
&\pi_{1*}(\mu_Y \otimes 1_Y) = \mu_X \otimes 1_Y =  
\pi^*_{2}(\mu_X \otimes 1_X).\\
\end{split}
\ee
\end{proof}

We set ${\Bbb V}_{X, S}:= H^1(X)\oplus \Q[S](-1)$. 

We say that Green functions $G_{\mu_Y}$ and $G_{\mu_X}$ 
are compatible if the constant in Lemma \ref{7.9.06.3ed} is zero. 
Choose a compatible collection of Green functions $G_{\mu_X}$ 
for the hyperbolic 
metrics on the modular curves $X$.   
We  get a {\it hyperbolic Hodge correlator map} on the modular curves:\footnote{Hyperbolic Hodge 
correlators are a special case 
of the Hodge correlators, specified by the choice of the metric.}
$$ 
{\rm Cor}^h_{{\cal H}, \mu_X}: 
{\cal C}{\rm T}({\Bbb V}_{X, S})  \lra H^2(X).
$$
The integral over the fundamental cycle of $X$ 
provides an isomorphism $H^2(X) \to \C$. So we may assume that for an individual 
curve $X$ the correlator is a complex number. However 
moving from $X$ to $Y$ we multiply it by ${\rm deg}\pi$. 

\begin{proposition} \label{7.9.06.3df} 
The hyperbolic Hodge correlator maps are 
compatible with natural projections $\pi: Y \to X$ 
of modular curves: 
\be \la{2.08.01.2}
{\rm Cor}^h_{{\cal H}, \mu_{Y}}{\cal C}\Bigl(\pi^* v_1 \otimes 
\ldots \otimes \pi^* v_{m} \Bigr) = 
\pi^*{\rm Cor}^h_{{\cal H}, \mu_{X}}{\cal C}\Bigl(v_1 \otimes 
\ldots \otimes v_{m} \Bigr)\in H^2(Y), \qquad v_i \in {\Bbb V}_{X, S}.
\ee
\end{proposition}

\begin{proof} Lemma \ref{7.9.06.3ed} 
plus the projection formula implies that, 
given $1$-forms $\omega_i $ on $X$, we have 
\be \la{2.08.01.1}
\int_{Y}\pi^*\omega_1\wedge \pi^*\omega_2 \cdot G_{\mu_Y}(y_1, y_2) = 
\int_{X}\omega_1\wedge \omega_2 \cdot \pi_{1*}G_{\mu_Y}(x_1, y_2) = 
  \int_{X}\omega_1\wedge \omega_2 \cdot \pi^*_{2}G_{\mu_X}(x_1, y_2).
\ee
Here the integrals are over $y_1$ and $x_1$. 
Similarly,  
we have 
\be \la{2.08.01.1a}
G_{\mu_Y}(\pi^*\{s\}, y_2) = \pi_2^*G_{\mu_X}(s, x_2). 
\ee
Therefore we have an identity  similar to (\ref{2.08.01.1}) where one or 
two of the decorating forms 
$\omega_i$'s are replaced by the decorating divisors $\pi^*\{s_i\}$ -- 
we understood (\ref{2.08.01.1}) as the Hodge correlator for a tree with one 
internal trivalent vertex.  

Let us use the induction on $m$. We understood the correlators are numbers 
via the isomorphism $H^2(X) \stackrel{\sim}{\to} \C$. 
The $m=3$ case follows from (\ref{2.08.01.1a}), using
$$
\int_{Y}\pi^*\eta_1\wedge \pi^*\eta_2 \cdot \pi^*\varphi_3 = {\rm deg}\pi 
\int_{X}\eta_1\wedge \eta_2 \cdot \varphi_3.
$$
For the induction step, 
take a plane trivalent tree $T$ decorated by 
$\pi^* v_1 \otimes \ldots \otimes \pi^* v_{m}$. 
Take two external edges sharing an internal vertex $w$ 
and decorated by $\pi^*v_i$ and $\pi^*v_{i+1}$. 
Applying (\ref{2.08.01.1}) to the integral over the 
curve $Y$ corresponding to $w$, 
we reduce the claim to the case of a tree with $m-1$ external vertices. 
\end{proof}

Let ${\cal M}$ be the universal modular curve, that is the projective limit of the tower of modular curves. 
It is defined over $\overline \Q$. 
Since Hecke operators split the Grothendieck motive  
$H^1({\cal M})$ into a direct sum of the cuspidal and Eisenstein parts, 
${\rm gr}^WH^1({\cal M})\otimes \overline \Q = 
H^1({\cal M})\otimes \overline \Q$. So one has  
\begin{equation} \label{sdds}
{\Bbb H}^{\vee}_{\cal M} = 
H^1({\cal M})\otimes \overline \Q = \bigoplus_M M^{\vee} 
\otimes V_M \bigoplus \Q(\widehat 
{\rm Cusps})_0(-1).
\end{equation}
Here the first sum is over all pure rank $2$, weight $-1$ motives $M$ over
$\Q$. The subscript $0$ in the second summand 
is used to denote the subspace of measures functions with zero mass. The $\overline \Q$-vector space $V_M$ 
is a representation of the finite adele group 
$GL_2({\Bbb A}^f_\Q)$. The tensor product of the 
Steinberg representation of $GL_2(\R)$ and $V_M$ is an automorphic 
representation of $GL_2({\Bbb A}_\Q)$ 
corresponding to the pure motive $M$ via the Langlands
correspondence. 

In the de Rham realization, we get an adelic description of 
$H_{DR}^1({\cal M}, \C)$:
$$
H_{DR}^1({\cal M}, \C) = \bigoplus_M (\Omega^1_M \oplus \overline \Omega^1_M) 
\otimes V_M \bigoplus \Q(\widehat 
{\rm Cusps})_0(-1).
$$
Here $\Omega^1_M$ is the $(1,0)$-part of the de Rham realization  
of the motive $M$. 
Proposition \ref{7.9.06.3df} implies that 
the hyperbolic Hodge correlator map provides 
a $GL_2({\Bbb A}_\Q^f)$-covariant map 
$$
{\rm Cor}^h_{{\cal H}}: 
{\cal C}{\rm T}(H_{DR}^1({\cal M}, \C))(1) \to \C.
$$ 
Indeed, it is $GL_2(\widehat \Z)$-covariant on the nose, 
and $GL_2(\Q)$-covariant thanks to the compatibility 
with the projections established in Proposition \ref{7.9.06.3df}. 
Since these two groups generate $GL_2({\Bbb A}_\Q^f)$, it descends to a map of $GL_2({\Bbb A}_\Q^f)$-coinvariants
\be \la{1.08.2.1q}
{\rm Cor}^h_{{\cal H}}: 
{\cal C}{\rm T}(H_{DR}^1({\cal M}, \C))(1)_{GL_2({\Bbb A}_\Q^f)} \lra \C.
\ee

The map (\ref{1.08.2.1q}) is a generalization of the Rankin-Selberg integrals: 
the latter appear for the triple tensor product of  copies of 
$H_{DR}^1({\cal M}, \C)$, as explained in \ref{sec11.1}.1. 
Restricting this map to the subspace 
$\overline {{\cal C}{\cal L}ie}^{\vee}(H_{DR}^1({\cal M}, \C))$
of the cyclic tensor product of $H_{DR}^1({\cal M}, \C)$ 
(Section \ref{9.4.2ref}) we get the periods 
of $\pi_1^{\rm nil}({\cal M})$.

\subsection{Motivic correlators  in towers} \la{sec11.2}
Let $\pi: X'\to X$ be a nontrivial map of projective curves, $S'$ and $S$ 
are subsets of 
points of $X'$ and $X$, and $\pi(S') = S$. 
We say that vectors $v_{s'}$ and $v_s$ at the points $s'\in S'$ and $s\in S$ 
are compatible with respect to the map $\pi$ if for local parameters 
$t'$ and $t$ such that $dt'$ and $dt$ are dual to the vectors $v'$ and $v$ 
the restriction of  $(t')^{m_{\pi}(s')}/\pi^*t$ to $s'$ equals $1$, 
where    $m_{\pi}(s')$ is the 
multiplicity of the map $\pi$ at the point $s'$. In particular, if 
$m_{\pi}(s')=1$, this means that $d\pi_*(v') =v$. 
Let us assume that for each $s \in S$ (respectively $s' \in S'$)
 there is a non-zero 
vector $v_s \in T_sX$ (respectively $v_{s'} \in T_{s'}X$), 
and that these vectors are compatible 
under the map $\pi$. 
Denote by ${\rm deg}(\pi)$ the degree of the map $\pi$. Recall 
$$
{\Bbb V}^{\vee}_{X, S}:= H^1(X) \oplus \Q[S](-1); \qquad  
{\Bbb V}^{\vee}_{X, S^*} \stackrel{\sim}{=}  {\rm gr}^WH^1(X-S). 
$$

There are maps 
\be \la{15:47}
\pi^*:{\Bbb V}^{\vee}_{X, S} \lra {\Bbb V}^{\vee}_{X', S'}, \quad 
\pi_*:{\Bbb V}^{\vee}_{X', S'} \lra {\Bbb V}^{\vee}_{X, S}, \quad 
\pi_*\pi^*= {\rm deg}(\pi).
\ee
The restriction of the map $\pi^*$ to $\Q[S](-1)$ is given by the following map: 
$$
\pi_S^*: \Q[S] \to \Q[S']; \quad \{s\} \to \sum_{s' \in \pi^{-1}(s)}m_{\pi}({s'})\{s'\}
$$
So ${\rm deg}\pi_S^*(\{s\}) = {\rm deg}\pi$. 
The maps $\pi^*$ and $\pi_S^*$ induce a map of vector spaces 
$$
\pi^*_{\otimes}:{\cal C}{\rm T}({\Bbb V}^{\vee}_{X, S})\otimes \Q[S]
\lra {\cal C}{\rm T}({\Bbb V}^{\vee}_{X', S'})\otimes \Q[S'].
$$

Motivic correlators require a choice of a tangential base point. 
One can consider the motivic correlator map obtained by using an arbitrary 
tangential base point. Then its source is the space 
${\cal C}{\rm T}({\Bbb V}^{\vee}_{X, S})\otimes \Q[S]$, where the second factor 
$\Q[S]$ bookkeeps the choice of the point $s\in S$ where we put the tangential base point.

\begin{lemma} \la{qwerf}
Assume that the tangential base points 
at $S'$ and $S$ are compatible with the map $\pi: X'-S' \to X-S$. 
Then there is a commutative diagram, with the horizontal rows 
provided by (\ref{der1}), and 
the vertical maps are induced by the map $\pi^*$: 
$$
\begin{array}{ccc}
{\cal C}{\rm T}({\Bbb V}^{\vee}_{X', S'})\otimes \Q[S']& 
\stackrel{{\rm Cor}_{X'- S'}}\lra & {\cal L}_{\rm Mot}\otimes H^2(X')\\
\uparrow \pi_{\otimes}^*&&\uparrow \\
{\cal C}{\rm T}({\Bbb V}^{\vee}_{X, S})\otimes \Q[S]&
\stackrel{{\rm Cor}_{X- S}}{\lra} & {\cal L}_{\rm Mot}\otimes H^2(X)
\end{array}
$$
\el

\begin{proof}  If $\pi$ is unramified at every point of $S'$ 
 it follows straight from the definition. 
The general case is deduced from this, say, by the specialization to a ramified point. 
\end{proof} 

Let 
${\rm tr}_X:H^2(X) \to \Q(-1)$ be the canonical 
map. Since
 ${\rm tr}_{X'}\pi^* = {\rm deg}(\pi){\rm tr}_X$, we get 
a commutative diagram
\begin{equation} \label{7.8.06.1}
\begin{array}{ccc}
{\cal C}{\rm T}({\Bbb V}^{\vee}_{X', S'})\otimes \Q[S'](1)
& \stackrel{{\rm Cor}_{X'- S'}}\lra & {\cal L}_{\rm Mot}\\
&&\\
~~~~~\uparrow \pi^*\otimes {\rm deg}(\pi)^{-1}\pi_{\otimes}^*&&\uparrow =\\
&&\\
{\cal C}{\rm T}({\Bbb V}^{\vee}_{X, S})\otimes \Q[S](1)
&\stackrel{{\rm Cor}_{X- S}}{\lra} & {\cal L}_{\rm Mot}
\end{array}
\end{equation}

\paragraph{Remark.} 
Denote by  $\pi_{\cal C}^*$ 
the left vertical map in (\ref{7.8.06.1}). It  
is not a map of Lie coalgebras. However it will become a 
Lie coalgebra map if we consider the quotient 
of the Lie coalgebra  ${\cal C}{\rm T}({\Bbb V}^{\vee}_{X', S'})(1)$ defined in Definition \ref{7.9.06.2}. 

\vskip 3mm
Now suppose that we have a projective system $(X_\alpha, S_\alpha)$ 
of curves $X_\alpha$ with finite subsets of 
points $S_\alpha$ on them and maps between them as above. 
Using the  normalized maps ${\rm deg}(\pi)^{-1}\pi^*_{\otimes}$, 
functions on the sets $S_\alpha$ in the limit 
turn into measures on the projective limit   
$\widehat S:= \lim S_\alpha$. We denote the space of measures by 
${\rm Meas}(\widehat S)$. 
So using (\ref{7.8.06.1}), we get a map of the corresponding inductive limit
\begin{equation} \label{7.8.06.1r}
\lim_{\lra} {\cal C}{\rm T}({\Bbb V}^{\vee}_{X_\alpha, S_\alpha})
\otimes {\rm Meas}(\widehat S)(1) \lra 
 {\cal L}_{\rm Mot}.
\end{equation}
Here ${\cal X}$ denotes the projective limit of the curves $X_\alpha-S_\alpha$. 
Abusing notation, 
we denote below the inductive limit in (\ref{7.8.06.1r}) by 
${\cal C}{\rm T}({\Bbb V}^{\vee}_{{\cal X}})\otimes {\rm Meas}(\widehat S)_0(1)$.

\paragraph{Examples.} 1. Let $E$ be an elliptic curve. 
The curves $E - E[N]$ and the isogenies 
$E - E[MN] ~\to~ E - E[N]$ form a tower. 
We assign to every missing point  on $E-E[N]$ 
the tangent vector $\frac{1}{N}v_\Delta$. We get 
a compatible family of tangential base points for the tower $\{E-E[N]\}$.  
Denote by $\widehat {E_{\rm tors}}$ the Tate module of $E$, that is the 
projective limit of the torsion groups of $E$. 
Therefore we arrive at the motivic correlator map, where the subscript $\widehat {E_{\rm tors}}$ refers to the coinvariants 
for the action of the group $\widehat {E_{\rm tors}}$: 
$$
{\rm Cor}_{\cal E}: \Bigl({\cal C}{\rm T}({\Bbb V}^{\vee}_{\cal E})\otimes 
{\rm Meas}(\widehat {E_{\rm tors}}(1))\Bigr)_{\widehat {E_{\rm tors}}} 
\lra {\cal L}_{\rm Mot}.
$$
There is a unique $\widehat {E_{\rm tors}}$-invariant measure $\mu_{\cal E}^0$ 
on the Tate module $\widehat {E_{\rm tors}}$. 
So we get a map 
$$
{\rm Cor}^0_{{\cal E}}: 
{\cal C}{\rm T}({\Bbb V}^{\vee}_{\cal E})(1)
\otimes \mu_{\cal E}^0
\lra {\cal L}_{\rm Mot}.
$$

2. Similarly, when $X = {\Bbb P}^1$ and 
$S = \{0\}\cup \{\infty\} \cup \mu_{l^\infty}$ we get the 
measures on $\Z_l(1)$.

\begin{definition} \label{7.9.06.2}
 ${\cal C}'{\rm T}({\Bbb V}^{\vee}_{X',S'})$ is the quotient of 
${\cal C}{\rm T}({\Bbb V}^{\vee}_{X',S'})$ 
by the subobject spanned by expressions
\begin{equation} \label{7.16.06.1}
{\cal C}\Bigl(\pi^* v_1 \otimes \ldots 
\otimes \pi^* v_{m-1}\otimes v_{m} \Bigr), \quad 
v_1, \ldots , v_{m-1} \in {\Bbb V}^{\vee}_{X,S}, 
\quad v_m \in  {\rm Ker}\pi_* \subset {\Bbb V}^{\vee}_{X',S'}.
\end{equation}
\end{definition}

\begin{lemma} \label{7.9.06.1t}
The correlator map ${\rm Cor}_{X-S}$ annihilates elements 
$$
{\cal C}\Bigl(\pi^* v_1 \otimes \ldots 
\otimes \pi^* v_{m-1}\otimes v_{m}\Bigr) \otimes \pi_S^*\{s\}(1) ~~ \mbox{such that} ~~\pi_*v_m=0.
$$
\end{lemma}

\begin{proof} We prove this in the $l$-adic realization. 
A proof for another realizations is similar. 
It is sufficient to prove the case when $\pi:X'-S'\to X-S$ is a 
Galois cover with a group $G$.  The action of  ${\rm Gal}(\overline F/F)$ 
commutes with the action of $G$. 
Thus 
\be
\begin{split}
&{\rm Cor}_{X'-S'}{\cal C}(\pi^*v_1 \otimes 
\ldots \otimes \pi^*v_{m-1}\otimes v_m) \otimes \pi_S^*\{s\} =\\
&{\rm Cor}_{X'-S'}\frac{1}{|G|}\sum_{g\in G}g^*{\cal C}(\pi^*v_1 \otimes 
\ldots \otimes \pi^*v_{m-1}\otimes v_m)\otimes \pi_S^*\{s\} = \\
&{\rm Cor}_{X'-S'}{\cal C}(\pi^*v_1 \otimes 
\ldots \otimes \pi^*v_{m-1}\otimes \pi_*v_m)\otimes \pi_S^*\{s\} = 0.
\end{split}
\ee
\end{proof}

Here is an analytic counterpart:
\begin{lemma} \label{7.9.06.3} Let $\pi:X'\to X$ and $\pi_*(\mu_{X'}) = \mu_X$.
Then 
$$
{\rm Cor}_{{\cal H}, \mu_{X'}}{\cal C}\Bigl(\pi^* v_1 \otimes 
\ldots \otimes \pi^* v_{m-1}\otimes v_{m} \Bigr) = 0 ~~\mbox{for}~~ 
m \geq 3 ~~\mbox{and}~~ \pi_*v_m=0.
$$
\end{lemma}

\begin{proof}  
We prove it by the induction on $m$. 
The $m=3$ case is clear from (\ref{2.08.01.1}). 
Take a plane trivalent tree $T$ decorated by 
$\pi^* v_1 \otimes \ldots \otimes \pi^* v_{m-1}\otimes v_{m}$. 
Take two external edges sharing an internal vertex $w$ and  
decorated by $\pi^*v_i$ 
and $\pi^*v_{i+1}$. Using  (\ref{2.08.01.1}) for the integral over the 
curve assigned to $w$ we reduce the claim 
to the case of a tree with $m-1$ external vertices. 
\end{proof}

\begin{lemma} \label{7.9.06.1}
${\cal C}'{\rm T}({\Bbb V}^{\vee}_{X',S'})(1)$ inherits a Lie coalgebra structure; 
$\pi_{\cal C}^*$ is a Lie coalgebra map. 
\end{lemma}

\begin{proof}   Let us check that elements (\ref{7.16.06.1}) span a coideal 
for the cobracket $\delta = 
\delta_{\rm Cas}+\delta_{S}$. 
A typical term of $\delta_{\rm Cas}(\ref{7.16.06.1})$ is $\sum_k{\cal C}(M_1 \otimes \alpha_k) \wedge 
{\cal C}(\alpha^{\vee}_k \otimes M_2)$, where $\sum_k \alpha_k\otimes \alpha^{\vee}_k$ is the Casimir element in 
$H^1(X') \otimes H^1(X')$, and ${\cal C}(M_1\otimes M_2)$ is the element (\ref{7.16.06.1}). There is a decomposition
$$
H^1(X') = \pi^*H^1(X) \oplus \pi^*{H^1(X)}^{\perp}. 
$$
The summands 
are orthogonal for the intersection 
pairing on $H^1(X')$. The map $\pi_*$ annihilates the second summand. 
We may assume that $v_m$ enters $M_2$. 
If $\alpha_k \in \pi^*H^1(X)$, then $\alpha_k^{\vee} \in \pi^*H^1(X)$, 
and ${\cal C}(\alpha^{\vee}_k \otimes M_2)$ 
is of type (\ref{7.16.06.1}). Otherwise 
$\alpha_k, \alpha_k^{\vee} \in \pi^*{H^1(X)}^\perp$, thus 
${\cal C}(\alpha_k \otimes M_1)$ is of type (\ref{7.16.06.1}). The 
$\delta_{S}$ term deals with the 
$S$-decorated vertices. Its typical term is 
$$
\sum_i {\cal C}(M_1  \otimes s_i) 
\wedge {\cal C}(s_i \otimes M_2) = 
\sum_i {\cal C}(M_1  \otimes (s_i-s_1)) \wedge {\cal C}(s_i \otimes M_2) +
{\cal C}(M_1  \otimes s_1)  
\wedge {\cal C}(\sum_i s_i \otimes M_2).
$$
 The left factor of the first summand, and the right factor
 of the second  are of type (\ref{7.16.06.1}). \end{proof}

\begin{corollary} \label{7.9.06.1s}
Set $\mu_S:= \frac{1}{|S|}\sum_{s\in S} \{s\}$. 
The map ${\rm Cor}_{X-S}$ gives rise  to a Lie coalgebra map. 
$$
{\rm Cor}^{\rm av}_{X-S}: {\cal C}'{\rm T}({\Bbb V}^{\vee}_{X,S})\otimes \mu_S \otimes 
H_2(X) \stackrel{}\lra  {\cal L}_{\rm Mot}. 
$$
\end{corollary}

\begin{proof} Follows from Lemmas \ref{7.9.06.1t} and \ref{7.9.06.1}. 
\end{proof}

\subsection{Motivic correlators for the universal modular curve} 
Recall the universal modular curve ${\cal M}$ and decomposition 
(\ref{sdds}). 
The tower of modular curves has a natural family of 
tangential base points at the cusps provided by  
exponents of the canonical parameter $z$ on the upper half plane. 
Precisely, if the stabilizer of a cusp is an index $N$ subgroup 
of the maximal unipotent subgroup in $SL_2(\Z)$, it is 
${\rm exp}(2\pi iz/N)$. They are compatible, in the sense of Section 
\ref{sec11.2}, with 
the natural projections of the modular curves. 
Going to the limit in the tower of modular curves, 
as explained in Section \ref{sec11.2}, 
 we arrive at the following picture.  
There is a Lie coalgebra in the $\otimes$-category of pure motives generated by $H^1(\overline {\cal M})$: 
$$ 
{\cal C}'{\rm T}({\Bbb V}^{\vee}_{\cal M})
\otimes {\rm Meas}(\widehat {\rm Cusps})_0(1).
$$ 
The group $GL_2(\widehat \Z)$ acts by  its 
automorphisms. The group $GL_2(\Q)$ acts by its automorphisms 
thanks to the compatibility with projections (Section \ref{sec11.2}). 
So the group  $GL_2({\Bbb A}^f_\Q)$ acts. 
Moreover, the canonical map of Lie
coalgebras 
$$
{\rm Cor}_{\cal M}: {\cal C}'{\rm T}({\Bbb V}^{\vee}_{\cal M})\otimes 
{\rm Meas}(\widehat {\rm Cusps})_0(1) \lra {\cal L}_{\rm Mot}  
$$
is  $GL_2({\Bbb A}^f_\Q)$-equivariant. 
(The group  acts trivially on the right).  
So it descends to a map of the coinvariants
$$
{\rm Cor}_{\cal M}: \Bigl({\cal C}'{\rm T}({\Bbb V}^{\vee}_{\cal M})\otimes 
{\rm Meas}(\widehat {\rm Cusps})_0(1)\Bigr)_{GL_2({\Bbb A}^f_\Q)} \lra {\cal L}_{\rm Mot}.  
$$
In the de Rham realization there is a commutative diagram

\begin{displaymath}
    \xymatrix{
    \Bigl({\cal C}'{\rm T}(H^1_{DR}({\cal M})))\otimes 
{\rm Meas}(\widehat {\rm Cusps}_0(1)\Bigr)_{GL_2({\Bbb A}^f_\Q)} ~~~~\ar[r]^{~~~~~~~~~~~~~~~~~~~~~~~~~~~{\rm Cor}_{\cal M}}  
\ar[dr]_{{\rm Cor}_{\cal H}}& 
{\cal L}_{\rm Mot}  \ar[d]^{{\cal P}}  \\
                                 &   
\C\\
                                 }
\end{displaymath}

where the diagonal arrow is the Hodge correlator map, 
and the vertical one is the period map. 

\begin{lemma} 
There is a unique $GL_2({\Bbb A}^f_\Q)$-covariant volume $1$ 
measure $\mu^0 \in {\rm Meas}(\widehat {\rm Cusps})$, 
 \el

\begin{proof}  Let $U$ 
be the upper triangular
 unipotent subgroup of $GL_2$ and $B$ the normalizing it Borel subgroup.  
Then  there are isomorphisms 
$$
\widehat {\rm Cusps} = 
GL_2(\widehat \Z)/U(\widehat \Z) = 
GL_2({\Bbb A}^f_\Q)/ B(\Q)U({\Bbb A}^f_\Q). 
$$ 
The first is obvious since the cusps on the modular curve $Y(N)$ are identified with the coset 
$GL_2(\Z/N\Z)/U(\Z/N\Z)$, the second follows from 
the decomposition 
$GL_2({\Bbb A}^f_\Q) = B(\Q)GL_2(\widehat \Z)$. 
A Haar measure on $GL_2({\Bbb A}^f_\Q)$ induces 
a covariant measure on the quotient. We normalize it 
by the volume one condition. 
\end{proof} 

Taking the covariant measure $\mu^0$ we are getting a map 
\be \la{zerocor}
{\rm Cor}^0_{\cal M}: {\cal C}'{\rm T}
({\Bbb V}^{\vee}_{\cal M})(1)_{GL_2({\Bbb A}^f_\Q)}\otimes 
\mu^0 \lra {\cal L}_{\rm Mot}. 
\ee

\begin{lemma} \label{7.6.06.1}
The motivic correlators for the classical Rankin-Selberg integrals 
(Section \ref{sec11.1ref}) are motivic ${\rm Ext}^1$'s:
$$
{\rm Cor}^0_{\cal M}{\cal C}(\{a_0\} \otimes \{a_1\}\otimes M^{\vee}_{\omega})
 \in {\rm Ext}^1_{{\rm Mot}}(\Q(0), M_\omega(1)).
$$
$$
{\rm Cor}^0_{\cal M}{\cal C}(\{a_0\} \otimes M^{\vee}_{\omega_1} \otimes M^{\vee}_{\omega_2})\in 
{\rm Ext}^1_{{\rm Mot}}(\Q(0), M_{\omega_1} \otimes M_{\omega_2} ).
$$
\end{lemma}

\begin{proof} The first and the second motivic correlators on the left 
live, respectively,  in the $M_\omega(1)$-isotypical component, and 
$M_{\omega_1} \otimes M_{\omega_2} $-isotypical component of the motivic Lie coalgebra. 
The first defines an element of ${\rm Ext}^1_{{\rm Mot}}(\Q(0), M_\omega(1))$ if and only if 
its coproduct there is zero. The second defines an element of 
${\rm Ext}^1_{{\rm Mot}}(\Q(0), M_\omega(1)\otimes M_{\omega_2} )$ if and only if 
its coproduct is zero. 

 Let us compute the second coproduct.  Set 
${\Bbb V}:= {\Bbb V}_{\overline {\cal M}}$. 
$$
\delta{\rm Cor}^0_{\cal M}{\cal C}(\{a_0\} \otimes M^{\vee}_{\omega_1}
\otimes M^{\vee}_{\omega_2}) =
 {\rm Cor}^0_{\cal M}{\cal C}(\{a_0\} \otimes {\Bbb V})
\wedge {\rm Cor}^0_{\cal M}{\cal C}({\Bbb H} 
\otimes M^{\vee}_{\omega_1}\otimes M^{\vee}_{\omega_2}). 
$$
One has 
${\rm Cor}_{\cal M}{\cal C}(\{a_0\} \otimes {\Bbb V})(1) \otimes \{s\} = a_0-s \in 
J_X \otimes \Q = {\rm Ext}^1_{\rm Mot}(\Q(0), {\Bbb V} )$. 
Since $s, a_0$ are cusps, it is zero by the Manin-Drinfeld theorem. 
Thus   
${\rm Cor}^0_{\cal M}{\cal C}(\{a_0\} \otimes {\Bbb V})(1)=0$, and hence 
the coproduct is zero. This implies the second claim. 
The proof of the first claim is similar. 
\end{proof}

For certain divisors $s_i$ these elements are not zero, and 
hypothetically generate the ${\rm Ext}^1$-groups. 

\vskip 3mm
How
 to describe the image of the motivic correlator map 
(\ref{zerocor})? 
Recall the hypothetical abelian category  ${\cal C}_{\cal M}$ 
of mixed motives generated by $H^1({\cal M})$.  
Let ${\cal P}_{\cal M}$ be the semi-simple category of pure motives over $\overline \Q$ 
 generated by $H^1({\cal M})$. Since $\overline \Q$  is of arithmetic dimension $1$, 
the motivic Lie algebra of the mixed category ${\cal C}_{\cal M}$ 
is conjectured to be a free Lie algebra in the category ${\cal P}_{\cal M}$. 
Let ${\cal L}_{\cal M}$ be the dual motivic Lie coalgebra. 
The map (\ref{zerocor})  lands in ${\cal L}_{\cal M}$: 
\be \la{zerocor123}
{\rm Cor}^0_{\cal M}: {\cal C}'{\rm T}
({\Bbb V}^{\vee}_{\cal M})(1)_{GL_2({\Bbb A}^f_\Q)}\otimes 
\mu^0 \lra {\cal L}_{\cal M}. 
\ee
The space of its cogenerators of type $N^{\vee}$, where $N$ is 
a pure motive of negative weight 
from ${\cal P}_{\cal M}$,  
should be isomorphic to 
${\rm Ext}^1(\Q(0), N) \bigotimes N^{\vee}$; 
Beilinson's conjectures predict its dimension. 

The map (\ref{zerocor123}) 
is not surjective for a very simple reason. 
The 
motivic correlators do not give 
non-trivial elements of  the Jacobian $J_{\cal M}\otimes \Q$ of the universal 
modular curve ${\cal M}$. 
Indeed, we get in the image of (\ref{zerocor123}) only the images of 
degree zero divisors supported at the cusps (Section \ref{hc8sec}.4.1), 
which are zero by the Manin-Drinfeld 
theorem. On the other hand the simplest component 
of  ${\cal L}_{\cal M}$ is given by 
${\rm Ext}^1(\Q(0), H_1({\cal M})) = J_{\cal M}\otimes \Q$. 
However Lemma \ref{7.6.06.1} and the Rankin-Selberg method 
suggest that we do get ${\rm Ext}^1(\Q(0), N)$ for 
 $N = H_1({\cal M})(1)$ and 
$N = H_1({\cal M}) \otimes H_1({\cal M})$.

\vskip 3mm
More specifically, let $V_{\omega_i}$ (resp. $M_{\omega_i}$) be the representation of 
$GL_2({\Bbb A}^f_{\Q})$ 
(resp. the weight $-1$ motive) corresponding to a Hecke eigenform $\omega_i$. 
It  follows from (\ref{zerocor123}) that the motivic correlator provides a map 
\begin{equation} \label{21:33}
{\rm Cor}^0_{\cal M}: {\cal C}'(\bigotimes_{i=1}^kM^{\vee}_{\omega_i} \otimes V_{\omega_i} )(1) 
\lra {\cal L}_{\otimes_{i=1}^kM_{\omega_i}(-1)} 
\bigotimes \otimes_{i=1}^kM^{\vee}_{\omega_i}(1).
\end{equation} 
It is essentially described by a map of vector spaces
$
{\cal C}'(\bigotimes_{i=1}^kV_{\omega_i} ) 
\lra {\cal L}_{\otimes_{i=1}^kM_{\omega_i}(-1)}.
$

\paragraph{Examples.} 
1. Consider the motivic correlator 
 corresponding to the tensor product of the motives of Hecke 
eigenforms $\omega_1, \ldots \omega_4$ as on Fig \ref{feyn32}. 
we conjecture that there are vectors $v_{\omega_i} \in V_{\omega_i}$ providing non-zero elements
$$
{\rm Cor}^0_{\cal M}{\cal C}(\bigotimes_{i=1}^4M^{\vee}_{\omega_i} \otimes v_{\omega_i} )(1) \in {\rm Ext}^1_{{\cal M}{\cal M}}(\Q(0), \otimes_{i=1}^4M_{\omega_i}(-1)).
$$

2. As discussed in Section 1.10, there are two types of cyclic elements, see  (\ref{4.3.10.1}) and 
Fig \ref{feyn33a},  
which should be related to $L(S^2M_{\omega}, 3)$.
Let us consider first the motivic correlator map 
assigned to the left one:
$$
{\rm Cor}^0_{\cal M}: 
{\cal C}\Bigl({\rm Meas}(\widehat 
{\rm Cusps})(-1)^{\otimes 2} 
\otimes S^2(M^{\vee}_{\omega} \otimes V_{\omega})(1)\Bigr) \lra {\cal L}_{\cal M}.
$$ 
By the Manin-Drinfeld theorem the components 
of the coproduct corresponding to the 
$s$-decorated vertices are zero. 
Furthermore, the only nontrivial  component of the coproduct 
is the one corresponding to the cut shown on Fig \ref{feyn33a} 
by a dotted line. It is given by  
\be \la{4.3.10.11}
\sum_i{\rm Cor}^0_{\cal M}(s_1 \otimes \omega_1 \otimes \psi_i(1)) \otimes 
{\rm Cor}^0_{\cal M}(\psi_i^\vee\otimes \omega_2 \otimes s_2  (1)) 
\ee
To formulate a criteria for this element to be zero, we proceed as  follows: 

i) Since $s_i$ are cuspidal divisors on a modular curve, 
each of the factors determines an element of ${\rm Ext}_{\rm Mot}^1$, 
e.g. ${\rm Cor}^0_{\cal M}(s_1 \otimes \omega_1 
\otimes \psi_i(1)) \in {\rm Ext}_{\rm Mot}^1(\Q(0), M_{\omega_1} \otimes M_{\psi_i})$. 

ii) We calculate the corresponding tensor product of the Hodge correlators, 
given explicitly by the Rankin-Selberg integrals:
\be \la{4.3.10.12}
\sum_i{\rm Cor}^0_{\cal H}(s_1 \otimes \omega_1 \otimes \psi_i(1)) \otimes 
{\rm Cor}^0_{\cal H}(\psi_i^\vee\otimes \omega_2 \otimes s_2  (1)) \in \C \otimes \C. 
\ee 
The conjectural injectivity of the regulator map implies that (\ref{4.3.10.12}) vanishes if and only if 
the coproduct (\ref{4.3.10.11}) is zero. 

The cyclic element provided by the right diagram on Fig \ref{feyn33a} is treated similarly. 
This completes the formulation of Conjecture \ref{9.30.13.1} relating $L({\rm Sym}^2f, 3)$ to the Hodge correlators. 
\begin{figure}[ht]
\centerline{\epsfbox{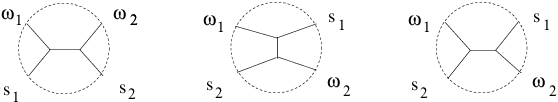}}
\caption{Correlator integrals of type $S^2M_\omega(1)$.}
\label{feyn33a}
\end{figure}

3. It would be interesting to generalize this picture, 
providing an automorphic description of a part of the motivic 
Galois group of the category of mixed motives over an 
arithmetic field $F$. Here is 
an example provided by the Drinfeld modular curves 
\cite{Dr2} 
in the function field case.

\paragraph{$l$-adic correlators for the universal Drinfeld modular curve.} 
Let $F$ be the function field of a curve $X$ over a finite
field. Choose a point $x$ on $X$ and take the Drinfeld
modular curve  
parametrising the rank $2$ elliptic 
modules corresponding to the affine curve $X-x$ \cite{Dr2}. 
Then there is the same kind of picture as above.  Namely, 
one should have a conjectural abelian 
category of mixed motives over $F$. One can define its $l$-adic version: 
a category of motivic  $l$-adic Galois representations of  the Galois group 
${\rm Gal}(\overline F/F)$, where $l \not = {\rm char}(F)$. 
An example is provided by  the $l$-adic  fundamental Lie algebra of the Drinfeld modular curve. 
It is a free pronilpotent Lie algebra generated by $H_1$ of the curve. 
There is a decomposition similar to 
(\ref{sdds}). 
Thus, just as above,  we arrive at the $l$-adic correlator map.

\section{Feynman integral for the Hodge correlators} \la{hc12sec}

\subsection{The Feynman integral} Let $\varphi$ be a smooth function 
on a complex curve $X(\C)$ with values in $N\times N$ complex matrices ${\rm Mat}_N(\C)$. We say that 
$\{\varphi\}$ is our space of  fields. 
Let us choose a cyclic word
\begin{equation} \label{ser}
W = {\cal C}\Bigl(
\{a_0\} \otimes \omega^{0}_{1} \otimes \ldots \otimes \omega^{0}_{n_0} \otimes 
\ldots 
\otimes 
\{a_m\} \otimes \omega^{m}_{1} \otimes \ldots \otimes 
\omega^{m}_{n_m}\Bigr) \in {\cal C}{\rm T}({\rm V}^{\vee}_{X,S})
\end{equation} 
where $a_i$ are points of $X(\C)$, not necessarily different, and 
$\omega^{i}_{j} \in  \Omega^1_X 
\oplus \overline \Omega^1_X$.   We are going to assign to $W$ a functional 
${\cal F}_W(\varphi)$ on the space of fields. 

There are the following elementary matrix-valued functionals on the space of fields:
\vskip 2mm
\noindent
(i) Each point $a \in X(\C)$ provides a functional
\begin{equation} \label{3.13.05.1}
{\cal F}_a(\varphi):= \varphi (a) \in {\rm Mat}_N(\C).
\end{equation}
(ii) A $1$-form $\omega$ on $X(\C)$ provides 
a functional
\begin{equation} \label{3.13.05.3}
{\cal F}_{\omega}(\varphi):= 
\int_{X(\C)}[\varphi, d^{\C}\varphi]\wedge \omega \in {\rm Mat}_N(\C).
\end{equation}
(iii) A pair of $1$-forms $(\omega_1, \omega_2)$ on $X(\C)$ provides 
a functional
\begin{equation} \label{3.13.05.2}
{\cal F}_{(\omega_1, \omega_{2})}(\varphi):= 
\int_{X(\C)}\varphi(x)\omega_1\wedge \omega_2 \in {\rm Mat}_N(\C).
\end{equation}

Let us choose a 
collection ${\cal P}$, possibly empty, of consecutive pairs 
$\{\omega^i_j, \omega^i_{j+1}\}$ of $1$-forms entering the cyclic word $W$. 
Take a field $\varphi$, and 
go along $W$, assigning to every element $\{a_i\}$ a matrix 
 ${\cal F}_{a_i}(\varphi)$, to every 
$1$-form $\omega$ which does not enter to ${\cal P}$ a matrix 
 ${\cal F}_{\omega}(\varphi)$, and to every pair of forms 
$\{\omega^i_j, \omega^i_{j+1}\}$ from the  
collection ${\cal P}$ a matrix 
 ${\cal F}_{(\omega^i_j, \omega^i_{j+1})}(\varphi)$. 
This way for a given collection ${\cal P}$  and a given field $\varphi$ we 
obtain a cyclic word in ${\rm Mat}_N(\C)$, so its  
trace is well defined, providing a  complex-valued 
functional, denoted by 
${\cal F}_{W, {\cal P}}(\varphi)$. 
The functional ${\cal F}_{W}(\varphi)$ is obtained by taking 
the sum over all possible collections ${\cal P}$:
$$
 {\cal F}_{W}(\varphi)= \sum_{{\cal P}}{\cal F}_{W, {\cal P}}(\varphi).
$$

{\bf Example}. Let $W = {\cal C}
\Bigl(\{s_0\} \otimes \{s_1\} \otimes \ldots \otimes \{s_m\}\Bigr)$. 
Then we recover formula (\ref{7.28.06.1}) for ${\cal F}_W(\varphi)$. 

\noindent

Given an integer $N$, we would like to define the correlator corresponding to $W$ via 
a Feynman integral
\begin{equation} \label{3.12.05.3}
{\rm Cor}_{X, N, h, \mu}(W):= 
\int {\cal F}_W(\varphi)
e^{iS(\varphi)}{\cal D}\varphi
\end{equation}
$$
S(\varphi):= \frac{1}{2\pi i}\int_{X(\C)} {\rm Tr}
\Bigl( \frac{1}{2}\partial  \varphi \wedge \overline \partial\varphi + 
\frac{1}{6}\hbar \cdot 
\varphi [\partial \varphi, \overline \partial \varphi]  \Bigr).
$$
Observe that $S(\varphi)$ is real: $\overline {S(\varphi)} = S(\overline \varphi)$. 
Unfortunately  formula (\ref{3.12.05.3}) does not have a precise mathematical sense. 
We understand it by postulating the perturbation series 
expansion with respect to a small parameter $\hbar$, 
using the standard Feynman rules, and then taking the leading 
term in the asymptotic expansion as $N \to \infty, \hbar=N^{-1/2}$. Before the implementation of this plan, let 
us  recall the finite dimensional situation.   

\subsection{Feynman rules} 

This subsection surveys  well-known results, mostly 
due to Feynman and t'Hooft,  
and is included for the convenience of the reader only. 

\paragraph{Feynman rules for finite-dimensional integrals.} We follow 
\cite{Wit}, Section 1.3. 
Let $V$ be a finite dimensional real vector space and $B(v,v) \in S^2(V^*)$ 
a positive definite quadratic bilinear symmetric form on $V$. 
Denote by $dv$ an invariant  Lebesgue measure on $V$ normalized by  
$$
\int_Ve^{-B(v,v)/2}dv=1.
$$
Let us consider a polynomial function on $V$ given by an expansion 
$$ 
Q(v) = \sum_{m\geq 1}g_mQ_m(v)/m!, \qquad Q_m(v)\in S^m(V^*)
$$ 
where $g_m$ are formal variables. Let $f_i(v) \in V^*$. Consider the following 
integral, the {\it correlator} of $f_1, \ldots, f_N$:
$$
\langle f_1, \ldots, f_N\rangle:= \int_Vf_1(v) \ldots 
f_N(v)e^{-B(v,v)/2+Q(v)}dv \in \C[[g_1,g_2,  \ldots]],
$$
where we expand $e^{Q(v)}$ into a series, and then calculate each term separately. 
This integral is computed using Feynman graphs as follows.  
Let ${\bf n}:=\{n_1, n_2, \ldots \}$ be any sequence of nonnegative 
integers which is eventually zero. Let $G(N, {\bf n})$ be the set of 
equivalence classes of graphs which have $N$ $1$-valent vertices labeled 
by $1, \ldots, N$, and $n_i$ unlabeled $i$-valent vertices, $i\geq 1$. 
The labeled vertices are called {\it external}, 
and unlabeled ones {\it internal} vertices. 

To each internal vertex $v$ of $\Gamma$ of valency $|v|$ we assign 
a tensor $Q_{|v|}$. To each external vertex labeled by $i$ we assign a vector $f_i$. Taking the product over all vertices of $\Gamma$, we get a vector 
$$
\prod_{i = 1}^Nf_{i} \prod_{\mbox{$v$}}Q_{|v|} 
$$ 
in the tensor product of 
$V^*$'s over the set of all flags of $G$. Further, let $B^{-1} \in S^2V$ 
be the inverse form to $B$. 
To each edge $e$ of the graph $\Gamma$ we assign a tensor $B^{-1}_e$ 
called the {\it propogator}. Then 
there is a vector in the the tensor product of 
$V$'s over the set of all flags, that is pairs (a vertex $v$, an edge $E$ incident to $v$),  of $G$: 
$$
\prod_{\mbox{edges $e$ of $\Gamma$}}B^{-1}_e. 
$$
 Contracting these two vectors 
 we get a complex number 
$$
F_{\Gamma}:= \left\langle \prod_{\mbox{vertices $v$ of $\Gamma$}}Q_{|v|}, 
\prod_{\mbox{edges $e$ of $\Gamma$}}B^{-1}_e\right\rangle \in \C
$$
where for external vertices $i$ we put $Q_{|i|}:=f_i$. 
\begin{theorem} \label{3.14.05.1} a) One has 
\begin{equation} \label{3.14.05.2}
\langle f_1, \ldots, f_N\rangle = \sum_n\prod_ig_i^{n_i}
\sum_{\Gamma \in G(N, {\bf n})}|{\rm Aut}(\Gamma)|^{-1}F_{\Gamma}(f_1, \ldots, f_n)
\end{equation} 
where ${\rm Aut}(\Gamma)$ denotes the group of automorphisms of 
$\Gamma$ which fix the external vertices. 

b) One has 
\begin{equation} \label{3.14.05.2a}
\log \langle f_1, \ldots, f_N\rangle = \sum_n\prod_ig_i^{n_i}
\sum_{\Gamma \in G_0(N, {\bf n})}|{\rm Aut}(\Gamma)|^{-1}F_{\Gamma}(f_1, \ldots, f_n)
\end{equation} 
where the sum now is over the subset $G_0(N, {\bf n}) \subset G(N, {\bf n})$ of connected graphs only.  
\end{theorem}

To deduce b) from a) observe that the sum in (\ref{3.14.05.2}) 
is over all, possibly disconnected, diagrams. 
The automorphisms of a graph are given by 
the automorphisms of the connected components, and 
the permutation groups $S_n$ of $n$ identical copies of connected components.
It remains to notice that 
$$
        \prod_i \sum_j {F_{\Gamma_i}^{n_j} \over n_j!} = \prod_i \exp(F_{\Gamma_i}) = \exp(\sum_i F_{\Gamma_i}).
$$

\paragraph{Feynman rules for the matrix model (\ref{3.12.05.3}).} 
Let us explain how to write an asymptotic expansion in $\hbar$ for the 
Feynman integral (\ref{3.12.05.3}).  The quadratic form in our case is
\begin{equation} \label{3.14.05.6}
B(\varphi)= (2\pi i)^{-1}\sum_{i,j=1}^N\frac{1}{2}
\int_{X(\C)}\partial \varphi^i_j\wedge \overline \partial \varphi^j_i 
 = -(2\pi i)^{-1}\sum_{i,j=1}^N
\int_{X(\C)}\varphi^i_j\cdot  \partial \overline \partial \varphi^j_i. 
\end{equation} 
The vector space $V$ is 
the infinite-dimensional space of ${\rm Mat}_N(\C)$-valued functions $\varphi$ on $X(\C)$. 

Let us assume first that $N=1$. 
Then the Laplacian
$\Delta = (2\pi i)^{-1}\overline \partial \partial $ has a one-dimensional kernel. 
So the inverse form $B^{-1}$ can be defined only on 
a complement to the kernel. 
Such a complement is described by a 
choice of a $2$-current $\mu$ on $X(\C)$ 
with a non-zero integral, which we normalize to be  $1$. 
It consists of functions $\varphi$ orthogonal to $\mu$. 
The  Green function $G_{\mu}(x,y)$ describes the bilinear form $B^{-1}$ 
on the complement 
as follows. The dual $V^*$ 
contains a dense subspace of smooth $2$-forms $\omega$ on $X(\C)$. 
The value of the bilinear form on 
$\omega_x\cdot \omega_y \in S^2V^*$ is $\int_{X(\C)^2}G_{\mu}(x,y)\omega_x\otimes 
\omega_y$.  
Its restriction to the subspace of $\omega$'s orthogonal to constants 
does not depend on the ambiguity in the definition of $G_\mu(x,y)$. 
\begin{figure}[ht]
\centerline{\epsfbox{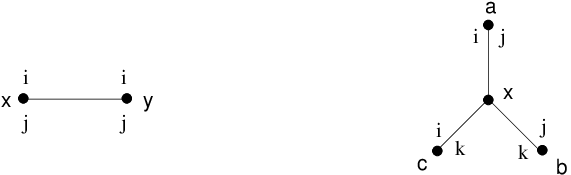}}
\caption{Feynman diagrams for the propagator and the vertex contribution}
\label{feyn11}
\end{figure}

Now let $N>1$. Let $e^i_j$ be an elementary matrix with 
the only non-zero entry, $1$, at the $(i,j)$ place. Then it is clear from 
(\ref{3.14.05.6}) that the propogator $B^{-1}$ is given by 
$$
B^{-1} = -G_{\mu}(x,y) \sum_{i,j=1}^Ne^i_j\otimes e_i^j.
$$
Its Feynman diagram is on the left of Fig \ref{feyn11}. 
Recall the cubical term of the action $S(\varphi)$: 
$$
\hbar\cdot\sum_{i,j, k=1}^N\int_{X(\C)}\Bigl(\varphi^i_j(x) \partial 
\varphi_k^j(x)  \wedge \overline 
\partial \varphi^k_i(x) - \varphi^i_j(x) \overline 
\partial \varphi_k^j(x) \wedge 
\partial \varphi^k_i(x)\Bigr). 
$$
Recall a convenient normalisation of the maps $\omega_n$ introduced in (\ref{5.16.06.10d}) 
and denoted by $\omega^*_n$. 

It follows that the contribution of the three valent vertex  shown 
 on the right of Fig \ref{feyn11} is
$$
-\hbar\cdot \int_{X(\C)}\Bigl(G_{\mu}(a,x) \partial G_{\mu}(b,x) \wedge \overline 
\partial G_{\mu}(c,x) - G_{\mu}(a,x) \overline \partial G_{\mu}(b,x) 
\wedge \partial G_{\mu}(c,x)\Bigr)  \sum_{i,j, k=1}^Ne^i_j\otimes e_k^j\otimes  e^k_i=
$$
\begin{equation} \label{3.18.05.1}
\hbar\cdot \int_{X(\C)}\omega^*_2\Bigl(G_{\mu}(a,x) \wedge  
G_{\mu}(b,x) \wedge G_{\mu}(c,x) \Bigr)\sum_{i,j, k=1}^N
e^i_j\otimes e_k^j\otimes  e^k_i.
\ee
Integral (\ref{3.18.05.1})
is convergent by Lemma 2.5. It  does not change if we add a constant to $G_{\mu}$.

\vskip 3mm

Recall that a ribbon graph is a graph with an additional structure: a cyclic 
order of the edges sharing each vertex. 
The contraction procedure for such propogator and vertex contributions implies that the non-zero 
contributions are given by ribbon graphs only  
(Fig \ref{feyn12}). Each ribbon graph $\Gamma$ contributes a Feynman integral 
entering with a certain weight $w_{\Gamma}$, 
calculated below. 
 
\begin{figure}[ht]
\centerline{\epsfbox{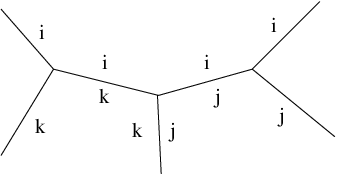}}
\caption{Non-zero 
contributions are given by ribbon graphs only.}
\label{feyn12}
\end{figure}

A ribbon  graph provides a surface $S_\Gamma$. 
It is obtained by taking a compact surface without boundary  
$\overline S_\Gamma$ and removing from it several discs. 
The graph $\Gamma$ is embedded into the surface $S_\Gamma$, 
so that the external vertices are on the boundary. It provides a decomposition 
of $\overline S_\Gamma - \Gamma$ 
into  polygons, called faces. 
Let $f_{\Gamma}$, $e_\Gamma$ and $v_\Gamma$  be the number of polygons, edges and 
vertices of $\Gamma$. So each Feynman diagram $\Gamma$ 
contributes an integral with the weight 
$
\omega_\Gamma = N^{f_{\Gamma}}\hbar^{v_{\Gamma} - (m+1)}.
$  

The Euler 
characteristic of the surface $\overline S_\Gamma$ is 
$
\chi_\Gamma = v_\Gamma+f_\Gamma-e_\Gamma.
$  
Let us shrink each boundary component of the surface $S_\Gamma$ into a point. 
The resulting surface is homeomorphic to $\overline S_\Gamma$. It has inside  a ribbon graph 
with all vertices of valency $3$ except the shrinked points, whose total valency is $m+1$.
Thus  $(m+1)+ 3v_\Gamma = 2e_\Gamma$. So 
$
\chi_\Gamma = f_\Gamma-v_\Gamma/2 - (m+1)/2
$.

\subsection{Hodge correlators and the matrix model} 
\begin{definition} \label{3.12.05.1} 
The correlator ${\rm Cor}_{X, \mu}(W)$ is the leading term 
of the 
asymptotics of ${\rm Cor}_{X, N, \mu}(W)$ as $N \to \infty$ and  $\hbar=N^{-1/2}$. 
\end{definition}

The number of ribbon graphs decorated by a given 
cyclic word $W$, even if we restrict our attention to connected diagrams, 
 is infinite: one could have a surface $S_\Gamma$ of arbitrary genus, 
and even if the genus is zero, one could have as many loops in $\Gamma$ as we want.
  However Lemma \ref{1.10.05.1w}  tells that 
all but finitely many of them give zero contribution to ${\rm Cor}_X(W)$. 
Precisely, only trees will contribute, and, of course, 
there are finitely many trees with given external vertices. 

\begin{lemma} \label{1.10.05.1w}
The logarithm $\log {\rm Cor}_{X, \mu}(W)$ of the Feynman integral 
correlator  is the sum 
of integrals assigned
to  connected  
Feynman diagrams, which are plane trees. 
\end{lemma}

\begin{proof} 
Since $\hbar=N^{-1/2}$, the finite-dimensional integral 
corresponding to a Feynman diagram $\Gamma$ enters 
with the weight 
$
\omega_\Gamma = N^{f_\Gamma-v_\Gamma/2+(m+1)/2} = N^{(m+1) + \chi_\Gamma}. 
$
Here $m$ is fixed by the cyclic word $W$. 
Thus the leading term of the asymptotics 
comes from the 
Feynman diagrams $\Gamma$ with the minimal possible value of $\chi_\Gamma$, that is plane trees. 
\end{proof}

\begin{theorem} \label{1.10.05.1}
The logarithm $\log {\rm Cor}_{X, \mu}(W)$ of the Feynman integral correlator
  equals, up to a sign computed below,  to the value of the Hodge correlator 
${\rm Cor}^*_{\cal H}(W)$ on $W$. 
\end{theorem}

\begin{proof} 
We need to show that the contribution of a $W$-decorated plane trivalent tree $T$  
to the Hodge correlator is the same as the finite dimensional Feynman integral $\varphi_T$ assigned 
by the Feynman rules to  the Feynman diagram 
corresponding to $T$. The contribution to $\varphi_T$ of a $3$-valent vertex $v$ of $T$ 
is $N^3$ times the integral 
\be \la{3.17.10.1}
-\int_{X(\C)}\Bigl(G_{\mu}(a,x) \partial G_{\mu}(b,x) \wedge \overline 
\partial G_{\mu}(c,x) + G_{\mu}(a,x) \overline \partial G_{\mu}(b,x) 
\wedge \partial G_{\mu}(c,x)\Bigr) =
\ee
$$
\int_{X(\C)}\omega^*_2\Bigl(G_{\mu}(a,x) \wedge  
G_{\mu}(b,x) \wedge G_{\mu}(c,x) \Bigr).
$$
Integral (\ref{3.17.10.1}) is skewsymmetric under the permutations of 
$(a,b,c)$. Indeed, it is evidently skewsymmetric under $b\leftrightarrow c$. 
To check skewsymmetry under $a\leftrightarrow b$, notice that the skewsymmetrisation leads to
 a complete derivative: 
$$
f(\partial g \wedge \overline \partial h) - 
f(\partial h\wedge \overline \partial g ) + g(\partial f \wedge \overline \partial h) 
-g(\partial h\wedge \overline \partial f) 
= \partial (fg\overline \partial h)+ 
\overline\partial (fg \partial h).
$$
Let us use  this symmetry to find a preferred presentation for 
the integrand of the Feynman integral $\varphi_T$ related to $T$. 

Recall that a flag at a vertex $v$ is a pair (the vertex $v$, an edge incident to $v$). 
So each internal vertex $v$ of $T$ gives rise to three flags. 
We say that among the three flags $(x,a), (x,b), (x,c)$ relevant to integral (\ref{3.17.10.1}) 
the flag $(x,a)$ is the {\it free flag at the vertex $x$} --
it contributes the Green function 
$G_\mu(a,x)$ rather then its derivative. 
Thanks to the symmetry, for any $3$-valent vertex $v$ of $T$,  
in the contribution to the Feynman integral $\varphi_T$ 
we can make any flag incident to $v$ to be the free flag.  
Choose an external edge of $T$, called the {\it root edge}. 
Then for each internal vertex there is a unique flag which is the closest to the root. 
We call it the {\it upper flag}. Using the symmetry of the integral (\ref{3.17.10.1}) one easily proves the following:

\bl Given a root edge of $T$, there is a unique presentation of the integrand for the integral $\varphi_T$
 for which the set of free flags is the set of upper flags.  
\el

On the other hand this contribution can be written as 
\be \la{Mcccccc}
2^{m-1} \int_{X(\C)^{m-1}}G_0 (\partial G_1 \wedge \overline \partial G_2) \wedge ... \wedge 
(\partial G_{2m-3} \wedge \overline \partial G_{2m-2}),  
\ee
where the integration is over the copies of $X$ assigned to the internal vertices of $T$, and 
$(G_{2k+1}, G_{2k+2})$ is 
the pair of Green functions 
assigned to the pair of non-upper edges of $k$-th vertex of $T$. 
We claim  that 
 $$
(\ref{Mcccccc}) = (-1)^{(m-1)(m-2)/2}\int_{X(\C)^{m-1}}\omega^*_{2m-2}(G_0, \ldots,  G_{2m-2}). 
$$
By the very definition, the integrands coincide up to a numerical coefficient. 
It is easy to see that the coefficient must be $\pm 1$.  
\end{proof}


\paragraph{A generalization of the Feynman integral.} 
Let ${\cal G}$ be a Lie algebra with an invariant symmetric 
bilinear form 
$Q: {\cal G} \otimes {\cal G}\lra \C$. Then there is a generalisation 
of the Lagrangian for ${\rm Mat}_N(\C)$ defined as follows. 
Consider ${\cal G}$-valued fields $\{\varphi\}$ on $X(\C)$. 
Then there is an action 
$$
S_{\cal G}(\varphi):= \frac{1}{2\pi i}\int_{X(\C)} 
Q\Bigl( \frac{1}{2}\partial\varphi \wedge \overline \partial  \varphi + 
\frac{1}{6}\hbar \cdot 
\varphi \otimes [\partial \varphi, \overline \partial \varphi]  \Bigr).
$$
Here the expression in parenthesis is a 
${\cal G}\otimes {\cal G}$-valued $2$-form on $X(\C)$, so applying 
the bilinear form $Q$ we get a $2$-form and integrate it over $X(\C)$. 
Formulas (\ref{3.13.05.1})-(\ref{3.13.05.3}) provide ${\cal G}$-valued functionals 
on the space of fields corresponding to a choice of a point $a$, a $1$-form, 
and a pair of $1$-forms on $X(\C)$. Repeating verbatim our construction we  
assign to a cyclic word $W$ a correlator

\begin{equation} \label{3.12.05.3x}
{\rm Cor}_{X, {\cal G}, h}(W):= 
\int {\cal F}_W(\varphi)
e^{iS_{\cal G}(\varphi)}{\cal D}\varphi.
\end{equation}
In the case when ${\cal G} = {\rm Mat}_N$, we recover the Feynman integral (\ref{3.12.05.3}).

The motion equation for the action is $\delta S_{\cal G}(\varphi)/\delta\varphi  =0$. 
Renormalising $\varphi \lms \varphi/\hbar$ to remove dependence on $\hbar$, we get   
$
\overline \partial\partial  \varphi  =  
[\overline \partial \varphi, \partial \varphi]  
$.

\end{document}